%% file: main.tex
\documentclass[
  paper = 17x24,
  language = english,
  acronym = nosymbol,
  acronymline = novertical,
  bibliography = combined,
  bibliographypart = all,
  titlesize = Huge,
  par = halfskip,
]{faupress}

\usepackage{amsmath}
\usepackage{amssymb}
\usepackage{blindtext}
\usepackage{physics}
\usepackage{bm}
\usepackage{relsize}

\usepackage{amsthm}
\theoremstyle{definition}
\newtheorem{theorem}{Theorem}
\newtheorem{definition}{Definition}
\usepackage{caption}
\usepackage{subcaption}
\usepackage{graphicx}
\usepackage{algorithmicx}
\usepackage{algorithm}
\usepackage{algpseudocode}
\usepackage{tikz}
\usepackage{booktabs}
\usepackage{backnaur}
\usepackage{hyperref}
\usepackage[titles]{tocloft}

\newcommand{\ps}[1]{\langle #1 \rangle}
\newcommand{\bps}[1]{ \ps{\bm{#1}} }

\colorlet{lightred}{red!33}
\colorlet{lightblue}{blue!33}

\numberwithin{equation}{chapter}
\newcounter{prodcounter}
\newenvironment{production}{%
\refstepcounter{prodcounter}

\begin{bnf}}
{\end{bnf}}

\setminted[python]{fontsize=\footnotesize, breaklines, frame=lines
, numbersep=7.5pt, numberblanklines=False}
\setminted[scala]{fontsize=\footnotesize, breaklines, frame=lines, numbersep=7.5pt, numberblanklines=False}
\usepackage{etoolbox}

\AtBeginEnvironment{tabular}{\small}

\title{Automating the Design of Multigrid Methods with Evolutionary Program Synthesis}
\subtitle{Automatisierung des Entwurfs von Mehrgitterverfahren mit Evolutionärer Programmsynthese}

\firstname{Jonas}
\lastname{Schmitt}
\origin{Forchheim}

\yearofpublication{2024}
\series{}
\volume{}
\doi{}
\isbn{}
\eisbn{}
\issn{}
\printinformation{%
}

\subject{Doktorarbeit}

\institute{Lehrstuhl für Informatik 10}
\supervisor{Prof. Dr. Harald Köstler}
\oralexam{14.12.2023}
\dean{Prof. Dr. Gerhard Wellein}
\reviewer{Prof. Dr. Harald Köstler \\ 
Prof. Dr. Penousal Machado \\
Prof. Dr. Dietmar Fey}

\begin{document}

% typeset the titlepage of the publisher
%\maketitle

% start the frontmatter (roman page numbering,
% lowercase roman sectioning numbering, no header)
\frontmatter
  
  % typesets the titlepage of the faculty
  \makefacultytitle

  % insert preface heading (has \label{ch:preface},
  % reference via \nameref{ch:preface})
  \begin{abstract}
    \input{contents/abstract.tex}
  \end{abstract}
  \begin{zusammenfassung}
    \input{contents/zusammenfassung.tex}
  \end{zusammenfassung}

  % typeset the table of contents
  \tableofcontents

% start the mainmatter (arabic page numbering [reset],
% arabic sectioning numbering, header)
\mainmatter

  % inserts introduction heading (hat \label{ch:intro},
  % reference via \autoref{ch:intro} or \ref{ch:intro})
\chapter{Introduction}
\input{contents/introduction.tex}
    %%% WRITE YOUR INTRODUCTION DIRECTLY HERE OR %%%%%%%%
    %%% INPUT AN EXTERNAL FILE WITH YOUR INTRODUCTION %%%
    %\input{contents/my_introduction.tex}

\chapter{Multigrid Methods for Solving Partial Differential Equations}
  %%% WRITE YOUR THESIS DIRECTLY HERE OR %%%%%%
  %%% INPUT EXTERNAL FILES WITH YOUR THESIS %%%
  \input{contents/discretization_of_pdes.tex}

  \input{contents/basic_iterative_methods.tex}

  \input{contents/multigrid_methods.tex}
\chapter{Formal Languages and Evolutionary Program Synthesis}
\label{chapter:formal-languages-and-gp}
  \input{contents/formal_languages.tex}
  \input{contents/genetic_programming.tex}
\chapter{A Formal Language for Expressing Multigrid Methods}
\label{chapter:multigrid-formal-language}
  \input{contents/multigrid_grammar.tex}

\chapter{Automated Multigrid Solver Design with EvoStencils -- Part 1: Core Implementation}
\label{chapter:evostencils-1}
  \input{contents/evostencils_part1.tex}
\chapter{Automated Multigrid Solver Design with EvoStencils -- Part 2: Generalization and Parallelization}
  \label{chapter:evostencils-2}
  \input{contents/evostencils_part2.tex}
\chapter{Experiments and Discussion}
\label{chapter:experiments}
  \input{contents/experiments.tex}

\chapter{Related Work and Conclusion}
  \input{contents/related_work.tex}
  \input{contents/conclusion.tex}

\appendix 
  %%% WRITE YOUR APPENDIX DIRECTLY HERE OR %%%%%%%%
  %%% INPUT AN EXTERNAL FILE WITH YOUR APPENDIX %%%
\input{contents/appendix.tex}

% start the backmatter (arabic page numbering [continued],
% no sectioning numbering, header)
\backmatter
  \faupressprintbibliography
   \faupressprintacronyms
  \begingroup
  \listofalgorithms
  \let\clearpage\relax
  \listoflistings
  \endgroup

\end{document}

%% file: contents/abstract.tex
Many of the most fundamental laws of nature can be formulated as partial differential equations (PDEs). 
Understanding these equations is, therefore, of exceptional importance for many branches of modern science and engineering. 
However, since the general solution of many PDEs is unknown, the efficient approximate solution of these equations is one of humanity's greatest challenges.
While multigrid represents one of the most effective methods for solving PDEs numerically, in many cases, the design of an efficient or at least working multigrid solver is an open problem.
This thesis demonstrates that grammar-guided genetic programming, an evolutionary program synthesis technique, can discover multigrid methods of unprecedented structure that achieve a high degree of efficiency and generalization.
For this purpose, we develop a novel context-free grammar that enables the automated generation of multigrid methods in a symbolically-manipulable formal language, based on which we can apply the same multigrid-based solver to problems of different sizes without having to adapt its internal structure.
Treating the automated design of an efficient multigrid method as a program synthesis task allows us to find novel sequences of multigrid operations, including the combination of different smoothing and coarse-grid correction steps on each level of the discretization hierarchy.
To prove the feasibility of this approach, we present its implementation in the form of the Python framework EvoStencils, which is freely available as open-source software.
This implementation comprises all steps from representing the algorithmic sequence of a multigrid method in the form of a directed acyclic graph of Python objects to its automatic generation and optimization using the capabilities of the code generation framework ExaStencils and the evolutionary computation library DEAP.
We furthermore describe how this implementation can be extended to yield multigrid methods that can efficiently solve multiple instances of the same PDE, thus achieving strong generalizability.
Even though generalization is one of the main goals of automated solver generation, artificial intelligence-based methods often fail to achieve it. 
While machine learning models have shown promise in replacing classical numerical solvers, they typically rely on a fixed-size neural network and can thus not easily be generalized to other problem sizes.
To speed up the evaluation of a large number of multigrid-based solvers, we derive a suitable distributed parallelization scheme based on the message-passing interface (MPI) that allows EvoStencils to leverage the computational power of modern clusters and supercomputers.
To investigate the effectiveness of our approach, we consider several different PDEs, including the indefinite Helmholtz equation, for which we obtain multigrid methods that achieve superior solving efficiency compared to classical multigrid cycles.
Moreover, some of the methods discovered with our approach are able to achieve convergence in the case of an extremely ill-conditioned Helmholtz problem, for which, at the same time, all known multigrid cycles fail to yield a converging solver.
Within our experiments, we also show that our implementation can be executed on recent clusters and supercomputers, such as SuperMUC-NG, currently one of Europe's largest supercomputing systems.
Finally, since our formal representation of multigrid methods can easily be translated into a human-readable format, we also perform an empirical analysis of the algorithmic features discovered with our evolutionary program synthesis approach.

%% file: contents/zusammenfassung.tex
Viele der grundlegendsten Naturgesetze können als partielle Differentialgleichungen (PDGs) formuliert werden. 
Das Verständnis dieser Gleichungen ist daher für viele Bereiche der modernen Wissenschaft und Technik von immenser Bedeutung. 
Da jedoch die allgemeine Lösung vieler PDGs unbekannt ist, stellt die effiziente Näherungslösung dieser Gleichungen eine der größten Herausforderungen der Menschheit dar.
Obwohl Mehrgitterverfahren eine der effektivsten Methoden zur numerischen Lösung von PDGs darstellen, ist der Entwurf eines effizienten oder zumindest funktionierenden Mehrgitterlösers in vielen Fällen ein offenes Problem.
In dieser Arbeit wird gezeigt, dass grammatikgeleitete genetische Programmierung, eine evolutionäre Programmsynthesetechnik, zur Entdeckung von Mehrgitterverfahren bisher unerreichter Struktur führen kann, welche zudem einen hohen Grad an Effizienz und Generalisierung erreichen.
Zu diesem Zweck entwickeln wir eine neuartige kontextfreie Grammatik, welche die automatisierte Generierung von Mehrgitterverfahren in einer symbolisch manipulierbaren formalen Sprache ermöglicht, auf deren Grundlage wir denselben mehrgitterbasierten Löser auf Probleme unterschiedlicher Größe anwenden können, ohne seine interne Struktur anpassen zu müssen.
Die Behandlung des automatisierten Entwurfs effizienter Mehrgitterverfahren als Programmsyntheseproblem erlaubt es uns neuartige Sequenzen von Mehrgitteroperationen zu finden, einschließlich der Kombination von verschiedenen Glättungs- und Grobgitterkorrekturschritten auf jeder Ebene der Diskretisierungshierarchie.
Um die Machbarkeit dieses Ansatzes zu beweisen, stellen wir seine Implementierung in Form des Python-Frameworks EvoStencils vor, das als Open-Source-Software frei verfügbar ist.
Diese Implementierung umfasst alle Schritte von der Darstellung der algorithmischen Sequenz eines Mehrgitterverfahrens in Form eines gerichteten azyklischen Graphen bestehend aus Python-Objekten bis hin zu seiner automatischen Generierung und Optimierung unter Verwendung der Fähigkeiten des ExaStencils-Frameworks zur Codegenerierung und der Bibliothek DEAP für die Implementierung evolutionärer Algorithmen.
Darüber hinaus beschreiben wir, wie diese Implementierung erweitert werden kann, um Mehrgittermethoden zu erhalten, die mehrere Instanzen derselben PDG effizient lösen können, wodurch eine starke Generalisierbarkeit erreicht werden kann.
Obwohl die Verallgemeinerung eines der Hauptziele bei der automatischen Generierung von Lösern ist, scheitern Methoden, die auf künstlicher Intelligenz basieren, oft daran diese zu erreichen. 
Zwar haben sich Modelle des maschinellen Lernens in einigen Fällen als vielversprechend bei der Ersetzung klassischer numerischer Löser erwiesen, doch basieren diese in der Regel auf einem neuronalen Netzwerk fester Größe und können daher nicht ohne weiteres auf andere Problemgrößen verallgemeinert werden.
Um die Evaluierung einer großen Anzahl von mehrgitterbasierten Lösern zu beschleunigen, leiten wir ein geeignetes verteiltes Parallelisierungsschema ab, das auf der Message-Passing-Schnittstelle (MPI) basiert und es EvoStencils ermöglicht, die Rechenleistung moderner Cluster und Supercomputer zu nutzen.
Um die Effektivität unseres Ansatzes zu untersuchen, betrachten wir verschiedene PDGs, darunter die indefinite Helmholtz-Gleichung, für die wir Mehrgitterverfahren erhalten, die im Vergleich zu klassischen Mehrgitterzyklen eine höhere Lösungseffizienz erreichen.
Darüber hinaus sind einige der mit unserem Ansatz entdeckten Methoden in der Lage, Konvergenz im Fall eines äußerst schlecht konditionierten Helmholtz-Problems zu erzielen, für das gleichzeitig alle bekannten Mehrgitterzyklen keinen konvergierenden Löser liefern.
Im Rahmen unserer Experimente zeigen wir auch, dass unsere Implementierung auf neueren Clustern und Supercomputern, wie SuperMUC-NG, einem der derzeit größten europäischen Hochleistungsrechner, ausgeführt werden kann.
Da unsere formale Darstellung von Mehrgitterverfahren leicht in ein für den Menschen lesbares Format übersetzt werden kann, führen wir schließlich eine empirische Analyse der algorithmischen Merkmale durch, welche unser evolutionären Programmsyntheseansatz hervorgebracht hat.

%% file: contents/introduction.tex
%TODO The following is about general algorithm design
Many of the most fundamental laws of nature can be formulated as partial differential equations (PDEs).
%Understanding these equations is, therefore, foundational for many branches of science and engineering.
%Since, for many PDEs, it is unknown whether a general solution exists, the efficient approximate solution of these equations represents one of humanity's greatest challenges.
Since the invention of modern computers, great efforts have been made to develop efficient frameworks and programming languages for solving these equations.
As a result of this effort, computer simulations nowadays represent an essential tool for researchers and engineers.
%However, leveraging the power of simulation-based methods, in many cases, requires designing solvers that achieve the highest possible degree of efficiency.
%Unfortunately, this task often not only requires a great deal of expertise and domain knowledge that only a limited group of people possesses, but often also uncountable working hours need to be invested for its accomplishment.
However, leveraging the power of simulation-based methods, in many cases, necessitates the use of PDE solvers that achieve the highest possible degree of efficiency.
This task often not only requires a great deal of expertise and domain knowledge that only a limited group of mathematical experts possesses, but often also a lot of effort needs to be invested for its accomplishment.
All this makes the manual design and implementation of efficient PDE solvers a difficult and labor-intensive endeavor.
%The automation of manual labor has always been one of the greatest incentives of the technical revolution since the first industrial revolution.
%However, in contrast to previous technological advancements, which were mostly concerned with freeing people from physical labor, the development of computing devices with ever-increasing power and speed has led to a point where the automation even of challenging cognitive tasks has started to come into reach.
In recent decades, the development of computing devices with ever-increasing power and speed has enabled the automation of increasingly challenging cognitive tasks.
Artificial intelligence (AI) methods have demonstrated super-human performance in numerous applications, such as image processing~\cite{krizhevsky2017imagenet}, game playing~\cite{schrittwieser2020mastering,reed2022generalist}, and natural language processing~\cite{brown2020language}.
Alongside the widespread success of AI methods in these domains, techniques for automated algorithm design have achieved breakthroughs in a number of cases, such as the design of SAT-solvers~\cite{khudabukhsh2016satenstein} and mixed-integer programming~\cite{hutter2010automated}.
In general, methods for automated algorithm design can be classified into \emph{top-down}, and \emph{bottom-up} approaches.
Top-down approaches, often also called algorithm configuration methods, aim to represent an algorithm design space as a finite list of global parameters.
Finding an optimal algorithm design thus corresponds to solving a combinatorial optimization problem, sometimes including continuous parameters, which can be tackled using classical black-box optimizers like evolutionary algorithms~\cite{back1996evolutionary} and bayesian optimization~\cite{frazier2018tutorial}.
However, a severe limitation of these approaches is that they do not allow modifying individual steps of an algorithm unless each of them is represented as a distinct global parameter.
Bottom-up design methods aim to overcome these limitations by considering the construction of an algorithm from its fundamental building blocks.
An algorithm is essentially a list of statements written in a formal language.
If we consider this language to be a programming language, the formulation of an algorithm is nothing else than writing a program in that language.
Therefore, the automated design of an optimal algorithm can be treated as a program synthesis task, which gives us the same degree of flexibility available in modern programming languages.
On the downside, bottom-up algorithm design requires the formulation of a programming language for expressing the individual steps of an algorithm as formal statements.
%Furthermore, since many algorithms already include a number of parameters facilitating the application of top-down design methods, which can often be formulated irrespective of the underlying algorithm structure~\cite{hutter2007proceedings,hutter2011sequential,lopez2016irace}.
Furthermore, the greater flexibility of bottom-up algorithm design means that the number of different algorithms considered is significantly larger than in the case of top-down methods.
While all these difficulties impede the widespread application of bottom-up algorithm design, recent works in the area of machine learning~\cite{real2020automl,co2021evolving,zz_ne1,zz_ne2} and matrix multiplication~\cite{fawzi2022discovering} demonstrate that these methods have the potential to discover completely novel algorithms in different domains, a feat that is not possible with classical top-down approaches.
In the future, we can expect these methods to become feasible in even more domains as, with the projected ongoing increase in computational power, the exploration of even larger algorithm design spaces comes into reach.
In contrast to the area of automated machine learning (AutoML), where the application of automated algorithm design and configuration methods has become an active field of research~\cite{ren2021comprehensive,hutter2019automated,elsken2019neural,he2021automl,schrodi2022towards}, the application of these methods to the design of PDE solvers is a largely unexplored research field.
This thesis aims to change this situation by introducing a novel framework for the automated design of multigrid methods, a class of numerical methods that offer the potential to solve many PDEs in an asymptotically optimal way.
Since multigrid methods can only achieve this property through the correct choice and composition of their individual operations, for many PDEs, the design of an optimal or even functioning multigrid-based solver is an open problem~\cite{trottenberg2000multigrid,ernst2012difficult}.
As a remedy, the works by Oosterlee et al.~\cite{oosterlee2003genetic}, Thekale et al.~\cite{thekale2010optimizing}, and Brown et al.~\cite{brown2021tuning} represent a first step towards the automated design of these methods.
However, the authors of these papers only consider a limited configuration space, which makes the discovery of completely novel algorithmic patterns impossible.
To overcome the inherent limitations of this approach, this thesis considers the task of constructing an optimal multigrid method from its basic components as a program synthesis task.
For this purpose, a novel formal language for the automated bottom-up design of multigrid methods will be introduced.
By levering the power of evolutionary computation, it will be demonstrated that this language enables the discovery of unique sequences of multigrid operations that can yield faster solvers than classical multigrid cycles.
In the following, a basic understanding of the fundamental theory of multigrid methods, formal languages, and evolutionary program synthesis will be established first to give the reader the necessary background for the main part of this thesis.
%In general, designing an algorithm based on a set of global parameters can be summarized under the term \emph{top-down}, as the choice of each parameter has a global effect on the method's behavior and thus might affect multiple computational steps simultaneously. 
%For instance, choosing a different value for the parameter $\mu$ in Algorithm~\ref{alg:multigrid-cycle} influences the number of recursive descents on every level and thus leads to a globally different type of multigrid cycle.
%In contrast, a bottom-up approach is characterized by the ability to change each computational step within an algorithm without affecting the behavior of any other.

%% file: contents/discretization_of_pdes.tex
\section{Discretization of Partial Differential Equations}\label{sec:discretization}
Many problems in science and engineering can be modeled as partial differential equations (PDEs)~\cite{folland2020introduction,evans2010partial}.
A PDE is an equation that contains functions of one or multiple variables together with their partial derivatives.
Consider, for instance, the equation
\begin{equation}
	-\alpha \nabla^2 u = f \quad \text{in} \; \Omega,
	\label{eq:heat-equation}
\end{equation}
where $u = u(\bm{x})$ and $f = f(\bm{x})$ are both functions with respect to the vector of space variables $\bm{x} = (x_1, x_2, \dots, x_n)^T$ and $\Omega \supset \mathbb{R}^d$.
Equation~\eqref{eq:heat-equation} describes the temperature distribution inside a medium whose thermal conductivity is determined by the coefficient $\alpha$ and which contains a heat source $f$.
Since Equation~\eqref{eq:heat-equation} is only satisfied in the interior of the domain $\Omega$, we, additionally, need to define a set of conditions at its boundaries.
These so-called \emph{boundary conditions} (BCs) can usually be classified into four different types:
\begin{description}
	\item[Dirichlet] $u(\bm{x}) = g(\bm{x})$
	\item[Neumann] $\frac{\partial}{\partial \vec{n}} u(\bm{x}) = 0$
	\item[Robin] $a u(\bm{x}) + b \frac{\partial}{\partial \vec{n}} u(\bm{x}) = g(\bm{x})$
	\item[Cauchy] $a u(\bm{x}) = g(\bm{x}), \; b \frac{\partial}{\partial \vec{n}} u(\bm{x}) = h(\bm{x})$.
\end{description}
Here $\frac{\partial}{\partial \vec{n}} u(\bm{x})$ denotes the partial derivative of $u$ with respect to the outwards-directed normal vector of the boundary.
The difference between Robin and Cauchy BCs is that in the former case, one condition is formulated as a weighted average of $u$ and its derivative in the normal direction, while in the latter, two conditions must be met individually.
While depending on the boundary conditions, an analytical solution for Equation~\eqref{eq:heat-equation} might exist, for many PDEs, such a solution has not been discovered, or its computation is infeasible.
A remedy is the application of so-called numerical methods that are based on approximating the solution of a given PDE on a discrete set of points. 
%\subsection{Grid Creation}
%Before computing a numerical approximation of the solution of a given PDE, we need to define a set of discrete points within the domain $\Omega$ at which we aim to obtain the solution.
Usually, these points are defined on a \emph{grid} or \emph{mesh} with a particular structure.
In general, a distinction is made between structured (or regular) and unstructured grids.
A structured grid is characterized by the uniform neighborhood of its grid points, which means that the number of neighbors is typically the same for each grid point.
In contrast, each point within an unstructured grid can have a varying number of neighbors.
Computations are usually easier to implement and more efficient on structured grids due to their regularity.
However, structured grids are often challenging to create on complicated and irregular domains.
On the other hand, unstructured grids offer a higher degree of flexibility and are also well-suited for the previously mentioned cases~\cite{knupp2020fundamentals}.
This thesis focuses on numerical methods that can be formulated on a hierarchy of structured grids, and, in the following, we restrict ourselves to this particular grid type.
Furthermore, while there exists a wide range of different PDEs, of which many also describe time-dependent phenomena, we only discuss discretization methods that can be applied in the spatial domain.
However, note that the techniques described in this thesis can, in principle, also be utilized for the numerical solution of time-dependent PDEs, for instance, by applying them within an implicit time-stepping scheme~\cite{ames2014numerical}.
One possibility to approximate a certain (time-independent) PDE on a structured grid is to compute the Taylor series expansion around each grid point, leading to the so-called \emph{finite difference method} (FDM).

\subsection{Finite Difference Method}
For the sake of simplicity, we consider the one-dimensional function $u(x)$.
To obtain an approximation for the derivatives of $u$, we can compute its Taylor expansion in the neighborhood of $x$ with a step size $h$, which yields
\begin{equation}
	u(x + h) = u(x) + h \dv{x} u(x) + \mathcal{O}(h^2).
\end{equation}
Assuming $h$ is sufficiently small, the first-order approximation 
\begin{equation}
	\dv{x} u(x) \approx \frac{u(x + h) -  u(x)}{h}
\end{equation}
is obtained.
Furthermore, we can derive an approximation for the second-order partial derivative $\dv[2]{x}$ by considering
\begin{equation}
	u(x + h) = u(x) + h \dv{x} u(x) + \frac{h^2}{2} \dv[2]{x} u(x) + \mathcal{O}(h^3),
	\label{eq:taylor-forward}
\end{equation}
\begin{equation}
	u(x - h) = u(x) - h \dv{x} u(x) + \frac{h^2}{2} \dv[2]{x} u(x) + \mathcal{O}(h^3).
	\label{eq:taylor-backward}
\end{equation}
Adding Equation~\eqref{eq:taylor-backward} to Equation~\eqref{eq:taylor-forward} then yields the second-order finite difference approximation
\begin{equation}
	 \dv[2]{x} u(x) \approx \frac{u(x + h) + u(x - h) - 2u(x)}{h^2}.
\end{equation}
Using the same technique similar approximation terms can be obtained for higher-dimensional functions and higher-order derivatives~\cite{strikwerda2004finite}.
While finite differences offer a simple and straightforward way to approximate a given PDE on a set of structured grid points, in many cases the underlying physical requirements and complex geometries necessitate the use of semi-structured or even unstructured grids together with more complicated discretization approaches such as the finite volume (FVM) and finite element method (FEM)~\cite{versteeg2007introduction,zienkiewicz2005finite}
\subsection{Model Problem}
To illustrate the approach described in the previous section, we consider the following model problem, which represents a two-dimensional version of Poisson's equation:
\begin{equation}
	\begin{split}
		-\frac{\partial^2}{\partial x^2} u(x,y) - \frac{\partial^2}{\partial y^2} u(x,y) = f(x, y) \quad & \forall x, y \in (0, 1) \\
		u(0, y) = u(x, 0) = u(1, y) = u(x, 1) = 0 \quad & \forall x, y \in (0, 1)
	\end{split}
	\label{eq:2D-poisson-model}
\end{equation}
Note that this equation corresponds to the two-dimensional steady-state heat equation with constant $\alpha = 1$ on the unit square $\Omega = ( 0, 1 )^2$.
We then discretize Equation~\eqref{eq:2D-poisson-model} using finite differences and a uniform step size $h$, which yields
\begin{equation}
	\begin{split}
		\frac{1}{h^2} (4 u_{i,j} - u_{i-1, j} - u_{i+1, j} - u_{i, j-1} - u_{i, j+1}) = f_{i, j} \quad & \forall i, j \in \{1, 2, \dots, n\} \\
		u_{0, j} = u_{i, 0} = u_{n+1, j} = u_{i, n+1} = 0 \quad & \forall i, j \in \{1, 2, \dots, n\},
	\end{split} 
	\label{eq:2D-poisson-model-discrete}
\end{equation}
with $u_{i,j} = u(ih, jh)$, $f_{i,j} = f(ih, jh)$ and $n = 1/h - 1$.
By defining a unique ordering of the grid points $u_{i, j}$, Equation~\eqref{eq:2D-poisson-model-discrete} can be formulated as a system of linear equation of type $A \bm{u} = \bm{f}$. 
For instance, setting $h = 0.25$, results in $n = 3$ and a total number of nine grid points.
%TODO draw grid
A natural ordering of the grid points then yields the system of linear equations
\begin{equation}
\frac{1}{h^2} \underbrace{ \begin{pmatrix}
\begin{array}{ccc|ccc|ccc}~4&-1&~0&-1&~0&~0&~0&~0&~0\\-1&~4&-1&~0&-1&~0&~0&~0&~0\\~0&-1&~4&~0&~0&-1&~0&~0&~0\\\hline -1&~0&~0&~4&-1&~0&-1&~0&~0\\~0&-1&~0&-1&~4&-1&~0&-1&~0\\~0&~0&-1&~0&-1&~4&~0&~0&-1\\\hline ~0&~0&~0&-1&~0&~0&~4&-1&~0\\~0&~0&~0&~0&-1&~0&-1&~4&-1\\~0&~0&~0&~0&~0&-1&~0&-1&~4\end{array}
	\end{pmatrix}}_{\textstyle A}
\underbrace{
	\begin{pmatrix}
	u_{11} \\ u_{12} \\ u_{13} \\ u_{21} \\ u_{22} \\ u_{23} \\ u_{31} \\ u_{32} \\ u_{33}
\end{pmatrix}}_{\textstyle{\bm{u}}} = \bm f,
\label{eq:2D-poisson-assembled-matrix}
\end{equation}
with a right-hand side $\bm f$. 
By incorporating the given Dirichlet BCs, which assume a constant value of zero at all boundary points, into the right-hand side $\bm{f}$, we obtain
\begin{equation}
\bm f = \begin{pmatrix}
		f_{11} + \frac{1}{h^2} (u_{10} + u_{01}) \\f_{12} + \frac{1}{h^2} u_{02} \\ f_{13} + \frac{1}{h^2} (u_{03} + u_{14})  \\ f_{21} + \frac{1}{h^2} u_{20} \\ f_{22} \\ f_{23} + \frac{1}{h^2} u_{24} \\ f_{31} + \frac{1}{h^2} (u_{30} + u_{41}) \\ f_{32} + \frac{1}{h^2} u_{42} \\ f_{33} + \frac{1}{h^2} (u_{34} + u_{43})
\end{pmatrix} = \begin{pmatrix}
f_{11} \\
f_{12} \\ 
f_{13} \\ 
f_{21} \\ 
f_{22} \\
f_{23} \\
f_{31} \\ 
f_{32} \\
f_{33}
\end{pmatrix}.
\label{eq:2D-poisson-assembled-rhs}
\end{equation}
Note that $A$ represents a sparse matrix, which means that only a minority of its entries are nonzero.
Moreover, since we already represent the solution of Equation~\eqref{eq:2D-poisson-model-discrete} on a regular grid, storing the resulting linear system in a matrix-vector format is inefficient.
If we assume that each operation required for solving a system of linear equations can be formulated as a matrix-vector or vector-vector operation, it is unnecessary to store the corresponding matrix explicitly.
Instead, we only have to store the computational pattern that corresponds to a matrix-vector multiplication performed on a vector containing the grid points in a fixed order.
On regular grids, each such pattern can be represented as a so-called \emph{stencil code}.

\subsection{Stencil Codes}
\label{subsec:stencil-codes}
In the previous section, we have already introduced the notion of a regular grid, which can be defined more formally as
\begin{equation}
	G_{h} = \left\{ \bm{x} : \bm{x} = \bm{i} \circ \bm{h}, \bm{i} \in \mathbb{N}^d \right\},
\end{equation}
where $d \in \mathbb{N}$ is the dimensionality of the problem, and $\bm{i} \circ \bm{h}$ represents the Hadamard (or element-wise) product between two vectors $\bm{i}$ and $\bm{h}$.\footnote{Even though when used as a subscript, the step size $h$ might represent a vector, we avoid the use of bold letters for better readability in all subsequent equations.}
A grid function or variable $u_h$ is then a mapping of the form
\begin{equation}
	\begin{split}
		& u_h : G_{h}\to \mathbb{C} \\
		& \bm{x} \to u_h(\bm{x}).
	\end{split}
\end{equation}
As in the previous section, one often makes use of the simplified notation
\begin{equation}
	\begin{split}
		 & u_i = u_{h_x}(i h_x) \\
		& u_{i,j} = u_{h_x, h_y}(i h_x, j h_y) \\
		& u_{i,j,k} = u_{h_x, h_y, h_z}(i h_x, j h_y, k h_z).
	\end{split}
\end{equation}
We then define the general stencil $S$ as a finite set of $m$ tuples of the following form:

\begin{equation}
	\begin{split}
			& S = \{(\bm{a}_k, b_k) \}_{k=1}^m = \{(\bm{a}_1, b_1),  (\bm{a}_2, b_2), \dots, (\bm{a}_m, b_m)\}, m \in \mathbb{N}
	\\ & \forall \, i, j, k \in \{1, 2, \dots, m \}: \,
	\bm{a}_k \in \mathbb{Z}^d \wedge \bm{a}_i \neq \bm{a}_j \; \text{if} \; i \neq j, \; b_k \in \mathbb{C}
	\end{split}
\label{eq:stencil-definition}
\end{equation}
Here, the left entry $\bm{a}_k$ of each tuple $(\bm{a}_k, b_k)$ denotes the \emph{offset} from the index of the current grid point and the left entry $b_k$ the respective \emph{weight} or \emph{value}.
From a mathematical point of view, stencils are indistinguishable from regular sets, and thus their special properties are simply derived from the way they are constructed.
Furthermore, note that a specific stencil instance is independent of the step size $\bm{h}$, and the interpretation of each individual offset $\bm{a}_k$ depends on the grid function to which the stencil is applied.
For one and two-dimensional problems, it is often more convenient to use the alternative notations 
\begin{equation}
	S_{h_x} = \begin{bmatrix}
		\cdots & s_{-1} & s_{0} & s_{1} & \cdots
	\end{bmatrix},
\end{equation}
and
\begin{equation}
	S_{h_x, h_y} = \begin{bmatrix}
		& \vdots & \vdots & \vdots & \\
		\cdots & s_{1,-1} & s_{1,0} & s_{1,1} & \cdots \\
		\cdots & s_{0,-1} & s_{0,0} & s_{0,1} & \cdots \\
		\cdots & s_{-1,-1} & s_{-1,0} & s_{-1,1} & \cdots \\
		& \vdots & \vdots & \vdots &
	\end{bmatrix}.
\end{equation}
We can also extend this notation to three dimensions by representing all components with the same offset as a separate two-dimensional stencil, such that
\begin{equation}
	S_{h_x, h_y, h_z} = 
	\begin{bmatrix}
		\cdots & S^{(k-1)} & S^{(k)} & S^{(k+1)} & \cdots 
	\end{bmatrix}.
\label{eq:3D-stencil-matrix-notation}
\end{equation}

Furthermore, we define the application of a stencil $S$ to a grid function $u_h(\bm x)$ with $\bm x \in G_h$:
\begin{equation}
	\begin{split}
		& S \cdot u_h(\bm{x}) = \sum_{k=1}^m b_k u_h({\bm x + \bm{a}_k} \circ \bm{h}) \quad 
		\text{with} \; \bm{x} \in G_h, m \in \mathbb{N} \\ & (\bm{a}_k, b_k) \in S \; \forall \, k \in \{ 1, 2, \dots, m \}
	\end{split}
\label{eq:stencil-application}
\end{equation}
For $d = 3$, we can define the application of a three-dimensional stencil $S_{h_x, h_y, h_z}$, based on Equation~\eqref{eq:3D-stencil-matrix-notation}, as a dot product of the form
\begin{equation}
	\begin{split}
	& S_{h_x, h_y, h_z} \cdot u_{i,j,k} = 	
	\begin{bmatrix}
	\cdots & S^{(k-1)} & S^{(k)} & S^{(k+1)} & \cdots 
	\end{bmatrix} \cdot
	\begin{bmatrix}
	\vdots \\ u_{i,j,k-1} \\ u_{i,j,k} \\ u_{i,j,k+1} \\ \vdots 
	\end{bmatrix} \\
	& = \cdots \, + S^{(k-1)} \cdot u_{i,j,k-1} + S^{(k)} \cdot u_{i,j,k} + S^{(k+1)} \cdot u_{i,j,k+1} + \, \cdots
	\end{split}
\end{equation}

As an example, we consider the following five-point stencil defined on a two-dimensional regular grid with uniform step size $h_x = h_y = h$:
\begin{equation}
	\begin{split}
		-\Delta_{h,h} = & \bigg\{ \left( \left( 0,0 \right), \frac{4}{h^2}\right), \left(\left(1,0\right), \frac{-1}{h^2}\right), \left(\left(-1,0\right), \frac{-1}{h^2}\right), \\ & \left(\left(0,1\right), \frac{-1}{h^2}\right), \left(\left(0,-1\right), \frac{-1}{h^2}\right) \bigg\}
	\end{split}
	\label{eq:five-point-stencil}
\end{equation}
Applying $-\Delta_{h,h}$ to a given point $u_{i,j}$ yields 
\begin{equation}
	-\Delta_{h,h} \cdot u_{i,j} = \frac{1}{h^2} \left(4 u_{i,j}  - u_{i-1,j} - u_{i+1,j} - u_{i,j-1} - u_{i,j-1}\right),
\end{equation}
which corresponds precisely to the left part of Equation~\eqref{eq:2D-poisson-model-discrete}.
From this example, we can conclude that applying a given stencil to a grid function $u_h(\bm{x})$ defined on a regular grid can always be considered as a sparse matrix-vector product, where the matrix is obtained by sorting the grid points according to a well-defined order.
We have already illustrated this fact using the example of Equation~\eqref{eq:2D-poisson-assembled-matrix}.
Therefore, in case it is possible to define each computational step for solving a given discretized PDE on a regular grid by means of a stencil application, we only need to obtain the respective stencil entries instead of assembling or even storing a complete matrix.
Furthermore, many operators defined for matrices can similarly be formulated as a stencil operation.
As a first elementary operation, we define the multiplication of a stencil with a scalar, which is simply achieved by multiplying the weight of each entry with the respective value:
% Scalar multiplication
\begin{equation}
	\alpha S = \{(\bm{a}_k, \alpha b_k) \}_{k=1}^m, m \in \mathbb{N}
\end{equation}

Next, we can formulate certain binary operations such as the addition, subtraction, and multiplication of two stencils $A$ and $B$.
Equation~\eqref{eq:stencil-addsub} shows a recursive definition of stencil addition and subtraction.
% Stencil addition and subtraction
\begin{equation}
	A \pm B = 
	\begin{cases}
		\{(\bm{a}, b\pm c ) \} \cup (\tilde{A} \pm \tilde{B}) & \text{if} \; A = 	\{(\bm{a}, b ) \} \cup \tilde{A} \\
		& \text{and} \; B = \{(\bm{a}, c ) \} \cup \tilde{B} \\
		
		\{(\bm{a}, b ) \} \cup (\tilde{A} \pm B) & \text{if} \; A = 	\{(\bm{a}, b ) \} \cup \tilde{A} \\
		& \text{and} \; \nexists  c \in \mathbb{C} : \{(\bm{a}, c ) \} \in B 
		\\
% TODO case probably unnecessary
%		\{(\bm{a}, c ) \} \cup (A \pm \tilde{B}) & \text{if} \; B = 	\{(\bm{a}, c ) \} \cup \tilde{B} \wedge
%		\{(\bm{a}, b ) \} \notin A \; \forall b \in \mathbb{C}
%		\\
		A & \text{if} \; B = \emptyset
		\\
		B & \text{else} 
		\\
	\end{cases}
\label{eq:stencil-addsub}
\end{equation}
Here, we simply add the weights of each entry with the same offset contained in both $A$ and $B$.
In case an offset is not contained in both stencils, the respective tuple is included unmodified.  
Based on Equation~\eqref{eq:stencil-addsub}, we can then provide a definition for the multiplication of two stencils $A$, and $B$, which can be seen as a recursive version of the iterative computation scheme described in~\cite{rittich2018extending}:
% Stencil multiplication
\begin{equation}
	A \cdot B = 
	\begin{cases}
		\{(\bm{a} + \bm{\hat{a}}, b \cdot \hat{b} ) \} \; + & \text{if} \; A = \{(\bm{a}, b ) \} \cup \tilde{A} \\
		\{(\bm{a}, b ) \} \cdot \tilde{B} + \tilde{A} \cdot B & 
		\text{and} \; B = \{(\bm{\hat{a}}, \hat{b} ) \} \cup \tilde{B} \\ & \\
		\emptyset & \text{else} %\text{if} \; A = \emptyset \vee B = \emptyset
	\end{cases}
\label{eq:stencil-mult}
\end{equation}
To compute the stencil product $A \cdot B$, for each pair of entries $\{(\bm{a}, b ) \} \in A$ and $\{(\bm{\hat{a}}, \hat{b} ) \} \in B$ the tuple $\{(\bm{a} + \bm{\hat{a}}, b \cdot \hat{b} ) \}$ needs to be formed.
We then define this computation recursively by first picking two entries $\{(\bm{a}, b ) \} \in A$ and $\{(\bm{\hat{a}}, \hat{b} ) \} \in B$ and combining them as described above.
Next, we multiply the entry $\{(\bm{a}, b ) \}$ chosen from $A$ with the remainder stencil $\tilde{B} = B \setminus \{(\bm{\hat{a}}, \hat{b} ) \}$, which means that we now obtain the combination of this entry with all remaining ones from $B$.
Finally, the process is continued recursively by computing the product of the remainder $\tilde{A} = A \setminus \{(\bm{a}, b ) \}$ with the original stencil $B$.
Since it is possible that the combination of different pairs of stencil offsets leads to the same result, we, additionally, need to accumulate the corresponding weights and combine them in a single tuple.
To sum up all tuples with matching offsets obtained within subcomputations, we employ stencil addition, as defined in Equation~\eqref{eq:stencil-addsub}.
\subsection{Systems of Partial Differential Equations}
\label{subsec:systems-of-pdes}
While so far, we have only considered a single partial differential equation, many phenomena can only be modeled in the form of a system of PDEs.
One of the simplest examples of such a system is the so-called biharmonic equation
\begin{equation}
	\begin{split}
		\nabla^2 u & = v  \\
		\nabla^2 v & = f.
	\end{split}
\label{eq:biharmonic-system}
\end{equation}
Note that this system is mathematically equivalent to the scalar equation 
\begin{equation}
	\nabla^4 u = f.
\end{equation}
By utilizing matrix-vector notation, we can reformulate Equation~\eqref{eq:biharmonic-system} as
\begin{equation}
	\underbrace{
	\begin{pmatrix}
		\nabla^2 & -1 \\
		0 & \nabla^2
	\end{pmatrix}}_{A}
\underbrace{ 
	\begin{pmatrix}
		u \\ v
	\end{pmatrix}
}_{\bm{u}}
=
\underbrace{
\begin{pmatrix}
	0 \\ f
\end{pmatrix}
}_{\bm{f}}.
\label{eq:biharmonic-system-matrix-formulation}
\end{equation}
Even though Equation~\eqref{eq:biharmonic-system-matrix-formulation} now includes partial derivatives of multiple variables, we can employ the same techniques to discretize each individual operator and obtain the corresponding system of linear equations. 
Consider the two-dimensional biharmonic system
\begin{equation}
	\begin{split}
		& \frac{\partial^2}{\partial x^2} u(x,y) + \frac{\partial^2}{\partial y^2} u(x,y) - v(x, y) = 0 \\
		& \frac{\partial^2}{\partial x^2} v(x,y) + \frac{\partial^2}{\partial y^2} v(x,y) = f(x, y) \quad \forall x, y \in (0, 1)^2.
	\end{split}
	\label{eq:2D-biharmonic-system}
\end{equation}
Discretizing Equation~\eqref{eq:2D-biharmonic-system} using finite differences with a uniform step size $h_x = h_y = h$ and rewriting the resulting equations in stencil form yields
\begin{equation}
		\begin{pmatrix}
			\Delta_{h, h} & -1 \\
			0 & \Delta_{h, h}
	\end{pmatrix}
		\begin{pmatrix}
			u_{i,j} \\ v_{i,j}
		\end{pmatrix}
	=
		\begin{pmatrix}
			0 \\ f_{i,j}
		\end{pmatrix} \quad
\forall i,j \in \{1, 2, \dots, n\},
	\label{eq:2D-biharmonic-system-stencil}
\end{equation}
where $\Delta_{h,h}$ is the two-dimensional five-point stencil defined in Equation~\eqref{eq:five-point-stencil}.
Equation~\eqref{eq:2D-biharmonic-system-stencil} represents a system of two linear equations which needs to be solved at every pair of grid points $u_{i,j}$ and $v_{i,j}$.
If we consider Equation~\eqref{eq:2D-biharmonic-system-stencil} at all grid points, a system of linear equations which corresponds to the discrete solution of Equation~\eqref{eq:2D-biharmonic-system} is obtained.

Before we conclude this section, it is important to mention that we have not yet discussed how each variable, in the given case $u_{i,j}$ and $v_{i,j}$, is placed on the grid.
While it is tempting to always place each variable at the same position within a grid, in certain applications, such as the incompressible Navier-Stokes equation, this can lead to unfavorable numerical properties.
As a consequence, often more complex grid placing strategies, such as so-called staggered discretizations, need to be used in practice~\cite{trottenberg2000multigrid}, which we, for the sake of brevity, do not further discuss here.
Also note that in order to solve Equation~\eqref{eq:2D-biharmonic-system}, or any other system of PDEs, we need to define suitable conditions at the boundaries of each domain on which a given variable, contained in one of the equations, is defined.
After briefly discussing numerical discretization methods that allow us to convert a given system of PDEs into a system of linear equations, we can next direct our attention to the efficient solution of these systems.
%These boundary conditions can be of the same type as described above.
%Again, for a more detailed treatment of boundary conditions for systems of PDEs, the reader is referred to %TODO insert suitable reference. 

%% file: contents/basic_iterative_methods.tex
\section{Basic Iterative Methods}
\label{sec:basic-iterative-methods}
In general, methods for solving systems of linear equations fall into two categories: Direct and iterative methods.
Direct methods are characterized by the fact that they are able to compute the exact solution of a linear system in a finite number of steps.
In contrast, iterative methods compute a series of approximations for the solution of the linear system.
Even though this series often converges to the exact solution, there is usually no guarantee that the approximations will ever reach the accuracy of a solution computed by a direct method.
Unfortunately, for many problems applying direct methods is infeasible due to their high computational complexity and memory storage requirements.
For instance, Gaussian elimination, in general, requires $\mathcal O(n^3)$ operations for solving a system of linear equations with $n$ unknowns.
Moreover, since the goal of Gaussian elimination is to transform a given matrix into upper-triangular form, the method is based on the direct manipulation of the input matrix and hence usually requires storing it explicitly.
In contrast, many iterative methods only require the computation of matrix-vector products and do not manipulate the input matrix directly.

As we have illustrated in Section~\ref{subsec:stencil-codes}, the discretization of many partial differential equations (PDEs) enables the representation of matrix-vector products as the application of a stencil code, whereas each stencil is directly derived from the discretization of a continuous differential operator.
If we assume that each stencil includes only a finite number of entries and at most one stencil per grid point needs to be stored, we can reduce the storage requirements to $\mathcal{O}(n)$, where $n$ is the total number of grid points.
If the same stencil applies to the whole domain, we even only need to store a single stencil, whose memory requirements are negligible compared to storing the grid itself.
Furthermore, due to the inevitable approximation of real numbers by floating-point numbers on computers, even the solution of a direct method is prone to numerical errors.
As a consequence, the exactness of direct methods can be undermined by these effects, and the approximations computed by an iterative method do not necessarily achieve a lower degree of numerical accuracy~\cite{higham2002accuracy}.
Assuming that we can compute an acceptable approximation using only a finite number of $m$ matrix-vector multiplications, the application of an iterative method reduces the computational complexity of solving a system of linear equations to $O(mn)$.
In the following, we first introduce basic iterative methods, such as the Jacobi and Gauss-Seidel method, and then derive fundamental statements about their convergence. 
%Finally, as multigrid methods represent the fundamental basis for this work, we discuss these methods in greater detail.   
\subsection{Jacobi and Gauss-Seidel} 
We begin our introduction of basic iterative methods by considering the general system of linear equations
\begin{equation}
	\underbrace{
	\begin{pmatrix}a_{11}&a_{12}&\cdots &a_{1n}\\a_{21}&a_{22}&\cdots &a_{2n}\\\vdots &\vdots &\ddots &\vdots \\a_{n1}&a_{n2}&\cdots &a_{nn}\end{pmatrix}}_{A}
\underbrace{\begin{pmatrix}
		x_1 \\ x_2 \\ \vdots \\ x_n
\end{pmatrix}}_{\bm{x}} = 
\underbrace{\begin{pmatrix}
		b_1 \\ b_2 \\ \vdots \\ b_n
\end{pmatrix}}_{\bm{b}},
	\label{eq:general-system-of-linear-equations}
\end{equation}
where $A \in \mathbb{C}^{n \times n}$ is the coefficient matrix, $\bm x \in \mathbb{C}^n$ the vector of unknowns and $\bm b \in \mathbb{C}^n$ the right-hand side.
At this point, we do not specify whether $A$ is represented as a dense/sparse matrix or a stencil as long as all the operations employed within subsequent steps are well-defined for a mathematical object of this type.
We can now rewrite Equation~\eqref{eq:general-system-of-linear-equations} to obtain
\begin{equation}
	\bm{x} = \bm{x} + \bm{b} - A \bm{x}.
	\label{eq:general-fixed-point}
\end{equation}
Which can be considered a fixed point of the form
\begin{equation}
	\bm x = f(\bm x).
\end{equation} 
Replacing $\bm x$ by $\bm{x}^{(k+1)}$ in the left and by $\bm{x}^{(k)}$ in the right part of Equation~\eqref{eq:general-fixed-point} yields the fixed-point iteration
\begin{equation}
	\bm{x}^{(k+1)} = \bm{x}^{(k)} + \bm b - A \bm{x}^{(k)}.
	\label{eq:richardson-iteration}
\end{equation}
Equation~\eqref{eq:richardson-iteration} is usually called the Richardson iteration and is the most basic form of an iterative method for solving a system of linear equations.
Here, the term $\bm{r}^{(k)} = \bm{b} - A \bm{x}^{(k)}$ represents the \emph{residual} or \emph{defect} in iteration $k$ of the method.
Next, we consider 
\begin{equation}
	M^{-1} A \bm{x} = M^{-1} \bm{b},
	\label{eq:general-preconditioned-system-of-linear-equations}
\end{equation}
which represents a modified version of Equation~\eqref{eq:general-system-of-linear-equations}, that is obtained by multiplying each side of the equation by the inverse of a matrix $M$.
By the rules of linear algebra, Equation~\eqref{eq:general-preconditioned-system-of-linear-equations} is equivalent to the original system and, therefore, has the same solution.
However, it can also be considered as a left-preconditioned version of Equation~\eqref{eq:general-system-of-linear-equations}, with $M$ as a preconditioner.
Considering again the fixed point of Equation~\eqref{eq:general-preconditioned-system-of-linear-equations} yields the iteration
\begin{equation}
	\bm{x}^{(k+1)} = \bm{x}^{(k)} + M^{-1}(\bm b - A \bm{x}^{(k)}),
	\label{eq:general-stationary-iterative-method}
\end{equation}
which represents the general form of a stationary iterative method. 
For instance, if we replace $M$ with the unit matrix $I$, we obtain the Richardson iteration.
Furthermore, setting $M = A$ allows us to compute the solution of the system in a single step since then
\begin{equation}
	\bm{x}^{(k+1)} = \bm{x}^{(k)} + A^{-1}(\bm b - A \bm{x}^{(k)}),
	\label{eq:one-step-iteration}
\end{equation}
which leads to $\bm{x}^{(k+1)} = A^{-1}\bm b$.
The result of this iteration $\bm{x}^{(k+1)}$ is thus independent of the choice of $\bm{x}^{(k)}$. 
If we insert $\bm{x} = A^{-1}\bm b$ into Equation~\eqref{eq:general-system-of-linear-equations}, it becomes apparent that this term always represents the correct solution of the respective system of linear equations.
Since the computation of the inverse of a general matrix $A$ is more expensive and numerically unstable than solving the system directly, the application of Equation~\eqref{eq:one-step-iteration} is, of course, impractical.
However, this derivation already provides us with an intuition about the choice of $M$.
In fact, the closer $M^{-1}$ is to the actual inverse of $A$, the faster we can expect the convergence of the respective stationary iterative method to be.
In contrast, the choice of a matrix $M$ that is easy to invert, with the Richardson iteration representing an extreme case with $I^{-1} = I$, leads to an iterative method that is easy to compute but which might suffer from slow convergence.
Before we introduce some basic iterative methods within the framework of Equation~\eqref{eq:general-stationary-iterative-method}, note that it can be reformulated as
\begin{equation}
	M (\bm{x}^{(k+1)} - \bm{x}^{(k)}) = \bm{b} - A \bm{x}^{(k)}. 
\end{equation}
By then defining $\bm{c}^{(k+1)}$ as the correction term
\begin{equation}
	\bm{c}^{(k+1)} = \bm{x}^{(k+1)} - \bm{x}^{(k)},
\end{equation}
we obtain a new system of linear equations
\begin{equation}
	M \bm{c}^{(k+1)} = \bm{b} - A \bm{x}^{(k)}. 
	\label{eq:general-stationary-iterative-method-system-formulation}
\end{equation}
The solution of this system $\bm{c}^{(k+1)}$ provides us with the approximate solution in step $i+1$ through the relation
\begin{equation}
	\bm{x}^{(k+1)} =  \bm{x}^{(k)} + \bm{c}^{(k+1)}.
\end{equation}
Therefore, instead of computing the inverse of $M$ in each iteration, we can solve the system of linear equations represented by Equation~\eqref{eq:general-stationary-iterative-method-system-formulation}.

Next, to derive a suitable choice for $M$ that leads to a faster convergence than Richardson's iteration, we define the splitting
\begin{equation}
	A = D - L - U,
\end{equation} 
where $D$ is the lower diagonal, $-L$ the lower triangular, and $-U$ the upper triangular part of $A$, respectively.
Both the Jacobi and Gauss-Seidel method are defined based on this splitting.
First of all, setting $M = D$, yields the Jacobi method
\begin{equation}
	\bm{x}^{(k+1)} = \bm{x}^{(k)} + D^{-1}(\bm b - A \bm{x}^{(k)}).
	\label{eq:jacobi-method}
\end{equation}
Since $D$ is a diagonal matrix, we can easily compute its inverse by inverting all of its diagonal entries.
Therefore, each iteration of the Jacobi method consists of a simple matrix-vector multiplication of the residual $\bm{r}^{(k)} = \bm{b} - A \bm{x}^{(k)}$ by the diagonal matrix $D^{-1}$.
To define the Jacobi method as a series of stencil operations, according to the derivation presented in Section~\ref{subsec:stencil-codes}, we need to extract the diagonal of a given stencil $S$.
Since a stencil entry with an offset of zero always refers to the current grid point, its value is equal to the corresponding diagonal entry of the matrix $A$.
We can therefore define the diagonal of $S$ as
\begin{equation}
	\text{diag}(S) = \begin{cases}
		\{(\bm{0}, b) \} & \text{if} \; (\bm 0, b) \in S \\
		\{(\bm{0}, 0) \} & \text{else}.
	\end{cases}
	\label{eq:stencil-diag}
\end{equation}
Based on this definition, we can also derive the inverse diagonal of a stencil as
\begin{equation}
	\text{diag-inv}(S) = \begin{cases}
		\{(\bm{0}, \frac{1}{b}) \} & \text{if} \; (\bm 0, b) \in S \\
		\{(\bm{0}, 0) \} & \text{else}.
	\end{cases}
	\label{eq:stencil-diag-inv}
\end{equation}
Because the Jacobi method, as defined in Equation~\eqref{eq:jacobi-method}, does often suffer from slow convergence, it is common to introduce a \emph{relaxation factor} or \emph{weight} $\omega$, which leads to the so-called weighted Jacobi method
\begin{equation}
	\bm{x}^{(k+1)} = \bm{x}^{(k)} + \omega D^{-1}(\bm b - A \bm{x}^{(k)}),
	\label{eq:weighted-jacobi-method}
\end{equation}
where $\omega$ is chosen from the interval $\left(0, 2\right)$.
Here, employing a value smaller than one is called \emph{underrelaxation}, while for a value larger than one, the term \emph{overrelaxation} is used.
Note that $\omega = 1$ leads to the original Jacobi method without any relaxation.

As a second variant of Equation~\eqref{eq:general-stationary-iterative-method}, we define the Gauss-Seidel method with $M = D - U$ which results in an iteration of the form
\begin{equation}
	(D - L) \bm{c}^{(k+1)} = \bm{b} - A \bm{x}^{(k)}, 
	\label{eq:gauss-seidel-method}
\end{equation}
where $\bm{c}^{(k+1)}$ again is the correction term $\bm{c}^{(k+1)} = \bm{x}^{(k+1)} - \bm{x}^{(k)}$.
Note that since $M = D - L$ does not represent a diagonal matrix, we can not easily compute its inverse but instead need to solve Equation~\eqref{eq:gauss-seidel-method} in every iteration of the method.
However, since $D - L$ is a lower triangular matrix, the solution of this system can be computed with a single step of back-substitution~\cite{saad2003iterative}. 
To express the Gauss-Seidel method as a sequence of stencil operations, we now additionally need to define a function that extracts the lower triangle of a given stencil.
As each diagonal entry corresponds to an offset of zero in each dimension, we can obtain the lower triangle by including only those entries with an offset lower than zero.
The resulting operation is defined in Equation~\eqref{eq:stencil-lower}.
\begin{equation}
	\text{lower}(S) = \begin{cases}
		\{(\bm{a}, b) \} \cup \text{lower}(\tilde{S}) & \text{if} \; \exists\, \bm a, b \; \text{with} \; \bm a < \bm 0 \\ & \text{and} \; S = (\bm a, b) \cup \tilde{S} \\
		\text{lower}(\tilde{S}) & \text{if} \; \exists\, \bm a, b \; \text{with} \; \bm a \geq \bm 0 \\ & \text{and} \; S = (\bm a, b) \cup \tilde{S} \\
		\emptyset & \text{else}
	\end{cases}
	\label{eq:stencil-lower}
\end{equation}

Similarly, we can also extract the upper triangle of a given stencil, which is formulated in Equation~\eqref{eq:stencil-upper}.
\begin{equation}
	\text{upper}(S) = \begin{cases}
		\{(\bm{a}, b) \} \cup \text{upper}(\tilde{S}) & \text{if} \; \exists\, \bm a, b \; \text{with} \; \bm a > \bm 0 \\ & \text{and} \; S = (\bm a, b) \cup \tilde{S} \\
		\text{upper}(\tilde{S}) & \text{if} \; \exists\, \bm a, b \; \text{with} \; \bm a \leq \bm 0 \\ & \text{and} \; S = (\bm a, b) \cup \tilde{S} \\
		\emptyset & \text{else}
	\end{cases}
	\label{eq:stencil-upper}
\end{equation}
While there are many more iterative methods that can fit into the framework of Equation~\eqref{eq:general-stationary-iterative-method}, the goal of this section is to introduce only those concepts necessary for a basic understanding of their functioning.
Therefore, we postpone the treatment of other and more advanced variants of the methods presented here to later chapters of this thesis.

\subsection{Convergence}
In contrast to direct methods, which compute the solution of a given system of linear equations in a fixed number of computational steps, iterative methods compute a series of approximations.
This series is said to converge when the difference between the actual solution and subsequent approximations approaches zero.
To quantify this behavior, we introduce a number of metrics used for evaluating the quality of a series of approximations computed with a specific iterative method.
First of all, assuming $\bm{x}^{(k)}$ is the $k$th approximation and $\bm{x}^{*}$ the correct solution of the system, we can define the error in iteration $k$ as
\begin{equation}
	\bm{e}^{(k)} = \bm{x}^{(k)} - \bm{x}^{*}.
\end{equation}
While the error gives us an immediate way to quantify the accuracy of an approximation, we usually do not know the correct solution of a given system of linear equations.
As a remedy, we instead consider the residual
\begin{equation}
	\bm{r}^{(k)} = \bm{b} - A \bm{x}^{(k)},
\end{equation}
which can always be computed irrespective of whether we know $\bm{x}^{*}$.
Note that for an error of zero, the residual vanishes as well.
Remember that the iterative methods introduced in the last section can all be considered as a fixed-point iteration of the form
\begin{equation*}
	\bm{x}^{(k+1)} = \bm{x}^{(k)} + M^{-1}(\bm b - A \bm{x}^{(k)}).
\end{equation*}
Assuming that $\bm{e}^{(k)} = \bm{0}$, we thus have $\bm{x}^{(k)} = \bm{x}^{*}$ and the above equation is reduced to
\begin{equation*}
	\bm{x}^{(k+1)} = \bm{x}^{*}.
\end{equation*}
The solution $\bm{x}^{*}$ of the system $A \bm{x} = \bm{b}$ is, therefore, a fixed point of each stationary iterative method.
However, note that this represents a necessary and not sufficient condition. 
Consider an arbitrary fixed-point $\bm{x}$ of Equation~\eqref{eq:general-stationary-iterative-method}
\begin{equation}
	\bm{x} = \bm{x} + M^{-1}(\bm b - A \bm{x}).
\end{equation}
Transforming this equation again into a system of linear equations yields
\begin{equation*}
	M^{-1} A \bm{x} = M^{-1}\bm b,
\end{equation*}
which means that each fixed-point $\bm{x}$ represents a solution of the preconditioned linear system~\eqref{eq:general-preconditioned-system-of-linear-equations} and hence also of the original system~\eqref{eq:general-system-of-linear-equations}.
If we assume that $A$ is a square, nonsingular matrix, the solution of each system of linear equations with $A$ as its coefficient matrix is unique. 
Therefore, it must be true that $\bm{x} = \bm{x}^{*}$, which means that the computed fixed point is equal to the correct solution of the linear system the method aims to solve.
However, the question remains to be answered under which conditions a sequence of the form of Equation~\eqref{eq:general-stationary-iterative-method} converges to a fixed point and how many iterations must be performed to achieve this goal.
For this purpose, we rewrite Equation~\eqref{eq:general-stationary-iterative-method} again to obtain
\begin{equation}
	\bm{x}^{(k+1)} = (I - M^{-1} A) \bm{x}^{(k)} + M^{-1}\bm b.
\end{equation}
Within this equation we can now set $\bm{x}^{(i)} = \bm{x}^{(0)}$ which yields
\begin{equation*}
	\bm{x}^{(1)} = (I - M^{-1} A) \bm{x}^{(0)} + M^{-1}\bm b,
\end{equation*}
and expand this sequence to the two-iteration series
\begin{equation*}
	\bm{x}^{(2)} = (I - M^{-1} A)^2 \bm{x}^{(0)} + (2I - M^{-1} A)M^{-1} \bm{b}.
\end{equation*}
By continuing this recursively for an arbitrary number of $k$ times, we see that the resulting equation will always be of the form
\begin{equation}
	\bm{x}^{(k)} = (I - M^{-1} A)^k \bm{x}^{(0)} + N^{(k)}\bm{b},
	\label{eq:stationary-iterative-method-series}
\end{equation}
which means that we can always separate the term $N^{(k)}\bm{b}$ from the rest of the equation.
Note that here the superscript $G^k$ denotes the $k$th power of $G$, while $N^{(k)}$ means that the matrix has been obtained through $k$ recursive substitutions of the respective iterate $\bm{x}^{(k)}$.
Since we have derived Equation~\eqref{eq:stationary-iterative-method-series} from Equation~\eqref{eq:general-stationary-iterative-method}, it still holds true that
\begin{equation}
	\bm{x}^{*} = (I - M^{-1} A)^k \bm{x}^{*} + N^{(k)}\bm{b},
\end{equation}
is a fixed point of this sequence and thus subtracting by $\bm{x}^{*}$ yields
\begin{equation}
	\bm{x}^{(k)} - \bm{x}^{*} = (\underbrace{I - M^{-1} A}_{G})^k (\bm{x}^{(0)} - \bm{x}^{*}).
	\label{eq:iteration-matrix-sequence}
\end{equation}
Here, $G = I - M^{-1} A$ is the so-called \emph{iteration matrix} of the considered stationary iterative method.
To reason about the convergence of this sequence, we introduce the spectral radius 
\begin{equation}
	\rho (A)=\max \limits_{1 \leq k \leq n} |\lambda _{k}|,
\end{equation}
where $A$ is a matrix of rank $n$ with the eigenvalues $\lambda_{1}, \dots, \lambda_{n}$.
Based on this definition, we can state the following theorem:
\begin{theorem}
$\lim \limits_{k \to  \infty} G^k = 0$ if and only if $\rho(G) < 1$~\cite{varga1962iterative,saad2003iterative}.
\label{theorem:general-convergence-result}
\end{theorem}

While Theorem~\ref{theorem:general-convergence-result} provides an answer whether the sequence~\eqref{eq:iteration-matrix-sequence} converges to zero, we still do not have any knowledge about the speed of this process.
For this purpose, we introduce the general \emph{convergence factor} $\rho$ of a sequence in the form of Equation~\eqref{eq:iteration-matrix-sequence} as
\begin{equation}
	\rho = \lim \limits_{k \to  \infty} \left( \max \limits_{x^{(0)} \in \mathbb{R}^{n}} \frac{\norm{\bm{x}^{(k)} - \bm{x}^{*}}}{\norm{\bm{x}^{(0)} - \bm{x}^{*}}} \right)^{\frac{1}{k}}.
\end{equation} 
Furthermore, the \emph{convergence rate} $\tau$ is defined as the inverse natural logarithm of the convergence factor
\begin{equation}
	\tau = -\ln \rho.
\end{equation}
We can now establish a link between the spectral radius of the iteration matrix and the convergence factor of the corresponding stationary iterative method by considering the following theorem:  
\begin{theorem}
	$\rho = \rho(G)$~\cite{varga1962iterative,saad2003iterative}.
	\label{theorem:convergence-factor}
\end{theorem}

As a result of Theorem~\ref{theorem:convergence-factor}, the spectral radius $\rho(G)$ gives a lower limit for the speed of convergence of every stationary iterative method with an iteration matrix $G = I - M^{-1} A$ that is independent of the choice of the initial vector $\bm{x}^{(0)}$.
Since all basic stationary iterative methods presented in this section are based on Equation~\eqref{eq:general-stationary-iterative-method}, these methods are easy to implement and analyze.
However, the application of these methods alone usually leads to slow convergence~\cite{briggs2000multigrid}.
While many different iterative methods for solving PDEs more efficiently have been proposed throughout the years~\cite{saad2003iterative}, the main focus of this thesis is on multigrid methods, which will be discussed in the following section.

%% file: contents/multigrid_methods.tex
\section{Multigrid Methods}\label{sec:multigrid-methods}
\begin{figure}[ht!]
	\centering
	\begin{subfigure}[t]{0.32\textwidth}
		\centering
		\includegraphics[width=\textwidth]{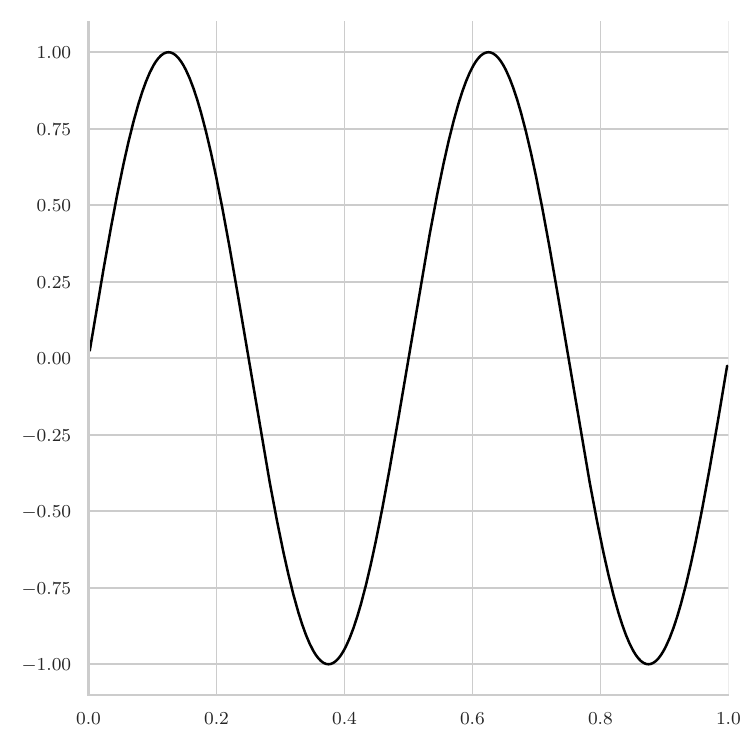}
		\includegraphics[width=\textwidth]{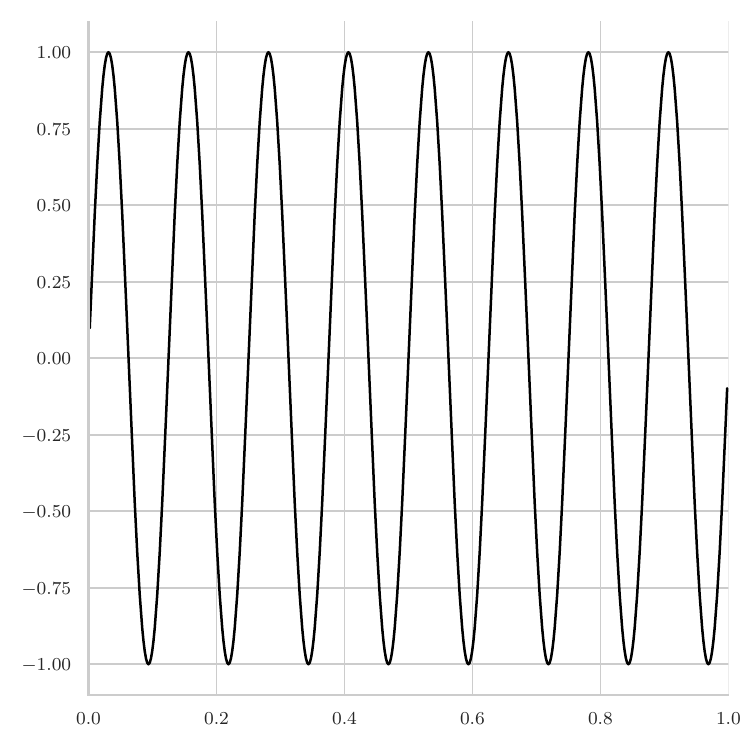}
		\includegraphics[width=\textwidth]{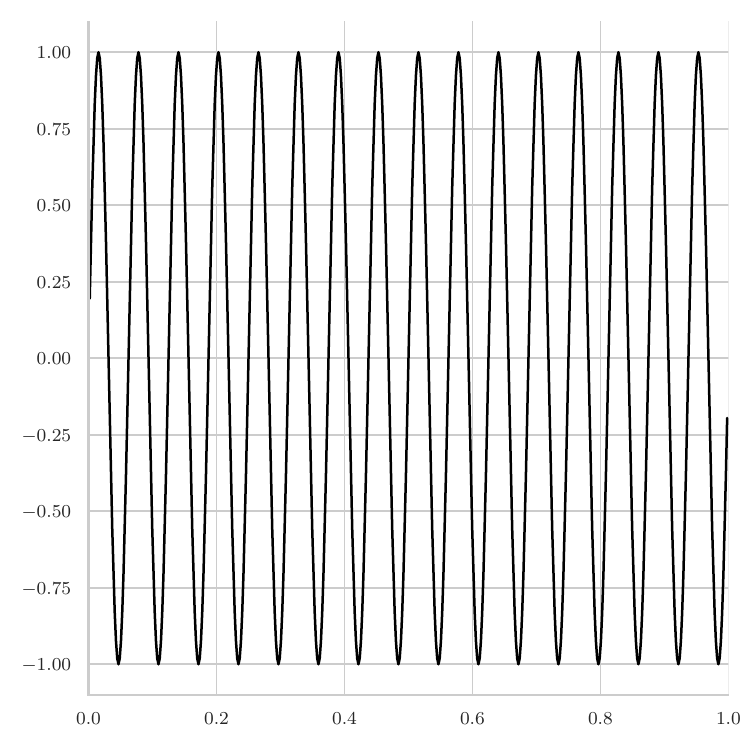}
		\caption{Initial error}
	\end{subfigure}
	\hfill
	\begin{subfigure}[t]{0.32\textwidth}
		\centering
		\includegraphics[width=\textwidth]{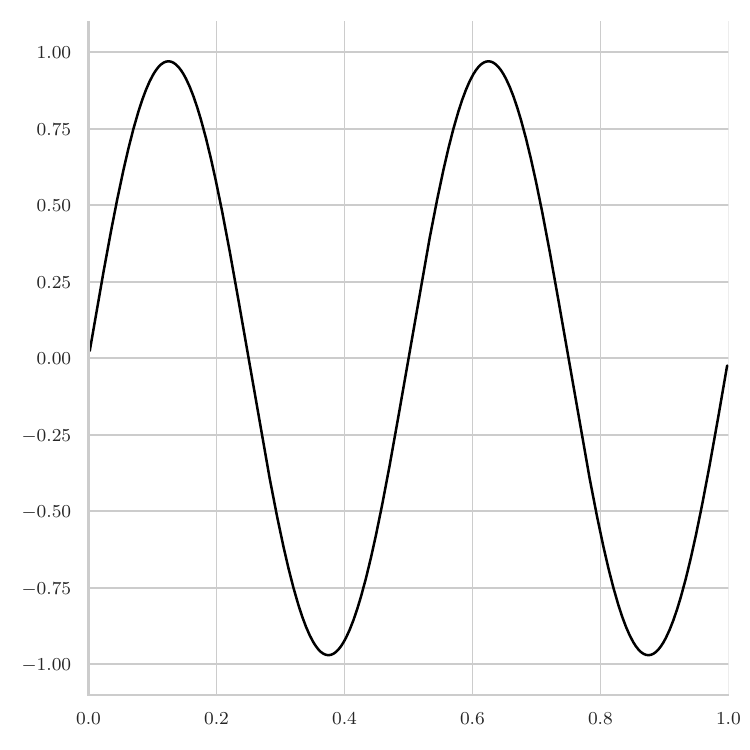}
		\includegraphics[width=\textwidth]{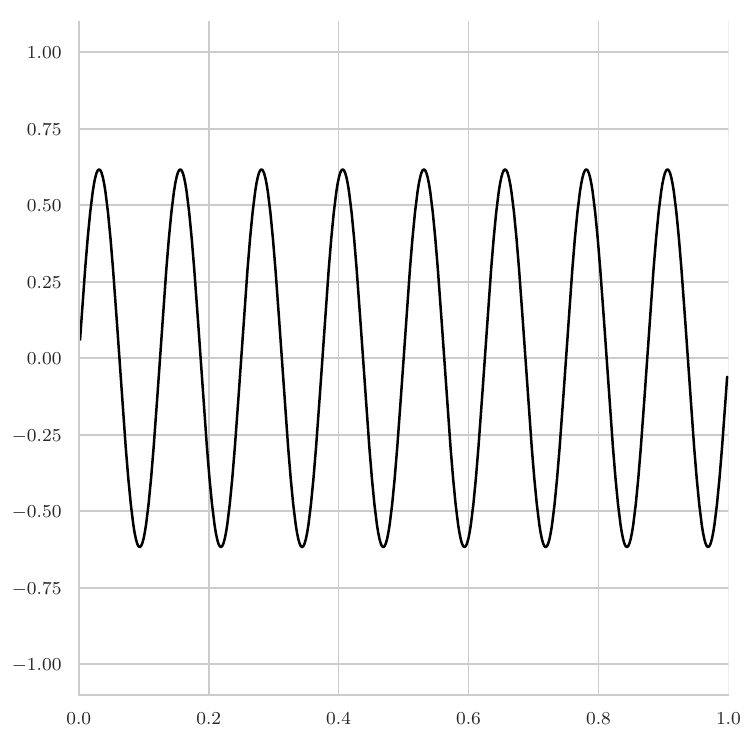}
		\includegraphics[width=\textwidth]{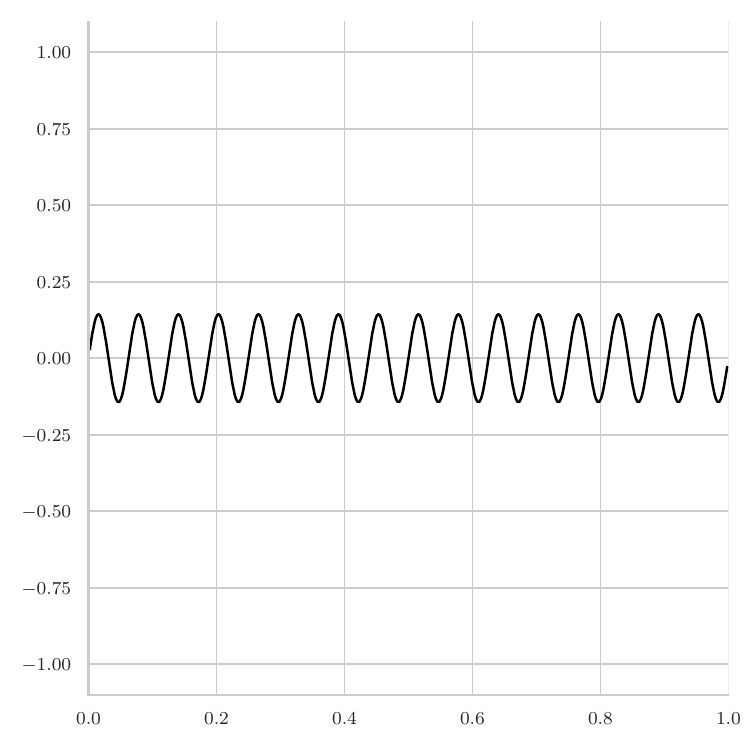}
	\caption{Error after applying Jacobi}
	\end{subfigure}
	\hfill
	\begin{subfigure}[t]{0.32\textwidth}
		\centering
		\includegraphics[width=\textwidth]{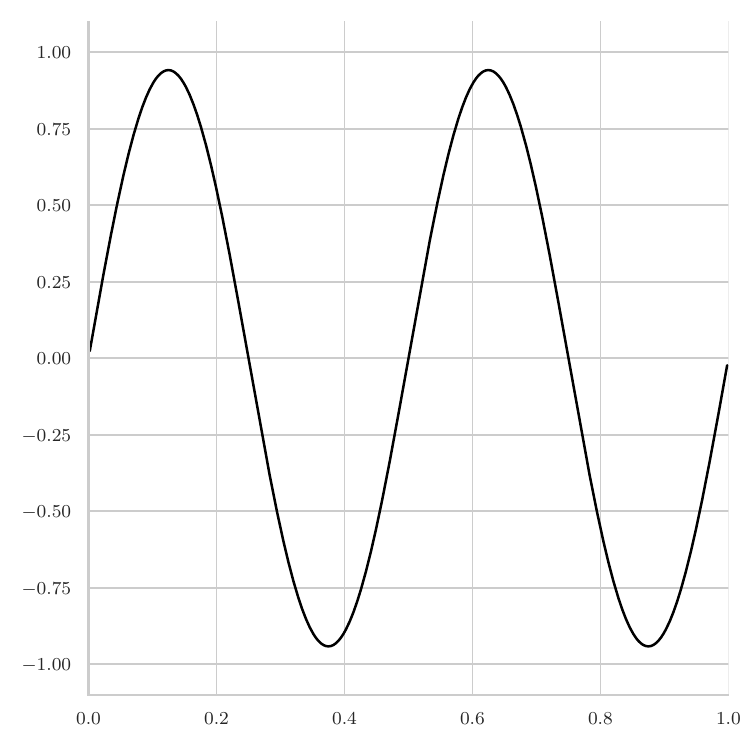}
		\includegraphics[width=\textwidth]{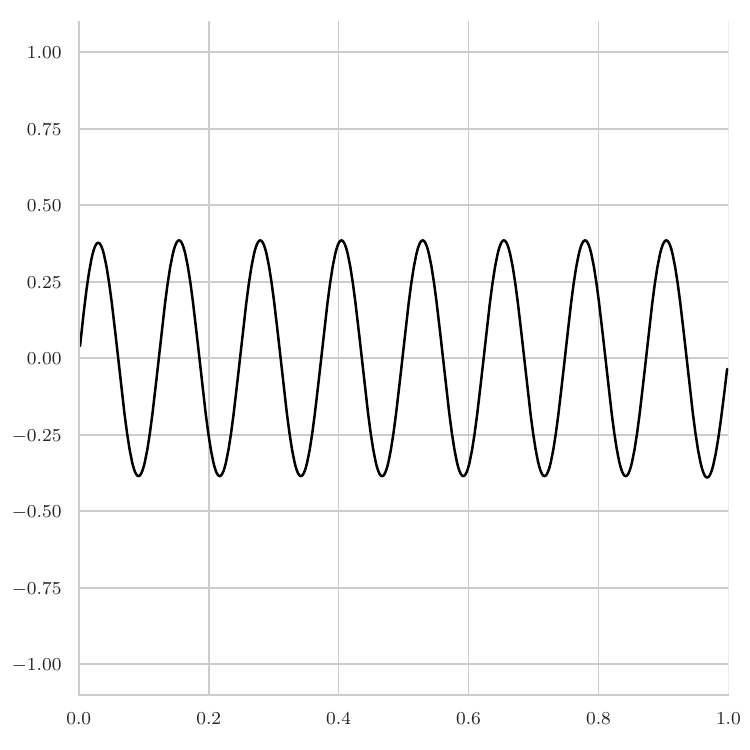}
		\includegraphics[width=\textwidth]{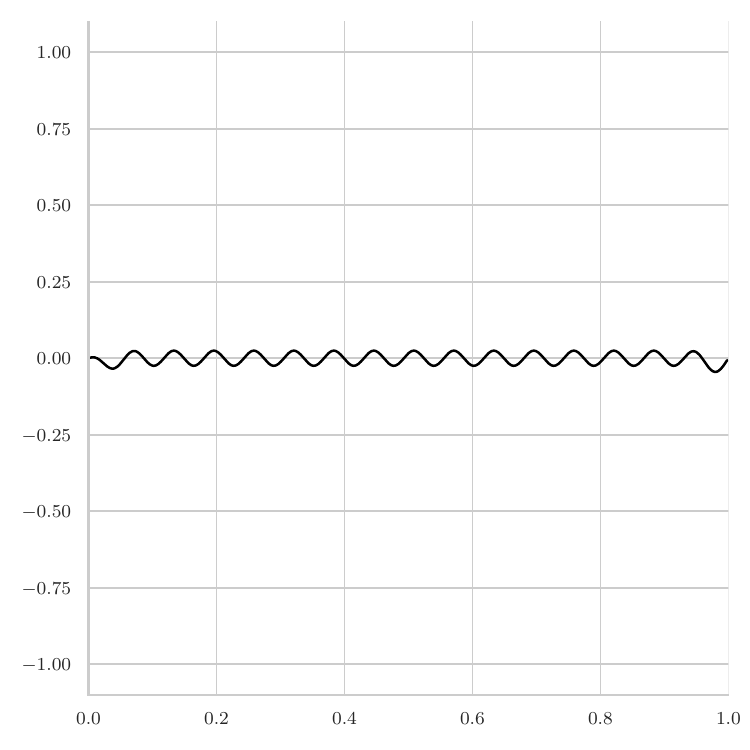}
	\caption{Error after applying GS}
	\end{subfigure}
	\caption[Different error components on a one-dimensional grid with step size $h = 2^{-9}$.]{Different error components on a one-dimensional grid with step size $h = 2^{-9}$ before and after applying 100 steps of the Jacobi or Gauss-Seidel (GS) method.}
	\label{fig:different-error-components}
\end{figure}

\begin{figure}
	\begin{subfigure}[t]{0.32\textwidth}
	\centering
		\includegraphics[width=\textwidth]{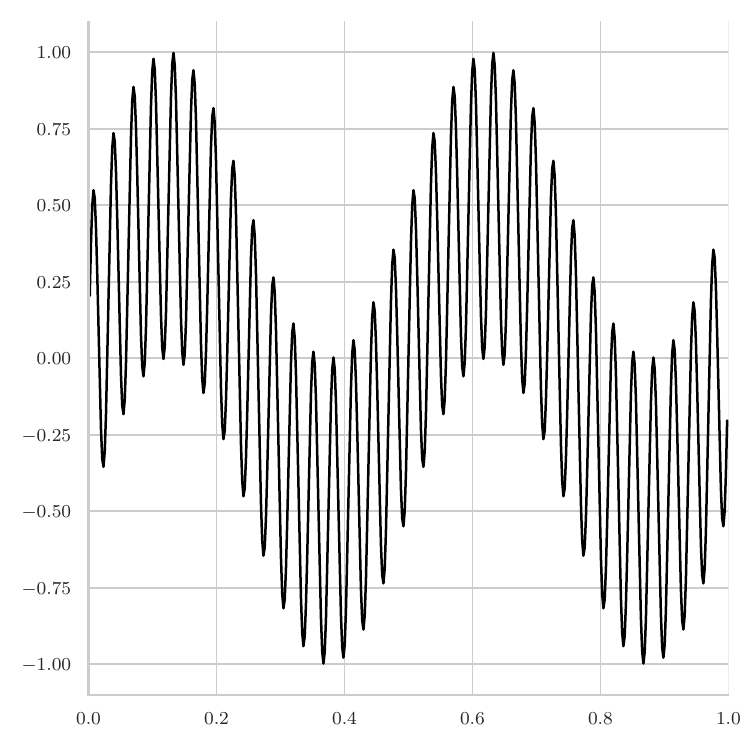}
		\caption{Initial error}
\end{subfigure}
\hfill
\begin{subfigure}[t]{0.32\textwidth}
	\centering
		\includegraphics[width=\textwidth]{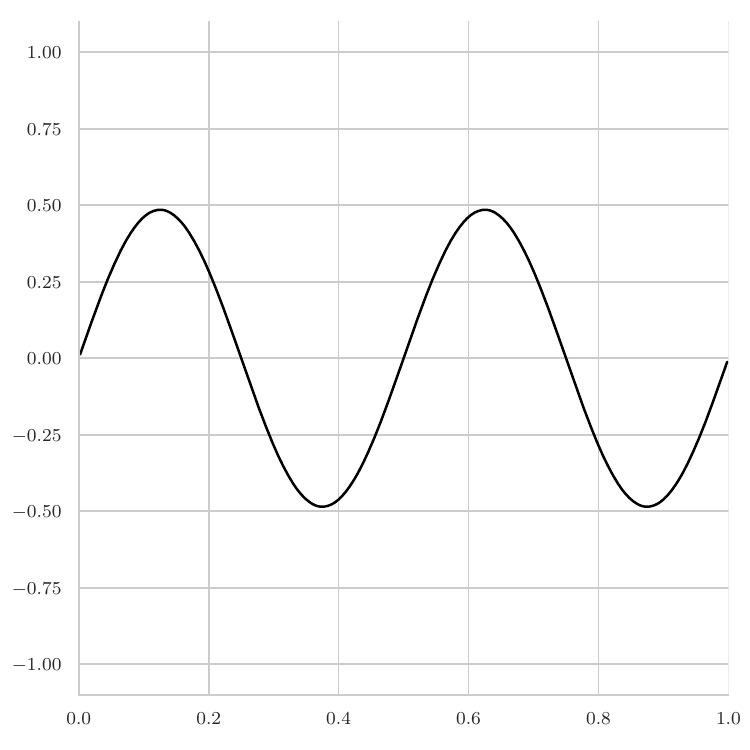}
		\caption{Error after applying Jacobi}
\end{subfigure}
\begin{subfigure}[t]{0.32\textwidth}
	\centering
	\includegraphics[width=\textwidth]{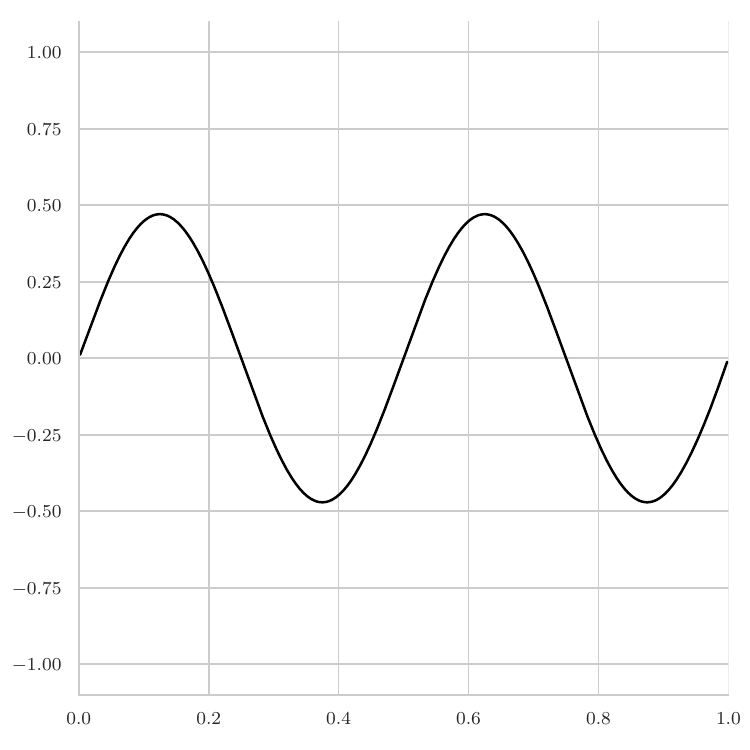}
	\caption{Error after applying GS}
\end{subfigure}
	\caption[Combination of two error components on a one-dimensional grid with step size $h = 2^{-9}$.]{Combination of two error components discretized on a one-dimensional grid with step size $h = 2^{-9}$ before and after applying 100 steps of the Jacobi or Gauss-Seidel (GS) method.}
\label{fig:combined-error}
\end{figure}

%Multigrid methods aim to accelerate the convergence of iterative methods by eliminating certain error components more efficiently on a coarser resolution of the same problem.
While basic iterative methods, such as Jacobi and Gauss-Seidel, are applicable to many PDE-based problems, in practice, their speed of convergence is often insufficient, which means that a large number of iterations is required until an acceptable approximation accuracy can be attained. 
The main reason for this behavior is that these methods are only efficient in the reduction of certain error components, while others remain mostly unaffected~\cite{briggs2000multigrid}.
This can be best understood by investigating their effect on oscillatory errors with different frequencies.
For this purpose, we consider the one-dimensional Laplace equation
\begin{equation}
		\begin{split}
			- \Delta u(x) & = 0 \quad \forall x \in (0, 1) \\
			u(0) & = u(1) = 0,
		\end{split}
		\label{eq:1D-laplace-model}
\end{equation}
which is discretized using the three-point stencil
\begin{equation}
	- \Delta_h = \frac{1}{h^2}\begin{bmatrix}
		-1 & 2 & -1
	\end{bmatrix}.
		\label{eq:1D-laplace-stencil}
\end{equation} 
Figure~\ref{fig:different-error-components} shows the impact of applying the Jacobi and Gauss-Seidel method to different periodic error components.
Here, the first column shows the initial error discretized on a grid with step size $h = 2^{-9}$ while the second and third include the remaining error after applying 100 Jacobi and Gauss-Seidel steps, respectively.
Note that the frequency of change increases from top to bottom, whereas the amplitude of the error is always the same.
As can be seen in the second and third column of the first row of Figure~\ref{fig:different-error-components}, the Jacobi and Gauss-Seidel methods do not yield a significant reduction of the low-frequency error components within 100 iterations.
In contrast, in the third row, which shows a highly-oscillating component, the application of 100 steps of the Jacobi method already reduces the initial error to less than one-fifth of its original value.
The same behavior can also be observed for the Gauss-Seidel method, whereby, compared to the Jacobi method, high-frequency error components are reduced even faster.
We can further illustrate the error reduction properties of basic iterative methods by investigating Figure~\ref{fig:combined-error}, which contains a combination of two error components with equal magnitude, one of them with low and the other one with high frequency.
Again, the first plot shows the initial error, while the second and third contain the reduced error after 100 iterations of Jacobi and Gauss-Seidel, respectively.
In accordance with our previous observations, the attained improvement achieved with both methods can be almost fully attributed to the reduction of the highly-oscillating component.
Because the remaining error is more smooth than initially, basic iterative methods are often called \emph{smoothers}, and their effectiveness is measured in terms of their capability to reduce the high-frequency components of a given error. 

\begin{figure}
	\begin{subfigure}[t]{0.32\textwidth}
	\centering
	\includegraphics[width=\textwidth]{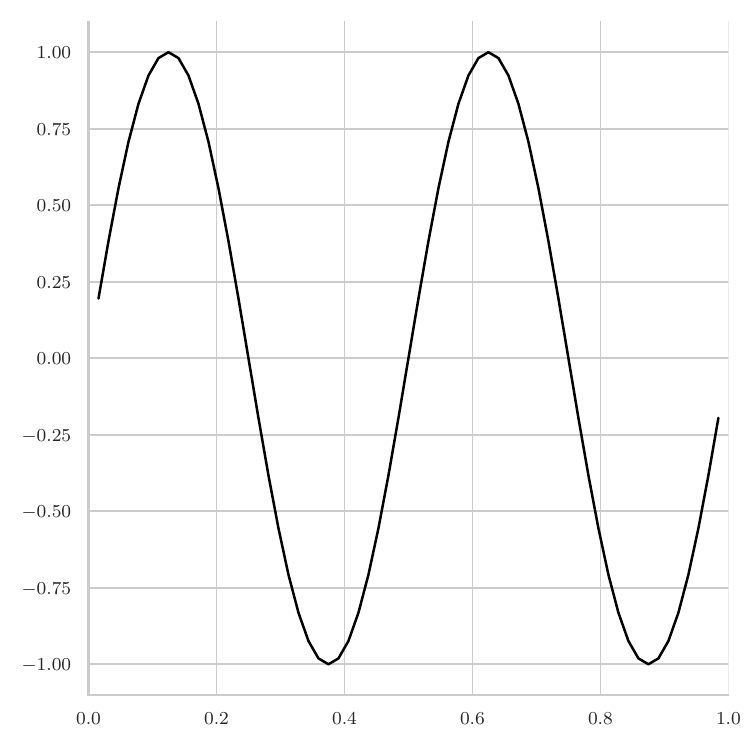}
	\caption{Initial error}
\end{subfigure}
\hfill
\begin{subfigure}[t]{0.32\textwidth}
	\centering
	\includegraphics[width=\textwidth]{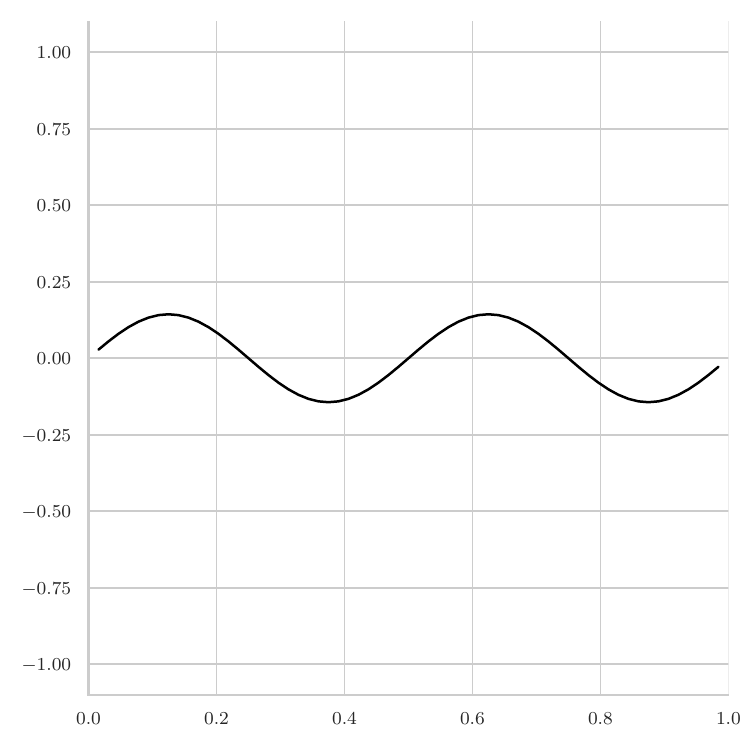}
	\caption{Error after applying Jacobi}
\end{subfigure}
	\hfill
	\begin{subfigure}[t]{0.32\textwidth}
		\centering
		\includegraphics[width=\textwidth]{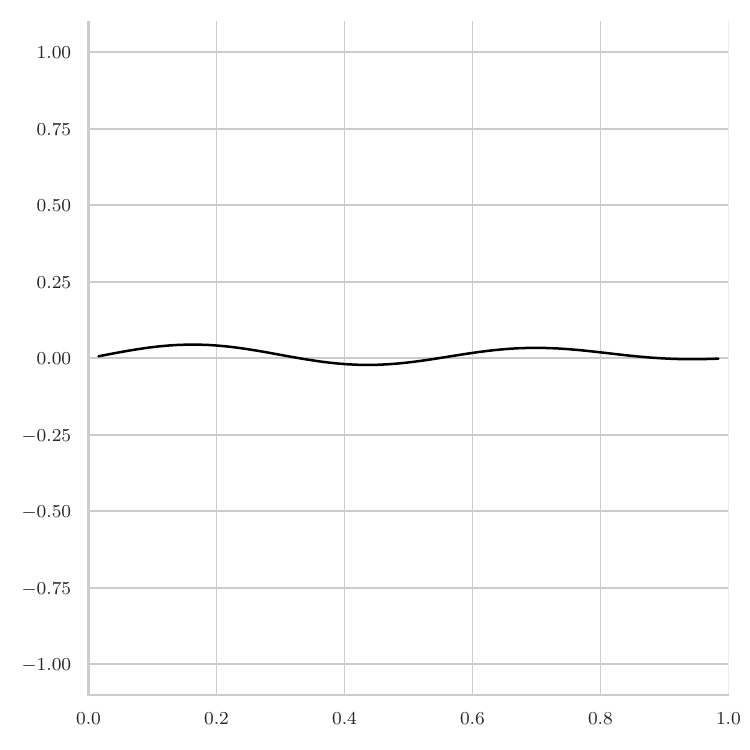}
		\caption{Error after applying GS}
	\end{subfigure}
	\caption[Low-frequency error component on a coarser one-dimensional grid with step size $h = 2^{-6}$.]{Low-frequency error component discretized on a coarser one-dimensional grid with step size $h = 2^{-6}$ before and after applying 100 steps of the Jacobi or Gauss-Seidel (GS) method.}
	\label{fig:low-frequency-error-component-coarse}
\end{figure}

Now observe what happens if we represent the same low-frequency error component shown in the first row of Figure~\ref{fig:different-error-components} on a grid with larger step size $h = 2^{-6}$ and thus a smaller number of grid points.
Because the number of (inner) grid points $n = 1/h - 1$ is inversely proportional to the step size, we call such a grid \emph{coarser}.
The resulting error reduction, again after 100 iterations of each method, is shown in Figure~\ref{fig:low-frequency-error-component-coarse}. 
As it can be seen, the amount of low-frequency error reduction is significantly higher for both methods than on the \emph{finer} grid with a step size of $h = 2^{-9}$.
While smoothers, such as Jacobi and Gauss-Seidel, are only effective in reducing the high-frequency components of a given error, we can overcome this limitation by representing the remaining error on a coarser grid.
Since, on this level, the remaining low-frequency components become more oscillatory, smoothing regains its effectiveness. 
\emph{Multigrid} methods extend this idea by recursively obtaining a coarser representation of the same problem, whose error can then be effectively reduced by employing only a few \emph{smoothing} iterations.
The result is then used to extinguish the remaining error on the next-higher level.
In the following, we introduce the basic components of multigrid methods, i.e., the smoothing, restriction, prolongation, and coarse-grid correction operations.
%Based on this definition, we then develop a formal language for the representation of multigrid solvers.
\subsection{Smoothing}
\label{sec:smoothing}
One of the central elements of multigrid methods is the utilization of a smoothing procedure for quickly reducing the oscillatory components of a given error.
We have already shown that the Jacobi and Gauss-Seidel method behave in such a way for the considered one-dimensional model problem.
To further improve the smoothing effectiveness of an iterative method, it is often beneficial to introduce an additional relaxation factor $\omega$, which yields the iteration 
\begin{equation}
	\bm{x}^{(k+1)} = \bm{x}^{(k)} + \omega M^{-1}(\bm b - A \bm{x}^{(k)}).
	\label{eq:general-weighted-stationary-iterative-method}
\end{equation}
Again, the weighted Jacobi or Gauss-Seidel method is obtained by replacing the matrix $M$ with the respective term.
\subsubsection{Red-Black Gauss-Seidel}
\label{sec:rb-gs}
While the Gauss-Seidel method often exhibits a superior smoothing property compared to the Jacobi method~\cite{briggs2000multigrid,trottenberg2000multigrid}, each iteration requires solving a lower triangular system of the form
\begin{equation*}
	(D - L) (\bm{x}^{(k+1)} - \bm{x}^{(k)}) = \bm{b} - A \bm{x}^{(k)},
\end{equation*}
with $D - L$ as the lower triangular part of the system matrix $A$.
Since $U = D - L - A$, we can rewrite this equation to obtain
\begin{equation}
	(D - L) \bm{x}^{(k+1)} = \bm{b} + U \bm{x}^{(k)},
\end{equation}  
which can then be solved using forward substitution by computing
\begin{equation}
	x_{i}^{(k+1)}={\frac {1}{a_{ii}}}\left(b_{i}-\sum _{j=1}^{i-1}a_{ij}x_{j}^{(k+1)}-\sum _{j=i+1}^{n}a_{ij}x_{j}^{(k)}\right),
	\label{eq:gauss-seidel-element-wise}
\end{equation}
where $x_{i}^{(k)}$ represents the $i$th element of the vector $\bm{x}^{(k)}$.
However, because $x_{i+1}^{(k)}$ depends on $x_{i}^{(k)}$ this computation can only be performed sequentially, which means that the individual components of $\bm{x}^{(k+1)}$ must be computed one after another. 
Modern compute architectures exhibit an ever-increasing degree of parallelism, and hence we must be able to perform all computations in parallel to fully utilize their capabilities.
Now consider the Jacobi method as defined by 
\begin{equation*}
	\bm{x}^{(k+1)} = \bm{x}^{(k)} + D^{-1}(\bm b - A \bm{x}^{(k)}).
\end{equation*}
Similar to the Gauss-Seidel method, we can rewrite this equation to obtain
\begin{equation}
	\bm{x}^{(k+1)} = D^{-1}(\bm b - (A - D)\bm{x}^{(k)}),
\end{equation}
which yields the following element-wise formulation of the Jacobi method:
\begin{equation}
x_{i}^{(k+1)}={\frac {1}{a_{ii}}}\left(b_{i}-\sum _{j\neq i}a_{ij}x_{j}^{(k)}\right).
	\label{eq:jacobi-element-wise}
\end{equation}
In contrast to Equation~\eqref{eq:gauss-seidel-element-wise} the computation of each subsequent element of the vector $\bm{x}^{(k+1)}$ does not depend on any previous one.
Consequently, the computation of each individual element of the new approximate solution $\bm{x}^{(k+1)}$ can be performed completely in parallel.

One possibility to enable a parallel Gauss-Seidel-like computation of the approximate solution $\bm{x}^{(k+1)}$ is to partition the grid points into multiple subsets.
The computation of each subset is then performed in a Jacobi-like fashion using the updated values from other subsets.
A common variant of this approach is the \emph{red-black} Gauss-Seidel (RB-GS) method.
Here, the grid points are assigned to two distinct subsets, where the first represents the red and the second one the black points.
We can define the RB-GS method in the following way:
\begin{equation}
	\begin{split}
		& \bm{x}^{(k+1/2)} = \bm{x}^{(k)} + P_R D^{-1} (\bm{b} - A \bm{x}^{(k)}) \\
		& \bm{x}^{(k+1)} = \bm{x}^{(k+1/2)} + P_B D^{-1} (\bm{b} - A \bm{x}^{(k+1/2)})
	\end{split}
\end{equation}
The entries of the matrices $P_R$ and $P_B$ are then defined as
\begin{equation}
	P_{R,ij} = \begin{cases}
	1 & \text{if} \; i = j \; \text{and} \; i,j \in R \\
	0 & \text{otherwise} 
	\end{cases}
\end{equation}
\begin{equation}
	P_{B,ij} = \begin{cases}
		1 & \text{if} \; i = j \; \text{and} \; i,j \in B \\
		0 & \text{otherwise},
	\end{cases}
\end{equation}
where $R$ and $B$ are the sets of grid indices that correspond to the red and black points, respectively.
For instance, an RB-GS method for the three-point stencil defined in Equation~\eqref{eq:1D-laplace-stencil} can be formulated as
\begin{equation}
\begin{split}
   	& x_{2i}^{(k+1/2)} = \frac {1}{2}\left(h^2 b_{2i} + x_{2i+1}^{(k)} + x_{2i-1}^{(k)}\right) \\
    & x_{2i+1}^{(k+1)} = \frac {1}{2}\left(h^2 b_{2i+1} + x_{2i+2}^{(k+1/2)} + x_{2i}^{(k+1/2)}\right),
\end{split}
\end{equation}
where each grid point with an even index belongs to the red and each one with an odd index to the black points.
Note that the update of each grid point is exclusively based on its neighbors, which have already been updated in the previous step of the method.
By always assigning neighboring points to different partitions, a similar effect can be achieved on any $d$-dimensional grid. 
For instance, a suitable partitioning for the discretized Laplace operator $\Delta_h$ is given by
\begin{equation}
		R = \{ \bm{i} : \bm{i} \in \mathbb{N}^d, \sum_{k=1}^d i_k \; \text{even} \}, \;
		B = \{ \bm{i} : \bm{i} \in \mathbb{N}^d, \sum_{k=1}^d i_k \; \text{odd} \}.
\end{equation}
In many cases, the resulting method has better smoothing properties than the Jacobi method without sacrificing much of its parallelism, as the computations on each partition can be performed concurrently~\cite{trottenberg2000multigrid}.
\subsubsection{Block Smoothing}
\label{subsec:block-smoothing}
So far, we have only discussed smoothers that operate in a pointwise manner.
For instance in the Jacobi method, Equation~\eqref{eq:jacobi-element-wise} is computed for each individual grid point.
The idea of block smoothing is to reorder the original system in such a way that the same operations can be defined on small subsets of grid points, which are usually chosen in the form of rectangular blocks of a particular size.
As a consequence, each scalar operation in the original pointwise method is replaced by a matrix or vector operation whose dimensionality corresponds to the chosen block size.

For example, we can rearrange our original system, as given by Equation~\eqref{eq:general-system-of-linear-equations}, in the following way:
\begin{equation}
\underbrace{
\begin{pmatrix}A_{11}&A_{12}&\cdots &A_{1m}\\A_{21}&A_{22}&\cdots &A_{2m}\\\vdots &\vdots &\ddots &\vdots \\A_{m1}&A_{m2}&\cdots &A_{mm}\end{pmatrix}}_{A}
\underbrace{
\begin{pmatrix}
\bm{x}_1 \\ \bm{x}_2 \\ \vdots \\ \bm{x_m} 
\end{pmatrix}}_{\bm{x}} =
\underbrace{
\begin{pmatrix}
	\bm{b}_1 \\ \bm{b}_2 \\ \vdots \\ \bm{b_m} 
\end{pmatrix}}_{\bm{b}}
\end{equation}
where $m = n / n_b$, if $n_b$ is the size of each block.
A block Jacobi method can then be defined as
\begin{equation}
	\bm{x}_{i}^{(k+1)}=A_{ii}^{-1}\left(\bm{b}_{i}-\sum _{j\neq i}A_{ij}\bm{x}_{j}^{(k)}\right).
	\label{eq:jacobi-block-wise}
\end{equation}
For instance, choosing a block size of two for the one-dimensional Laplace equation yields
\begin{equation}
	\begin{pmatrix}
		2 & -1 \\
		-1 & 2
	\end{pmatrix}
	\begin{pmatrix}
		x_{j}^{(k+1)} \\ x_{j+1}^{(k+1)} 
	\end{pmatrix}
= 	-  \begin{pmatrix}
	0 & 0 \\
	-1 & 0
\end{pmatrix} 	
\begin{pmatrix}
x_{j+2}^{(k)} \\ x_{j+3}^{(k)}
\end{pmatrix} -
\begin{pmatrix}
	0 & -1 \\
	0 & 0
\end{pmatrix} 	
\begin{pmatrix}
	x_{j-2}^{(k)} \\ x_{j-1}^{(k)} 
\end{pmatrix},
\end{equation}
where $j = n_b (i - 1) + 1$.
This method can be defined in a similar way using our previously introduced stencil notation
\begin{equation}
\begin{split}
	& \begin{pmatrix}
		\left[0 \quad 2 \quad -1 \right] \cdot x_{j}^{(k+1)} \\ \left[ -1 \quad 2 \quad 0 \right] \cdot x_{j+1}^{(k+1)} 
	\end{pmatrix}
	= \\ - & 
	\begin{pmatrix}
		\left[ 0 \right] \cdot x_{j+2}^{(k)} \\ \left[-1 \quad 0 \quad 0 \right] \cdot x_{j + 3}^{(k)}
	\end{pmatrix} -
	\begin{pmatrix}
		\left[0 \quad 0 \quad -1 \right] \cdot x_{j-2}^{(k)} \\ \left[ 0 \right] \cdot x_{j-1}^{(k)} 
	\end{pmatrix}.
\end{split}
\end{equation}
In both cases, we obtain a system of two linear equations
\begin{equation}
	\begin{pmatrix}
		2 x_{j}^{(k+1)} - x_{j+1}^{(k+1)} \\ 2 x_{j+1}^{(k+1)} - x_{j}^{(k+1)} 
	\end{pmatrix}
	=
	\begin{pmatrix}
		x_{j - 1}^{(k)} \\ x_{j + 2}^{(k)}
	\end{pmatrix}.
\end{equation}
which must be solved for each block, for instance, using a direct solver such as Gaussian elimination.
As in the pointwise Jacobi method, each block can be solved independently, and hence operations on different blocks can be performed in parallel. 
Considering the fact that a direct solver, in general, requires $\mathcal{O}(n_b^3)$ operations to solve each block, the overall computational cost of applying a block smoother can be estimated with $\mathcal{O}(n_b^3 \cdot n / n_b) = \mathcal{O}(n_b^2 n)$.
Note that choosing $n_b = n$ means that we treat the whole matrix $A$ as a single block and thus solve the original system using Gaussian elimination.
In contrast, the choice of $n_b = 1$ corresponds to the pointwise Jacobi method, which can be computed with a constant number of operations per grid point.
While we here only provide a brief introduction to block smoothers, it must be mentioned that it is also possible to define block variants for other smoothers, such as the Gauss-Seidel and red-black Gauss-Seidel method.
Finally, it is also possible to define overlapping block smoothers, which means that multiple blocks contain the same grid point as an unknown~\cite{trottenberg2000multigrid}.
\subsection{The Coarse-Grid Correction Scheme}
\begin{figure}[t]
	\begin{subfigure}[b]{0.32\textwidth}
		\centering
		\includegraphics[width=\textwidth]{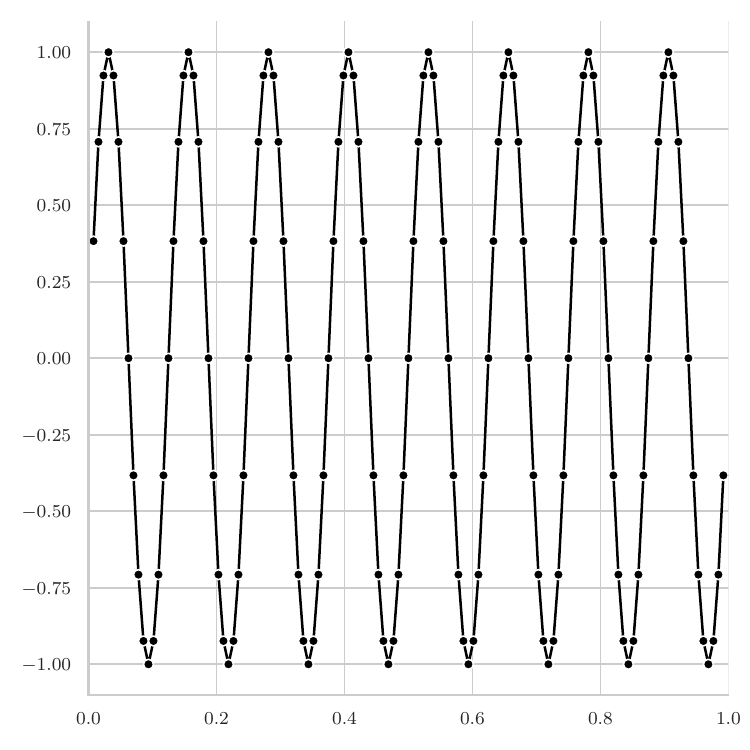}
	\end{subfigure}
	\hfill
	\begin{subfigure}[b]{0.32\textwidth}
		\centering
		\includegraphics[width=\textwidth]{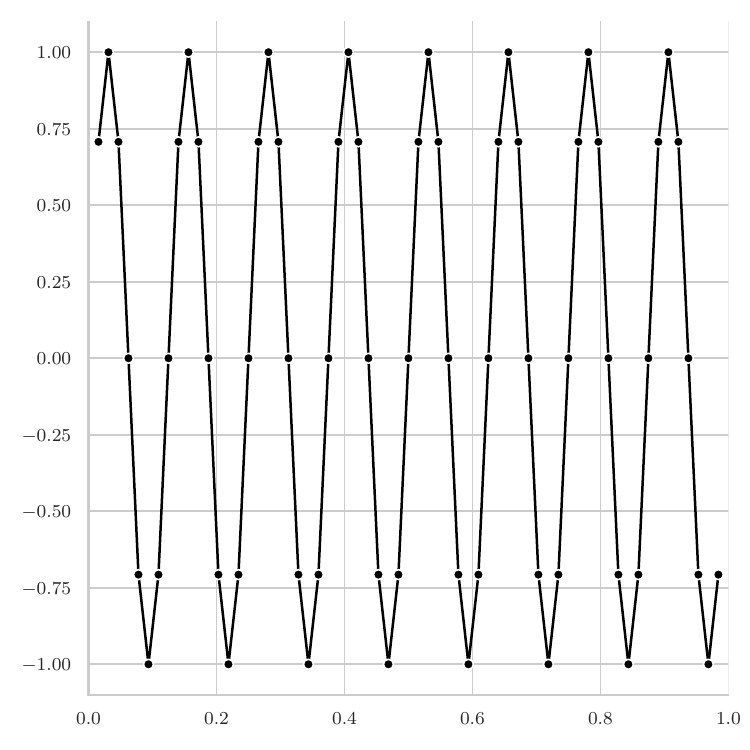}
	\end{subfigure}
	\hfill
	\begin{subfigure}[b]{0.32\textwidth}
		\centering
		\includegraphics[width=\textwidth]{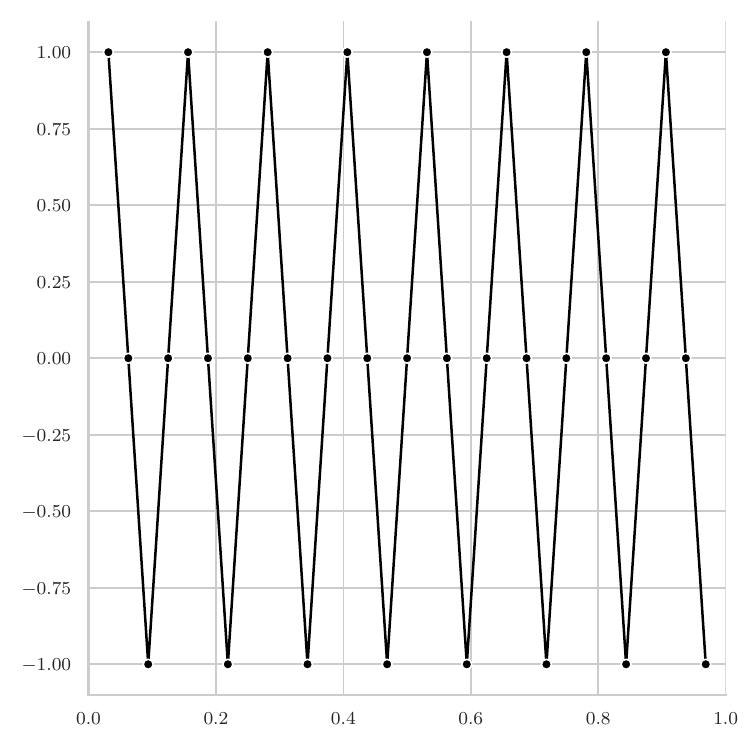}
	\end{subfigure}
		\begin{subfigure}[b]{0.32\textwidth}
		\centering
		\includegraphics[width=\textwidth]{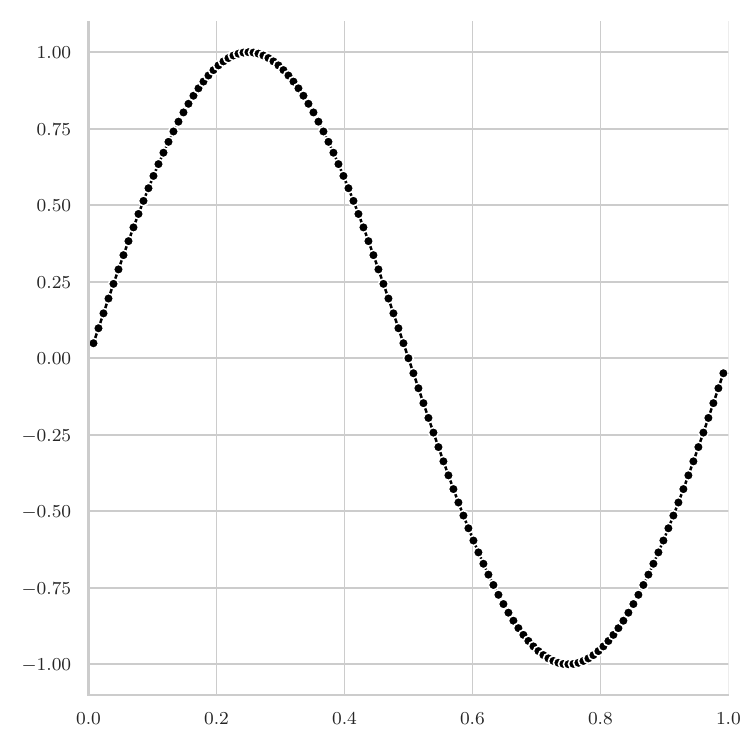}
		\caption{$h = 2^{-7}$}
	\end{subfigure}
	\hfill
	\begin{subfigure}[b]{0.32\textwidth}
		\centering
		\includegraphics[width=\textwidth]{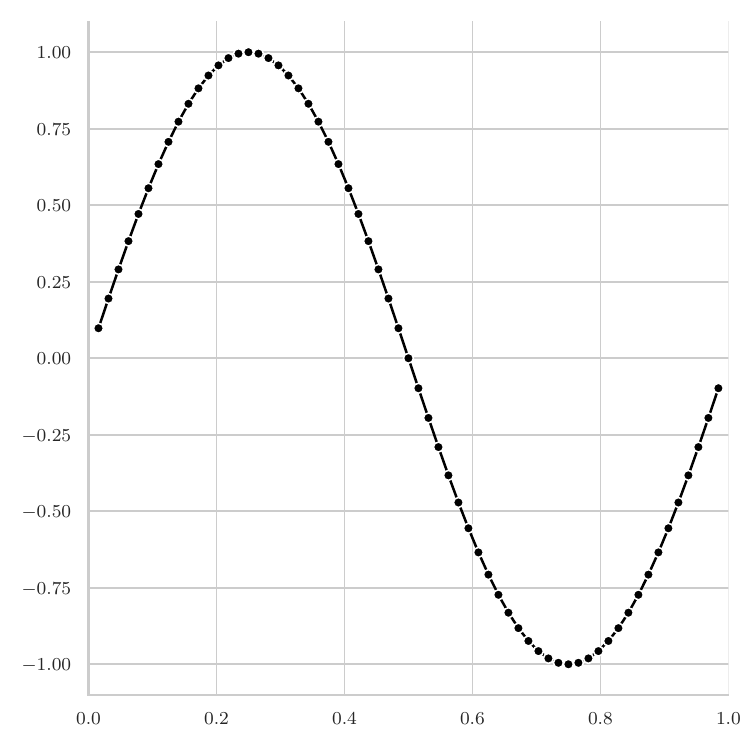}
		\caption{$h = 2^{-6}$}
	\end{subfigure}
	\hfill
	\begin{subfigure}[b]{0.32\textwidth}
		\centering
		\includegraphics[width=\textwidth]{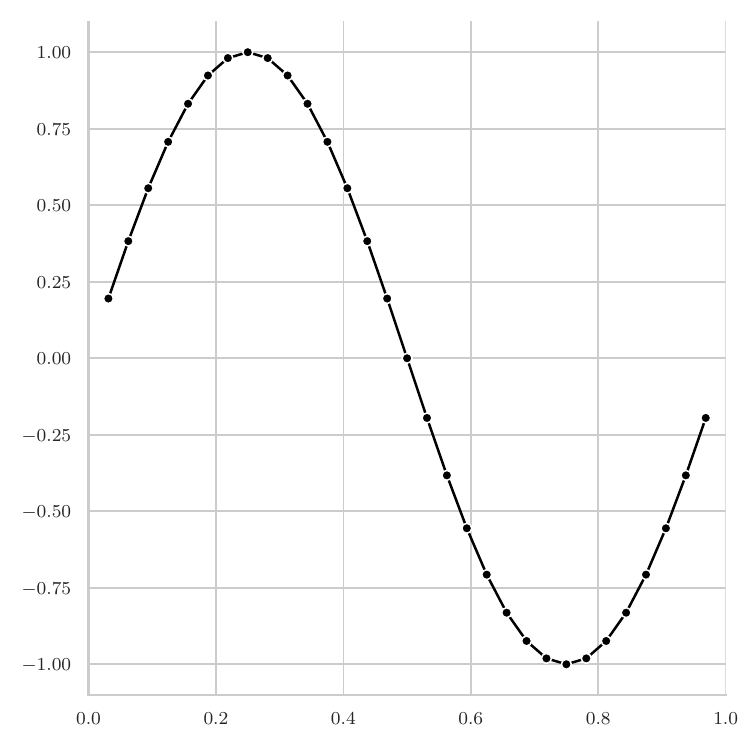}
		\caption{$h = 2^{-5}$}
	\end{subfigure}
	\caption[Oscillatory and smooth error components discretized on a hierarchy of grids with decreasing step size]{Oscillatory and smooth error components discretized on a hierarchy of grids with decreasing step size.}
	\label{fig:error-on-multiple-levels}
\end{figure}

The core idea behind multigrid methods is to reduce the oscillatory components of an error by computing an approximation of the same problem on a coarser grid.
As we have illustrated in Figure~\ref{fig:low-frequency-error-component-coarse}, these components can then be efficiently reduced using the same smoothing techniques already employed on the fine grid.
Before we can define this approach algorithmically, note that the system of linear equations
\begin{equation}
	A_h \bm{x}_h = \bm{b}_h
\end{equation}
can be reformulated as
\begin{equation}
	A_h \left(\bm{x}_h - \bm {x}^{(0)}_h\right) = \bm{b}_h - A_h \bm{x}^{(0)}_h,
\end{equation}
with an arbitrary-chosen initial guess $\bm{x}^{(0)}_h$.
Here, the subscript in $A_h$ and $\bm{x}_h$ indicates a discretization with step size $\bm{h}$.
Therefore, each entry of the vector $\bm{x}_h$ represents a grid function value $u_h(\bm{i} \circ \bm{h})$ at the respective position.
Furthermore, introducing the error $\bm{e}_h = \bm{x}_h - \bm {x}^{(0)}_h$, yields the equation
\begin{equation}
	A_h \bm{e}_h = \bm{b}_h - A_h \bm{x}^{(0)}_h.
	\label{eq:linear-system-error-equation}
\end{equation}
Based on the solution $\bm{e}_h$ of this system, which depends on $\bm{x}^{(0)}_h$, we obtain the solution of the original system by computing
\begin{equation}
	\bm{x}_h = \bm{x}^{(0)}_h + \bm{e}_h.
\end{equation}
Now, assume there exist two inter-grid operators, $I_h^{2h}$ and $I_{2h}^h$, that enable us to compute the approximation of a given vector $\bm{x}_h$ on a coarser and finer grid, respectively, such that
\begin{equation}
	\bm{x}_{2h} \approx I_h^{2h} \bm{x}_{h}, \;
	\bm{x}_{h} \approx I_{2h}^{h} \bm{x}_{2h}. 
\end{equation}
In general, this approximation can obviously never be exact.
However, we have already observed that if a certain error exclusively consists of smooth components, we can approximate them on a coarser grid without a significant loss of accuracy.
This is illustrated in Figure~\ref{fig:error-on-multiple-levels}, which shows the discretization of two error components with different frequencies on one-dimensional uniform grids with decreasing step size.
Here, the first error component, which possesses a higher frequency, can not be accurately represented on the coarsest grid with step size $h = 2^{-5}$.
In contrast, the slope of the second error component, which is relatively smooth, is still clearly visible on this grid. 
Assuming that an error $\bm{e}_{2h}$ on the coarse grid consists exclusively of smooth components and thus $\bm{e}_{h} \approx I_{2h}^{h} \bm{e}_{2h}$, we can define a coarse-grid correction 
\begin{equation}
	\bm{x}^{(k+1)}_h = \bm{x}^{(k)}_h + I_{2h}^h \bm{e}_{2h}.
\end{equation} 
Now the question remains to be answered how we can compute an approximation for the error $\bm{e}_{2h}$ on the coarse grid.
As we have pointed out above, Equation~\eqref{eq:linear-system-error-equation} is equivalent to the original system.
We can, therefore, make use of our previously defined inter-grid transfer operators to define the error equation on a coarser grid with step size $2\bm{h}$
\begin{equation}
	\underbrace{I_{h}^{2h} A_h I_{2h}^h}_{A_{2h}} \bm{e}_{2h} = I_{h}^{2h} \left(\bm{b}_h - A_h \bm{x}^{(0)}_h\right).
	\label{eq:coarse-grid-error-equation}
\end{equation}
Note that we have again made use of the fact that in case $\bm{e}_{2h}$ is smooth, $\bm{e}_{h} \approx I_{2h}^{h} \bm{e}_{2h}$ is a reasonably accurate approximation.
Note that in Equation~\eqref{eq:coarse-grid-error-equation}, the coarse operator $A_{2h}$ is directly obtained from the original operator $A_{2h}$, which is called \emph{Galerkin coarsening}.
However, in cases where the operator $A_h$ can be directly discretized on a coarser grid with step size $2\bm{h}$, it is often also possible to use the resulting operator $A_{2h}$ instead.
By bringing all these components together, we can formulate the two-level method shown in Algorithm~\ref{alg:two-grid-method}.
\begin{algorithm}
	\caption{Two-Grid Method}
	\label{alg:two-grid-method}
	\begin{algorithmic}
		\State{Smooth on $A_h \bm{x}_h = \bm{b}_h$ to obtain an approximation $\bm{x}_h$}
		\State{Compute the residual $\bm{r}_h = \bm{b}_h - A_h \bm{x}_h$}
		\State{\hskip1em Obtain $A_{2h}$ through Galerkin coarsening or rediscretization}
		\State{\hskip1em Restrict the residual $\bm{b}_{2h} = I_h^{2h} \bm{r}_h$}
		\State{\hskip1em Solve $A_{2h} \bm{x}_{2h} = \bm{b}_{2h}$ for $\bm{x}_{2h}$}
		\State{Perform the correction $\bm{x}_h = \bm{x}_h + I_{2h}^h \bm{x}_{2h}$}
		\State{Smooth again on $A_h \bm{x}_h = \bm{b}_h$ now using $\bm{x}_h$} as an initial guess
	\end{algorithmic}
\end{algorithm}
Note that while on the coarse grid, we are actually solving the error equation, it has been redefined as $A_{2h} \bm{x}_{2h} = \bm{b}_{2h}$.
Based on this notation, we could similarly define a two-level method starting from a grid with step size $2\bm{h}$.
However, to put this approach into practice, we need to choose an initial guess for the error equation on this level.
For this purpose, note that one step of smoothing applied to Equation~\eqref{eq:coarse-grid-error-equation} with an initial guess of zero corresponds to
\begin{equation}
	\bm{e}_{2h}^{(1)} = M_{2h}^{-1} I_{h}^{2h} \left(\bm{b}_h - A_h \bm{x}^{(0)}_h\right).
	\label{eq:initial-coarse-grid-relaxation}
\end{equation}
Since the error $\bm{e}_{2h}^{(k+1)}$ in step $k+1$ is defined as
\begin{equation*}
	\bm{e}_{2h}^{(k+1)} = \bm{x}_{2h}^{(k+1)} - \bm{x}_{2h}^{(k)},
\end{equation*}
and assuming that $\bm{x}_{h}^{(0)}$ is smooth, we can choose $\bm{x}_{2h}^{(0)} = I_{2h}^{h} \bm{x}_{h}^{(0)}$ and thus obtain the iteration
\begin{equation}
	\bm{x}_{2h}^{(1)} = I_{h}^{2h} \bm{x}_{h}^{(0)} + M_{2h}^{-1} ( I_{h}^{2h} \bm{b}_h - \underbrace{I_{h}^{2h} A_h I_{2h}^{h}}_{A_{2h}} I_{h}^{2h} \bm{x}_{h}^{(0)} ).
\end{equation}
Performing one smoothing step with an initial guess of zero on the coarse-grid error equation is thus similar to smoothing on the equation
\begin{equation}
	I_{h}^{2h} A_h I_{2h}^h \bm{x}_{2h} = I_{h}^{2h} \bm{b}_h,
\end{equation}
with an initial guess of $\bm{x}_{2h}^{(0)} = I_{h}^{2h} \bm{x}_{h}^{(0)}$.
Because smoothing aims to remove the remaining oscillatory components of the error that have been transferred from the fine grid, applying it in the form of Equation~\eqref{eq:initial-coarse-grid-relaxation} precisely serves this purpose.

\subsection{Restriction and Prolongation}
\label{subsec:restriction-and-prolongation}
%TODO fix stencil notation in this section
Before we can define an actual multigrid method, which recursively applies the techniques introduced in the last section to a hierarchy of discretizations consisting of more than two levels, we need to define suitable inter-grid operators, $I_{h}^{2h}$ and $I_{2h}^{h}$.
Here, the \emph{restriction} operator $I_{h}^{2h}$ is supposed to yield an accurate approximation of the current residual on a coarser grid, while the goal of the \emph{prolongation} operator $I_{2h}^{h}$ is to transfer the computed solution of the error equation back to a finer grid.
In both cases, we are, first and foremost, interested in preserving the low-frequency components, as all the others can already be quickly reduced by smoothing.
In general, the implementation of these operators depends on the chosen method of discretization. 
Since, as mentioned in Section~\ref{sec:discretization}, this work focuses on the discretizations of PDEs on regular grids, we do not consider inter-grid operators defined on other grid types, such as unstructured grids.
For a detailed treatment of these cases, the reader is referred to~\cite{trottenberg2000multigrid,ruge1987algebraic,stuben2001introduction}.
As a first step towards defining a suitable restriction operator, note that on a regular grid, the set of coarse-grid points is always contained in the set of fine-grid points.
Therefore, the easiest way to define such an operator is to simply carry over the values present at the respective points of the fine grid over to the coarse grid, which leads to the so-called \emph{injection} restriction operator.
While this approach might lead to a functioning multigrid method in case the residual is sufficiently smooth, it neglects the information contained within all fine grid points that do not coincide with one on the coarse grid.
In most cases, it is beneficial to incorporate this information into the coarse grid by computing a weighted average over all neighboring fine grid points.
This idea leads to the so-called \emph{full-weighting} and \emph{half-weighting} restriction operators.
Using our previously defined stencil notation, the one-dimensional full-weighting restriction operator is given by
\begin{equation}
	I_{h_x}^{2 h_x} = \frac{1}{4}\{((-1), 1), ((0), 2), ((1), 1)\}_{h_x}^{2h_x},
\end{equation} 
or equivalently using the matrix notation
\begin{equation}
	I_{h_x}^{2h_x} =  \frac{1}{4} \begin{bmatrix}
			1 & 2 & 1
		\end{bmatrix}_{h_x}^{2h_x}.
	\label{eq:full-weighting-restriction}
\end{equation} 
Since the application of this stencil yields a grid function of different dimensionality than the one to which it is applied, we additionally include the respective step sizes as a sub- and superscript.
If we treat Equation~\eqref{eq:full-weighting-restriction} as a row vector, we can define the two- and three-dimensional full-weighting restriction operator as a tensor product with the corresponding column vector, such that
\begin{equation}
	I^{2h_x, 2h_y}_{h_x, h_y} = \frac{1}{4} \begin{bmatrix}
		1 \\ 2 \\ 1
	\end{bmatrix}_{h_y}^{2h_y} \otimes \frac{1}{4} \begin{bmatrix}
		1 & 2 & 1
	\end{bmatrix}_{h_x}^{2h_x} =
\frac{1}{16} 
\begin{bmatrix}
	1 & 2 & 1 \\
	2 & 4 & 2 \\
	1 & 2 & 1
\end{bmatrix}^{2h_x, 2h_y}_{h_x, h_y}
\end{equation} 
\begin{equation}
\begin{split}
	& I^{2h_x, 2h_y, 2h_z}_{h_x, h_y, h_z} = \frac{1}{4} \begin{bmatrix}
		1 & 2 & 1
	\end{bmatrix}_{h_z}^{2h_z} \otimes 
	\frac{1}{16} 
	\begin{bmatrix}
		1 & 2 & 1 \\
		2 & 4 & 2 \\
		1 & 2 & 1
	\end{bmatrix}^{2h_x, 2h_y}_{h_x, h_y} \\
& = \frac{1}{64} \begin{bmatrix}
\begin{bmatrix}
	1 & 2 & 1 \\
	2 & 4 & 2 \\
	1 & 2 & 1
\end{bmatrix} &	\begin{bmatrix}
2 & 4 & 2 \\
4 & 8 & 4 \\
2 & 4 & 2
\end{bmatrix} &
\begin{bmatrix}
	1 & 2 & 1 \\
	2 & 4 & 2 \\
	1 & 2 & 1
\end{bmatrix}
\end{bmatrix}^{2h_x, 2h_y, 2h_z}_{h_x, h_y, h_z}
\end{split}
\end{equation}
A second restriction operator based on the idea of averaging neighboring fine grid points is the half-weighting restriction operator, whose two- and three-dimensional versions can be defined as
\begin{equation}
	I^{2h_x,2h_y}_{h_x, h_y} = \frac{1}{8}
	\begin{bmatrix}
		0 & 1 & 0 \\
		1 & 4 & 1 \\
		0 & 1 & 0
	\end{bmatrix}^{2h_x, 2h_y}_{h_x, h_y}
\end{equation} 
\begin{equation}
	I^{2h_x, 2h_y, 2h_z}_{h_x, h_y, h_z} = 
\frac{1}{12} \begin{bmatrix}
	\begin{bmatrix}
		0& 0 & 0 \\
		0 & 1 & 0 \\
		0& 0 & 0
	\end{bmatrix}
 &		\begin{bmatrix}
 	0 & 1 & 0 \\
 	1 & 6 & 1 \\
 	0 & 1 & 0
 \end{bmatrix} &
	\begin{bmatrix}
		0& 0 & 0 \\
		0 & 1 & 0 \\
		0& 0 & 0
	\end{bmatrix}
\end{bmatrix}^{2h_x, 2h_y, 2h_z}_{h_x, h_y, h_z}
\end{equation} 
% \begin{equation}
% \begin{split}
% 	& I^{2h_x, 2h_y, 2h_z}_{h_x, h_y, h_z} = 
% \frac{1}{4} \begin{bmatrix}
% 	1 & 2 & 1
% \end{bmatrix}_{h_z}^{2h_z}
% \otimes 
% \frac{1}{8}
% 	\begin{bmatrix}
% 	0 & 1 & 0 \\
% 	1 & 4 & 1 \\
% 	0 & 1 & 0 
% \end{bmatrix}^{2h_x,2h_y}_{h_x, h_y} \\
% & =
% \frac{1}{32} \begin{bmatrix}
% 	\begin{bmatrix}
% 		0& 1 & 0 \\
% 		1 & 4 & 1 \\
% 		0& 1 & 0
% 	\end{bmatrix}
%  &		\begin{bmatrix}
%  	0 & 2 & 0 \\
%  	2 & 6 & 2 \\
%  	0 & 2 & 0
%  \end{bmatrix} &
% 	\begin{bmatrix}
% 	0& 1 & 0 \\
% 	1 & 4 & 1 \\
% 	0& 1 & 0
% \end{bmatrix}
% \end{bmatrix}^{2h_x, 2h_y, 2h_z}_{h_x, h_y, h_z}
% \end{split}
% \end{equation} 
Note that in contrast to our original definition of the stencil application in Equation~\eqref{eq:stencil-application}, the restriction stencils presented here only have to be applied to each point on the fine grid that coincides with a coarse-grid point.
We can thus replace it with the following slightly adapted definition of stencil application
\begin{equation}
	\begin{split}
		& I_{h}^{2h} \cdot u_h(\bm{x}) = \sum_{k=1}^m b_k u_h({\bm x + \bm{a}_k} \circ \bm{h}) \quad 
		\text{with} \; \bm{x} \in G_{2h}, m \in \mathbb{N} \\ & (\bm{a}_k, b_k) \in I_{h}^{2h} \; \forall \, k \in \{ 1, 2, \dots, m \},
	\end{split}
	\label{eq:stencil-restriction-application}
\end{equation}
where $u_h(\bm{x})$ with $\bm{x} \in G_{2h} \supset G_h$ represents an arbitrary point on the fine grid for which there exists a unique coarse-grid point defined at the same spatial position within the computational domain.
Note that this is an immediate consequence of the fact that we have defined the set of coarse-grid points as a subset of the fine-grid points.

On the other hand, for prolongation, our goal is to define an operator that transfers an approximation of the error computed on a certain discretization level to a grid of higher resolution.
Therefore, we must be able to restore a larger number of grid points based on the given values on the coarse grid, which can be regarded as a distribution process.
The application of this operator to a given coarse-grid point $u_{2h}(\bm{x})$ with $\bm{x} \in G_{2h}$ can be defined as
\begin{equation}
	I_{2h}^{h} \cdot u_{2h}(\bm{x}) \rightarrow
	\begin{cases}
		& \forall \, (\bm{a}_k, b_k) \in I_{2h}^{h} \; \text{with} \; k \in \{ 1, 2, \dots, m \} \; \wedge \; \bm{x} \in G_{2h} : \\
		& u_{h}(\bm{x} + \bm{a}_k \circ \bm{h}) = u_{h}(\bm{x} + \bm{a}_k \circ \bm{h}) + b_k u_{2h}(\bm{x}), 
	\end{cases}
	\label{eq:stencil-prolongation application}
\end{equation}
where we assume that initially $\forall u_h(\bm{x}) \; \text{with} \; \bm{x} \in \mathcal G_h : u_h(\bm{x}) = 0$.
A common choice for one-dimensional prolongation is the linear interpolation operator~\cite{trottenberg2000multigrid}, which can be defined as
\begin{equation}
	I_{2h_x}^{h_x} =  \frac{1}{2} \begin{bmatrix}
		1 & 2 & 1
	\end{bmatrix}_{2h_x}^{h_x}.
	\label{eq:linear-interpolation}
\end{equation}
Similar to full-weighting restriction, we can derive two- and three-dimensional versions of this operator as tensor products of the form
\begin{equation}
	I_{2h_x, 2h_y}^{h_x, h_y} = \frac{1}{2} \begin{bmatrix}
		1 \\ 2 \\ 1
	\end{bmatrix}_{2h_y}^{h_y} \otimes \frac{1}{2} \begin{bmatrix}
		1 & 2 & 1
	\end{bmatrix}_{2h_x}^{h_x} =
	\frac{1}{4} 
	\begin{bmatrix}
		1 & 2 & 1 \\
		2 & 4 & 2 \\
		1 & 2 & 1
	\end{bmatrix}_{2h_x, 2h_y}^{h_x, h_y}
\end{equation}
\begin{equation}
	\begin{split}
		I_{2h_x, 2h_y, 2h_z}^{h_x, h_y, h_z} = & \frac{1}{2} \begin{bmatrix}
			1 & 2 & 1
		\end{bmatrix}_{2h_z}^{h_z} \otimes 
		\frac{1}{4} 
		\begin{bmatrix}
			1 & 2 & 1 \\
			2 & 4 & 2 \\
			1 & 2 & 1
		\end{bmatrix}_{2h_x,2h_y}^{h_x,h_y} \\
		= & \frac{1}{8} \begin{bmatrix}
			\begin{bmatrix}
				1 & 2 & 1 \\
				2 & 4 & 2 \\
				1 & 2 & 1
			\end{bmatrix}&	\begin{bmatrix}
				2 & 4 & 2 \\
				4 & 8 & 4 \\
				2 & 4 & 2
			\end{bmatrix} &
			\begin{bmatrix}
				1 & 2 & 1 \\
				2 & 4 & 2 \\
				1 & 2 & 1
			\end{bmatrix}
		\end{bmatrix}_{2h_x,2h_y,2h_z}^{h_x,h_y,h_z}
	\end{split}
\end{equation}
Finally, we want to emphasize that while the prolongation and restriction operators presented here represent a common choice for regular grids with uniform step sizes, for instance, the discretization of PDEs with varying coefficients often necessitates the use of more complex operators, such as~\cite{dendy1982black}.

\subsection{The Multigrid Cycle}\label{sec:multigrid-cycles}
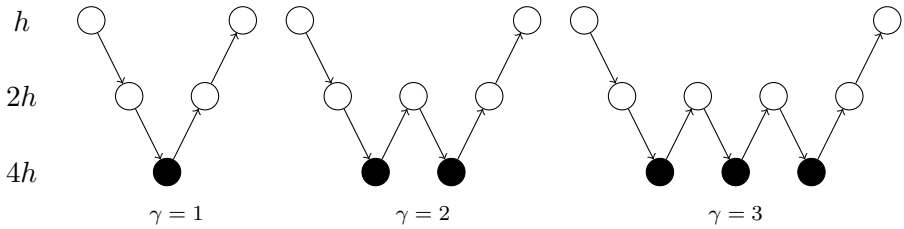
\begin{figure}
	\captionsetup{justification=centering}
   \begin{subfigure}[t!]{0.075\textwidth}
		\begin{tikzpicture}
			\node   (h) at (-0.75, 4){$h$};
			\node   (2h) at (-0.75, 3){$2h$};
			\node   (4h) at (-0.75, 2){$4h$};
		\end{tikzpicture}
		\subcaption*{\phantom{$\gamma = 1$}}
	\end{subfigure}
	\begin{subfigure}[t!]{0.21\textwidth}
		\begin{tikzpicture}
			\node	(a) at (0,4) [draw, circle,scale=1] {};
			\node	(b) at (0.5,3) [draw, circle,scale=1] {};
			\node	(c) at (1,2) [draw,circle,fill=black,scale=1] {};
			\node	(d) at (1.5,3) [draw, circle,scale=1] {};
			\node	(e) at (2,4) [draw, circle, scale=1] {};
			\draw 
			(a) edge[->] (b) 
			(b) edge[->] (c)
			(c) edge[->] (d)
			(d) edge[->] (e)   
			;
		\end{tikzpicture}
		\subcaption*{$\gamma = 1$}
	\end{subfigure}
	\begin{subfigure}[t!]{0.29\textwidth}
		\begin{tikzpicture}
			\node	(a) at (0,4) [draw, circle,scale=1] {};
			\node	(b) at (0.5,3) [draw, circle,scale=1] {};
			\node	(c) at (1,2) [draw, circle,fill=black,scale=1] {};
			\node	(d) at (1.5,3) [draw, circle,scale=1] {};
			\node	(e) at (2,2) [draw, circle, fill=black, scale=1] {};
			\node	(f) at (2.5,3) [draw, circle, scale=1] {};
			\node	(g) at (3,4) [draw, circle, scale=1] {};
			
			\draw 
			(a) edge[->] (b) 
			(b) edge[->] (c)
			(c) edge[->] (d)
			(d) edge[->] (e)   
			(e) edge[->] (f)
			(f) edge[->] (g)
			;
		\end{tikzpicture}
		\subcaption*{$\gamma = 2$}
	\end{subfigure}
	\begin{subfigure}[t!]{0.35\textwidth}
		\begin{tikzpicture}
			\node	(a) at (0,4) [draw, circle,scale=1] {};
			\node	(b) at (0.5,3) [draw, circle,scale=1] {};
			\node	(c) at (1,2) [draw, circle,fill=black,scale=1] {};
			\node	(d) at (1.5,3) [draw, circle,scale=1] {};
			\node	(e) at (2,2) [draw, circle, fill=black, scale=1] {};
			\node	(f) at (2.5,3) [draw, circle, scale=1] {};
			\node	(g) at (3,2) [draw, circle, fill=black,scale=1] {};
			\node	(h) at (3.5,3) [draw, circle, scale=1] {};	
			\node	(i) at (4,4) [draw, circle, scale=1] {};	
			\draw 
			(a) edge[->] (b) 
			(b) edge[->] (c)
			(c) edge[->] (d)
			(d) edge[->] (e)   
			(e) edge[->] (f)
			(f) edge[->] (g)
			(g) edge[->] (h)
			(h) edge[->] (i)
			;
		\end{tikzpicture}
		\subcaption*{$\gamma = 3$}
	\end{subfigure}
	\caption{Three-grid cycles ($k = 3$).}
	\label{fig:three-grid-cycles}
\end{figure}

\begin{figure}
	\captionsetup{justification=centering}
	\begin{subfigure}[t!]{0.075\textwidth}
		\begin{tikzpicture}
			\node   (h) at (-0.75, 4){$h$};
			\node   (2h) at (-0.75, 3){$2h$};
			\node   (4h) at (-0.75, 2){$4h$};
			\node   (8h) at (-0.75, 1){$8h$};
		\end{tikzpicture}
		\subcaption*{\phantom{$\gamma = 1$}}
	\end{subfigure}
	\begin{subfigure}[t!]{0.295\textwidth}
		\begin{tikzpicture}
			\node	(a) at (0,4) [draw, circle,scale=1] {};
			\node	(b) at (0.5,3) [draw, circle,scale=1] {};
			\node	(c) at (1,2) [draw, circle,scale=1] {};
			\node	(d) at (1.5,1) [draw, circle,fill=black, scale=1] {};
			\node	(e) at (2,2) [draw, circle, scale=1] {};
			\node	(f) at (2.5,3) [draw, circle,scale=1] {};
			\node	(g) at (3,4) [draw, circle,scale=1] {};
			\draw 
			(a) edge[->] (b) 
			(b) edge[->] (c)
			(c) edge[->] (d)
			(d) edge[->] (e)   
			(e) edge[->] (f)
			(f) edge[->] (g)
			
			;
		\end{tikzpicture}
		\subcaption*{$\gamma = 1$}
	\end{subfigure}
	\begin{subfigure}[t!]{0.59\textwidth}
		\begin{tikzpicture}
			\node	(a) at (0,4) [draw, circle,scale=1] {};
			\node	(b) at (0.5,3) [draw, circle,scale=1] {};
			\node	(c) at (1,2) [draw, circle,scale=1] {};
			\node	(d) at (1.5,1) [draw, circle,fill=black, scale=1] {};
			\node	(e) at (2,2) [draw, circle, scale=1] {};
			\node	(f) at (2.5,1) [draw, circle,fill=black,scale=1] {};
			\node	(g) at (3,2) [draw, circle,scale=1] {};
			\node	(h) at (3.5,3) [draw, circle,scale=1] {};
			\node	(i) at (4,2) [draw, circle,scale=1] {};
			\node	(j) at (4.5,1) [draw, circle,fill=black, scale=1] {};
			\node	(k) at (5,2) [draw, circle, scale=1] {};
			\node	(l) at (5.5,1) [draw, circle,fill=black,scale=1] {};
			\node	(m) at (6,2) [draw, circle,scale=1] {};
			\node	(n) at (6.5,3) [draw, circle, scale=1] {};
			\node	(o) at (7,4) [draw, circle, scale=1] {};
			
			\draw 
			(a) edge[->] (b) 
			(b) edge[->] (c)
			(c) edge[->] (d)
			(d) edge[->] (e)   
			(e) edge[->] (f)
			(f) edge[->] (g)
			(g) edge[->] (h)
			(h) edge[->] (i)
			(i) edge[->] (j)
			(j) edge[->] (k)
			(k) edge[->] (l)
			(l) edge[->] (m)
			(m) edge[->] (n)
			(n) edge[->] (o)
			;
		\end{tikzpicture}
		\subcaption*{$\gamma = 2$}
	\end{subfigure}
	\caption{Four-grid cycles ($k = 4$).}
	\label{fig:four-grid-cycles}
\end{figure}
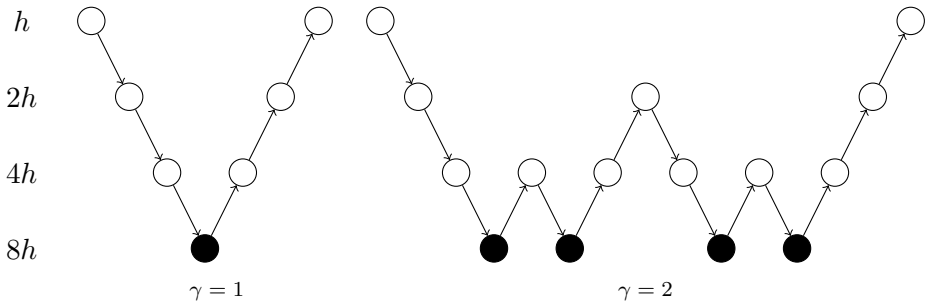

\begin{figure}
	\begin{subfigure}[t]{0.4\textwidth}
		\begin{tikzpicture}
			\node   (h) at (-0.75, 4){$h$};
			\node   (2h) at (-0.75, 3){$2h$};
			\node   (4h) at (-0.75, 2){$4h$};
			\node   (8h) at (-0.75, 1){\phantom{$8h$}};
			\node	(a) at (0,4) [draw, circle,scale=1] {};
			\node	(b) at (0.5,3) [draw, circle,scale=1] {};
			\node	(c) at (1,2) [draw, circle,fill=black,scale=1] {};
			\node	(d) at (1.5,3) [draw, circle,scale=1] {};
			\node	(e) at (2,2) [draw, circle, fill=black, scale=1] {};
			\node	(f) at (2.5,3) [draw, circle, scale=1] {};
			\node	(g) at (3,4) [draw, circle, scale=1] {};
			
			\draw 
			(a) edge[->] (b) 
			(b) edge[->] (c)
			(c) edge[->] (d)
			(d) edge[->] (e)   
			(e) edge[->] (f)
			(f) edge[->] (g)
			;
		\end{tikzpicture}
  		\subcaption*{}
	\end{subfigure}
	\begin{subfigure}[t]{0.59\textwidth}
		\begin{tikzpicture}
				\node   (h) at (-0.75, 4){$h$};
				\node   (2h) at (-0.75, 3){$2h$};
				\node   (4h) at (-0.75, 2){$4h$};
				\node   (8h) at (-0.75, 1){$8h$};
			\node	(a) at (0,4) [draw, circle,scale=1] {};
			\node	(b) at (0.5,3) [draw, circle,scale=1] {};
			\node	(c) at (1,2) [draw, circle,scale=1] {};
			\node	(d) at (1.5,1) [draw, circle,fill=black, scale=1] {};
			\node	(e) at (2,2) [draw, circle, scale=1] {};
			\node	(f) at (2.5,1) [draw, circle,fill=black,scale=1] {};
			\node	(g) at (3,2) [draw, circle,scale=1] {};
			\node	(h) at (3.5,3) [draw, circle,scale=1] {};
			\node	(i) at (4,2) [draw, circle,scale=1] {};
			\node	(j) at (4.5,1) [draw, circle,fill=black, scale=1] {};
			\node	(k) at (5,2) [draw, circle, scale=1] {};
			\node	(l) at (5.5,3) [draw, circle,scale=1] {};
			\node	(m) at (6,4) [draw, circle,scale=1] {};
			
			\draw 
			(a) edge[->] (b) 
			(b) edge[->] (c)
			(c) edge[->] (d)
			(d) edge[->] (e)   
			(e) edge[->] (f)
			(f) edge[->] (g)
			(g) edge[->] (h)
			(h) edge[->] (i)
			(i) edge[->] (j)
			(j) edge[->] (k)
			(k) edge[->] (l)
			(l) edge[->] (m)
			;
		\end{tikzpicture}
    		\subcaption*{}
	\end{subfigure}
	\caption{Computational pattern of the F-cycle with a different number of coarsening steps.}
	\label{fig:f-cycle}
\end{figure}
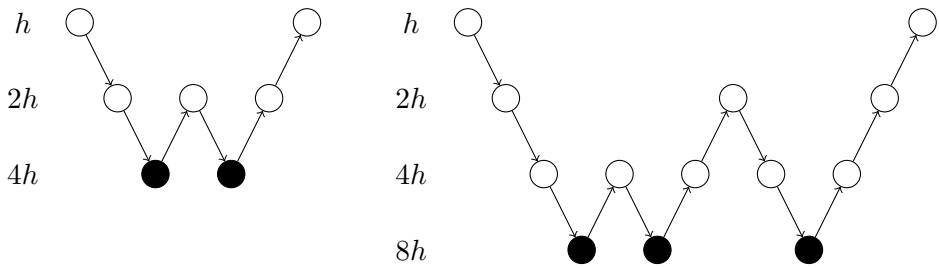

By putting this all together, we can now implement a recursive version of Algorithm~\ref{alg:two-grid-method} that allows us to perform an arbitrary number of coarsening steps until the resulting system of linear equations is small enough to be solved directly.
The corresponding implementation of a \emph{multigrid cycle} in the form of the function \textsc{mg-cycle} is shown in Algorithm~\ref{alg:multigrid-cycle}.
\begin{algorithm}[ht]
	\caption{Multigrid Cycle}
	\label{alg:multigrid-cycle}
	\begin{algorithmic}
		\Function{mg-cycle}{$k$, $\gamma$, $\bm{x}_h$, $A_h$, $\bm{b}_h$, $\nu_1$, $\nu_2$, $\omega$}
		\For{$i = 1, \dots, \nu_1$}
		\State{$\bm{x}_h = \bm{x}_h + \omega M_h^{-1} \left( \bm{b}_h - A_h \bm{x}_h \right)$ where $A_h = M_h + N_h$}
		\EndFor
		\State{$\bm{r}_h = \bm{b}_h - A_h \bm{x}_h$}
		\State{$\bm{b}_{2h} = I_h^{2h} \bm{r}_h$}
		\If{$k = 1$}
		\State{ Solve $A_{2h} \bm{x}_{2h} = \bm{b}_{2h}$ for $\bm{x}_{2h}$}
		\Else
		\State{$\bm{x}_{2h} = 0$}
		\For{$i = 1, \dots, \gamma$}
		\State{$\bm{x}_{2h} =  \textsc{mg-cycle}(\bm{x}_{2h}, A_{2h}, \bm{b}_{2h}, k-1, \gamma, \nu_1, \nu_2, \omega)$}
		\EndFor
		\EndIf
		\State{$\bm{x}_h = \bm{x}_h + I_{2h}^h \bm{x}_{2h}$}
		\For{$i = 1, \dots, \nu_2$}
		\State{$\bm{x}_h = \bm{x}_h + \omega M_h^{-1} \left( \bm{b}_h - A_h \bm{x}_h \right)$ where $A_h = M_h + N_h$}
		\EndFor
		\State \Return{$\bm{x}_h$}
		\EndFunction
	\end{algorithmic}
\end{algorithm}
Note that this function has a number of additional parameters compared to our original two-grid method.
First of all, $k$ defines the number of discretization levels and thus determines how many recursive coarsening steps need to be performed until the respective system of linear equations is solved directly.
Furthermore, since sometimes a single recursive application of this function is not sufficient to obtain a reasonably accurate approximation on the coarse grid, additional coarse-grid correction steps can be performed, as controlled by the parameter $\gamma$.
Finally, the parameters $\nu_1$ and $\nu_2$ specify the number of smoothing steps before and after the coarse-grid correction is performed, respectively.
Note that when applying multiple sweeps of smoothing or coarse-grid correction, the initial guess is replaced by the approximation obtained in the previous step.
While, in principle, the parameters of \textsc{mg-cycle} can be freely chosen, one usually classifies multigrid cycles according to the choice of $\gamma$, as each value yields a distinct computational pattern.
For instance, choosing $\gamma = 1$ means that only a single recursive descent is performed on each discretization level.
Figure~\ref{fig:three-grid-cycles} and~\ref{fig:four-grid-cycles} illustrate the algorithmic structure resulting from different values of $\gamma$ on a hierarchy of three and four grids, respectively.
Each white node corresponds to one or multiple steps of smoothing on the respective level, while a black node implies that the resulting error equation is solved directly.
Since, as it can be seen in Figure~\ref{fig:three-grid-cycles}, the choice of $\gamma = 1$ results in a V-shaped computational pattern, the corresponding multigrid method is usually called a V-cycle.
Similarly, the computational pattern of a method with $\gamma = 2$ on a three-grid hierarchy resembles the letter W, and hence the resulting method is called a W-cycle.
As it can be seen in Figure~\ref{fig:four-grid-cycles}, the amount of computational work within a W-cycle drastically increases with the number of coarsening steps.
While applying a V-cycle always results in significantly fewer computations, in many cases, a single coarse-grid correction step is not sufficient to reduce the low-frequency components of the initial error on the finest grid~\cite{trottenberg2000multigrid}.
Due to the resulting drastic increase in the number of computations, values of $\gamma$ larger than two are usually impractical for multigrid methods.
While Algorithm~\ref{alg:multigrid-cycle} enables the realization of different multigrid methods based on choosing different values for the parameters $k$, $\gamma$, $\nu_1$, $\nu_2$ and $\omega$, the structural composition of the resulting methods is limited.
For instance, in Algorithm~\ref{alg:multigrid-cycle}, it is assumed that the same value of $\gamma$ is chosen on each level, which restricts the set of possible multigrid cycles to those portrayed in Figure~\ref{fig:three-grid-cycles} and~\ref{fig:four-grid-cycles}.
One way to overcome this limitation is to combine different cycle types in a single method.
For instance, combining W- and V-cycles on each level results in a so-called F-cycle, whose algorithmic structure is illustrated in Figure~\ref{fig:f-cycle}.
Except for $k = 3$, in which case the F-cycle and W-cycle are equivalent, this method represents a compromise between a pure V-cycle and W-cycle.
While the composition of a multigrid method from different cycle types represents an additional degree of freedom, the individual cycles are still derived from the rules of Algorithm~\ref{alg:multigrid-cycle} and are thus fully described by the aforementioned parameters.
Since, according to the classical formulation of a multigrid method as in~\cite{brandt1977multi,hackbusch2013multi,trottenberg2000multigrid,briggs2000multigrid}, these parameters are considered global, adapting one of their values uniformly changes the method's behavior on each level. 
While for many applications, this approach has demonstrated to yield efficient multigrid methods~\cite{trottenberg2000multigrid}, there are still cases where multigrid could not yet achieve its full potential~\cite{ernst2012difficult,benzi2005numerical}.
Even though, to date, there exists a rich amount of research on the optimization of multigrid cycles based on a fixed set of global parameters, the variation of its components on an individual level has not been considered yet.
Such an approach would grant the flexibility to generate multigrid cycles that consist of varying restriction and coarse-grid correction steps, each with a different combination of smoothers and relaxation factors.
In this work, we aim to realize this vision by expressing each multigrid cycle as a finite sequence of instructions whose order is restricted by the rules of linear algebra and multigrid theory.
To specify these rules in a mathematically formal way, we make use of programming language theory, such that the design of an efficient multigrid method can be treated as a program synthesis task.
In the next section, we, therefore, provide a brief overview of programming language theory and introduce the concept of formal languages and grammars.
We then utilize this formalism together with the theoretical background presented in this chapter to derive a grammar-based representation of multigrid cycles that enables us to automate the design of these methods using evolutionary program synthesis techniques.

%% file: contents/formal_languages.tex
\section{Formal Languages and Grammars}
\label{sec:formal-languages}
In the previous section, we have introduced the theoretical background necessary to understand the basic functioning of multigrid methods.
Since the goal of this thesis is to automate the design of these methods, we next direct our attention to the area of formal languages and grammars, which will allow us to treat this problem as a program synthesis task.
A formal language represents a fundamental concept of computer science that enables the expression of arbitrarily-complex computational systems in an unambiguous and well-defined manner.
Similar to most natural languages, a formal language is characterized by its \emph{alphabet}. 
The alphabet of a formal language is a non-empty set of symbols $\Sigma$.
Based on $\Sigma$ \emph{strings} can be created, which represent finite sequences of symbols.
For instance, the alphabet $\Sigma = \left\{a, b\right\}$ contains only the symbols $a$ and $b$, which makes both $aabbb$ and $abaab$ strings on the alphabet $\Sigma$.
The length of a string, denoted by $|\cdot|$, is equal to the number of symbols contained in it. 
A concatenation of two strings 
\begin{equation}
	\begin{split}
		v & = a_1 a_2 \cdots a_n \\
		w & = b_1 b_2 \cdots b_m
	\end{split}
\end{equation}
is achieved by appending the second string at the end of the first string which yields
\begin{equation}
	vw = a_1 a_2 \cdots a_n b_1 b_2 \cdots b_m.
\end{equation}
The length of the concatenated string is thus equal to the sum of the lengths of the individual strings, which means that
\begin{equation}
|vw| = |v| + |w|.
\end{equation}
Finally, we define the empty string $\lambda$ as
\begin{equation}
	\begin{split}
		& |\lambda| = 0 \\
		& v \lambda = \lambda v = v,
	\end{split}
\end{equation}
where $v$ is a string on an arbitrary alphabet.
Assume $\Sigma$ is an alphabet, then $\Sigma^*$ is the set of strings obtained by concatenating zero or an arbitrary number of symbols from $\Sigma$.
Consequently, $\Sigma^*$ always contains the empty string $\lambda$.
To exclude $\lambda$ from the set, we define $\Sigma^+ = \Sigma^* \setminus \lambda$.
Since the number of strings that can be created from an alphabet through concatenation is infinite, both $\Sigma^+$ and $\Sigma^*$ represent infinite sets.
For example, if we again define $\Sigma = \left\{a, b\right\}$, then
\begin{equation*}
	\Sigma^{*} = \left\{\lambda, a, b, aa, ab, ba, bb, aaa, aab, aba, \dots \right\}.
\end{equation*} 
In general, for a given alphabet the \emph{language} $L$ is defined as a subset of $\Sigma^*$
\begin{equation}
	L \subset \Sigma^*.
	\label{eq:language-basic-definition}
\end{equation}
However, in practice, we usually want to define a language $L_G$ that represents a specific subset of $\Sigma^*$.
One way to achieve this is to specify a list of rules that generate $L_G$.
These rules can be defined in the form of a \emph{grammar} $G$.
\begin{definition}[Grammar]
	\begin{equation}
		G = \left(V, T, S, P \right),
	\end{equation}
	where $V$ is a finite set of \emph{variables},
	$T$ a finite set of \emph{terminals},
	$S \in V$ is the \emph{start} variable and 
	$P$ a finite set of \emph{productions} or \emph{production rules}.
	We also assume that $V \neq \emptyset$, $T \neq \emptyset$ and $V \cap T = \emptyset$.
\end{definition}
The core of a grammar is the definition of the productions $P$, which are usually specified as a list of mappings
\begin{equation*}
	x \to y,
	\label{eq:unrestricted-production}
\end{equation*}
where $x \in \left(V \cup T\right)^+$ and $y \in \left(V \cup T\right)^*$.
Note that alternatively, the operator
\begin{equation*}
	x \bnfpo y
\end{equation*}
is also often used to denote a production.
Each production $x \to y$ can then be applied in the following way.
Given a string $u$ of the form 
\begin{equation*}
	u = vxw,
\end{equation*}
we can replace $x$ with $y$ to \emph{derive} a new string
\begin{equation}
	u' = vyw,
\end{equation}
which is usually written as $u \Rightarrow u'$.
This process can then be continued by applying a sequence of derivations chosen arbitrarily from the set of available productions $P$.
Within this sequence, each production can be applied whenever its conditions on the left-hand side of the rule are fulfilled, i.e., the respective pattern occurs anywhere within the current string.
Also, note that there is no limit on how many times a certain production can be applied.
Each string $u$ that can be derived by applying an arbitrary sequence of productions starting from $S$
\begin{equation*}
	S \Rightarrow \dots \Rightarrow u,
\end{equation*}
is then an element of the language $L_G$.
Assuming that 
\begin{equation*}
	S \overset{*}{\Rightarrow} u
\end{equation*}
represents the application of an arbitrary sequence of productions, we can define the language $L_G$ generated by a grammar $G$ as follows:
\begin{definition}[Language]\label{def:language}
	Let $G = \left\{V, T, S, P\right\}$ be a grammar, then
	\begin{equation}
		L_G = \left\{u \in T^* : S \overset{*}{\Rightarrow} u\right\}
	\end{equation}
is the language generated by $G$.
\end{definition}
Assuming $u \in L_G$ and 
\begin{equation}
	S \Rightarrow u_1 \Rightarrow u_2 \Rightarrow \dots \Rightarrow u_n \Rightarrow u
\end{equation}
is a \emph{derivation} of $u$, then the strings $S, u_1, u_2, \dots, u_n$ are called its \emph{sequential forms}.
\subsection{The Chomsky Hierarchy}
\label{sec:chomsky-hierarchy}
So far, we have not yet imposed any restrictions on the individual components of a grammar $G$.
We call such a grammar, whose productions are of the general form of Equation~\eqref{eq:unrestricted-production}, \emph{unrestricted}.
It can be shown that any language generated by an unrestricted grammar is recursively enumerable, which means that there exists a Turing machine capable of enumerating all strings contained in that language~\cite{linz2006introduction}.
As a consequence, one can prove that unrestricted grammars are equally powerful as Turing machines and hence both computational models are equivalent.
Since we are only interested in languages that can be processed by any Turing-complete system on a modern computer, there is no use in considering languages that do not fall into this category.
While Turing machines and unrestricted grammars both represent a universal model of computation, it can be useful to consider grammars that impose certain restrictions on their productions and thus simplify both the derivation as well as the manipulation of strings.
\emph{Context-sensitive} grammars (CSG) represent the first step in this direction.
\begin{definition}[Context-Sensitive Grammar]
A grammar $G = \left\{V, T, S, P\right\}$ is context-sensitive if all productions can be written as
\begin{equation}
	xAy \to xuy,
\end{equation}
where $A \in V$ and $x, y, u \in \left(V \cup T\right)^*$.
\label{def:context-sensitive-grammar}
\end{definition}
This means that a certain production $A \to u$ can only be performed in the \emph{context} of $x$ and $y$ thus leading to the term context-sensitive.
While unrestricted grammars can be described by a Turing machine, the equivalent model of computation for a CSG is the \emph{linear bounded automaton}, which can be described as a Turing machine whose tape is linearly bounded by the length of the string~\cite{linz2006introduction}.
Furthermore, by prohibiting any context on the left-hand side of a production one arrives at the class of \emph{context-free} grammars (CFG).
\begin{definition}[Context-Free Grammar]
	A grammar $G = \left\{V, T, S, P\right\}$ is context-free if all productions are of the form
	\begin{equation}
		A \to u,
	\end{equation}
	where $A \in V$ and $u \in \left(V \cup T\right)^*$.
\end{definition}
The equivalent model of computation for a CFG is the \emph{pushdown automaton}~\cite{linz2006introduction}.  
CFGs have a number of desirable properties that explain their widespread use in computer science, especially in the theory of programming languages~\cite{pierce2002types}.
For instance, each string generated by a CFG, as well as the derivation of each string, can be represented as a tree.
Moreover, it is possible to check whether an arbitrary string of length $n$ is contained in the language generated by a CFG using only $\mathcal{O}(n^3)$ steps.
Note that the class of CFGs is contained in the class of CSGs, since we can simply replace each production 
\begin{equation*}
	A \to u
\end{equation*}
of a CFG by one of the form
\begin{equation*}
	\lambda A \lambda \to \lambda u \lambda.
\end{equation*} 
Since $\lambda \in \left(V \cup T\right)^*$, the resulting grammar meets the requirements of Definition~\ref{def:context-sensitive-grammar}.
Finally, the introduction of further restrictions on the productions leads to the class of regular grammars (RGs).
\begin{definition}[Regular Grammar]
	A grammar $G = \left\{V, T, S, P\right\}$ is regular if all productions are of the form
	\begin{equation}
		\begin{split}
			A & \to xB \\
			A & \to x,
		\end{split}
	\end{equation}
in which case $G$ is said to be \emph{left-linear} or of the form
	\begin{equation}
	\begin{split}
		A & \to Bx \\
		A & \to x,
	\end{split}
	\end{equation}
in which case $G$ is said to be \emph{right-linear}. In both cases, $A, B \in V$ and $x \in T^*$.
\end{definition}
Since both $xB \in \left(V \cup T\right)^*$ and $Bx \in \left(V \cup T\right)^*$, it is obvious that each regular grammar can also be considered as a CFG.
Each language generated by an RG is a regular language and can thus be defined in the form of a regular expression. 
RGs, therefore, possess the same expressive power as \emph{finite-state machines}~\cite{linz2006introduction}.
With the definition of RGs, we have arrived at the bottom of the so-called \emph{Chomsky hierarchy} of formal grammars and languages, which is given by the following set of relations
\begin{equation}
	\mathcal{G} \supset \mathcal{G}_{CS} \supset \mathcal{G}_{CF} \supset \mathcal{G}_{R}, 
\end{equation} 
where $\mathcal{G}$ is the set of all unrestricted (type-0), $\mathcal{G}_{CS}$ the set of all context-sensitive (type-1), $\mathcal{G}_{CF}$ the set of all context-free (type-2) and $\mathcal{G}_{R}$ the set of all regular (type-3) grammars.
The same is true for the corresponding languages
\begin{equation}
	\mathcal{L} \supset \mathcal{L}_{CS} \supset \mathcal{L}_{CF} \supset \mathcal{L}_{R}. 
\end{equation}
Consequently, each level in this hierarchy corresponds to a particular class of grammars introduced in this section, where each subsequent level introduces further restrictions on the productions but, on the other hand, enables the use of faster and more efficient algorithms for manipulating the strings contained in its generated language.
After we have now introduced a class of formal languages for representing strings of symbols in a structured way, as well as those grammars able to generate them, we will next discuss how we can represent the construction and manipulation of such strings as an optimization problem, whose solution we aim to approximate using evolutionary search methods.

%% file: contents/genetic_programming.tex
\section{Evolutionary Program Synthesis}
\label{sec:gggp}
Evolutionary program synthesis can be described as the automation of program generation through the use of evolutionary algorithms, which is a class of search algorithms for solving optimization problems based on the principle of natural evolution.
The class of evolutionary algorithms that is concerned with finding optimal programs is often summarized under the term \emph{Genetic Programming} (GP), which was first proposed by John Koza~\cite{koza1994genetic}.
In general, a program $p$ can be considered as a mapping from the set of inputs $\mathcal{I}$ to the corresponding set of outputs $\mathcal{O}$
\begin{equation}
	p : \mathcal{I} \to \mathcal{O}.
	\label{eq:gp-program}
\end{equation}
However, since not all input-output pairs are known in advance, this mapping is usually given in the form of a finite set of training cases.
The goal of GP is then to find a program that correctly computes the correct output for each input contained in the training set.
Since there is often not a unique program that satisfies this condition, usually a number of additional constraints are applied to assess the quality of each correct program.
For instance, in practice, one is often interested in finding the shortest or fastest program that is able to pass all training cases.
Based on a program's effectiveness in solving all training cases under the given constraints, GP assigns a \emph{fitness} value to each program, which is then treated as an \emph{individual} within a \emph{population} of programs.
In evolutionary computation, the population is the set of individuals currently considered within the search and, therefore, spans a subspace within the space of solutions for the given search problem.
Each subsequent step of the search is then performed by generating a new population based on the previous one through \emph{mutation} or recombination (often called \emph{crossover}) of the individuals contained in the current one.
Usually, the candidates for mutation and crossover are sampled from the current population based on the fitness of each individual.
If we repeat this procedure for a number of $n$ steps, we arrive at a final population $P_n$, which can then additionally be evaluated on a validation set.
The resulting search method is summarized in Algorithm~\ref{alg:genetic-programming}, which gives an overview of the general structure of a GP method.
\begin{algorithm}[t]
	\caption{Genetic Programming}
	\label{alg:genetic-programming}
	\begin{algorithmic}[1] % The number tells where the line numbering should start
		\State \textbf{Randomly generate} an initial population $P_0$ of programs
		\State \textbf{Evaluate} $P_0$ on the \textbf{training set} 
		\For{$i := 1, \dots, n$}
		\State \textbf{Select} a subset of individuals $M_i \subset P_{i-1}$ based on their \textbf{fitness}
		\State \textbf{Generate} new programs $C_i$ based on $M_i$ with \textbf{mutation} and \textbf{crossover}
		\State \textbf{Evaluate} $C_i$ on the \textbf{training set} 
		\State \textbf{Select} $P_{i}$ from $C_i \cup P_{i-1}$
		\EndFor
		\State \textbf{Evaluate} the final population $P_{n}$ on a \textbf{validation set}  to obtain the best overall program
	\end{algorithmic}
\end{algorithm}
However, we have not yet defined the individual operations, such as the generation of an initial population and the creation of new individuals through mutation and crossover.
Since all these operations are based on manipulating the internal structure of a given program, the first step towards the implementation of a GP method is the choice of a suitable program representation.
Note that this representation does not necessarily need to be equal to the target language in which the actual program is supposed to be implemented but rather needs to define a unique mapping that enables its automatic generation.
This process is usually called \emph{genotype} to \emph{phenotype} mapping, where the genotype refers to the internal representation used within the GP method while the phenotype represents the actual program implemented on the target machine.
One of the most widely used genotype representations is tree-based GP, where each program is internally represented as a tree of expressions~\cite{koza1994genetic,poli2008field}.
\subsection{Representation}\label{sec:gggp-representation}
\begin{figure}
	\begin{subfigure}{0.58\textwidth}
 \centering
		\includegraphics[scale=0.9]{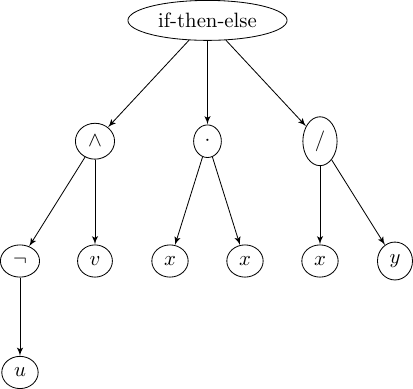}
	\end{subfigure}
	\begin{subfigure}{0.41\textwidth}
  \centering
		\includegraphics[scale=0.9]{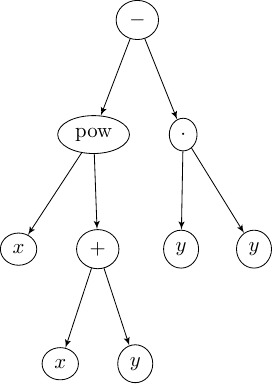}
	\end{subfigure}
 \caption[Expression trees of the functions $f_1$ and $f_2$~\eqref{eq:gp-example-functions}.]{Expression trees of the functions $f_1$ and $f_2$~\eqref{eq:gp-example-functions}, where we have simplified the if-then-else construct to a single ternary operator.}
 \label{fig:gp-expression-tree-examples}
\end{figure}

\begin{figure}[h]
	\begin{subfigure}{0.58\textwidth}
 \centering
	\includegraphics[scale=0.46]{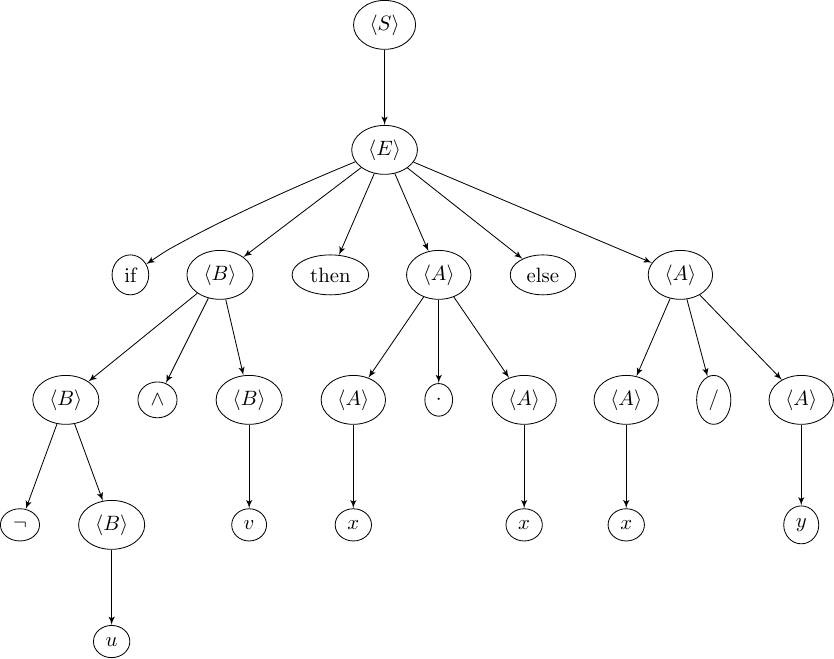}
	\end{subfigure}
	\begin{subfigure}{0.41\textwidth}
 \centering
	\includegraphics[scale=0.46]{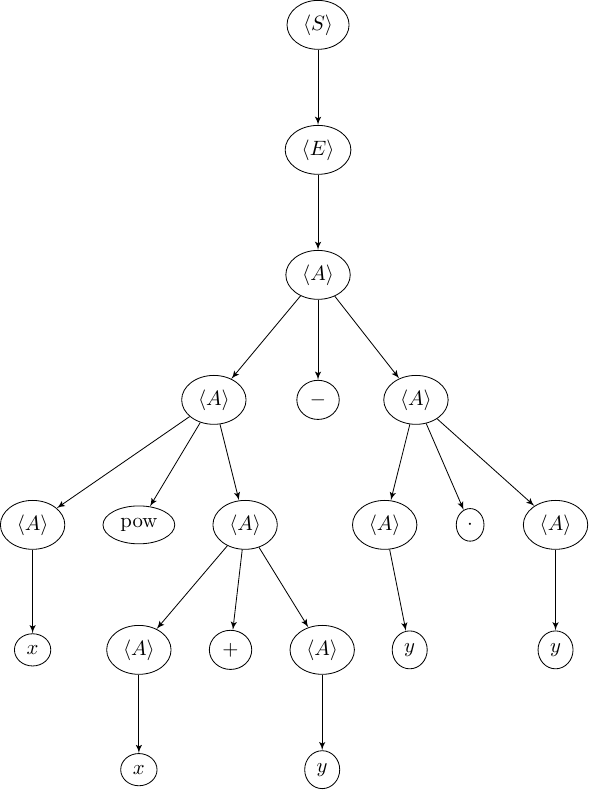}
	\end{subfigure}
	\caption{Derivation trees of the functions $f_1$ and $f_2$~\eqref{eq:gp-example-functions}.}
	\label{fig:gp-derivation-tree-examples}
\end{figure}

In contrast to other evolutionary algorithms, which represent the solution to an optimization problem as an array of discrete or continuous numbers~\cite{back1997handbook}, in evolutionary program synthesis, we are concerned with the automated discovery of programs. Therefore, each discovered solution corresponds to an executable program.
%Before we consider these operations in detail, we first need to define a procedure to generate program expression trees in a structured way. 
To illustrate this, we consider the following example grammar $G$ with the productions
\begin{equation}
	\begin{split}
		\ps S \; \; \bnfpo & \; \; \ps E \\
		\ps E \; \; \bnfpo & \; \; \text{if} \; \ps B \; \text{then} \; \ps E \; \text{else} \; \ps E \; | \; \ps A \\
		\ps A \; \; \bnfpo & \; \; -\ps A \; | \; (\ps A + \ps A) \; | \; (\ps A - \ps A) \; | \\
		 & \; \; (\ps A \cdot \ps A) \; | \; (\ps A / \ps A) \; | \ps{A}^{\ps{A}} \; | \; x \; | \; y \\  
		\ps B \; \; \bnfpo & \; \;  \neg \ps B \; | \; (\ps B \wedge \ps B) \; | \; (\ps B \vee \ps B) \; | \; u \; | \; v,
	\end{split}
\label{eq:gp-example-grammar}
\end{equation}
where $V = \{\ps S, \ps E, \ps A, \ps B\}$ is the set of variables and $S = \ps S$ the start variable.
Semantically the symbols $x$, $y$ represent numbers while the symbols $u$, $v$ correspond to boolean values, i.e. $\top$ or $\bot$.
Note that $G$ is context-free since the left-hand side of each production exclusively consists of a single variable.
If we treat $x$, $y$, $u$, and $v$ as an input, each expression generated by $G$ can be considered as a quaternary function $f(x,y,u,v) \in L_{G}$, where $L_G$ is the language generated by $G$ according to Definition~\ref{def:language}.
For example, the functions
\begin{equation}
	\begin{split}
		f_1(x,y,u,v) & = \text{if} \; (\neg u \wedge v) \; \text{then} \; (x \cdot x) \; \text{else} \; (x / y) \\
		f_2(x,y,u,v) & = x^{(x + y)} - (y \cdot y)
	\end{split}
\label{eq:gp-example-functions}
\end{equation} can both be generated by $G$ and are thus included in $L_G$.
By applying the rules of arithmetic and boolean algebra, we can compute the result of each such function based on a given list of inputs.
Therefore, assuming both $x$ and $y$ represent real-valued numbers, we can evaluate the correctness of the mapping 
\begin{equation*}
    f: x, y, u, v \to \mathbb{R}
\end{equation*}
for a variable number of problem instances.
Next, to define operations for generating arbitrary functions $f \in L_G$, we need to choose a suitable data structure for its representation.
While a program's structure depends on the programming language it is written in, in many cases, it is possible to represent it as a tree of program expressions.
This is especially true for languages from the Lisp family, which were originally used by John Koza. 
His work can be considered the first implementation of GP~\cite{koza1994genetic}.
In tree-based GP, all operations are performed on program expression trees, and thus both mutation and crossover must be designed with respect to this representation.
%TODO: Probably OK to leave out
%Expression trees can be created in a straightforward way after rewriting each expression in postfix notation.
%To obtain the corresponding tree, the expression is traversed sequentially by putting each encountered symbol, either in the form of a terminal or an operator that has already been processed, on a stack.
%For each new operator encountered we then first retrieve the required operands from the stack based on which the tree is created in a bottom-up manner.
Figure~\ref{fig:gp-expression-tree-examples} shows the corresponding expression trees for $f_1$ and $f_2$.
Based on a given expression tree, we can easily restore the original expression by traversing the tree recursively in a top-down manner while generating the corresponding expression for each node as soon as each of its operands has been recursively processed or if it is available as a terminal symbol.
While expression trees can be easily generated from a given program and offer an efficient way to evaluate it, they possess an inherent limitation, which is the absence of type information.
For a better understanding of this limitation, we again consider the context-free grammar shown in Equation~\eqref{eq:gp-example-grammar}.
Here, all productions on the variable $\ps B$ only result in the generation of boolean expressions, while all productions starting from $\ps A$ exclusively lead to arithmetic expressions.
As a consequence, we are not allowed to intermix expressions that are derived from $\ps A$ and $\ps B$.
However, this information is not contained in any of the expressions generated by $G$ and is thus also missing in the corresponding tree.
As the generation of each new expression is performed based on the productions of $G$, it is guaranteed to be valid.
However, the main step of GP is the creation of new individuals based on the ones contained in the population, either through the recombination of two individuals (crossover) or by altering certain parts of an individual (mutation). 
Both operations require us to investigate which nodes and subtrees of an expression tree can be safely replaced by an alternative branch without violating the type constraints imposed by the productions of our grammar.
While it is possible to deal with this problem by simply evicting those individuals that violate the type constraints from the population, depending on the number of constraints, this will be the case for a high percentage of individuals, rendering this method extremely inefficient.
A different and usually more efficient approach is to annotate each node within an expression tree with additional type information, which leads to the concept of strongly-typed GP, first proposed by David Montana~\cite{montana1995strongly}.
Strongly-typed GP represents a viable solution to the problem of retaining type correctness, but it requires us to transform the implicit type information contained within the productions of our grammar into explicit type annotations at the nodes of each expression tree.
As an alternative, we can instead utilize the information encoded in the productions of a grammar $G$ by considering the \emph{derivations} of an expression $e \in L_G$, as defined by
\begin{equation*}
	S \Rightarrow e_1 \Rightarrow e_2 \Rightarrow \dots \Rightarrow e_n \Rightarrow e.
\end{equation*}
%Furthermore, if $G$ is regular, its derivations can be further simplified to finite state machine, where each state corresponds to a sequential form.
Consider again the grammar whose productions are shown in Equation~\eqref{eq:gp-example-grammar} and the example functions $f_1$ and $f_2$ as defined in Equation~\eqref{eq:gp-example-functions}.
% TODO Produce always leftmost symbol
Based on these productions, we can formulate derivations that lead to the expressions of both functions, where
\begin{equation}
	\begin{aligned}
		\bps{S} & \Rightarrow \bps{E} \Rightarrow \text{if} \; \bps B \; \text{then} \; \ps E \; \text{else} \; \ps E 
		\\ & \Rightarrow \text{if} \; (\bps B \wedge \ps B) \; \text{then} \; \ps E \; \text{else} \; \ps E
		\\ & \Rightarrow \text{if} \; (\neg \bps B \wedge \ps B) \; \text{then} \; \ps E \; \text{else} \; \ps E 
		\\ & \Rightarrow \text{if} \; (\neg u \wedge \bps B) \; \text{then} \; \ps E \; \text{else} \; \ps E 
		\\ & \Rightarrow \text{if} \; (\neg u \wedge v) \; \text{then} \; \bps E \; \text{else} \; \ps E 	
		\\ & \Rightarrow \text{if} \; (\neg u \wedge v) \; \text{then} \; \bps A \; \text{else} \; \ps E 	
		\\ & \Rightarrow \text{if} \; (\neg u \wedge v) \; \text{then} \; (\bps A \cdot \ps A) \; \text{else} \; \ps E 	 
		\\ & \Rightarrow \text{if} \; (\neg u \wedge v) \; \text{then} \; (x \cdot \bps A) \; \text{else} \; \ps E 
		\\ & \Rightarrow \text{if} \; (\neg u \wedge v) \; \text{then} \; (x \cdot x) \; \text{else} \; \bps E
		\\ & \Rightarrow \text{if} \; (\neg u \wedge v) \; \text{then} \; (x \cdot x) \; \text{else} \; \bps A
		\\ & \Rightarrow \text{if} \; (\neg u \wedge v) \; \text{then} \; (x \cdot x) \; \text{else} \; (\bps A / \ps A)
		\\ & \Rightarrow \text{if} \; (\neg u \wedge v) \; \text{then} \; (x \cdot x) \; \text{else} \; (x / \bps A)
		\\ & \Rightarrow \text{if} \; (\neg u \wedge v) \; \text{then} \; (x \cdot x) \; \text{else} \; (x / y) = f_1(x,y,u,v)
	\end{aligned}
\label{eq:gp-derivation-f1}
\end{equation}
represents a valid derivation for $f_1$ while the same is true for the derivation
\begin{equation}
\begin{aligned}
	\bps{S} & \Rightarrow \bps E \Rightarrow \bps A \Rightarrow \bps A - \ps A \Rightarrow \bps{A}^{\ps A} - \ps A
	\Rightarrow x^{\bps A} - \ps A 	
	\\ & \Rightarrow x^{(\bps A + \ps A)} - \ps A \Rightarrow x^{(x + \bps A)} - \ps A 
	\\ & \Rightarrow x^{(x + y)} - \bps A \Rightarrow x^{(x + y)} - (\bps A \cdot \ps A) 
	\\ & \Rightarrow x^{(x + y)} - (y \cdot \bps A) \Rightarrow x^{x + y} - (y \cdot y) = f_2(x,y,u,v)
\end{aligned}
\end{equation}
with respect to $f_2$. 
Note that while each derivation contains the same sequential forms, the order in which they occur is not fixed and depends on the application order of the productions.
Here, the variable (written in bold font) that occurs leftmost within each sequential form is transformed next.
If a grammar is context-free, which is the case for the example grammar shown in Equation~\eqref{eq:gp-example-grammar}, each of its possible derivations can be represented as a tree, where the root node is always the start variable $S$, and each vertex of the tree corresponds to a transition between sequential forms within the derivation~\cite{linz2006introduction}.
More formally, we can define a derivation tree as follows:
\begin{definition}[Derivation Tree]\label{def:derivation-tree}
	Let $G = \left\{V, T, S, P\right\}$ be a context-free grammar, then a derivation tree is an ordered tree with the following properties
	\begin{enumerate}
		\item The start variable $S$ is the root of the tree.
		\item If $a$ is a leaf, then $a \in T \cup \{\lambda \}$.
		\item If $A$ is an interior node, then $A \in V$.
		\item If $A$ is an interior node and its children are $a_1, a_2, \dots, a_n$, then $P$ must contain a production of the form
		\begin{equation*}
			A \to a_1 a_2 \dots a_n.
		\end{equation*}  
	\end{enumerate}
\end{definition}
Note that multiple derivations can refer to the same tree, as the order in which the productions are performed is not specified within its structure.
%TODO discuss ambiguity yes or not? 
%In general, a context-free grammar $G$ is said to \emph{ambiguous} if there exists an expression $e \in L_G$ that has at least two distinct derivation trees~\cite{linz2006introduction}.
%With respect to the terminology of evolutionary algorithms, this means that there exist multiple genotype representations that map to the same phenotype, which results in an unnecessarily bloated search space.
Again, as an example, the derivation trees for $f_1$ and $f_2$ are shown in Figure~\ref{fig:gp-derivation-tree-examples}.
If we consider the structure of these trees, the advantage of this representation becomes obvious.
In contrast to Figure~\ref{fig:gp-expression-tree-examples}, a derivation tree stores the variable on the left-hand side of each production as a parent node of the respective subtree.
Each variable encodes which production can be applied at a particular point within a derivation.
Consequently, by only allowing the replacement of a certain subtree by those productions whose left-hand side matches with its parent node, we can ensure that the resulting tree will always correspond to an expression contained in $L_G$. 
However, while derivation trees contain additional information, they are still structurally close to the expression trees originally introduced in this section, which means that a certain structural change within a derivation tree will have a comparably strong effect on the corresponding phenotypic expression.
In the theory of evolutionary computation, this property is referred to as \emph{locality}.
%In fact, we can easily obtain the corresponding expression tree by applying a bottom-up replacement of each variable node by a vertex that represents a combination of its direct leafs.
To describe the use of derivation trees as a genotype representation within GP, the term tree-based \emph{grammar-guided genetic-programming} (G3P) has been introduced~\cite{mckay2010grammar,whigham1995grammatically}.
Since tree-based G3P combines the advantage of ensuring automatic type correctness with a high locality between genotype and phenotype, it represents a widely-used GP variant~\cite{mckay2010grammar}.
Furthermore, within the last decades, \emph{grammatical evolution} (GE) has been established as a valid alternative~\cite{oneill2001grammatical,zz_ge1,zz_ge2}, where the genotype is encoded as a dynamically-sized array of integers. 
Each entry of this array then corresponds to the choice of a particular production.
GE offers the advantage of having a simpler typesafe genotype representation, which, however, comes at the price of requiring additional steps for the genotype-to-phenotype mapping.
Furthermore, a small change in the genotype of a GE program does not necessarily result in a small phenotypic change since the interpretation of each array entry depends on its predecessors.
Therefore, changing a single value within the genotype can affect the interpretation of all subsequent genes.
For instance, consider the encoding
\begin{equation}
	\begin{split}
		\ps S \; \; \bnfpo & \; \; \underbrace{\ps E}_{0} \\
		\ps E \; \; \bnfpo & \; \; \underbrace{\text{if} \; \ps B \; \text{then} \; \ps E \; \text{else} \; \ps E}_{0} \; | \; \underbrace{\ps A}_{1} \\
		\ps A \; \; \bnfpo & \; \; \underbrace{-\ps A}_0 \; | \; \underbrace{(\ps A + \ps A)}_{1} \; | \; \underbrace{(\ps A - \ps A)}_2 \; |   \\  
		&  \; \; \underbrace{(\ps A \cdot \ps A)}_{3} \; | \; \underbrace{(\ps A / \ps A)}_{4} \; | \; \underbrace{\ps{A}^{\ps A}}_{5} \; | \; \underbrace{x}_{6} \; | \; \underbrace{y}_{7} \\
		\ps B \; \; \bnfpo & \; \; \underbrace{\neg \ps B}_{0} \; | \; \underbrace{(\ps B \wedge \ps B)}_{1} \; | \; \underbrace{(\ps B \vee \ps B)}_{2} \; | \; \underbrace{u}_3 \; | \; \underbrace{v}_4
	\end{split}
	\label{eq:gp-example-grammar-ge-encoding}
\end{equation}
of our example grammar, based on which we can represent $f_1$ as an array consisting of the following integers:
\begin{equation*}
	\left[ 0, 0, 1, 0, 3, 4, 1, 3, 6, 6, 1, 4, 6, 7 \right].
\end{equation*}
Since we are always processing the leftmost variable first, each array entry to a sequential form in the derivation shown in Equation~\eqref{eq:gp-derivation-f1}.
While in this case, we have chosen each number from the available range of productions, in practice, we can always map any non-negative integer to this range by applying a modulo $n$ operation, where $n$ is the number of available productions for each variable.
We can now investigate how changing a single array entry will affect the interpretation of all subsequent ones.
For instance, changing the second entry to a value of one results in the choice of the production 
\begin{equation*}
	\ps E \Rightarrow \ps A,
\end{equation*}
instead of
\begin{equation*}
	\ps E \Rightarrow \text{if} \; (\ps B \wedge \ps B) \; \text{then} \; \ps E \; \text{else} \; \ps E,
\end{equation*}
which means that all subsequent entries are now interpreted as productions starting from $\ps A$.
This example clearly illustrates that even slight changes in the genotype of a GE representation can have a dramatic effect on its phenotypic expression.
However, it needs to be mentioned that newer GE variants mitigate this issue by learning a statistical distribution about the effect of each production on the fitness of an individual~\cite{zz_ge1,zz_ge2}.
Since tree-based G3P and GE both represent a different tradeoff between the simplicity and locality of their genotype representation, in general, it is impossible to predict which variant will lead to a better outcome.
Furthermore, we have to mention that, in addition to tree-based GP and GE, numerous other variants have been proposed since the original invention of GP~\cite{poli2008field}, the most prominent ones being linear GP~\cite{brameier2007linear} and cartesian GP~\cite{miller2008cartesian}. 
However, since the implementation presented in this thesis is based on tree-based G3P, we will exclusively focus on this variant in our subsequent treatment of mutation, crossover, and fitness evaluation.
%TODO Probably too much, since GE is not considered further
%Furthermore, if at a certain position within the array the corresponding sequential form contains already no more variables, the remainder is usually ignored or removed.
%Now we can investigate the effect of changing a single number within this array on the resulting derivation.
%For instance changing the second entry leads to 
%\begin{equation}
%	\left[ 0, 1, 1, 0, 3, 4, 1, 3, 6, 7, 1, 4, 6, 7 \right],
%\end{equation}
%which can be translated to the following derivation
%\begin{equation}
%	\begin{aligned}
%		0 \mod 1 = 0: \quad & \bps{S} \Rightarrow \bps E \\
%		1 \mod 2 = 1: \quad & \Rightarrow \bps A \\
%		1 \mod 8 = 1: \quad & \Rightarrow \bps A + \ps A \\ 
%		0 \mod 8 = 0: \quad & \Rightarrow -\bps A + \ps A \\ 
%		3 \mod 8 = 3: \quad & \Rightarrow -(\bps A \cdot \ps A) + \ps A \\ 
%		\dots & 
%	\end{aligned}
%\end{equation}
\subsection{Initialization}
\label{sec:gggp-initialization}
\begin{figure}
	\centering
	\includegraphics[scale=0.46]{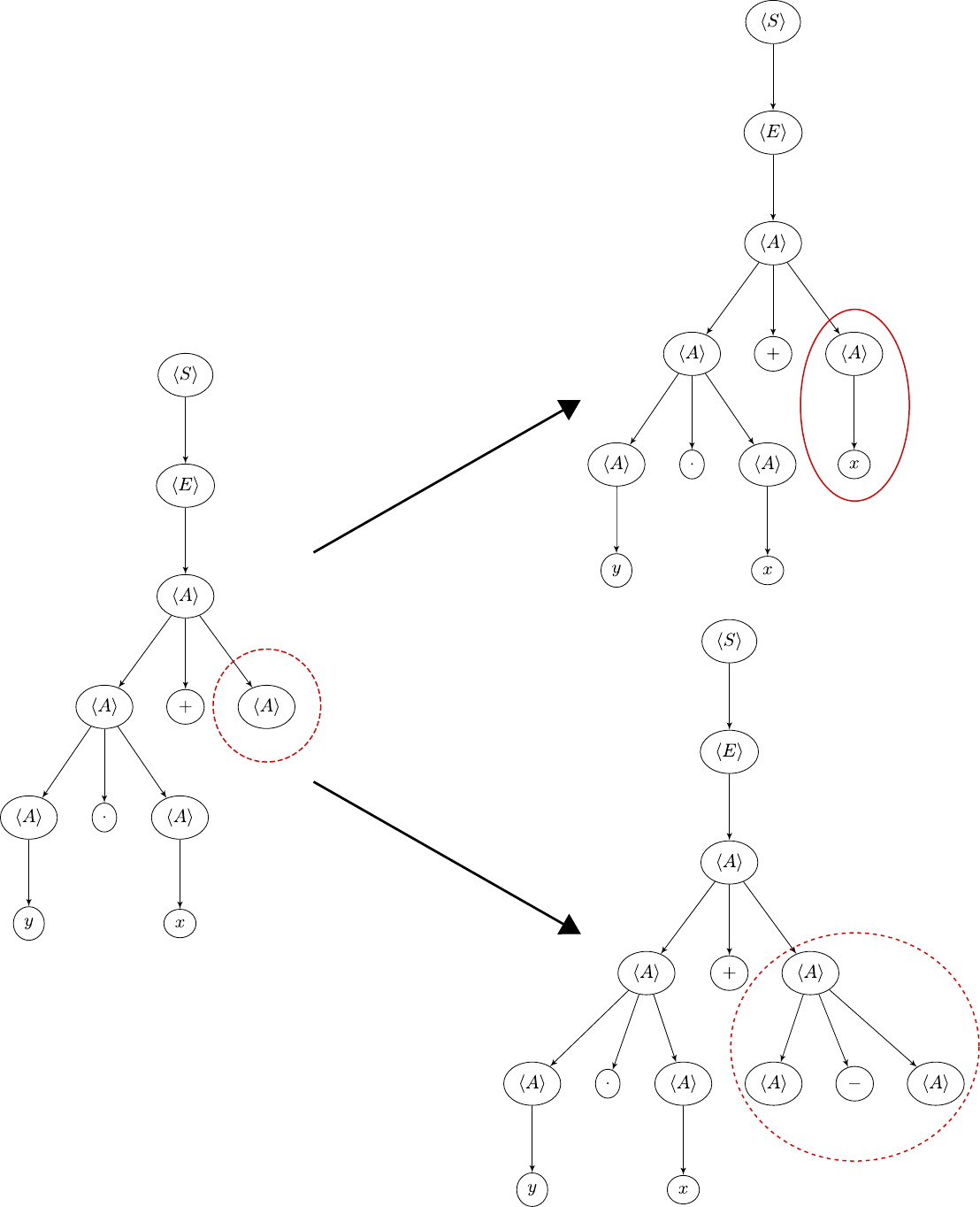}
	\caption[Options for extending a derivation tree]{Options for extending a derivation tree - An unfinished branch consisting of a variable can either be extended by applying a production that includes at least one variable (lower-right tree), or it is ended by inserting a terminal symbol as a leaf (upper-right tree).}
	\label{fig:gp-tree-growing-options}
\end{figure}

After choosing a suitable genotype representation, the next ingredient of an evolutionary program synthesis approach is the definition of suitable operators for the creation of new programs, either from scratch or based on an existing population.
Due to the stochastic nature of evolutionary algorithms, these operations are usually performed with a certain degree of randomness.
As shown in Algorithm~\ref{alg:genetic-programming}, the first step within each evolutionary search method is the initialization of the population. 
In some cases, it is possible to choose a certain proportion of the initial population from a set of promising individuals, either generated in previous experiments on a similar problem or hand-picked by a human domain expert, which is called \emph{seeding}.
However, seeding the population with too many individuals may introduce an excessive bias into the search, and as a result, the descendants of the seeded individuals might quickly take over the population~\cite{poli2008field}.
This risk can be mitigated by including a sufficiently high number of randomly-generated individuals in the initial population. 
In tree-based G3P, each individual is represented as a grammar derivation tree.
Therefore, to generate an individual from scratch, starting with $S$ as a root node, we can successively extend a tree by choosing a production for each node that corresponds to a variable until all leaves consist exclusively of terminals.
Note that within the choice of each production, one has to decide between growing the tree further, which corresponds to applying a production that contains at least one variable, or between cutting the current branch off by choosing one that generates exclusively terminal symbols.
Consider, for instance, the tree shown in Figure~\ref{fig:gp-tree-growing-options}, which has been generated by the example grammar formulated in Equation~\eqref{eq:gp-example-grammar} and includes an unfinished branch with the variable $\ps{A}$.
The productions available for the variable $\ps{A}$ can be classified into two categories, which we call \emph{terminal} and \emph{non-terminal productions}.
Terminal productions are those that generate subexpressions exclusively consisting of terminal symbols, leading to the creation of one or multiple leaves and thus ending the growth of the tree at this particular point.
This is illustrated in the upper-right tree in Figure~\ref{fig:gp-tree-growing-options}.
In contrast, non-terminal productions include all those that generate expressions with at least one variable.
The consequence of applying such a production is that the growth of the tree continues at the respective branch based on the productions available for the generated variables.
This is illustrated in the lower-right tree shown in Figure~\ref{fig:gp-tree-growing-options}. 
We can, therefore, utilize these different types of productions to control the shape of the derivation trees generated during the initialization of the population.
Two commonly used initialization operators that are based on this observation are the \emph{full} and \emph{grow} operators~\cite{poli2008field}.  
The goal of both operators is the generation of a tree that satisfies a certain \emph{depth} limit, which is defined as the maximum over the length of all paths from the root node to any leaf.
As the name suggests, the \emph{full} operator aims to generate a tree with even depth throughout all branches, which means that every path to a leaf has roughly the same length.
This is, for instance, achieved by applying non-terminal productions within each subtree until the necessary depth is reached.
In contrast, the \emph{grow} operator does not impose any constraints on the depth of each path within the tree but instead only requires the longest path to fulfill this condition, which means that both terminal and non-terminal productions are allowed until the depth limit is reached.
Note that, in general, a grammar is not required to include both terminal and non-terminal productions for each variable, and hence, the \emph{full} operator can only be applied to a subset of all possible grammars.
In contrast, the \emph{grow} operator is more widely applicable since it only requires a single path within the tree to reach a certain depth. 
It, however, has the disadvantage that the average shape of the generated trees strongly depends on the ratio of the terminal to non-terminal productions and the arity of the latter, i.e., the number of variables in the generated expressions.
One way to mitigate this problem is to adapt the probability of choosing a production from both categories accordingly.
Finally, another possibility to initialize a population is to employ both operators, \emph{full} and \emph{grow}, on a certain proportion of the individuals using a variety of different depth limits. 
If the proportion of individuals is roughly 50 \% for both operators, the initialization is called \emph{ramped half-and-half}~\cite{poli2008field,koza1994genetic}.
Furthermore, while both initialization operators presented here can also be applied to other tree-based GP variants, alternative methods specifically designed for G3P have been proposed in~\cite{garcia2006initialization,criado2020grammatically}.
These operators aim to exploit the structure of a given grammar, for instance, in order to generate a population that is more uniformly distributed across the search space~\cite{criado2020grammatically}.

\subsection{Fitness Evaluation and Selection}
\label{sec:gggp-evaluation-and-selection}
Since G3P is a program optimization technique, each individual corresponds to a program in the form of Equation~\eqref{eq:gp-program}.
However, as the search itself is performed on the genotype representation of each individual, as described in the last section, this representation first needs to be translated to the phenotype, i.e., a program in the target language.
In many program synthesis tasks, the grammar employed within the search is either equivalent or represents a subset of that of the target programming language.
While in such cases, a direct correspondence between the derivation and expression tree of the target programming language exists, there is also the possibility to employ multiple layers of intermediate representations before an actual program is generated, which necessitates the use of multiple code generation steps.
After an executable program is available in the target language, the final step within fitness evaluation is to measure its performance in a number of problem instances, which is then condensed to one or multiple measures of an individual's fitness.
For instance, one simple way to assess a program's quality is to represent its performance in each individual problem instance as a separate optimization objective.
An alternative possibility is to represent different aspects of a program's quality in the form of distinct objectives, each of which is determined from a combination of its performance in multiple or even all considered problem instances.
With respect to fitness evaluation, GP algorithms can be categorized into single- and multi-objective variants.
Similar to other evolutionary algorithms, GP treats fitness evaluation as a black box, which means that its inner workings are independent of the computational details of fitness evaluation.
The search is then performed exclusively based on the resulting fitness landscape~\cite{pitzer2012comprehensive}.
However, as fitness evaluation often represents the computationally most expensive part of an evolutionary algorithm, it drastically impacts to what extent the given search space can be explored.
Since each individual represents a distinct program that can be executed on a unique set of problem instances, each evaluation can be performed independently, which facilitates its concurrent execution on a multi-processor system.
For this purpose, parallel~\cite{sudholt2015parallel} and distributed~\cite{gong2015distributed} evolutionary algorithm variants, which can be executed on recent manycore architectures and computer clusters, have been proposed.
Due to the inherent similarity of the internal structure of Algorithm~\ref{alg:genetic-programming} to other evolutionary algorithms, these approaches are equally applicable to GP-based search methods.

After a single- or multi-objective fitness value has been assigned to each individual, the question remains how the search should be progressed by selecting a number of individuals for the application of the two main evolutionary operators, recombination and mutation.
While initialization aims to generate a population that is evenly spread over the whole search space, the purpose of mutation and recombination is to create novel individuals based on an existing population that are located in a more promising subspace.  
For this purpose, evolutionary algorithms usually employ an additional selection phase to identify candidates in the current population based on which these operations are performed.
Note that within this process, certain individuals can be selected multiple times.
Common operators are fitness-proportionate selection, where individuals are randomly chosen according to a probability equal to their relative fitness value~\cite{lipowski2012roulette}, and tournament selection, where individuals are compared in a tournament-like fashion to identify the fittest among them~\cite{fang2010review}.
In general, selection is independent of the genotype representation of an individual but depends solely on the objective function of the target optimization problem.
Consequently, a selection operator is agnostic of the type of evolutionary algorithm it is applied to. 
For this reason, we only briefly discuss selection here, and for a more complete treatment of this operation, the reader is referred to one of the following publications:~\cite{back1997handbook,beyer2002evolution,goldberg1991comparative}.
To deal with multi-objective optimization problems, several selection operators have been proposed, about which an overview can be found in~\cite{coello2007evolutionary,deb2011multi,deb2015multi}.
Many of these operators are based on the idea of establishing a domination hierarchy between the individuals of the population, whereby an individual is said to dominate another one if it possesses superior fitness in all objectives.
One way to determine the dominance relation between all individuals in the population is to employ a non-dominated sorting procedure, which has been, for instance, proposed in~\cite{deb2002fast,deb2013evolutionary}.
The actual selection is then performed using a dominance-based tournament selection, while additional diversity metrics can be incorporated to decide between non-dominating individuals~\cite{coello2007evolutionary}.
Another recently proposed multi-objective selection operator that is especially geared towards program synthesis problems where the fitness evaluation is performed on multiple problem instances is lexicase selection.
In lexicase selection, individuals are selected according to their performance in a randomly chosen ordering of the considered problem instances, whereby decreasing precedence is given to each subsequent instance.
For a more detailed description of this operator, the reader is referred to~\cite{helmuth2014solving,la2016epsilon}.
 
\subsection{Mutation and Recombination}
\label{sec:gggp-mutation-and-recombination}
\begin{figure}
	\centering
	\includegraphics[scale=0.46]{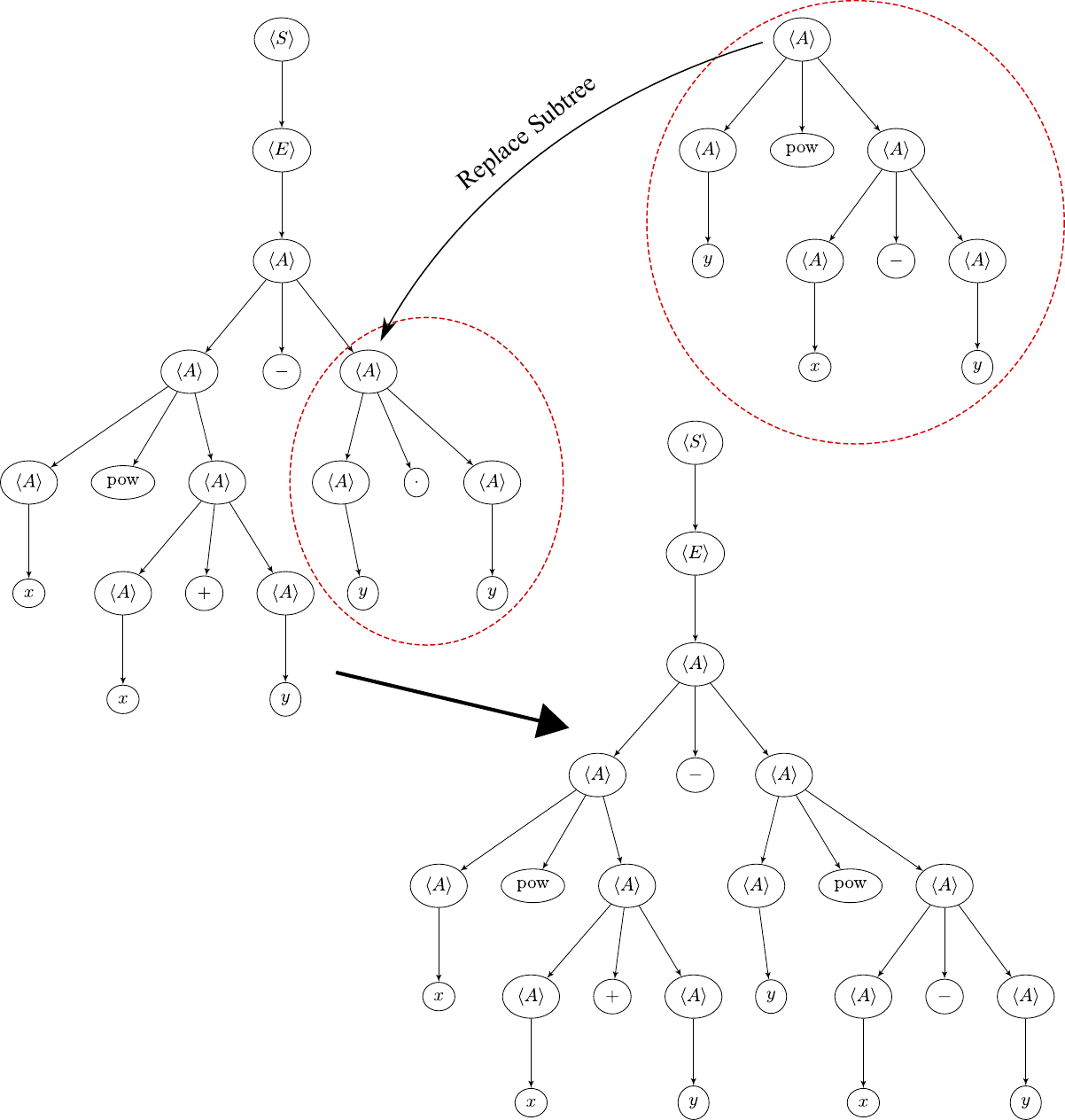}
	\caption{Subtree replacement applied to the derivation tree of the function $f_2$~\eqref{eq:gp-example-functions}.}
	\label{fig:gp-replacement-mutation}
\end{figure}

After selecting a number of candidate individuals from the population, mutation and recombination represent two complementary options to extend the current population toward more promising regions of the search space.
\paragraph{Mutation}
The term mutation is usually referred to as the process of altering the genotype of a given individual, usually in a randomized way, with the purpose of creating an individual that, hopefully, has higher fitness than its predecessor.
Therefore, mutation is usually applied to one individual at a time.
Since the possibilities of altering an individual depend on its genotype representation, mutation operators need to be customized to the type of evolutionary algorithm employed.
For this purpose, different mutation operators have been proposed for tree-based genotype representations~\cite{poli2008field,koza1994genetic}.
The most commonly used tree-mutation operator is \emph{subtree replacement}, which replaces a certain branch with a randomly-generated subtree.
This subtree can, for instance, be created using a similar tree-generation operator as within the population initialization~\cite{poli2008field}.
While subtree replacement was originally defined for classical GP, it can be adapted to G3P by introducing the additional constraint that branches with a specific variable as their root node can only be replaced by subtrees with the same root node.
The resulting mutation operator is illustrated in Figure~\ref{fig:gp-replacement-mutation}, where it is applied to the derivation tree of the function $f_2$ formulated in Equation~\eqref{eq:gp-example-functions}.

\begin{figure}
	\centering
	\includegraphics[scale=0.46]{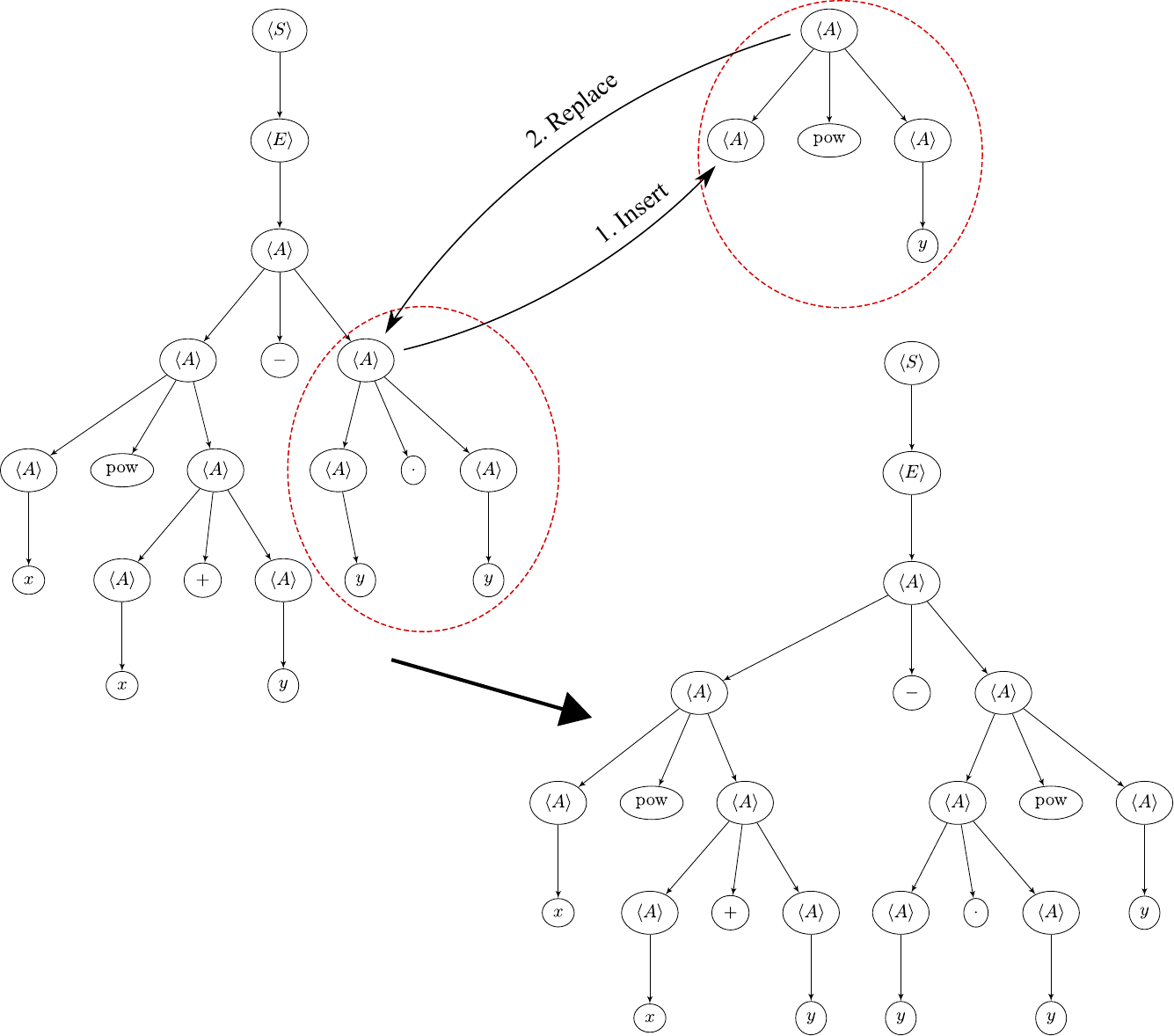}
	\caption{Subtree insertion applied to the derivation tree of the function $f_2$~\eqref{eq:gp-example-functions}.}
	\label{fig:gp-insertion-mutation}
\end{figure}

A second mutation operator, which can be directly derived from subtree replacement, is \emph{insertion}.
In subtree mutation, the original branch is replaced and consequently removed entirely from the original tree, whereas insertion only results in the creation of a partial subtree to which the original subtree is attached.
This operation is illustrated in Figure~\ref{fig:gp-insertion-mutation}.
In contrast to subtree replacement, whose applicability in the context of G3P is independent of the grammar, insertion requires the generation of a subtree that includes the root node of the original branch as a child node.
While both subtree replacement and insertion have the potential to induce drastic changes in the structure of the original tree, a less intrusive way to alter a given tree is to replace only a single node with one of the same arity, which is called \emph{point mutation}.
In the context of G3P, this operation corresponds to replacing a certain production with a different one that is available for the same variable.
However, to be valid, the newly inserted production needs to generate a sequence of child nodes in which the same variables occur as in the original sequence.
As a consequence, point mutation can only be applied to the subset of productions that fulfills this condition, which limits its applicability within G3P.

\begin{figure}
    \centering
	\includegraphics[scale=0.46]{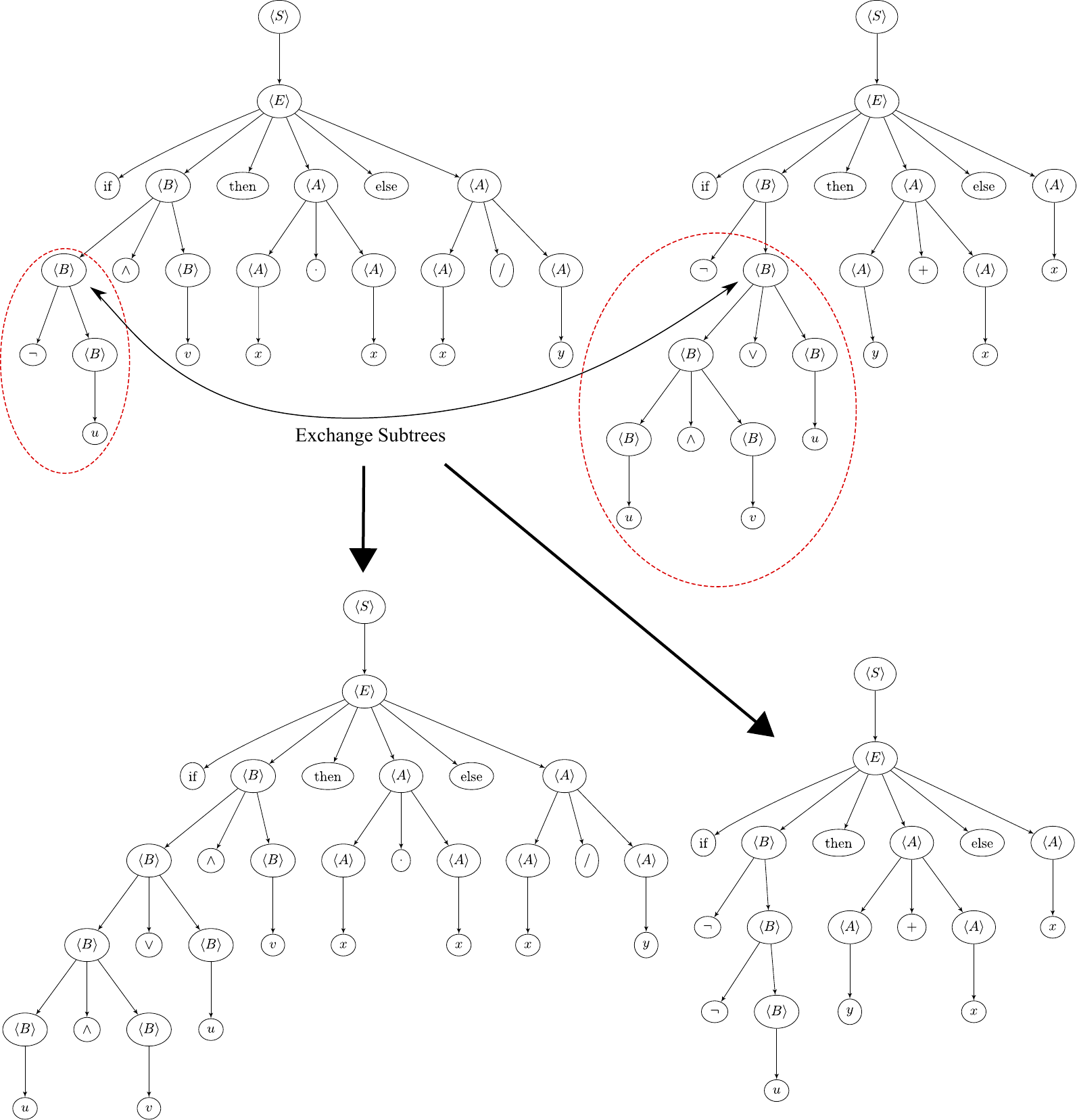}
	\caption{Subtree crossover between two derivation trees, where the first one corresponds to the function $f_1$~\eqref{eq:gp-example-functions}.}
	\label{fig:gp-subtree-crossover}
\end{figure}

\paragraph{Recombination}
While the goal of mutation is to introduce new genetic information into the population, recombination, often also called crossover, aims to combine the genotypes of existing individuals in novel ways to construct individuals with improved fitness compared to their predecessors.
In general, recombination operators can be classified into \emph{homologous} and \emph{non-homologous} ones, whereby those from the former category preserve the position of the individual elements of the genome.
A straightforward way to perform recombination in tree-based GP is to exchange subtrees between two individuals at a certain crossover point, which is called \emph{subtree crossover}.
Again, to apply this operator within G3P, the crossover point must be chosen in a way that the root nodes of the exchanged branches correspond to the same variable.
The resulting operation is illustrated in Figure~\ref{fig:gp-subtree-crossover}, where it is applied to the derivation tree of the example function $f_1$ formulated in Equation~\eqref{eq:gp-example-functions}.  
As subtree crossover permits the exchange of arbitrary branches with matching root nodes between two individuals, this operator is non-homologous.
To obtain a similarly functioning homologous recombination operator, the \emph{one-point crossover} operator has been proposed~\cite{poli1998schema}.
One-point crossover first aligns both trees by traversing them recursively until there occurs a mismatch between the arity of two non-terminal nodes of both trees at the same position.
The crossover point is then chosen from all positions in both trees where such a mismatch occurs while the subtree exchange is performed in a similar way as in subtree crossover.
While this operator preserves the relative position of each node within both trees, only a limited number of crossover points can be selected.
For instance, the subtree exchange shown in Figure~\ref{fig:gp-subtree-crossover} is not a valid one-point crossover operation, as the first arity mismatch already occurs above the chosen crossover point within the respective branches of both trees.
In general, the applicability of one-point crossover within G3P depends on the number of productions available for each variable.
In particular, if most productions generate only a single variable, the average derivation tree possesses a small branching factor, and thus, there is only a limited number of possible crossover points.
Finally, it must be mentioned that, besides the mutation and crossover operators discussed here, a variety of alternative ones have been proposed since the invention of GP, of which an overview can be found in~\cite{poli2008field}.
However, these operators are designed for general tree-based GP systems and hence do not take the underlying grammatical representation into account. 
As a remedy, in~\cite{couchet2007crossover}, mutation and crossover operators that are tailored toward grammar-based representations have been proposed.

After creating a number of novel individuals by means of mutation and crossover, the remaining question that needs to be addressed is how the search should be progressed by forming a new population based on the previous one and the individuals created through mutation and crossover, as it is shown in Line 7 of Algorithm~\ref{alg:genetic-programming}.
Here one possibility is to completely replace the previous generation with new individuals.
However, this approach incurs the danger of evicting individuals that are located in promising areas of the search space and thus might negatively affect the outcome of the search.
As a remedy, many evolutionary algorithms employ \emph{elitism}, which means that individuals from the previous population are allowed to pass over to the next one in case their fitness is not surpassed by enough newly-created ones.
While the application of elitism reduces the occurrence of fitness regressions from one generation to the next, it comprises the danger of getting stuck in local optima.
Therefore, deciding upon the amount of elitism always represents a compromise between the loss of promising individuals and the danger of premature convergence to local optima.
Options to mitigate this risk and prevent locally-optimal individuals from quickly overtaking the population are the restriction of elitism to only a subset of the population or the introduction of additional niching or diversity criteria into the selection process for the next population. 
The latter is especially common in multi-objective evolutionary algorithms~\cite{coello2007evolutionary}.

%% file: contents/multigrid_grammar.tex
After establishing a theoretical foundation on both multigrid methods and evolutionary program synthesis, we can now focus on developing a formal language for expressing arbitrarily-structured multigrid solvers in a generalized way.
As we have seen in Section~\ref{sec:multigrid-cycles}, according to the classical formulation of multigrid, each of these methods belongs to a particular family of cycles.
Each cycle possesses a distinct computational structure that stems from the number of recursive descents performed on each level of the method, as determined by the parameter $\gamma$ in Algorithm~\ref{alg:multigrid-cycle}.
For instance, a V-cycle is characterized by exactly one recursive descent per level.
Furthermore, each classical multigrid cycle employs a fixed number of smoothing steps per discretization level, which is determined by the parameters $\nu_1$ and $\nu_2$.
While the representation of a multigrid method as a recursive cycle yields a formally simple and easily parameterizable algorithmic formulation, it also enforces unnecessary restrictions on its structure.
Consider, for instance, the multigrid method shown in Figure~\ref{fig:non-traditional-multigrid-cycle}.
\begin{figure}
		\begin{tikzpicture}
            \node   (h) at (-0.75, 4){$h$};
			\node   (2h) at (-0.75, 3){$2h$};
			\node   (4h) at (-0.75, 2){$4h$};
			\node   (8h) at (-0.75, 1){$8h$};
			\node	(a) at (0,4) [draw, circle,scale=1.0] {};
			\node	(b) at (0.5,3) [draw, circle,scale=1.0] {};
			\node	(c) at (1,2) [draw, circle,scale=1.0] {};
			\node	(d) at (1.5,3) [draw, circle, scale=1.0] {};
			\node	(e) at (2,2) [draw, circle, scale=1.0] {};
			\node	(f) at (2.5,1) [draw, circle,scale=1.0,fill=black] {};
			\node	(g) at (3,2) [draw, circle,scale=1.0] {};
			\node	(h) at (3.5,3) [draw, circle,scale=1.0] {};
			\node	(i) at (4,2) [draw, circle,scale=1.0] {};
			\node	(j) at (4.5,3) [draw, circle,scale=1.0] {};
			\node	(k) at (5,2) [draw, circle,scale=1.0] {};
			\node	(l) at (5.5,1) [draw, circle,scale=1.0,fill=black] {};
			\node	(m) at (6,2) [draw, circle,scale=1.0] {};
			\node	(n) at (6.5,3) [draw, circle,scale=1.0] {};
			\node	(o) at (7,4) [draw, circle,scale=1.0] {};
			\node	(p) at (7.5,3) [draw, circle,scale=1.0] {};
			\node	(q) at (8,4) [draw, circle,scale=1.0] {};
			\node	(r) at (8.5,3) [draw, circle,scale=1.0] {};
			\node	(s) at (9,2) [draw, circle,scale=1.0] {};
			\node	(t) at (9.5,3) [draw, circle,scale=1.0] {};
			\node	(u) at (10,4) [draw, circle,scale=1.0] {};
			\draw 
			(a) edge[->] (b) 
			(b) edge[->] (c)
			(c) edge[->] (d)
			(d) edge[->] (e)   
			(e) edge[->] (f)
			(f) edge[->] (g)
			(g) edge[->] (h)
			(h) edge[->] (i)
			(i) edge[->] (j)
			(j) edge[->] (k)
			(k) edge[->] (l)
			(l) edge[->] (m)
			(m) edge[->] (n)
			(n) edge[->] (o)
			(o) edge[->] (p)
			(p) edge[->] (q)
			(q) edge[->] (r)
			(r) edge[->] (s)
			(s) edge[->] (t)
			(t) edge[->] (u)
			;
		\end{tikzpicture}
	\caption{Example for a non-classical multigrid method.}
	\label{fig:non-traditional-multigrid-cycle}
\end{figure}
While this method reaches the coarsest level twice and hence, at first sight, looks similar to a W-cycle, it employs a unique pattern of computations that is completely different from any known multigrid cycle.
As a consequence, this method can not be represented within the classical framework of multigrid cycles, as formulated in Algorithm~\ref{alg:multigrid-cycle}.
To overcome the limitations of this formulation by constructing multigrid methods with a unique sequence of computations on each discretization level, a new formal language is needed.
The first step toward the development of such a language is to find a way to represent the state within each step of a multigrid method.
Based on this representation, we can then define transition rules between the individual states, which will allow us to design multigrid methods of novel structure.
\section{Multigrid States}
\label{sec:multigrid-states}
In order to derive a representation of the state of a multigrid method, we need to reconsider our original formulation of a multigrid cycle in Algorithm~\ref{alg:multigrid-cycle}.
For the sake of simplicity, in the following, symbols that correspond to vectors are written in regular font.
In case the mathematical interpretation of a certain lowercase symbol is ambiguous, its meaning will be explicitly stated. 
With the exception of the coarsest level, where the only allowed operation is the application of the coarse-grid solver, we can identify three elementary multigrid operations that can be performed within a cycle.
\begin{definition}[Elementary Multigrid Operations]
\label{def:elementary-multigrid-operations}
On each level except for the lowest, the following three operations can be performed within a multigrid cycle.
\begin{itemize}
	\item \textbf{Smoothing}: Reduce the oscillatory error components of the approximate solution $x_h$ on the current level. 
	\begin{equation}
		x_h = x_h + \omega M_h^{-1} \left( b_h - A_h x_h \right) \; \text{where} \; A_h = M_h + N_h
	\end{equation}
	\item \textbf{Coarsening}: A coarse problem is obtained by restricting the residual.
	\begin{equation}
    \begin{split}
     	x_{2h} & = 0 \\
		b_{2h} & = I_h^{2h} (b_h - A_h x_h)
    \end{split}
	\end{equation}
	\item \textbf{Coarse-Grid Correction}: Prolongate a correction $x_{2h}$ obtained on a coarser grid to reduce the low-frequency error components of the approximate solution $x_h$.
	\begin{equation}
		x_h = x_h + I_{2h}^h x_{2h}
	\end{equation}
\end{itemize}
\end{definition}
Except for the operators $A_h$, $I_h^{2h}$ and $I_{2h}^h$ the result of each of these operations exclusively depends on the current value of the approximate solution on subsequent levels, i.e., $x_{h}$ and $x_{2h}$, and the right-hand side $b_h$.
However, in contrast to the coarse-grid correction step, which makes use of the current approximate solution on the coarse grid, both smoothing as well as coarsening requires us first to compute the residual $r_h = b_h - A_h x_h$, which can be considered as an intermediate step.
While this differentiation is not strictly necessary, it leads to simpler expressions since both operations can then be split into two steps.
Furthermore, as in practice, the error can usually not be computed directly, the residual is often the only available metric to investigate whether a multigrid iteration has achieved a certain amount of error reduction and hence has to be computed repeatedly.
%The residual is then either restricted and assigned to the right-hand side $b_{2h}$ or employed to reduce the oscillatory error components of the approximate solution by computing a coarse-grid correction term of the form $\omega A_h^{-1} r_h$.
To derive a general representation for the state of a multigrid method based on our previous observations, we next consider the sequence of operations shown in Algorithm~\ref{alg:example-three-grid-method}, which corresponds to a three-grid V-cycle that performs one step of underrelaxed Jacobi smoothing on the second finest level.
\begin{algorithm}[t]
	\begin{algorithmic}[1]
		\State $x_{h} = x_{h}^{0}$
		\State $r_{h} = b_{h} - A_h x_{h} $
		\State $ x_{2h} = 0$
		\State $ b_{2h} = I_{h}^{2h} r_{h}$
		\State $ r_{2h} = b_{2h} - A_{2h} x_{2h}$
		\State $ x_{2h} = x_{2h} + I_{4h}^{2h} A_{4h}^{-1} I_{2h}^{4h} r_{2h}$
		\State $ r_{2h} = b_{2h} - A_{2h} x_{2h}$
		\State $ x_{2h} = x_{2h} + 0.6 \cdot D_{2h}^{-1} r_{2h}$
		\State $x_{h} = x_{h}  + I_{2h}^h x_{2h}$
	\end{algorithmic}
\caption{Example of a Three-Grid V-Cycle}
\label{alg:example-three-grid-method}
\end{algorithm}
In each step of this sequence, either the approximate solution, right-hand side, or residual is updated on a particular level.
While in practice, each of these variables corresponds to a data structure containing numerical values, our goal is to represent the algorithmic structure of a multigrid method in the form of symbolic expressions.
In this case, the value of each variable is determined by the expression that computes its value.
Updating a variable, therefore, corresponds to assigning a new expression to the corresponding symbol.
For this purpose, we consider the tuple
\begin{equation*}
	Z_h = (x_h, b_h, r_h)
\end{equation*} 
on each level with step size $h$, which then refers to the expression of each of the corresponding three variables.
Starting from the first line, we can progressively update the contents of this tuple with the expression given in each line, whereby each occurrence of one of the three symbols is replaced by the expression currently contained in the respective entry of the tuple.
Figure~\ref{fig:example-tree-grid-method-states} shows the content of the state tuple for each line of Algorithm~\ref{alg:example-three-grid-method}.
\begin{figure}
	\begin{equation*}
		\begin{array}{l l}
			\hline
			\bm{Z_h} & \bm{(x_h, \, b_h, \, r_h)}  \\
			1: &  x_{h}^0, \, b_h, \, \lambda \\
			2: &  x_{h}^0, \, b_h, \, b_{h} - A_h x_{h}^0 \\ \hline
			\bm{Z_{2h}} &  \bm{(x_{2h}, \, b_{2h}, \, r_{2h})} \\
			3: &  0, \, \lambda, \, \lambda \\
			4: &  0, \, I_{h}^{2h}\underbrace{(b_{h} - A_h x_{h}^0)}_{r_{h}}, \, \lambda \\
			5: &  0, \, I_{h}^{2h}(b_{h} - A_h x_{h}^0), \,\underbrace{I_{h}^{2h}(b_{h} - A_h x_{h}^0)}_{b_{2h}} - A_{2h} 0 \\
			6: & 0 + I_{4h}^{2h} A_{4h}^{-1} I_{2h}^{4h} (\underbrace{I_{h}^{2h}(b_{h} - A_h x_{h}^0) - A_{2h} 0}_{r_{2h}}), \, I_{h}^{2h}(b_{h} - A_h x_{h}^0), \, \lambda\\
			7: & 0 + I_{4h}^{2h} A_{4h}^{-1} I_{2h}^{4h} (I_{h}^{2h}(b_{h} - A_h x_{h}^0) - A_{2h} 0), \, I_{h}^{2h}(b_{h} - A_h x_{h}^0), \\ 
			&  \underbrace{I_{h}^{2h}(b_{h} - A_h x_{h}^0)}_{b_{2h}} - A_{2h} (\underbrace{0 + I_{4h}^{2h} A_{4h}^{-1} I_{2h}^{4h} (I_{h}^{2h}(b_{h} - A_h x_{h}^0) - A_{2h} 0)}_{x_{2h}}) \\
			8: &   (\underbrace{0 + I_{4h}^{2h} A_{4h}^{-1} I_{2h}^{4h} (I_{h}^{2h}(b_{h} - A_h x_{h}^0) - A_{2h} 0)}_{x_{2h}}) + 0.6 \cdot D_{2h}^{-1} \cdot \\ 
			& (\underbrace{I_{h}^{2h}(b_{h} - A_h x_{h}^0) - A_{2h} (0 + I_{4h}^{2h} A_{4h}^{-1} I_{2h}^{4h} (I_{h}^{2h}(b_{h} - A_h x_{h}^0) - A_{2h} 0))}_{r_{2h}}), \\ 
			& I_{h}^{2h}(b_{h} - A_h x_{h}^0), \, \lambda \\ \hline 
			\bm{Z_h} & \bm{(x_h, \, b_h, \, r_h)}  \\
			9: & x_{h}^0 + I_{2h}^h ((0 + I_{4h}^{2h} A_{4h}^{-1} I_{2h}^{4h} (I_{h}^{2h}(b_{h} - A_h x_{h}^0) - A_{2h} 0)) + 0.6 \cdot D_{2h}^{-1} \cdot \\ 
			& (I_{h}^{2h}(b_{h} - A_h x_{h}^0) - A_{2h} (0 + I_{4h}^{2h} A_{4h}^{-1} I_{2h}^{4h} (I_{h}^{2h}(b_{h} - A_h x_{h}^0) - A_{2h} 0)))), \\ 
			&  b_h, \, \lambda \\
			\hline
		\end{array}
	\end{equation*}
	\caption{State tuple in each step of Algorithm~\ref{alg:example-three-grid-method}.}
	\label{fig:example-tree-grid-method-states}
\end{figure}
As introduced in Section~\ref{sec:formal-languages}, the empty symbol $\lambda$ denotes that a certain component of the tuple is unspecified.
After the last step of the sequence (Line~9), the first component $x_h$ of the tuple combines all computational steps of the method in a single expression.

While Figure~\ref{fig:example-tree-grid-method-states} illustrates that the ternary tuple $Z_h$ contains all relevant information of a multigrid method's current state on a certain level, we have not yet discussed how we can transition between the states of two subsequent levels.
First, we consider the case of transitioning from a level with step size $h$ to the next coarser grid with step size $2h$, which is shown in Line~3~and~4 of Figure~\ref{fig:example-tree-grid-method-states}.
Here, the tuple
\begin{equation*}
	Z_{2h} = (0, \, I_{h}^{2h}(b_{h} - A_h x_{h}^0), \, \lambda)
\end{equation*} 
corresponds to the coarse-grid error equation 
\begin{equation*}
	A_{2h} x_{2h} = I_{h}^{2h}(b_{h} - A_h x_{h}^0),
\end{equation*}
whose solution is to be approximated, starting with an initial guess of zero.
Therefore, all necessary information for the creation of this state is obtained from the next higher level in the form of the restricted residual 
\begin{equation*}
    I_h^{2h} r_h = I_h^{2h} (b_{h} - A_h x_{h}^0).
\end{equation*}
On the other hand, consider the transition from a coarse grid back to the next finer grid in the form of the coarse-grid correction in Line~9.
In this case, we require both the approximate solution $x_{2h}$, computed on the coarse grid, as well as the previous values of the first two entries of the fine-grid state tuple, i.e., $x_h$ and $b_h$.
While for coarsening, we can neglect all previous states on the coarser grid, as their information is explicitly contained in the residual, for a coarse-grid correction, the previous state on the fine grid needs to be restored.
Note that for a multigrid method that operates on a hierarchy of discretizations of even larger depth than the example shown in Algorithm~\ref{alg:example-three-grid-method}, this process must be carried out recursively for each bottom-up transition.
To resolve this issue, we need to extend our original formulation of a multigrid state with a fourth component that has the purpose of preserving the current state on the next-higher level.
While at the topmost level, this component is always empty, whenever coarsening is performed, the state of the current level is included.
For instance, in Line~4 of Figure~\ref{fig:example-tree-grid-method-states} we have to extend the given ternary tuple 
\begin{equation*}
Z_{2h} = (0, \, I_{h}^{2h}(b_{h} - A_h x_{h}^0), \, \lambda)
\end{equation*}
by including the state 
\begin{equation*}
Z_h = (x_{h}^0, \, b_h, \, \lambda, \, \lambda) 
\end{equation*} 
as an additional fourth entry. 
As a consequence, all required information for restoring the previous fine-grid state values is included in the quaternary tuple 
\begin{equation*}
	Z_{2h} = (0, \, I_{h}^{2h}(b_{h} - A_h x_{h}^0), \, \lambda, \, Z_h).
\end{equation*}
Note that since $Z_h$ represents the current state on the finest grid, its fourth component is empty, while otherwise, it would refer to the previous state of the next higher level in the discretization hierarchy.
To assess the feasibility of this approach, we have to check whether all components required to construct the coarse-grid correction expression, as in Line~9 of Figure~\ref{fig:example-tree-grid-method-states}, are available within the current state.
Since in addition to $x_{2h}$, both the approximate solution $x_h$ and right-hand side $b_h$ are now contained in the fourth component of the coarse-grid state tuple, this expression can be assembled in a straightforward manner.
At this point, note that a coarse-grid correction on the lowest level represents a special case. 
Since the only allowed operation on this level is the application of the coarse-grid solver, which is denoted by multiplication with the inverse of the system matrix, it is represented as a single operation given by the expression 
\begin{equation*}
	x_{2h} = x_{2h} + I_{4h}^{2h} A_{4h}^{-1} I_{2h}^{4h} r_{2h}.
\end{equation*}
After resolving the issue of restoring previous states during a coarse-grid correction, we can thus provide a complete formal definition of the state of a multigrid method.
\begin{definition}[Multigrid State]
\label{def:multigrid-state}
The state of a multigrid method for solving the equation $A_{h} x_{h} = b_{h}$ on a grid with step size $h$ is given by the quaternary tuple
\begin{equation}
	Z_{h} = \left( x_{h}, b_{h}, c_{h}, Z_{h/2}\right), 
\end{equation}
where
\begin{itemize}
	\item $x_{h}$ is an expression for computing an approximate solution of the above equation,
	\item $b_{h}$ is an expression for computing the right-hand side of the above equation,
	\item $c_{h}$ is a correction term for improving the accuracy of the approximate solution,
	\item $Z_{h/2}$ is either the current state on the next finer grid with a step size $h/2$ or the empty symbol $\lambda$ in case $h$ represents the topmost level in the given hierarchy of discretizations.
\end{itemize}
\end{definition}
Note that in Definition~\ref{def:multigrid-state}, the third component $c_{h}$ of a multigrid state no longer has to refer to the residual but now represents an expression that corresponds to a general correction term.
While in the classical multigrid formulation, as presented in Section~\ref{sec:multigrid-methods}, all smoothing expressions are directly derived from the residual, this is not necessarily the case for all smoothers, such as distributive smoothing~\cite{trottenberg2000multigrid}.
Furthermore, the general notion of a correction term enables us to represent intermediate expressions that do not specifically refer to the current residual within a multigrid state.
For instance,
\begin{equation*}
    c_h = \omega M_h^{-1} \left( b_h - A_h x_h \right)
\end{equation*}
contains the complete expression for correcting the current approximate solution within smoothing.
%Finally, note that the computation of the coarse-grid correction does not require the value of the previous fine-grid residual.
%Therefore, the corresponding expression does not need to be restored and the residual can be omitted in the respective state tuple when performing the fine-to-coarse grid transition.
\section{State Transition Functions}
\label{sec:multigrid-state-transitions}
After developing a representation for the state of a multigrid method, the next step towards a formal language for representing arbitrary sequences of multigrid operations,
such as the one shown in Figure~\ref{fig:non-traditional-multigrid-cycle}, is the derivation of a list of rules that describes all possible transitions between states.
For this purpose, we first need to investigate which operations are allowed within a multigrid method.
From a mathematical point of view, all operations performed within a multigrid method exclusively consist of either matrix-vector multiplications or vector additions and subtractions.
First of all, the computation of the residual $r_h = b_h - A_h x_h$ represents a fundamental operation, based on which an approximation of the solution of the target system $A_h x_h = b_h$ is iteratively improved.
The latter is either achieved through smoothing or by computing a coarse-grid correction.
The \textsc{residual} function shown in Algorithm~\ref{alg:state-transition-residual} implements the corresponding state transition for computing the residual on a particular level within the discretization hierarchy with step size $h$\footnote{While so far, the letter $h$ did always refer to the spacing of the finest grid, in this case, it represents the grid spacing on an arbitrary level.}.
The residual expression is assembled based on the system matrix $A_h$ and the given state $Z_h = (x_h, b_h, \lambda, Z_{h/2})$, which contains the current approximate solution $x_h$ and right-hand side $b_h$.
The resulting expression is then included in $Z_h$ as a correction term $c_h$, which is returned at the end of the function.
\begin{algorithm}[t]
	\begin{algorithmic}
		\Function{residual}{$A_h$, $(x_h, b_h, \lambda, Z_{h/2})$}
		\State $c_h \gets b_h - A_h x_h$
		\State return $(x_h, b_h, c_h, Z_{h/2})$
		\EndFunction
	\end{algorithmic}
\caption{Residual Computation}
\label{alg:state-transition-residual}
\end{algorithm}
After constructing an initial correction $c_h$ from the residual expression $b_h - A_h x_h$, we can next apply an operator $B_h$ to this term, either as part of a smoothing expression or in the form of a restriction that yields the right-hand side $b_h$ of the coarse-grid error equation.
Since both operations can be formulated as a matrix-vector multiplication, we consider this operator application as an intermediate step that is implemented in the form of the \textsc{apply} function shown in Algorithm~\ref{alg:state-transition-apply}.
This function applies the operator $B_h$ to the current correction term. 
The resulting expression then serves as a new correction term within the returned state.
\begin{algorithm}[h]
	\begin{algorithmic}
		\Function{apply}{$B_h$, $(x_h, b_h, c_h, Z_{h/2})$}
		\State return $(x_h, b_h, B_h\cdot c_h, Z_{h/2})$
		\EndFunction
	\end{algorithmic}
\caption{Operator Application}
\label{alg:state-transition-apply}
\end{algorithm}
In general, the choice of the operator $B_h$ leads to the two elementary operations available within a multigrid method, i.e., smoothing and solving the given problem on a coarser grid.
For instance, $B_h = D_h^{-1}$ leads to the expression
\begin{equation*}
	c_h = D_h^{-1} (b_h - A_h x_h),
\end{equation*}
which corresponds to one step of the Jacobi method, as formulated in Equation~\eqref{eq:jacobi-method}.
In contrast, choosing $B_h$ as the restriction operator $I_h^{2h}$ leads to
\begin{equation*}
	c_{2h} = I_{h}^{2h} (b_h - A_h x_h),
\end{equation*}
which then serves as the right-hand side $b_{2h}$ of the corresponding coarse-grid error equation.
In each of the two cases, the correction term obtained from the \textsc{apply} function serves a different purpose.
In the first case, we must be able to generate an expression that computes an improved approximate solution $x_h$ by applying an update in the form of the correction term $c_h$ to the previous value of $x_h$. 
This behavior is implemented in the function \textsc{update}\footnote{Note that in previous publications~\cite{schmitt2020constructing,schmitt2021evostencils} we have used the name \textsc{iterate} for this function. However, as the name \textsc{update} reflects its meaning in a more concrete way, we have chosen to rename it.} shown in Algorithm~\ref{alg:state-transition-update}, which returns a state tuple that contains the expression for computing an updated approximate solution as its first entry.
\begin{algorithm}[t]
	\begin{algorithmic}
		\Function{update}{$\omega$, $P$, ($x_h$, $b_h$, $c_h$, $Z_{h/2}$)}
			\State $x_h \gets x_h + \omega \cdot c_h$ with $P$
			\State return ($x_h$, $b_h$, $\lambda$, $Z_{h/2}$) 
		\EndFunction
	\end{algorithmic}
 \caption{Approximate Solution Update}
\label{alg:state-transition-update}
%\caption{Iteration}
\end{algorithm}
In addition to the actual correction term, this function also includes a relaxation factor $\omega$ and a partitioning operator $P$, which enable the formulation of an underrelaxed, overrelaxed, and colored version of each operation, as described in Section~\ref{sec:smoothing}.
The second possibility of applying an operator to the current residual is in the form of the restriction expression $c_{2h} = I_{h}^{2h} r_h$.
As denoted by the subscript, the result of this operation is an expression on a coarser grid of step size $2h$.
To construct the coarse-grid error equation 
\begin{equation*}
A_{2h} x_{2h} = I_{h}^{2h} r_h,
\end{equation*} 
we have to assign this expression to the right-hand side $b_{2h}$ that is contained in the corresponding state tuple.
Furthermore, as discussed in the last section, before we can transition to a state on the next coarser grid, we have to store the previous state in the fourth entry of the resulting tuple.
The resulting state transition is implemented in the function \textsc{coarsening} shown in Algorithm~\ref{alg:state-transition-cycle}, whose name indicates that the computation proceeds on a lower level and thus a coarser grid.
\begin{algorithm}[h]
	\begin{algorithmic}
	\Function{coarsening}{$A_{2h}$, $x_{2h}^0$, ($x_h$, $b_{h}$, $c_{2h}$, $Z_{h/2}$)}
		\State $x_{2h} \gets x_{2h}^0$ 
		\State $b_{2h} \gets c_{2h}$
		\State $c_{2h} \gets b_{2h} - A_{2h} x_{2h}$ 
		\State $Z_h \gets$ ($x_{h}$, $b_{h}$, $\lambda$, $Z_{h/2}$)
		\State return ($x_{2h}$, $b_{2h}$, $c_{2h}$, $Z_h$)
	\EndFunction
	\end{algorithmic}
 \caption{Coarsening}
\label{alg:state-transition-cycle}
\end{algorithm}
Within this function, we already construct an expression for the initial residual
\begin{equation*}
	r_{2h} = b_{2h} - A_{2h} x_{2h}^0,
\end{equation*} 
and include it as a correction term in the newly created state tuple, where we assume that the initial approximate solution $x^0_{2h}$ has already been set to zero.
The reasoning behind this is that performing a coarse-grid correction with $x^0_{2h}$ is not expected to lead to any improvement, and thus computing the residual always represents the first computational step on the coarse grid.
Finally, after deriving an expression $x_{2h}$ for computing the approximate solution of the error equation on the coarse grid, we have to transfer this term back to the fine grid, where it can then be applied as a correction term to improve the approximate solution on this level.
We, therefore, first need to restore the previous state on the next finer level, including the current approximate solution, right-hand side, and the preceding state.
To then obtain the coarse-grid correction expression, we need to apply the prolongation operator to the current approximate solution on the coarse grid.
The resulting state transition is implemented in the function \textsc{cgc} shown in Algorithm~\ref{alg:state-transition-cgc}.
Note that the complete coarse-grid correction step is obtained through subsequent application of the \textsc{update} function, which updates the approximate solution with the correction term returned by the \textsc{cgc} function.
\begin{algorithm}[t]
	\begin{algorithmic}
		\Function{cgc}{$I_{2h}^{h}$, $(x_{2h}, b_{2h}, \lambda, Z_{h})$}
		\State ($x_h$, $b_{h}$, $\lambda$, $Z_{h/2}$) $\gets Z_{h}$
		\State $c_h \gets I_{2h}^{h} \cdot x_{2h}$
		\State return ($x_h$, $b_{h}$, $c_h$, $Z_{h/2}$)
		\EndFunction
	\end{algorithmic}
\caption{Coarse-Grid Correction}
 \label{alg:state-transition-cgc}
\end{algorithm}
%TODO describe coarse-grid solver A_{4h}^{-1}
As it has been described in this section, the application of each of these functions corresponds to a particular transition between multigrid states, and hence, for any multigrid method, a sequence of function applications can be derived.
Furthermore, note that the evaluation of this sequence leads to a state whose first component $x_h$ contains an expression that corresponds to the stepwise execution of the method.
The functions described in this section, therefore, can be considered as a \emph{language} for the formal representation of multigrid methods. 
The main components of this language are summarized in Algorithm~\ref{table:grammar-semantics}.
\begin{algorithm}[ht!]
	\caption{Multigrid State Transition Functions}
	\label{table:grammar-semantics}
	\begin{algorithmic}
	\Function{residual}{$A_h$, $(x_h, b_h, \lambda, Z_{h/2})$}
	\State $c_h \gets b_h - A_h x_h$
	\State return $(x_h, b_h, c_h, Z_{h/2})$
	\EndFunction
	\State
	\Function{apply}{$B_h$, $(x_h, b_h, c_h, Z_{h/2})$}
	\State return $(x_h, b_h, B_h\cdot c_h, Z_{h/2})$
	\EndFunction
	\State
	\Function{update}{$\omega$, $P$, ($x_h$, $b_h$, $c_h$, $Z_{h/2}$)}
	\State $x_h \gets x_h + \omega \cdot c_h$ with $P$
	\State return ($x_h$, $b_h$, $\lambda$, $Z_{h/2}$) 
	\EndFunction
	\State
	\Function{coarsening}{$A_{2h}$, $x_{2h}^0$, ($x_h$, $b_{h}$, $c_{2h}$, $Z_{h/2}$)}
	\State $x_{2h} \gets x_{2h}^0$ 
	\State $b_{2h} \gets c_{2h}$
	\State $c_{2h} \gets b_{2h} - A_{2h} x_{2h}$ 
	\State $Z_h \gets$ ($x_{h}$, $b_{h}$, $\lambda$, $Z_{h/2}$)
	\State return ($x_{2h}$, $b_{2h}$, $c_{2h}$, $Z_h$)
	\EndFunction
	\State
	\Function{cgc}{$I_{2h}^{h}$, $(x_{2h}, b_{2h}, \lambda, Z_{h})$}
	\State ($x_h$, $f_{h}$, $c_h$, $Z_{h/2}$) $\gets Z_{h}$
	\State $c_h \gets I_{2h}^{h} \cdot x_{2h}$
	\State return ($x_h$, $f_{h}$, $c_h$, $Z_{h/2}$)
	\EndFunction
	\end{algorithmic}
\end{algorithm}

Now let us revisit one of the main goals of this language: The description of multigrid methods that execute a unique sequence of operations on each level, which can not be formulated in the classical framework of multigrid cycles.
To achieve this goal, we have ensured that each possible state transition within a multigrid method is described by the application of a specific sequence of well-defined functions. 
At the same time, we have carefully avoided the inclusion of transitions spanning over multiple states, as this would reduce the expressiveness of our language, and thus similar to Algorithm~\ref{alg:multigrid-cycle}, restrict it to only a subset of all possible multigrid methods.
The remaining step within the development of a formal system for the construction of arbitrary sequences of multigrid operations is thus the derivation of a formal grammar that generates the corresponding strings contained in our multigrid language, as it has been described in Section~\ref{sec:formal-languages}.

\section{A Novel Family of Multigrid Grammars}
\label{sec:multigrid-grammar}
In Section~\ref{sec:formal-languages}, we have already introduced the general grammar $G$ as the quaternary tuple 
\begin{equation*}
	G = \left(V, T, S, P \right),
\end{equation*}
where $V$ is the set of variables, $T$ the set of terminal symbols, $S \in V$ the start variable and $P$ the set of productions.
However, before we start defining the individual components of a grammar, we need to decide which constraints we want to apply within its productions.
As it has been discussed in Section~\ref{sec:chomsky-hierarchy}, the Chomsky hierarchy defines four grammatical levels, where each subsequent level introduces additional restrictions.
While an unrestricted or type-0 grammar represents the most general model of computation, it also leads to the greatest complexity.
In particular, enumerating all strings generated by an unrestricted grammar requires a Turing machine, which represents the most general model of computation available.
The main goal of this work is to enable the automated design of multigrid methods using evolutionary program synthesis techniques.
As we have discussed in Section~\ref{sec:gggp}, the application of tree-based grammar-guided genetic programming (G3P) requires each derivation of a grammar to be representable as a tree, which is then called a derivation tree.
The family of grammars that fulfills this property corresponds to the type-2 category of the Chomsky hierarchy.
Type-2 grammars are characterized by the constraint that only a single variable is allowed to be placed on the left-hand side of each production, which means that the applicability of each production is independent of the context in which the variable is embedded. 
These grammar are thus said to be context-free.
Since each non-terminal node of a derivation tree corresponds to a variable of the corresponding context-free grammar (CFG), all operations that are required within a G3P method can be defined in a straightforward manner, which has already been shown in Section~\ref{sec:gggp-initialization}~and~\ref{sec:gggp-mutation-and-recombination}.
In the following, we will derive a family of CFGs for the construction of multigrid methods on discretization hierarchies of variable depth that is based on the state representation and transition rules defined in Section~\ref{sec:multigrid-states}~and~\ref{sec:multigrid-state-transitions}.
Therefore, starting from the initial state 
\begin{equation}
	Z^0_h = (x_h^0, b_h, \lambda, \lambda),
 \label{eq:initial-mg-state}
\end{equation}
we need to define a production for each possible operation within a multigrid method.
Each of these productions needs to be expressed as a sequence of transitions that yields the correct state after the respective operation.
In Section~\ref{sec:multigrid-state-transitions}, we have made the distinction between states where a new approximate solution $x_h$ is computed and those with the purpose of constructing a correction term $c_h$ based on the residual expression $r_h = b_h - A_h x_h$.
To define the respective productions on every level of a given discretization hierarchy, we need to recursively determine how the final state of a multigrid method can be reached.
Note that in this state, all computational steps are combined in $x_h$ as a single expression.
We represent the final state of a multigrid method as the variable $\ps{s_h}$ and thus obtain the initial production 
\begin{production}
	\bnfprod{$S$} {
		\bnfpn{$s_h$}
	}
\end{production}Next, we have to consider the different possibilities of deriving a new expression for the approximate solution $x_{h}$ on the finest level.
According to Definition~\ref{def:elementary-multigrid-operations}, $x_{h}$ can either be updated by means of smoothing or by applying a coarse-grid correction.
In the case of smoothing, the approximate solution is updated with a correction term obtained by multiplying the residual expression with an operator.
This operation can be expressed as a combination of the state transition functions $\textsc{update}$ and $\textsc{apply}$ yielding the production
\begin{production}
	\bnfprod{$s_h$} {
		\bnfts{\textnormal{\textsc{update}}}(\bnfts{$\omega$}, \bnfsp \bnfpn{$P$}, \bnfsp \bnfts{\textnormal{\textsc{apply}}}(\bnfpn{$B_h$}, \bnfsp \bnfpn{$c_h$}))
	}.
\label{prod:smoothing}
\end{production}Here, the current state, containing the previously computed residual, is represented by the variable $\ps{c_h}$.
Different smoothers can be realized by expressing the generation of the operator and partitioning variable, i.e., $\ps{B_h}$ and $\ps{P}$, as two separate productions
\begin{production}
	\bnfprod{$B_h$} {
		\bnfts{\textnormal{\textsc{inverse}}}(\bnfts{$M_h$}) \bnfsp \bnfts{\textnormal{with}} \bnfsp \bnfts{$A_{h} = M_{h} + N_{h}$}
	}\label{prod:smoothing-operator} \\
  \refstepcounter{prodcounter}
 	\bnfprod{$P$} {
		\bnfts{\textnormal{\textsc{partitioning}}} \bnfor \bnfes
	},\label{prod:partitioning}
\end{production}where the empty symbol $\lambda$ indicates that no partitioning is used. 
On the topmost level, a correction term is always obtained by computing the residual based on the current approximate solution, which leads to the production
\begin{production}
	\bnfprod{$c_h$} {
		\bnfts{\textnormal{\textsc{residual}}}(\bnfts{$A_h$}, \bnfsp \bnfpn{$s_h$}) 
	}.
\label{prod:residual}
\end{production}Multiple smoothing steps can thus be realized by means of a repeated application of the Productions~\eqref{prod:residual}~and~\eqref{prod:smoothing}.
In contrast to smoothing, a coarse-grid correction updates the current approximate solution with a correction term that is constructed based on an approximate solution for the error equation on the next lower level in the discretization hierarchy.
We can express this operation as a combination of the \textsc{update} and \textsc{cgc} transition functions, which yields the production
\begin{production}
	\bnfprod{$s_h$} {
		\bnfts{\textnormal{\textsc{update}}}(\bnfts{$\omega$}, \bnfsp \bnfes, \bnfsp \bnfts{\textnormal{\textsc{cgc}}}(\bnfts{$I_{2h}^h$}, \bnfsp \bnfpn{$s_{2h}$}))
	}.
\label{prod:coarse-grid-correction}
\end{production}Similar to Production~\eqref{prod:smoothing}, the variable $\ps{s_{2h}}$ refers to the previous state passed as a second argument to the \textsc{cgc} function, which means that this state is already the result of a sequence of productions.
However, since a partitioned computation of the coarse-grid correction does not yield any benefits, the second argument of the \textsc{update} function is intentionally left empty.
Furthermore, note that while we have encoded the prolongation operator as a terminal symbol, it would equally be possible to specify its choice in the form of a separate production, as in the case of the smoothing operator $\ps{B_h}$.
Since the right-hand side of Production~\eqref{prod:coarse-grid-correction} contains the variable $\ps{s_{2h}}$, which corresponds to the state on a coarser grid with step size $2h$, we next have to consider the different possibilities to transition to this state.
First of all, as multigrid methods recursively apply the same operations throughout the complete hierarchy of discretizations, we can define the Productions~\eqref{prod:smoothing}, \eqref{prod:smoothing-operator}, \eqref{prod:residual} and \eqref{prod:coarse-grid-correction} on the corresponding level which yields
\begin{production}
	\bnfprod{$s_{2h}$} {
		\bnfts{\textnormal{\textsc{update}}}(\bnfts{$\omega$}, \bnfsp \bnfpn{$P$}, \bnfsp \bnfts{\textnormal{\textsc{apply}}}(\bnfpn{$B_{2h}$}, \bnfsp \bnfpn{$c_{2h}$})) \bnfor
	} \\
 \refstepcounter{prodcounter}
	\bnfmore {
		\bnfts{\textnormal{\textsc{update}}}(\bnfts{$\omega$}, \bnfsp \bnfes, \bnfsp \bnfts{\textnormal{\textsc{cgc}}}(\bnfts{$I_{4h}^{2h}$}, \bnfsp \bnfpn{$s_{4h}$}))
		%\bnfts{\textnormal{\textsc{update}}}(\bnfts{\textnormal{\textsc{apply}}}(\bnfts{$I_{4h}^{2h}$}, \bnfsp \bnfpn{$c_{4h}$}), \bnfsp \bnfts{$\omega$}, \bnfsp \bnfes)
	}\label{prod:coarse-coarse-grid-correction} \\
  \refstepcounter{prodcounter}
	\bnfprod{$c_{2h}$} {
		\bnfts{\textnormal{\textsc{residual}}}(\bnfts{$A_{2h}$}, \bnfsp \bnfpn{$s_{2h}$})
	} \\
  \refstepcounter{prodcounter}
	\bnfprod{$B_{2h}$} {
		\bnfts{\textnormal{\textsc{inverse}}}(\bnfts{$M_{2h}$}) \bnfsp \bnfts{\textnormal{with}} \bnfsp \bnfts{$A_{2h} = M_{2h} + N_{2h}$}.
	}
\end{production}Note that the right-hand side of Production~\eqref{prod:coarse-coarse-grid-correction} now contains the variable $\ps{s_{4h}}$ defined on an even coarser grid with step size $4h$, based on which we can proceed in a similar fashion to obtain a multigrid method that operates on a hierarchy of discretizations with even larger depth.
In contrast to a coarse-grid correction, which represents a transition to a state on the next higher level, we have not yet considered the opposite direction: The construction of the error equation on a coarser grid based on the state on the current level.
While the state on the finest grid is always derived from the original system of linear equations that the method aims to solve, we can apply the state transition function $\textsc{coarsening}$ in combination with $\textsc{apply}$ to obtain a new state on the next lower level that includes the initial residual as a correction term, which leads to the production
\begin{production}
\bnfprod{$c_{2h}$} {
	\bnfts{\textnormal{\textsc{coarsening}}}(\bnfts{$A_{2h}$}, \bnfsp \bnfts{$x^0_{2h}$}, \bnfsp \bnfts{\textnormal{\textsc{apply}}}(\bnfts{$I_h^{2h}$}, \bnfsp \bnfpn{$c_h$}))
}.
\label{prod:coarsening}
\end{production}With the definition of Production~\eqref{prod:coarsening}, we are now capable of generating state transitions both in a top-down as well as a bottom-up direction within the given hierarchy of discretizations.
Finally, the remaining step for the formulation of complete grammar is the definition of the initial problem on the finest grid and the application of the coarse-grid solver on the coarsest grid.
For the former, note that the initial problem corresponds to the state $Z_h^0$, as given by Equation~\eqref{eq:initial-mg-state}.
Since this state represents the starting point of each possible sequence of multigrid operations, we have to include it as an additional production of the variable $\ps{s_h}$
\begin{production}
	\bnfprod{$s_h$} {
	 (\bnfts{$x_h^0$}, \bnfts{$b_h$}, \bnfes, \bnfes)
	}.
\label{prod:initial-solution}
\end{production}If we now consider the structure of all the productions defined so far, it becomes obvious that, with the exception of the topmost level, only the productions~\eqref{prod:smoothing-operator} and~\eqref{prod:partitioning} generate an expression that consists exclusively of terminal symbols.
Furthermore, note that all other productions generate exactly a single variable of type $\ps{s_H}$ or $\ps{c_H}$, where $H$ is the step size on an arbitrary level in the given hierarchy of discretizations, and that starting from each of these variables there exists a sequence of productions that leads to an expression containing the variable $\ps{s_h}$.
Based on the latter, we can then apply the Production~\eqref{prod:initial-solution} to terminate the derivation process.
Since each sequence of productions also starts with $\ps{s_h}$, all derivations of our grammar are of the form
\begin{equation*}
	\ps{S} \Rightarrow \ps{s_h} \Rightarrow \dots \Rightarrow u \, \ps{s_h} \, v \Rightarrow \, u (x_h^0, b_h, \lambda, \lambda) \, v,
\end{equation*}
where $u  \, (x_h^0, b_h, \lambda, \lambda) \, v$ is an arbitrary word in the language generated by our multigrid grammar.
This word then represents a multigrid method whose computational structure is determined through the order of productions in the corresponding sequence of function applications.
As a final step, we have to define one or multiple productions that correspond to the application of the coarse-grid solver on the lowest level, which allows us to improve the accuracy of the approximate solution on the above level.
By utilizing a combination of the \textsc{apply} and \textsc{update} functions, we can express this operation in the form of the two production
\begin{production}
	\bnfprod{$c_{2H}$} {
		\bnfts{\textnormal{\textsc{apply}}}(\bnfts{$A^{-1}_{H}$}, \bnfsp \bnfts{\textnormal{\textsc{apply}}}(\bnfts{$I_{H}^{2H}$}, \bnfsp\bnfpn{$c_{H}$}))
  \label{prod:coarse-grid-solver}
	}  \\
  \refstepcounter{prodcounter}
	\bnfprod{$s_H$}{
		\bnfts{\textnormal{\textsc{update}}}(\bnfts{$\omega$}, \bnfsp \bnfes, \bnfsp \bnfts{\textnormal{\textsc{apply}}}(\bnfts{$I_{2H}^{H}$}, \bnfsp \bnfpn{$c_{2H}$})).
  \label{prod:coarse-grid-solver-correction}
	}
\end{production}where $H$ is the step size on the second lowest level.
Having now assembled all necessary components, we can finally specify the complete grammar for generating multigrid methods on a given hierarchy of discretizations.
As an example, we consider a five-grid hierarchy with a uniform coarsening factor of two and a step size of $h$ on the finest grid, which means that on the coarsest grid, a step size of $16h$ is obtained.
Therefore, we first define the set of variables 
\begin{equation}\tag{4.6}
\begin{split}
	V = \{ & \ps{S}, \ps{s_h}, \ps{c_h}, \ps{B_h}, \ps{s_{2h}}, \ps{c_{2h}}, \ps{B_{2h}}, \ps{s_{4h}}, \ps{c_{4h}}, \ps{B_{4h}}, \\
	& \ps{s_{8h}}, \ps{c_{8h}}, \ps{B_{8h}}, \ps{c_{16h}}, \ps{P} \},
\end{split}
\label{eq:multigrid-grammar-variables}
\end{equation}
for each of which we can then specify the list of productions available on the respective level. 
Algorithm~\ref{table:multigrid-grammar} shows the resulting productions for the complete discretization hierarchy.
\begin{algorithm}
	\caption{Productions for Generating a Five-Grid Method}
	\label{table:multigrid-grammar}
	\begin{bnf*}
		\bnfprod{$S$} {
			\bnfpn{$s_h$}
		} \\
		\bnfprod{$s_h$} {
			\bnfts{\textnormal{\textsc{update}}}(\bnfts{$\omega$}, \bnfsp \bnfpn{$P$}, \bnfsp \bnfts{\textnormal{\textsc{apply}}}(\bnfpn{$B_h$}, \bnfsp \bnfpn{$c_h$})) \bnfor
		} \\
		\bnfmore {
			\bnfts{\textnormal{\textsc{update}}}(\bnfts{$\omega$}, \bnfsp \bnfes, \bnfsp \bnfts{\textnormal{\textsc{cgc}}}(\bnfts{$I_{2h}^h$}, \bnfsp \bnfpn{$s_{2h}$})) \bnfor (\bnfts{$x_h^0$}, \bnfts{$b_h$},\bnfes, \bnfes)
		} \\
		\bnfprod{$c_h$} {
			\bnfts{\textnormal{\textsc{residual}}}(\bnfts{$A_h$}, \bnfsp \bnfpn{$s_h$}) 
		} \\
		\bnfprod{$B_h$} {
			\bnfts{\textnormal{\textsc{inverse}}}(\bnfts{$M_h$}) \bnfsp \bnfts{\textnormal{with}} \bnfsp \bnfts{$A_{h} = M_{h} + N_{h}$}
		} \\ \\
		\bnfprod{$s_{2h}$} {
			\bnfts{\textnormal{\textsc{update}}}(\bnfts{$\omega$}, \bnfsp \bnfpn{$P$}, \bnfsp \bnfts{\textnormal{\textsc{apply}}}(\bnfpn{$B_{2h}$}, \bnfsp \bnfpn{$c_{2h}$})) \bnfor
		} \\
		\bnfmore {
			\bnfts{\textnormal{\textsc{update}}}(\bnfts{$\omega$}, \bnfsp \bnfes, \bnfsp \bnfts{\textnormal{\textsc{cgc}}}(\bnfts{$I_{4h}^{2h}$}, \bnfsp \bnfpn{$s_{4h}$}))
			%\bnfts{\textnormal{\textsc{update}}}(\bnfts{\textnormal{\textsc{apply}}}(\bnfts{$I_{4h}^{2h}$}, \bnfsp \bnfpn{$c_{4h}$}), \bnfsp \bnfts{$\omega$}, \bnfsp \bnfes)
		} \\
		\bnfprod{$c_{2h}$} {
			\bnfts{\textnormal{\textsc{residual}}}(\bnfts{$A_{2h}$}, \bnfsp \bnfpn{$s_{2h}$}) \bnfor 
		} \\
        \bnfmore{
        \bnfts{\textnormal{\textsc{coarsening}}}(\bnfts{$A_{2h}$}, \bnfsp \bnfts{$x^0_{2h}$}, \bnfsp \bnfts{\textnormal{\textsc{apply}}}(\bnfts{$I_h^{2h}$}, \bnfsp \bnfpn{$c_h$}))
        } \\
		\bnfprod{$B_{2h}$} {
			\bnfts{\textnormal{\textsc{inverse}}}(\bnfts{$M_{2h}$}) \bnfsp \bnfts{\textnormal{with}} \bnfsp \bnfts{$A_{2h} = M_{2h} + N_{2h}$}
		} \\ \\
		\bnfprod{$s_{4h}$} {
			\bnfts{\textnormal{\textsc{update}}}(\bnfts{$\omega$}, \bnfsp \bnfpn{$P$}, \bnfsp \bnfts{\textnormal{\textsc{apply}}}(\bnfpn{$B_{4h}$}, \bnfsp \bnfpn{$c_{4h}$})) \bnfor
		} \\
		\bnfmore {
			\bnfts{\textnormal{\textsc{update}}}(\bnfts{$\omega$}, \bnfsp \bnfes, \bnfsp \bnfts{\textnormal{\textsc{cgc}}}(\bnfts{$I_{8h}^{4h}$}, \bnfsp \bnfpn{$s_{8h}$}))
			%\bnfts{\textnormal{\textsc{update}}}(\bnfts{\textnormal{\textsc{apply}}}(\bnfts{$I_{4h}^{2h}$}, \bnfsp \bnfpn{$c_{4h}$}), \bnfsp \bnfts{$\omega$}, \bnfsp \bnfes)
		} \\
		\bnfprod{$c_{4h}$} {
			\bnfts{\textnormal{\textsc{residual}}}(\bnfts{$A_{4h}$}, \bnfsp \bnfpn{$s_{4h}$}) \bnfor 
		} \\
          \bnfmore{
        \bnfts{\textnormal{\textsc{coarsening}}}(\bnfts{$A_{4h}$}, \bnfsp \bnfts{$x^0_{4h}$}, \bnfsp \bnfts{\textnormal{\textsc{apply}}}(\bnfts{$I_{2h}^{4h}$}, \bnfsp \bnfpn{$c_{2h}$}))
        } \\
		\bnfprod{$B_{4h}$} {
			\bnfts{\textnormal{\textsc{inverse}}}(\bnfts{$M_{4h}$}) \bnfsp \bnfts{\textnormal{with}} \bnfsp \bnfts{$A_{4h} = M_{4h} + N_{4h}$}
		} \\ \\
		\bnfprod{$s_{8h}$} {
			\bnfts{\textnormal{\textsc{update}}}(\bnfts{$\omega$}, \bnfsp \bnfpn{$P$}, \bnfsp \bnfts{\textnormal{\textsc{apply}}}(\bnfpn{$B_{8h}$}, \bnfsp \bnfpn{$c_{8h}$})) \bnfor
		} \\
		\bnfmore {
			\bnfts{\textnormal{\textsc{update}}}(\bnfts{$\omega$}, \bnfsp \bnfes, \bnfsp \bnfts{\textnormal{\textsc{apply}}}(\bnfts{$I_{16h}^{8h}$}, \bnfsp \bnfpn{$c_{16h}$}))
		} \\
		\bnfprod{$c_{8h}$} {
			\bnfts{\textnormal{\textsc{residual}}}(\bnfts{$A_{8h}$}, \bnfsp \bnfpn{$s_{8h}$}) \bnfor 
		} \\ 
          \bnfmore{
        \bnfts{\textnormal{\textsc{coarsening}}}(\bnfts{$A_{8h}$}, \bnfsp \bnfts{$x^0_{8h}$}, \bnfsp \bnfts{\textnormal{\textsc{apply}}}(\bnfts{$I_{4h}^{8h}$}, \bnfsp \bnfpn{$c_{4h}$}))
        } \\
		\bnfprod{$B_{8h}$} {
			\bnfts{\textnormal{\textsc{inverse}}}(\bnfts{$M_{8h}$}) \bnfsp \bnfts{\textnormal{with}} \bnfsp \bnfts{$A_{8h} = M_{8h} + N_{8h}$}
		} \\ \\
		\bnfprod{${c}_{16h}$} {
			\bnfts{\textnormal{\textsc{apply}}}(\bnfts{$A^{-1}_{16h}$}, \bnfsp \bnfts{\textnormal{\textsc{apply}}}(\bnfts{$I_{8h}^{16h}$}, \bnfsp \bnfpn{$c_{8h}$}))
		} \\
		\bnfprod{$P$} {
			\bnfts{\textnormal{\textsc{partitioning}}} \bnfor \bnfes
		}
	\end{bnf*}
\end{algorithm}
For the sake of simplicity, we omit the specification of all terminal symbols here, but the reader can safely assume that each symbol that occurs on the right-hand side of a production and is not contained in the set of variables is a terminal symbol.
Note that, as we have discussed above, the same three productions are repeated on each level of the method, while on the topmost level, the generation of the initial state, defined in Production~\eqref{prod:initial-solution}, and on the lowest level the application of the coarse-grid solver, defined in Production~\eqref{prod:coarse-grid-solver} and~\eqref{prod:coarse-grid-solver-correction}, must be included additionally.
While we have used a five-grid hierarchy as an example for defining the complete grammar, it is possible to specify a similar grammar on a discretization hierarchy with a different number of coarsening steps.
Furthermore, we can extend the grammar to also encompass multigrid methods that employ a fewer number of coarsening steps, i.e., two-, three- and four-grid methods.
For this purpose, we only have to include additional productions for the application of the coarse-grid solver on each subsequent level, which in the case of Algorithm~\ref{table:multigrid-grammar} leads to
\begin{bnf*}
	\bnfprod{$c_{2h}$} {
		\bnfts{\textnormal{\textsc{apply}}}(\bnfts{$A^{-1}_{2h}$}, \bnfsp \bnfts{\textnormal{\textsc{apply}}}(\bnfts{$I_{h}^{2h}$}, \bnfsp\bnfpn{$c_{h}$}))
	} \\
	\bnfprod{$s_{h}$}{
		\bnfts{\textnormal{\textsc{update}}}(\bnfts{$\omega$}, \bnfsp \bnfes, \bnfsp \bnfts{\textnormal{\textsc{apply}}}(\bnfts{$I_{2h}^{h}$}, \bnfsp \bnfpn{$c_{2h}$}))
	} \\ \\
	\bnfprod{$c_{4h}$} {
		\bnfts{\textnormal{\textsc{apply}}}(\bnfts{$A^{-1}_{4h}$}, \bnfsp \bnfts{\textnormal{\textsc{apply}}}(\bnfts{$I_{2h}^{4h}$}, \bnfsp\bnfpn{$c_{2h}$}))
	} \\
	\bnfprod{$s_{2h}$}{
		\bnfts{\textnormal{\textsc{update}}}(\bnfts{$\omega$}, \bnfsp \bnfes, \bnfsp \bnfts{\textnormal{\textsc{apply}}}(\bnfts{$I_{4h}^{2h}$}, \bnfsp \bnfpn{$c_{4h}$}))
	} \\ \\
	\bnfprod{$c_{8h}$} {
		\bnfts{\textnormal{\textsc{apply}}}(\bnfts{$A^{-1}_{8h}$}, \bnfsp \bnfts{\textnormal{\textsc{apply}}}(\bnfts{$I_{4h}^{8h}$}, \bnfsp\bnfpn{$c_{4h}$}))
	} \\
	\bnfprod{$s_{4h}$}{
		\bnfts{\textnormal{\textsc{update}}}(\bnfts{$\omega$}, \bnfsp \bnfes, \bnfsp \bnfts{\textnormal{\textsc{apply}}}(\bnfts{$I_{8h}^{4h}$}, \bnfsp \bnfpn{$c_{8h}$}))
	}.
\end{bnf*}
\subsection{Grammar-Based Algorithm Generation}
\label{sec:grammar-based-algorithm-generation}
\begin{figure}
	\centering
	\includegraphics[width=\textwidth]{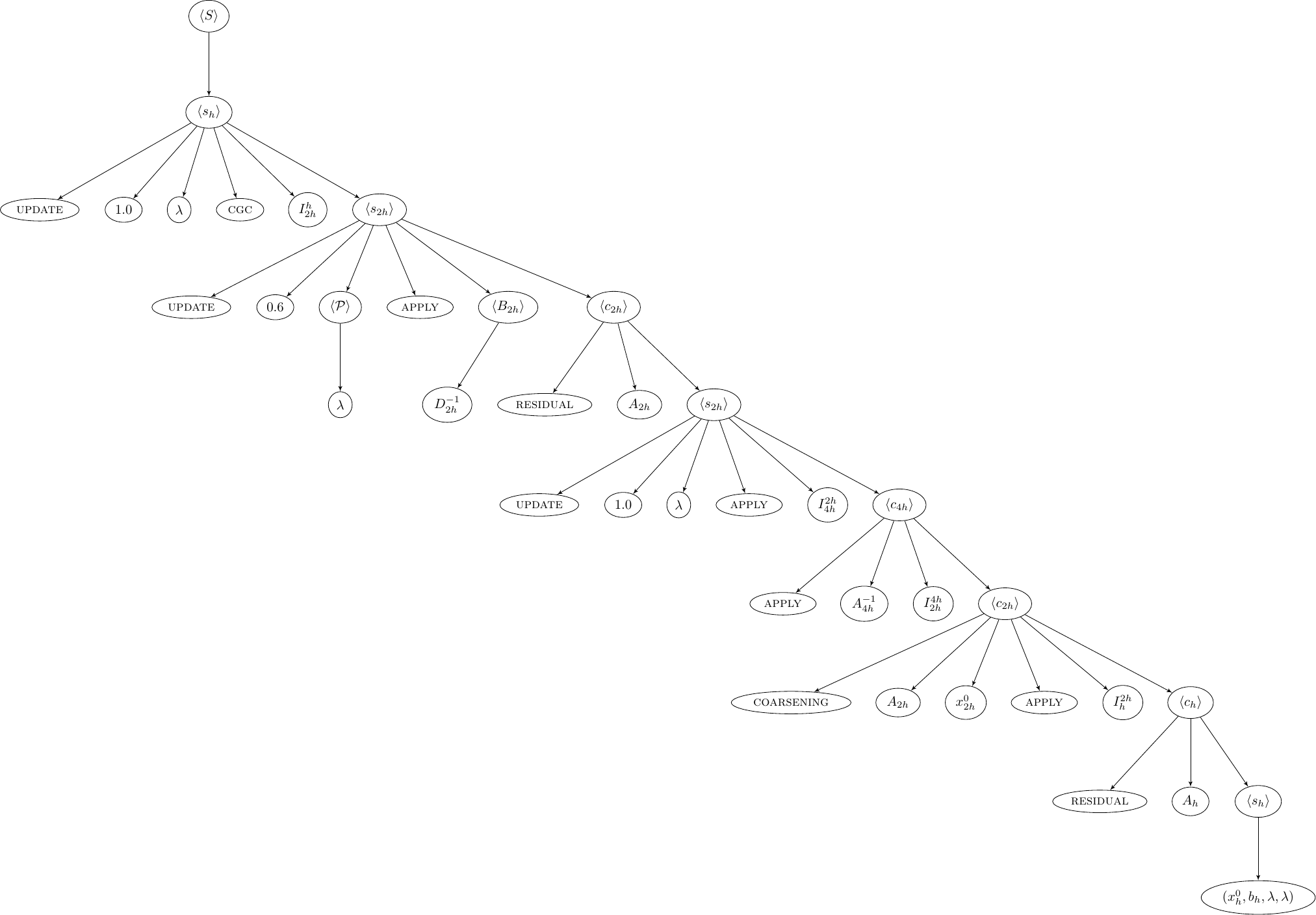}
	\caption{Grammar derivation tree for the three-grid method formulated in Algorithm~\ref{alg:example-three-grid-method}.}
	\label{fig:example-three-grid-method-derivation-tree}
\end{figure}

Up to this point, we have derived a formal language for representing multigrid methods as a sequence of state transitions, where each transition corresponds to a particular computational step within the method.
Furthermore, we have demonstrated that it is possible to define a CFG for generating representations of arbitrarily-structured multigrid methods in that language.
As we have already discussed in Section~\ref{sec:gggp-representation}, every derivation of a CFG can be represented as a tree.
Consider, for instance, the three-grid method shown in Algorithm~\ref{alg:example-three-grid-method}.
We can express this method as a sequence of the productions defined in Algorithm~\ref{table:multigrid-grammar}, which yields the derivation tree shown in Figure~\ref{fig:example-three-grid-method-derivation-tree}.
Note that each inner tree node corresponds to a symbol contained in the set of variables $V$ of our grammar, as shown in Equation~\eqref{eq:multigrid-grammar-variables}, while all leaf nodes refer to terminals.
Since in the context of Algorithm~\ref{table:multigrid-grammar}, the application order of each function is well-defined, we have omitted all parentheses in the tree for the sake of simplicity.
%TODO evtl darauf eingehen, dass der CGS auf dem entsprechenden level spezifiziert sein muss

\begin{figure}
	\captionsetup{justification=centering}
	\begin{subfigure}[b]{0.49\textwidth}\centering
		\includegraphics[height=0.925\textheight]{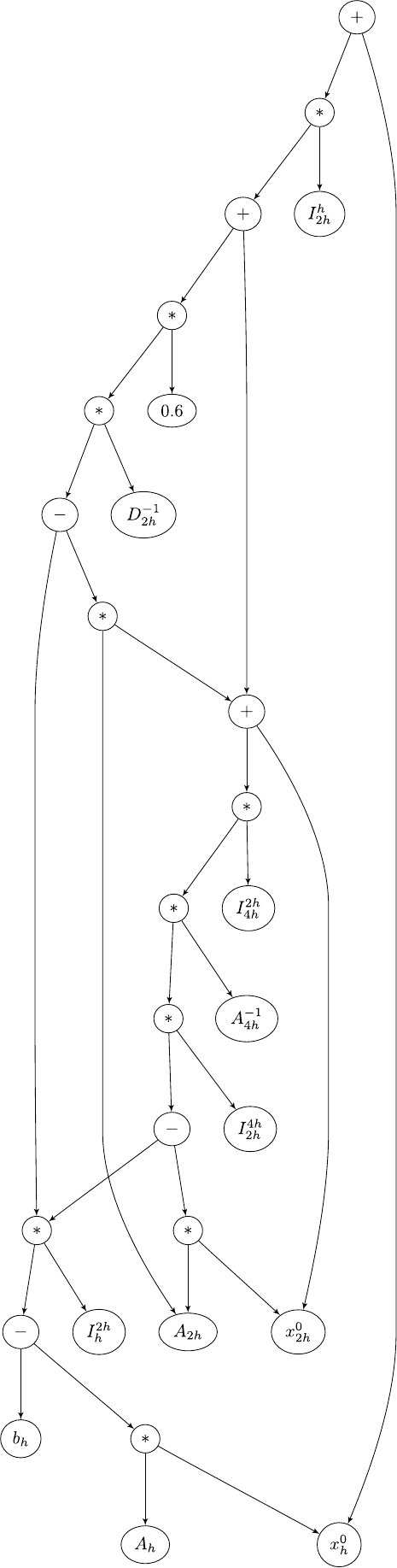}
		\caption{Original graph}
		\label{fig:example-three-grid-method-computational-graph-original}
	\end{subfigure}
	\begin{subfigure}[b]{0.49\textwidth}
 	\centering
		\includegraphics[height=0.925\textheight]{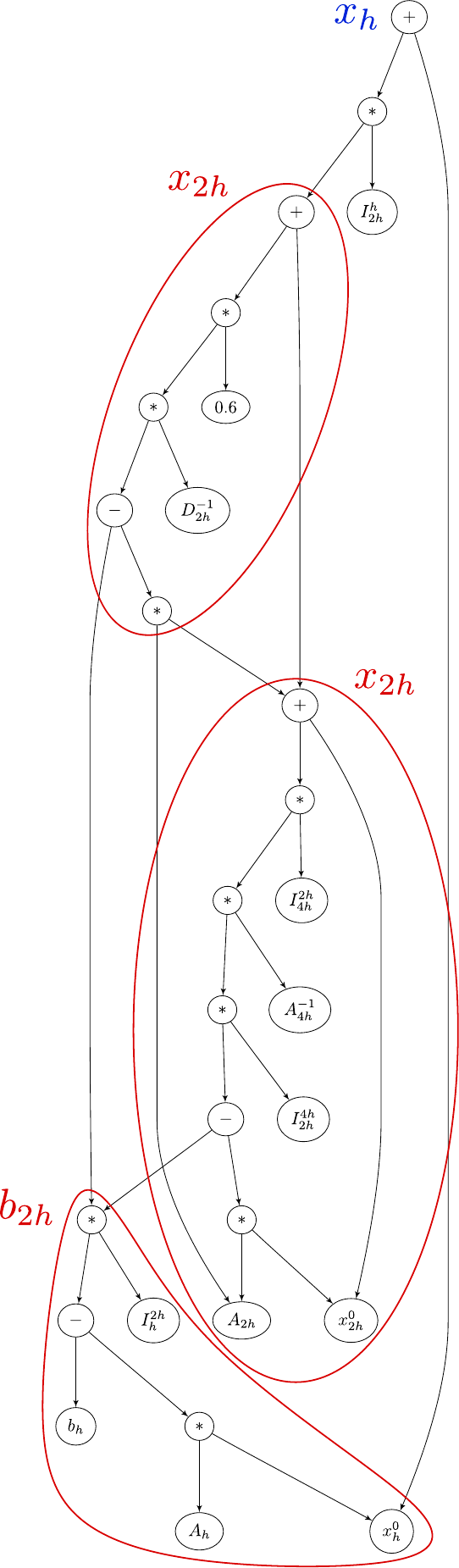}
		\caption{Annotated graph}
		\label{fig:example-three-grid-method-computational-graph-annotated}
	\end{subfigure}
	\caption{Computational graph for the three-grid method defined in Algorithm~\ref{alg:example-three-grid-method}.}
	\label{fig:example-three-grid-method-computational-graph}
\end{figure}

Until now, we have been exclusively concerned with the problem of developing a formal representation for multigrid methods that is based on their algorithmic formulation.
However, within G3P a derivation tree is usually generated by applying a sequence of productions in a randomly-chosen order, either to alter an already existing tree or to create a completely new one from scratch.
We, therefore, have to consider the task of deriving a method's algorithmic formulation from its representation as a derivation tree and the semantical information contained within its nodes.
As we have seen in Section~\ref{sec:multigrid-state-transitions}, the application of each transition function generated by our grammar returns a new state $Z_H$ consisting of a tuple $\left( x_{H}, b_{H}, c_{H}, Z_{H/2}\right)$, where $H$ is the grid spacing on the current level.
Starting from the initial state tuple $\left(x_{h}^0, b_{h}, \lambda, \lambda\right)$, we can hence apply the respective sequence of transition functions in a bottom-up manner until we arrive at the final state $\left(x_{h}, b_{h}, \lambda, \lambda\right)$.
As Figure~\ref{fig:example-tree-grid-method-states} demonstrates, the first component $x_{h}$ of this tuple then combines all computational steps of the corresponding multigrid method in a single expression
\begin{equation}\tag{4.7}
	\begin{split}
		x_h = \; & x_{h}^0 + I_{2h}^h ((x_{2h}^0 + I_{4h}^{2h} A_{4h}^{-1} I_{2h}^{4h} (I_{h}^{2h}(b_{h} - A_h x_{h}^0) - A_{2h} x_{2h}^0)) \\
		& + 0.6 \cdot D_{2h}^{-1} (I_{h}^{2h}(b_{h} - A_h x_{h}^0) - A_{2h} \\
		& \cdot (x_{2h}^0 + I_{4h}^{2h} A_{4h}^{-1} I_{2h}^{4h} (I_{h}^{2h}(b_{h} - A_h x_{h}^0) - A_{2h} x_{2h}^0)))).
		\label{eq:example-three-grid-method-expression}
	\end{split}
\end{equation}
While this expression includes all information about the computational structure of the method, certain terms occur multiple times, which might lead the redundant computations.
We can remove this redundancy by interpreting Equation~\eqref{eq:example-three-grid-method-expression} as a computational graph with directed edges, where each node either corresponds to an arithmetic operation or a predefined symbol, such as $x^0_h$, $b_h$ and $A_h$.
The resulting directed graph is shown in Figure~\ref{fig:example-three-grid-method-computational-graph}.
Again each node within this graph either corresponds to a predefined symbol or an operation on vectors and matrices, where in the case of the latter, there exists a direct edge to each node that refers to one of its arguments.  
Therefore, whenever a node serves as an argument to more than one operation, multiple edges are directed toward it.
Based on this representation, we can now define a redundancy-free sequence of computations that corresponds to the given multigrid method.
The most straightforward way to achieve this for an arbitrary graph is to simply introduce a temporary value for each non-leaf node.
However, note that in the graph shown in Figure~\ref{fig:example-three-grid-method-computational-graph}, every node with multiple incoming edges refers to a term for computing an approximate solution or right-hand side on the respective level.
To verify this assumption, consider again the three elementary multigrid operations listed in Definition~\ref{def:elementary-multigrid-operations}.
While a correction term is always discarded after its application, the current approximate solution is needed in the computation of the residual as well as whenever its value is updated, either through smoothing or by applying a coarse-grid correction.
Note that the same is true for the right-hand side, which is required whenever a new residual is computed on the respective level.
This is illustrated in Figure~\ref{fig:example-three-grid-method-computational-graph-annotated}, where we have annotated all subgraphs that correspond to the computation of an approximate solution or right-hand side.
Based on this graph, we can now easily obtain an algorithmic representation for a multigrid solver while avoiding redundant computations.
For this purpose, we first determine which node does not have any incoming edge.
Since this node computes the final approximate solution, we mark it as $x_h$ in Figure~\ref{fig:example-three-grid-method-computational-graph-annotated}.
The expression that needs to be assigned to this variable is then generated by recursively visiting its children.
Since the graph is traversed starting from the last operation, the algorithmic instructions are generated in reverse order. 
To avoid redundant computations, whenever a node is reached that either corresponds to an approximate solution or a right-hand side, we have to ensure that the respective subgraph is only generated once and assigned to a separate variable, which can then be referenced in all subsequent expressions.
Finally, the process of algorithm generation is carried out until only nodes without any outgoing edges are reached.
For Figure~\ref{fig:example-three-grid-method-computational-graph-annotated}, we thus obtain the algorithmic representation shown in Algorithm~\ref{alg:example-three-grid-method-generated}.
Note that this formulation is mathematically equivalent to the three-grid method shown in Algorithm~\ref{alg:example-three-grid-method}, where the only difference is that certain intermediate expressions, such as the residual, are not assigned to a variable.
\begin{algorithm}
	\begin{algorithmic}[1]
		\State $ b_{2h} = I_{h}^{2h} \left(b_{h} - A_h x_{h}^0 \right)$
		\State $ x_{2h} = x^0_{2h} + I_{4h}^{2h} A_{4h}^{-1} I_{2h}^{4h} \left(b_{2h} - A_{2h} x^0_{2h}\right)$
		\State $ x_{2h} = x_{2h} + 0.6 \cdot D_{2h}^{-1} \left(b_{2h} - A_{2h} x_{2h}\right)$
		\State $x_{h} = x_{h}  + I_{2h}^h x_{2h}$
	\end{algorithmic}
	\caption{Example of a Three-Grid V-Cycle (Generated)}
	\label{alg:example-three-grid-method-generated}
\end{algorithm}
Since we can apply the same approach to any computational graph derived from a derivation tree of our multigrid grammar, we arrive at the following three-step procedure of algorithm generation:
\begin{enumerate}
	\item Construct a multigrid method in the form of its grammar derivation tree.
	\item Transform the derivation tree to a graph-based intermediate representation.
	\item Generate expressions for each approximate solution and right-hand side through recursive graph traversal.
\end{enumerate}
With the formulation of this approach, we have now outlined all necessary steps from the grammar-based generation of a multigrid method to its translation into an algorithmic representation, based on which these methods can be employed as numerical solvers.
However, we have not yet discussed how the individual steps of our approach, as described in this section, can be carried out in an automatic way on a modern computer.
In the remaining part of this thesis, we will close this gap by presenting the implementation of our evolutionary program synthesis framework \emph{EvoStencils}.
Since the EvoStencils framework builds substantially on the techniques presented in Section~\ref{sec:gggp}, it is, however, necessary to first evaluate the possibility of identifying the optimal multigrid method for a given problem using a mere brute-force approach.
%TODO rewrite this part, as it is a bit unclear
For this purpose, we have to estimate the size of the search space spanned by the rules of the grammar defined by Equation~\eqref{eq:multigrid-grammar-variables} and Algorithm~\ref{table:multigrid-grammar}, i.e., the number of different individuals that can be generated based on its productions.

\subsection{Search Space Estimation}
\label{sec:search-space-estimation}
\begin{figure}
	\centering
	\includegraphics[width=\textwidth]{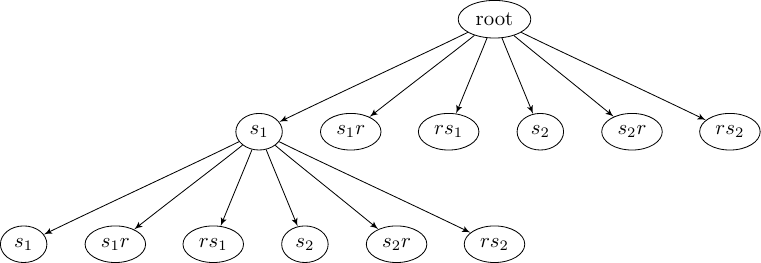}
	\caption[Simplified search space]{Simplified search space -- Each node in the decision tree has the same number of children.}
	\label{fig:decision-tree}
\end{figure}

If we treat the grammar-based design of multigrid methods as a search problem, the size of the search space spanned by the underlying grammar corresponds to the number of different methods that can be constructed based on its productions.
For a small search space, it can be feasible to simply enumerate and evaluate all possible solutions. 
Such an approach has the advantage that the best solution can be determined in a deterministic manner while heuristic search methods, such as those presented in Section~\ref{sec:gggp}, are, in general, not guaranteed to find it.
We consider the following simplified model estimation, which acts as a lower bound of the size of the search space of a five-grid hierarchy, as defined by the productions in Algorithm~\ref{table:multigrid-grammar}. 
As we have discussed in Section~\ref{sec:multigrid-methods}, multigrid methods need to combine smoothing and coarse-grid correction steps to effectively reduce both the high- and low-frequency components of the error.
While in practice often multiple smoothing steps are used, the application of a coarse-grid correction without any smoothing is rarely effective on any level~\cite{trottenberg2000multigrid}.
To simplify the search space estimation, we, therefore, introduce the additional constraint that, on any level, the total number of smoothing steps is always greater or equal to the number of coarse-grid corrections.
We, furthermore, enforce that coarsening and coarse-grid corrections can only be performed after smoothing, with the exception of the topmost level, where we additionally allow coarsening to be applied as the first operation.
Since the original system is solved on the finest grid, finding an optimal sequence of operations on this level is indespensible for the efficiency of a multigrid method.
We, therefore, want to enforce fewer constraints on the order of the individual operations.
Using the symbols $s$, $r$, and $p$ for smoothing, coarsening, and coarse-grid correction\footnote{For the sake of simplicity, we use the starting letters of restriction and prolongation as abbreviations for coarsening and coarse-grid correction.}, respectively, the following operations can be defined in each step:
\begin{equation}\tag{4.8}
\begin{split}
	s | rs | sr & \quad \text{for} \; l = l_{max} \\
	A^{-1} & \quad \text{for} \; l = l_{min} \\
	s | sr | sp & \quad \text{otherwise}.
\end{split}
\end{equation}
Consequently, with the exception of the coarsest level, where the only allowed operation is the application of the coarse-grid solver, there always exist three different orderings of the available operations.
In addition, each step includes a smoother, which additionally needs to be chosen.
Assuming the number of available smoothers is $n$, which includes different types of smoothers, partitionings, and relaxation factors, this results in a total number of $3 n$ options in each step of the method.
While none of these constraints is present in our original grammar, these simplifications enable us to derive an expression in closed form, which then acts as a lower bound for the actual size of the search space.
We can derive such an expression as a finite sum of polynomials by considering the total number of smoothing steps within a multigrid method.
Here we know that for each step, we need to choose from $3 n$ options independent of the level on which the operations are performed.
As a consequence, if our method incorporates $i$ smoothing steps, the total number of options is 
\begin{equation}\tag{4.9}
	N_{n,i} = (3 n)^i
	\label{eq:simplified-number-of-options}
\end{equation}
We can illustrate this expression by considering the corresponding decision tree, which is shown in Figure~\ref{fig:decision-tree} for $n = 2$.
Each node in this tree corresponds to a particular choice of operations, which includes the application of either $s_1$ or $s_2$ as a smoother, optionally accompanied by coarsening or a coarse-grid correction.
Therefore, the depth of the tree is equal to the total number of smoothing steps $k$ applied within the multigrid method.
Since each leaf of the tree corresponds to a structurally-different multigrid method, the total number of leaves is equal to $N_{2,i} = 6^i$.
However, since the optimal amount of smoothing is not known beforehand, we must consider a varying number of smoothing steps.
Now note that each decision tree that corresponds to a certain number of smoothing steps $i$ includes all those with fewer smoothing steps as a subtree, where each of these subtrees can be obtained by performing a top-down breadth-first traversal until the depth of the tree is equal to $i$.
Therefore, the total number of different multigrid methods with at least $i_{min}$ but at most $i_{max}$ smoothing steps is equal to the total number of leaf nodes of all subtrees of a decision tree of depth $i_{max}$ that have at least depth $i_{min}$.
Since the number of leaf nodes of a decision tree is given by Equation~\eqref{eq:simplified-number-of-options}, the total number of different multigrid methods that employ between $i_{min}$ and $i_{max}$ smoothing steps is equal to
\begin{equation}\tag{4.10}
	N(n, i_{min}, i_{max}) = \sum_{i = i_{min}}^{i_{max}} N_{n,i} = \sum_{i = i_{min}}^{i_{max}} (3 n)^i.
	\label{eq:search-space-estimation}
\end{equation}
We can now finally make use of this formula to compute a lower bound for the total number of different five-grid methods that can be generated based on our grammar.
For this purpose, we need to choose a reasonable interval for the minimum and maximum number of smoothing steps.
In Algorithm~\ref{alg:multigrid-cycle}, the number of smoothing steps per level depends on the parameters $\nu_1$, $\nu_2$, and $\gamma$, which set the number of pre-smoothing steps, post-smoothing steps and recursive descents per level, respectively.
Therefore, the total number of smoothing steps is given by
\begin{equation}\tag{4.11}
	\sum_{i = 0}^{k - 2} \gamma^i (\nu_1 + \nu_2),
\end{equation}
where $k$ is the number of levels of the discretization hierarchy, as in Algorithm~\ref{alg:multigrid-cycle}.
%TODO only true for W-cycle, fix formula!
%TODO check if this has just been forgotten
Now let us consider how different choices of $\nu_1$, $\nu_2$, and $\gamma$ lead to a different total number of smoothing steps for a five-grid method ($k = 5$):
\begin{itemize}
	\item V(1,0): $\gamma = 1, \nu_1 = 1, \nu_2 = 0 \Rightarrow 4 \; \text{smoothing steps}$
	\item V(3,3): $\gamma = 1, \nu_1 = 3, \nu_2 = 3 \Rightarrow 24 \; \text{smoothing steps}$
	\item W(1,0): $\gamma = 2, \nu_1 = 1, \nu_2 = 0 \Rightarrow 15 \; \text{smoothing steps}$
	\item W(1,1): $\gamma = 2, \nu_1 = 1, \nu_2 = 0 \Rightarrow 30 \; \text{smoothing steps}$
	\item W(3,3): $\gamma = 2, \nu_1 = 3, \nu_2 = 3 \Rightarrow 90 \; \text{smoothing steps}$
\end{itemize}
Setting $i_{min} = 4$ and $i_{max} = 30$ as an upper limit for the number of smoothing steps in Equation~\eqref{eq:search-space-estimation} represents a reasonable compromise. 
It allows us to generate multigrid methods that perform a similar amount of smoothing as the majority of commonly used V-cycles, while methods comparable to less expensive W-cycle variants can also be obtained.
Finally, the only remaining parameter is the number of different methods available for each smoothing step.
To obtain a conservative estimate, we set $n = 2$, which means that in each step, we can only choose from two alternatives, for example, the Jacobi and red-black Gauss-Seidel method.
Based on Equation~\eqref{eq:search-space-estimation} we hence obtain 
\begin{equation}\tag{4.12}
	N(2, 4, 30) = \sum_{i = 4}^{30} 6^i \approx 2.65 \cdot 10^{24}
\end{equation}
as an estimated lower bound for the search space spanned by the productions defined in Algorithm~\ref{table:multigrid-grammar}.
If we now assume that a modern multi-core CPU is capable of evaluating each multigrid method on average in one millisecond,
even a supercomputer consisting of one trillion ($10^{12}$) such processors would require more than eight years to evaluate all methods contained in this search space.
In practice, it is usually beneficial to consider a higher number of different smoothers and relaxation factors, which results in an even larger search space.
As a consequence, while a pure brute-force-approach that aims to evaluate all possible multigrid methods may be feasible for the constrained parameter space that results from the classical formulation in Algorithm~\ref{alg:multigrid-cycle}, a grammar-based design requires the utilization of heuristic search method, such as those presented in Section~\ref{sec:gggp}.
\paragraph{Final Remarks}
In the remainder of this thesis, we will describe our implementation of grammar-guided genetic programming (G3P) in the Python programming language and how it can be utilized for the automated grammar-based design of multigrid methods.
%We will, furthermore, demonstrate EvoStencils' effectiveness by evolving efficient multigrid methods for a number of different partial differential equations (PDEs).
We, therefore, assume a basic familiarity with Python due to its widespread use in both academia and industry.
If this should not be the case, the reader is advised to consult one of the many available resources for learning Python\footnote{Python:~\url{https://www.python.org/doc/}}.
We, furthermore, expect the reader to have a basic knowledge of programming and parallelization techniques for recent multiprocessor and cluster systems, of which an overview can be found in~\cite{sterling2017high,hager2010introduction}.
Finally, while hardware-specific code optimization is not a focus of this thesis, the reader should have some basic understanding of the architectural features of modern computers.
%Along with this, we also assume a familiarity with programming techniques and tools for performing a shared and distributed memory parallelization on these systems, such as OpenMP\footnote{OpenMP:~\url{https://www.openmp.org/}} and MPI\footnote{MPI:~\url{https://www.mpi-forum.org/}}. 
%\begin{figure}
%	\captionsetup{justification=centering}
%	\begin{subfigure}{0.1\textwidth}
%		\begin{tikzpicture}
%			\node   (h) at (-0.75, 4){$h$};
%			\node   (2h) at (-0.75, 3){$2h$};
%			\node   (4h) at (-0.75, 2){$4h$};
%			\node   (8h) at (-0.75, 1){$8h$};
%			\node   (16h) at (-0.75, 0){$8h$};
%		\end{tikzpicture}
%	\end{subfigure}
%	\begin{subfigure}{0.9\textwidth}
%		\begin{tikzpicture}
%			\node	(a) at (0,4) [draw, circle,fill=blue,scale=0.8] {};
%			\node	(b) at (0.5,3) [draw, circle, fill=blue, scale=0.8] {};
%			\node	(c) at (1,2) [draw, circle, fill=blue, scale=0.8] {};
%			\node	(d) at (1.5,1) [draw, circle, fill=blue, scale=0.8] {};
%			\node	(e) at (2,0) [draw, circle, ,fill=black, scale=0.8] {};
%			\node	(f) at (2.5,1) [draw, circle,scale=0.8] {};
%			\node	(g) at (3,2) [draw, circle,scale=0.8] {};
%			\node	(h) at (3.5,3) [draw, circle, scale=0.8] {};
%			\node	(i) at (4,4) [draw, circle,scale=0.8] {};
%			\draw 
%			(a) edge[->] (b) 
%			(b) edge[->] (c)
%			(c) edge[->] (d)
%			(d) edge[->] (e)   
%			(e) edge[->] (f)
%			(f) edge[->] (g)
%			(g) edge[->] (h)
%			(h) edge[->] (i)
%			;
%		\end{tikzpicture}
%	\end{subfigure}
%	
%	\caption{Five-grid cycle with $\nu_1 = 1$ and $\nu_2 = 0$. Smoothing is only performed in the steps corresponding to blue-colored nodes.}
%	\label{fig:five-grid-V-cycle}
%\end{figure}

%% file: contents/evostencils_part1.tex
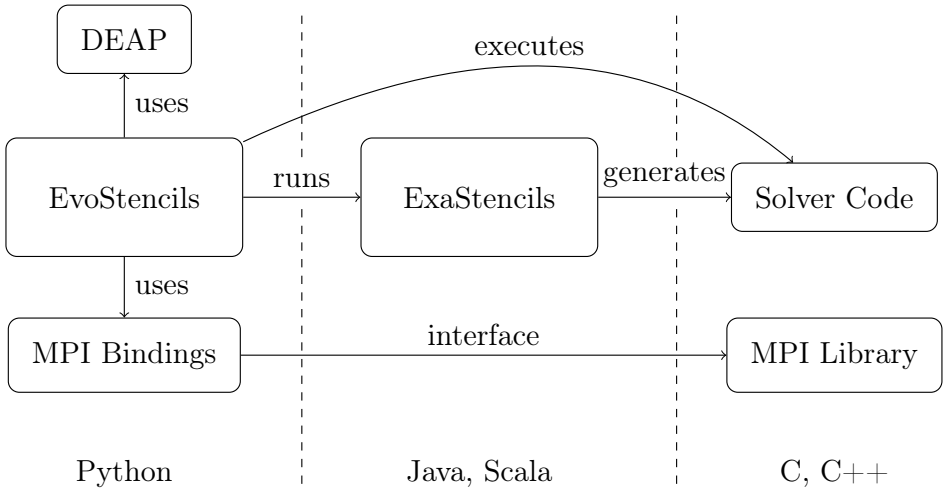
\begin{figure}
	\resizebox{\columnwidth}{!}{%
		\begin{tikzpicture}
			%\draw [help lines] (-10,-10) grid (10,10);
			\node[draw, minimum width=3cm, minimum height=1.5cm, rounded corners] (evo) at (0,0) {EvoStencils};
			\node[draw, inner sep=3mm, rounded corners] (bindings) at (0, -2) {MPI Bindings};
			\node[draw, inner sep=3mm, rounded corners] (deap) at (0, 2) {DEAP};
			\node[draw, minimum width=3cm, minimum height=1.5cm, rounded corners] (exa) at (4.5, 0) {ExaStencils};
			\node[draw, inner sep=3mm, rounded corners] (code) at (9, 0) {Solver Code};
			\node[draw, inner sep=3mm, rounded corners] (mpi) at (9, -2) {MPI Library};
			\draw[dashed] (2.25, 2.3) -- (2.25,0.5);
			\draw[dashed] (2.25, -0.15) -- (2.25,-3.7);
			\draw[dashed] (7, 2.3) -- (7,0.5);
			\draw[dashed] (7, -0.15) -- (7,-3.7);
			\node (python) at (0, -3.5) {Python};
			\node (java) at (4.5, -3.5) {Java, Scala};
			\node (c) at (9, -3.5) {C, C++};
			\draw[->] (evo)-- node[anchor=west] {uses} (deap);
			\draw[->] (evo)-- node[anchor=south]{runs} (exa);
			\draw[->] (exa)--node[anchor=south] {generates} (code);
			\draw[->] (evo) to [out=25,in=140] node[anchor=south] {executes} (code);
			\draw[->] (bindings)--node[anchor=south] {interface}(mpi);
			\draw[->] (evo)--node[anchor=west] {uses}(bindings);
			%\draw[->] (mpi)--(code);
		\end{tikzpicture}
	}\caption{Software Architecture of EvoStencils.}
	\label{fig:evostencils-architecture}
\end{figure}

In the first part of this thesis, we have established a theoretical foundation for multigrid methods, formal languages, and evolutionary program synthesis.
In Chapter~\ref{chapter:multigrid-formal-language}, building on this foundation, we have then developed a novel formal language and grammar for the automatic generation of multigrid methods. 
While we have already demonstrated the capability of this approach to modify each individual step of a multigrid method, we could not yet demonstrate its benefits compared to classical multigrid cycles, such as V-, F-, and W-cycles.
We aim to achieve this goal with the implementation of \emph{EvoStencils}, a prototypical Python framework for the automated grammar-based design of multigrid methods, which we have made available as open-source software\footnote{EvoStencils: \url{https://github.com/jonas-schmitt/evostencils}}.
Using this framework, we will demonstrate the discovery of multigrid methods that are able to solve certain PDE-based problems faster than all classical multigrid cycles.
Before we discuss EvoStencils' features and their implementation in Python, we want to provide an overview of its general workflow and software architecture.
Here, we distinguish between EvoStencils' core implementation and functionality that builds upon external libraries.
Since we have expressed the rules for constructing a multigrid method in the form of a context-free grammar, we can apply the evolutionary program synthesis techniques presented in Chapter~\ref{chapter:formal-languages-and-gp} without significant adaption.
For this purpose, we employ the widely-used evolutionary computation framework DEAP\footnote{DEAP: \url{https://github.com/deap/deap}}~\cite{rainville2012deap}, which enables us to implement grammar-guided genetic programming (G3P) in a modular way.
However, to realize this approach, we need to evaluate each multigrid method obtained through G3P in an automatic and reproducible manner.
As we have seen in Section~\ref{sec:grammar-based-algorithm-generation} the 
evaluation of the sequence of state transitions contained in a particular derivation tree produces a computational graph in the form of Figure~\ref{fig:example-three-grid-method-computational-graph}.
This graph can then be translated to an algorithmic representation similar to Algorithm~\ref{alg:example-three-grid-method-generated}.
%While a domain expert could manually implement the corresponding multigrid solver based on this representation using a numerical software package, our framework has to perform the evaluation of each method in an automatic way without requiring any human intervention.
Recently, code-generation techniques that only require the specification of a numerical solver in a high-level domain-specific language (DSL) have become increasingly powerful~\cite{kostler2020code}.
An example of this approach is the ExaStencils framework~\cite{lengauer2020exastencils,lengauer2014exastencils}, which has been specifically designed for the automatic generation of fast and scalable implementations of multigrid-based solvers specified in a tailored DSL called ExaSlang~\cite{schmitt2014exaslang,schmitt2016systems,kuckuk2016automatic}.
ExaSlang represents a multigrid method as a sequence of high-level operations while granting the user the flexibility to apply further optimizations through the addition of code transformations and lower-level statements.
To evaluate a given solver obtained from a grammar-based representation, we emit its corresponding algorithmic formulation as an ExaSlang specification, based on which we then generate a scalable C++ implementation using the capabilities of the ExaStencils framework.
The resulting program can then be executed on a number of test cases in order to measure its desired performance characteristics.
Finally, note that the execution of an evolutionary program synthesis method requires the evaluation of a large number of different programs, each representing a unique numerical solver.
Depending on the problem that one aims to solve, it can be infeasible to run this method on a single compute node, which necessitates a multi-node parallelization.
The message-passing interface (MPI)~\cite{walker1996mpi} provides a unified interface for performing parallel computations on a distributed system that is supported by the majority of available supercomputing devices.
While MPI was originally designed for the traditional scientific computing languages Fortran and C, it has recently been made available within Python~\cite{dalcin2021mpi4py}. 
With the addition of MPI as a distributed computing backend, we arrive at the high-level view of EvoStencils' software architecture shown in Figure~\ref{fig:evostencils-architecture}.
In the following, we will discuss the individual parts of this architecture in more detail, starting with the core implementation of EvoStencils, which does not depend on any of the other tools and libraries mentioned here.
As a first step, we will outline the implementation of an intermediate representation (IR) for multigrid methods that can be generated in a straightforward manner based on a given derivation tree.
This IR will then act as a basis for all subsequent steps of solver generation and evaluation.

\section{Intermediate Representation}
\label{sec:intermediate-representation}
Before we derive an IR for each component of a multigrid method, note that each of them needs to be defined with respect to the chosen discretization.
%Note that algebraic multigrid methods~\cite{stuben2001introduction,ruge1987algebraic}, which are not considered in this work, represent an exception to this, as they operate on sparse matrix data structures.
In Section~\ref{sec:discretization}, we have already made the assumption of discretizing the underlying PDE on a hierarchy of structured grids.
To identify a grid within this hierarchy, certain information is required, which we store in a \emph{Grid} data structure whose implementation is shown in Listing~\ref{code:ir:grid}.
In general, a structured grid is defined by its \emph{spacing} ($h$) in each dimension and its \emph{size}, i.e., the number of grid points.
In addition, we also include the grid's \emph{level} to identify it within the discretization hierarchy.
Note that in case the grid is uniform, we only need to store a single value for its spacing in each dimension, while otherwise, a value would need to be stored for each pair of grid points.
As the problems considered in this work are all solved on a hierarchy of uniform grids, we focus on this particular case.
\begin{listing}
\begin{minted}{python}

class Grid:
    def __init__(self, size, spacing, level):
        assert len(size) == len(spacing), 
        "Dimensions of the size and step size must match"
        self._size = size
        self._spacing = spacing
        self._level = level

    @property
    def size(self):
        return self._size

    @property
    def spacing(self):
        return self._spacing

    @property
    def level(self):
        return self._level

    @property
    def dimension(self):
        return len(self.size)
\end{minted}
	\caption{IR -- Grid Data Structure}
	\label{code:ir:grid}
\end{listing}
After defining a data structure that provides all relevant information about a specific grid within the discretization hierarchy, we can start defining expressions that operate on this data structure.
In Listing~\ref{code:ir:abstact-base-class}, a common abstract base class is provided from which all subsequent expression classes are derived.
\begin{listing}
\begin{minted}{python}
import abc

class Expression(abc.ABC):

    @property
    @abc.abstractmethod
    def shape(self):
        pass

    @property
    @abc.abstractmethod
    def grid(self):
        pass
\end{minted}
	\caption{IR -- Abstract Expression}
	\label{code:ir:abstact-base-class}
\end{listing}
In addition to the already mentioned grid data structure, this class also defines a \emph{shape} for each expression.
From a mathematical point of view, each expression within a multigrid method either computes a matrix or a vector whose shape can be derived recursively from its operands.
This shape is thus defined as a pair $(r, c)$, whose first entry $r$ corresponds to the number of rows and the second $c$ to the number of columns of the corresponding matrix.
Note that a vector of size $n$ can simply be considered a matrix of shape $(n,1)$.
Based on this abstract class, we can distinguish between predefined entities, such as the system matrix and right-hand side, and expressions that refer to the mathematical operations defined within a multigrid method.
First of all, Listing~\ref{code:ir:entity} contains the implementation of the \emph{Entity} base class.
\begin{listing}
\begin{minted}{python}
class Entity(Expression):

    def __init__(self, name, grid, shape):
        self._name = name
        self._grid = grid
        self._shape = shape
        super().__init__()

    @property
    def name(self):
        return self._name

    @property
    def grid(self):
        return self._grid

    @property
    def shape(self):
        return self._shape
\end{minted}
	\caption{IR -- Entity}
	\label{code:ir:entity}
\end{listing}
In addition to the previously mentioned attributes, this class is also given a \emph{name} to identify the respective entity.
We can then further define classes for representing an approximate solution, right-hand side, and operator, which are shown in the Listings~\ref{code:ir:approximation}~and~\ref{code:ir:operator}.
\begin{listing}
\begin{minted}{python}
from operator import mul as operator_mul

class Approximation(Entity):
    def __init__(self, name, grid):
        shape = (reduce(operator_mul, grid.size), 1)
        super().__init__(name, grid, shape)

    @property
    def predecessor(self):
        return None

class RightHandSide(Approximation):
    pass
\end{minted}
	\caption{IR -- Approximate Solution and Right-Hand Side}
	\label{code:ir:approximation}
\end{listing}
The shape of each of these entities is determined by computing the product of the grid size over all dimensions.
From a mathematical point of view, the \emph{Approximation} and \emph{RightHandSide} classes both represent vectors.
We, therefore, only implement the former and make use of inheritance to avoid unnecessary code duplication in the case of the latter.
In addition to the attributes defined in its parent class, an \emph{Approximation} includes a \emph{predecessor} attribute.
However, we postpone the discussion of this attribute's purpose until we implement the actual grammar.
While the two Python classes shown in Listing~\ref{code:ir:approximation} correspond to the solution and right-hand side of a discretized PDE, Listing~\ref{code:ir:operator} represents its operator, given as one or multiple stencil codes.
\begin{listing}
\begin{minted}{python}
from operator import mul

class Operator(Entity):
    def __init__(self, name, grid, stencil_generator=None):
        n = reduce(mul, grid.size)
        shape = (n, n)
        self._stencil_generator = stencil_generator
        super().__init__(name, grid, shape)

    @property
    def stencil_generator(self):
        return self._stencil_generator

    def generate_stencil(self):
        if self._stencil_generator is None:
            return None
        return self._stencil_generator.generate_stencil(self._grid)
\end{minted}
	\caption{IR -- Operator}
	\label{code:ir:operator}
\end{listing}
In Section~\ref{subsec:stencil-codes}, we have already introduced a mathematical notation for stencil codes in the form of Equation~\eqref{eq:stencil-definition}, which can be implemented in a straightforward manner leading to the \emph{Stencil} class shown in Listing~\ref{code:ir:stencil}.
%TODO introduce stencil implementation here
\begin{listing}
\begin{minted}{python}
class Stencil:
    def __init__(self, entries, dimension=None):
        self._dimension = dimension
        self._entries = tuple(entries)

    @property
    def entries(self):
        return self._entries

    @property
    def dimension(self):
        if self._dimension is None:
            return len(self.entries[0][0])
        else:
            return self._dimension

    @property
    def number_of_entries(self):
        return len(self.entries)
\end{minted}
	\caption{IR -- Stencil}
	\label{code:ir:stencil}
\end{listing}
The \emph{entries} attribute of this class directly corresponds to our definition of a stencil, whereby we use a \emph{tuple} object to represent this mathematical entity in the Python programming language.
Each entry $e_i$ of this object is thus given as
\begin{equation}
	e_i = \left(\bm{a}_i, b_i \right),
\end{equation} 
where $\bm{a}_i$ is the offset from the current grid point in each dimension, given as an array of integer values, and $b_i$ the stencil value that corresponds to each offset, stored as a floating point number.
Based on this class, we can then provide an implementation for each stencil operation defined in Section~\ref{subsec:stencil-codes}.
%Furthermore, in accordance with Section~\ref{subsec:systems-of-pdes} and~\ref{subsec:block-smoothing}, we can extend this implementation to systems of PDEs and block smoothers, which will be briefly covered in the next Chapter of this thesis. %TODO insert ref to next chapter here or remove sentence
Now note that in our implementation of the \emph{Operator} class in Listing~\ref{code:ir:operator}, the stencil is not included directly, but instead, we provide a so-called \emph{stencil generator}.
A stencil generator is a function that returns the discretization of an operator on a particular grid in the form of a stencil.
For instance, the finite-difference discretization of the Laplace operator on a two-dimensional uniform grid leads to the stencil 
\begin{equation*}
	\begin{split}
		\Delta_{h,h} = & \; \big\{ \left( \left( 0,0 \right), 4 / h^2 \right), \left(\left(1,0\right), -1/h^2\right), \left(\left(-1,0\right), -1 / h^2\right), \\ & \left(\left(0,1\right), -1/h^2\right), \left(\left(0,-1\right), -1/h^2\right) \big\}_{h,h}
	\end{split}
\end{equation*}
in which the value at each offset depends on the grid spacing $h$.
Furthermore, in certain cases, it is possible to derive the higher-dimensional version of a stencil from its lower-dimensional counterpart, as we have shown in Section~\ref{subsec:restriction-and-prolongation} for the considered prolongation and restriction operators.
Therefore, instead of storing a unique stencil for each individual operator instance, we can parametrize its generation based on the features of the applied discretization by including a reference to the respective generator function.
Finally, as both prolongation and restriction operators transfer information between adjacent grids within a hierarchy of discretizations, they can be considered as a special case of a general operator, which, for instance, leads to a different value of the \emph{shape} attribute.
For this purpose, the class \emph{InterGridOperator}, from which all restriction and prolongation operators are derived, extends the \emph{Operator} class with the required functionality.
For the sake of brevity, the implementation of this class and the respective \emph{Restriction} and \emph{Prolongation} subclasses can be found in Section~\ref{appendix:ir} of the appendix.

After discussing the implementation of the different entities based on which a multigrid method is built, we next shift our attention to the implementation of IR classes for representing the expressions that constitute its computational structure. 
For this purpose, we first provide basic classes for representing unary and binary expressions, which are shown in Listing~\ref{code:ir:unary-expression} and~\ref{code:ir:binary-expression}.
\begin{listing}
\begin{minted}{python}
class UnaryExpression(Expression):

    def __init__(self, operand):
        self._operand = operand
        self._shape = operand.shape
        super().__init__()

    @property
    def operand(self):
        return self._operand

    @property
    def shape(self):
        return self._shape

    @property
    def grid(self):
        return self.operand.grid
\end{minted}
	\caption{IR -- Unary Expression}
	\label{code:ir:unary-expression}
\end{listing}
\begin{listing}
\begin{minted}{python}
class BinaryExpression(Expression):

    def __init__(self, operand1, operand2):
        self._operand1 = operand1
        self._operand2 = operand2
        super().__init__()

    @property
    def operand1(self):
        return self._operand1

    @property
    def operand2(self):
        return self._operand2

    @property
    def shape(self):
        raise NotImplementedError("Shape undefined in binary expression")

    @property
    def grid(self):
        return self.operand1.grid
\end{minted}
	\caption{IR -- Binary Expressions}
	\label{code:ir:binary-expression}
\end{listing}
In both cases, all necessary properties are obtained from the operands of the expression in a recursive manner.
However, since the result of a binary expression depends on the type of operation, we raise an error if this method is not implemented in one of the derived classes. 
To illustrate that the majority of multigrid operations can already be expressed based on these two classes, we have included a number of specific examples in Section~\ref{appendix:ir} of the appendix.
As discussed in Section~\ref{sec:grammar-based-algorithm-generation}, we aim to represent the computational structure of each multigrid method as a redundancy-free directed graph. 
While the previously defined base classes allow us to represent the arithmetic expressions that occur within the correction terms of a multigrid method, there are two operations that require special treatment.
As we have already seen in Figure~\ref{fig:example-three-grid-method-computational-graph}, it is necessary to access previously computed intermediate results on multiple occasions within a multigrid method.
In particular, each time a coarse-grid correction is performed, we have to restore the previous approximate solution and right-hand side on the respective level.
Furthermore, whenever we compute a new residual, the current expression of both the right-hand side and the approximate solution is required.
For this purpose, we implement the classes \emph{Residual} and \emph{Cycle}, which allow us to establish additional references to previously defined expressions within a multigrid method.
Each of these references then corresponds to one of the subgraphs in Figure~\ref{fig:example-three-grid-method-computational-graph} that possess root nodes with multiple incoming edges, i.e., those subgraphs that have been annotated in Figure~\ref{fig:example-three-grid-method-computational-graph-annotated}.
We will later see that the complete graph of a multigrid method can be constructed based on these references.
Listing~\ref{code:ir:residual} shows the implementation of the \emph{Residual} class, which contains references to the system operator $A_H$, the current approximate solution $x_H$ and right-hand side $b_H$, where $H$ is the grid spacing on the current level.
Based on these components, we can easily construct the corresponding residual expression $r_H = b_H - A_H x_H$.
\begin{listing}
\begin{minted}{python}
class Residual(Expression):
    def __init__(self, operator, approximation, rhs):
        self._operator = operator
        self._approximation = approximation
        self._rhs = rhs
        super().__init__()

    @property
    def shape(self):
        return self.rhs.shape

    @property
    def grid(self):
        return self.rhs.grid

    @property
    def operator(self):
        return self._operator

    @property
    def approximation(self):
        return self._approximation

    @property
    def rhs(self):
        return self._rhs
\end{minted}
	\caption{IR -- Residual}
	\label{code:ir:residual}
\end{listing}
As a final step in the implementation of our intermediate representation, the \emph{Cycle} class can be found in Listing~\ref{code:ir:cycle}.
\begin{listing}
\begin{minted}{python}
class Cycle(Expression):
    def __init__(self, approximation, rhs, correction=None, relaxation_factor=1.0, partitioning=None, predecessor=None):
        self.approximation = approximation
        self.rhs = rhs
        self.correction = correction
        self.relaxation_factor = relaxation_factor
        self.partitioning = partitioning
        self.predecessor = predecessor
        super().__init__()

    @property
    def shape(self):
        return self.approximation.shape

    @property
    def grid(self):
        return self.approximation.grid
\end{minted}
	\caption{IR -- Cycle}
	\label{code:ir:cycle}
\end{listing}
This class represents the execution of a complete multigrid cycle on a particular level.
It thus computes a new value for the approximate solution $x_H$ on a particular level with spacing $H$, i.e.,
\begin{equation*}
	x_H = x_H + \omega c_H \; \text{with} \; P,
\end{equation*}
where $c_H$ is a correction term, $\omega$ the relaxation factor and $P$ a partitioning.
Note that $x_H$ and $c_H$ contain all previous computations performed within the given cycle.
To make the right-hand side available in subsequent steps of the method, such as for the computation of the residual, the data structure includes an additional reference to the corresponding object.
Additionally, a reference to the previous state on the next-higher level is needed to restore the expression for the approximate solution and right-hand side after applying a coarse-grid correction.
To better understand the purpose of the \emph{Residual} and \emph{Cycle} class, consider the example shown in Listing~\ref{code:ir:example.py}, which demonstrates the construction of a computational graph based on the IR described in this section.
\begin{listing}
\begin{minted}{python}
l = 10 # level
h = 1/2**l # fine-grid spacing
n_h = 1/h - 1
H = 2*h # coarse-grid spacing
n_H = 1/H - 1

# Generate the fine and coarse grid
grid = Grid((n_h, n_h), (h, h), l)
coarse_grid = Grid((n_H, n_H), (H, H), l-1)

# Define entities on the fine grid
x_h = Approximation("x_h", grid)
b_h = RightHandSide("b_h", grid)
A_h = Operator("A_h", grid)
r_h = Residual(A_h, x_h, b_h)
# Create a temporary Cycle object without a correction term
x_h = Cycle(x_h, b_h)

# Define the entities on coarse-grid
A_H = Operator("A_H", coarse_grid)
x_H = Approximation("x_H", coarse_grid)
b_H = Multiplication(Restriction("I_hH", grid, coarse_grid), r_h)
r_H = Residual(A_H, x_H, b_H)

# Define a complete Cycle on the coarse level
x_H = Cycle(x_H, b_H, r_H, predecessor=x_h)
# Restore the state on the fine level
x_h = x_H.predecessor
# Replace the old correction term with the computed coarse-grid correction
x_h.correction = Multiplication(Prolongation("I_Hh", grid, coarse_grid), x_H)
\end{minted}
	\caption{Example Usage of the Intermediate Representation}
	\label{code:ir:example.py}
\end{listing}
Starting on the original problem on the finest grid, we first store references to the initial approximate solution and right-hand side in a \emph{Cycle} object, which itself is included as a \emph{predecessor} reference in the subsequently created coarse-grid \emph{Cycle} object.
In order to apply the latter as a coarse-grid correction on the finest grid, the original fine-grid \emph{Cycle} object is restored, and its \emph{correction} variable is replaced by the respective expression, which is obtained by applying the prolongation operator to the approximate solution that has been previously computed on the coarse grid.
As this example demonstrates, our IR enables the assembly of the computational graph in a stepwise manner.
However, to automate the process of generating a multigrid method from its grammar-based representation, we must be able to translate any derivation of the grammar into a corresponding IR object.
For this purpose, we utilize the functionality of the evolutionary computation framework DEAP~\cite{rainville2012deap}.
However, before we proceed with this task, we want to address some final remarks about the IR presented in this section.
The main purpose of the implementation presented here is to uniquely represent the computational structure of a multigrid method in the form of a redundancy-free directed graph.
Therefore, to enable the construction of such a graph in a step-wise manner based on the formal system described in Chapter~\ref{chapter:multigrid-formal-language}, each \emph{Cycle} node needs to include additional references to the current approximate solution and right-hand side.
While this information is important for the construction of the graph, it could later be discarded by replacing each such node with the respective arithmetic expression for computing an updated approximate solution, as it is shown in the \textsc{update} function in Section~\ref{sec:multigrid-state-transitions}.
However, preserving these additional references has the advantage of being able to easily traverse the sequence of \emph{Cycle} objects within each graph.
This not only lets us determine the computational structure of a multigrid method by traversing the corresponding graph data structure but also facilitates the identification of potential errors in the implementation.
Consequently, we represent each newly computed approximate solution as a \emph{Cycle} object, which is then referenced in subsequent computational steps of the method.
Note that even though the translation of a graph-based representation to an algorithmic one later requires us to transform each of these objects into an expression for updating the approximate solution, this operation can be performed while traversing the graph and thus does not induce a significant overhead within the process of algorithm generation.

\section{Grammar Generation}
According to Section~\ref{sec:multigrid-grammar}, our family of context-free grammars (CFGs) for the generation of multigrid methods consists of three components, a set of terminals, variables, and productions, while additionally, we have to choose a starting symbol $\ps{S}$ from the set of variables.
In Algorithm~\ref{table:grammar-semantics}, we have defined the semantics of each state transition function occurring within the productions listed in Algorithm~\ref{table:multigrid-grammar}.
While each instance of our grammar has to be formulated on a specific grid hierarchy, we have already mentioned that it is possible to define a structurally-equivalent grammar on a different hierarchy of discretizations with the same number of coarsening steps.
By treating the number of coarsening steps as a parameter, we can hence automate the process of grammar generation for different problems and discretizations.
For this purpose, we first need to generate the set of terminals that is defined on each level of the hierarchy, which we encapsulate in the class \emph{Terminals} shown in Listing~\ref{code:grammar:terminals}.
\begin{listing}
\begin{minted}{python}
from evostencils.ir import partitioning

class Terminals:
    def __init__(self, x_h, A_h, A_2h, restriction_operators, prolongation_operators, CGS_2h, relaxation_factor_interval, partitionings=None):
        self.A_h = A_h
        self.A_2h = A_2h
        self.x_h = x_h
        self.prolongation_operators = prolongation_operators
        self.restriction_operators = restriction_operators
        self.CGS_2h = CGS_2h
        self.relaxation_factor_interval = relaxation_factor_interval
        self.no_partitioning = partitioning.Single
        self.partitionings = partitionings

    @property
    def grid(self):
        return self.A_h.grid

    @property
    def coarse_grid(self):
        return self.A_2h.grid
\end{minted}
	\caption{Terminals Data Structure}
	\label{code:grammar:terminals}
\end{listing}
Note that this class comprises a few notable differences compared to our grammar formulation in Section~\ref{sec:multigrid-grammar}.
First of all, since each smoother considered in this work is derived from the system operator, it can be generated automatically within the grammar and, thus, does not need to be included explicitly as a terminal.
Also, while so far we have abstractly represented the application of the coarse-grid solver in the form of its multiplication with the inverse $A^{-1}_H$ on the coarsest level with spacing $H$, the coarse-grid solver itself can also be considered as a degree of freedom, and hence might be provided by the user.
The implementation of a Python class for representing the coarse-grid solver can be found in Listing~\ref{code:ir:coarse-grid-solver}.
\begin{listing}
\begin{minted}{python}
class CoarseGridSolver(Entity):
    def __init__(self, name, operator, expression=None):
        shape = operator.shape
        self._operator = operator
        self._expression = expression
        super().__init__(name, operator.grid, shape)

    @property
    def operator(self):
        return self._operator

    @property
    def expression(self):
        return self._expression
\end{minted}
	\caption{IR -- Coarse-Grid Solver}
	\label{code:ir:coarse-grid-solver}
\end{listing}
Note that each object of this class may additionally contain an expression that represents a complete multigrid method whose finest grid coincides with the coarsest grid of the original method.
This enables the construction of multigrid methods in a hierarchical manner, which means that after obtaining a multigrid method on a certain hierarchy of discretizations, we can employ it as a coarse-grid solver for another multigrid method formulated on top of our original discretization hierarchy.

\subsection{State Transition Functions}
\label{sec:evostencils:state-transition-functions}
As a second step, based on the \emph{Terminals} class, we can implement the state transition functions defined in Section~\ref{sec:multigrid-state-transitions} to construct the IR object of a particular grammar derivation tree.
Listing~\ref{code:grammar:basic-functions} contains the implementation of each of the five state transition functions defined in Algorithm~\ref{table:grammar-semantics}.
\begin{listing}
\begin{minted}{python}
def residual(state):
    x_h, b_h = state
    return base.Cycle(x_h, b_h, base.Residual(terminals.A_h, approximation, rhs), predecessor=approximation.predecessor)

def apply(A_h, cycle):
    cycle.correction = base.Multiplication(A_h, cycle.correction)
    return cycle

def update(relaxation_factor_index, partitioning, cycle):
    relaxation_factor = terminals.relaxation_factor_interval[relaxation_factor_index]
    cycle.relaxation_factor = relaxation_factor
    cycle.partitioning = partitioning
    x_h, b_h = cycle, cycle.rhs
    return x_h, b_h

def coarsening(A_2h, x_2h, cycle):
    r_2h = base.Residual(A_2h, x_2h, cycle.correction)
    new_cycle = base.Cycle(x_2h, cycle.correction, r_2h)
    new_cycle.predecessor = cycle
    return new_cycle

def coarse_grid_correction(P_2h, state):
    cycle = state[0]
    correction = base.Multiplication(P_2h, cycle)
    cycle.predecessor.correction = correction
    return cycle.predecessor
\end{minted}
	\caption{Basic State Transition Functions}
	\label{code:grammar:basic-functions}
\end{listing}
While each of these implementations is semantically equivalent to the corresponding function definition, certain adaptions are required due to the properties of the \emph{Cycle} class.
In Section~\ref{sec:intermediate-representation}, we have already discussed the advantages of storing the current state of a multigrid method directly within each \emph{Cycle} object.
As a consequence, the application of each state transition function either alters a given \emph{Cycle} object or returns a new object of this type.
The \emph{residual} function, however, represents an exception to this rule since it expects a \emph{state} variable containing an \emph{Approximation} and \emph{RightHandSide} object as its argument.
Note that according to Algorithm~\ref{table:multigrid-grammar}, every derivation ends with computing the residual based on the initial approximate solution $x_h^0$ and right-hand side $b_h$, given in the form of an \emph{Approximation} and \emph{RightHandSide} object, respectively.
In this case, we, therefore, need to pass the initial state 
\begin{equation*}
	Z_h^0 = (x_h^0, b_h, \lambda, \lambda)
\end{equation*} explicitly to the \emph{residual} function.
Since the third and fourth entry of this state is empty, in Python, a binary \emph{tuple} is sufficient to store the respective \emph{Approximation} and \emph{RightHandSide} object.
Even though in all subsequent computations, the first entry of this \emph{tuple} then consists of a \emph{Cycle} object, for the sake of simplicity, the right-hand side is always included as a second entry.
This decision allows us to employ the same \emph{residual} function within all grammar productions.
As a consequence, also the \emph{update} function needs to be adapted accordingly, such that the respective binary state, in the form of a \emph{Cycle} object for computing a new approximate solution, and the current right-hand side, is returned.
Finally, since the \emph{state} argument of the \emph{coarse\_grid\_correction} function results from an application of the \emph{update} function, it also represents a binary \emph{tuple}.
We, therefore, need to extract the first entry of this \emph{tuple} to obtain the corresponding \emph{Cycle} object.
The second main difference is in the implementation of the \emph{update} function compared to its original definition in Algorithm~\ref{table:grammar-semantics} function.
Here we represent the relaxation factor as an index within a uniformly-sampled interval that is included in the respective \emph{Terminal} object.
While we could explicitly store the relaxation factor as a floating point number, its representation accuracy then depends on the underlying floating point format.
Assume we want to encode a certain derivation tree in a string-based format.
This format then later needs to be decoded in a different environment to restore the original information.
If each relaxation factor is stored as a floating point number, we need to ensure that each number is stored with the same accuracy in both environments.
In contrast, an index can always be accurately represented in the form of a single positive integer value.

While Listing~\ref{code:grammar:basic-functions} provides us with the basic functionality to generate an IR representation of arbitrarily-structured multigrid methods, the majority of the productions defined in Algorithm~\ref{table:multigrid-grammar} consists of a combination of two different state transition functions.
For instance, smoothing is performed through the consecutive application of \textsc{apply} and \textsc{update}.
We can simplify the process of grammar generation by identifying each possible combination of state transitions, each of which can then be implemented in a single Python function.
In Definition~\ref{def:elementary-multigrid-operations}, we have already identified the three elementary multigrid operations \emph{smoothing}, \emph{coarsening}, and \emph{coarse-grid correction}.
In Algorithm~\ref{table:multigrid-grammar}, each of these three operations is defined as a combination of two state transition functions.
First, consider Production~\eqref{prod:smoothing}, which is similarly defined on each level and corresponds to the \emph{smoothing} operation in Definition~\ref{def:elementary-multigrid-operations}.
Each of the resulting productions corresponds to the application of an operator
\begin{equation*}
	B_{H} = \left( M_{H} \right)^{-1},
\end{equation*}
where $M_{H}$ is defined based on the splitting $A_{H} = M_{H} + N_{H}$ on a grid with spacing $H$.
While in Algorithm~\ref{table:multigrid-grammar} $M_{H}$ is provided as a terminal symbol, in practice, it can usually be derived from the system operator $A_{H}$.
Therefore, we can implement its generation within the function \emph{generate\_smoother}, which returns a similar splitting for each operator provided as an argument.
The generation of the actual smoothing expression is then performed in the function \emph{smoothing}, which is shown in Listing~\ref{code:grammar:smoothing}.
\begin{listing}
\begin{minted}{python}
def smoothing(relaxation_factor_index, partitioning, generate_smoother, cycle):
    assert isinstance(cycle.correction, base.Residual), 'Invalid production: expected residual'
    smoothing_operator = generate_smoother(cycle.correction.A_h)
    cycle = apply(base.Inverse(smoothing_operator), cycle)
    return update(relaxation_factor_index, partitioning, cycle)
\end{minted}
	\caption{State Transition -- Smoothing}
	\label{code:grammar:smoothing}
\end{listing}
Similar to Production~\eqref{prod:smoothing}, the implementation of this function consists of a combination of \textsc{apply} and \textsc{update}, whereby we generate $\ps{B_{h}}$ using the function \emph{generate\_smoother}.
To realize different smoothers, we thus only need to provide a generator function for each of them.
For instance, Listing~\ref{code:grammar:jacobi} shows how a Jacobi-based smoother can be implemented.
\begin{listing}
\begin{minted}{python}
def generate_jacobi_operator(operator: base.Operator):
    return base.Diagonal(operator)

def jacobi(relaxation_factor_index, partitioning, cycle):
    return smoothing(relaxation_factor_index, partitioning, generate_jacobi_operator, cycle)
\end{minted}
	\caption{Example for Generating Jacobi-Based Smoothers}
	\label{code:grammar:jacobi}
\end{listing}
The second basic multigrid operation \emph{coarsening} is realized in a similar way by combining the functions \textsc{apply} and \textsc{coarsening}.
After the former applies a restriction operator to the previously generated residual expression, the latter creates a new \emph{Cycle} object on the next lower level using the restricted residual as a right-hand side.
Next, we need to provide an implementation for Production~\eqref{prod:coarse-grid-correction}, which corresponds to the \emph{coarse-grid correction} operation in Definition~\ref{def:elementary-multigrid-operations}.
Note that in contrast to the state transition function \textsc{cgc}, this operation additionally updates the current approximate solution with the computed correction.
We can realize this behavior by additionally applying the \textsc{update} function, which generates the respective expression based on the given correction term.
Listing~\ref{code:grammar:inter-grid-operations} shows the implementation of Production~\eqref{prod:coarsening} and~\ref{prod:coarse-grid-correction} in the form of the Python functions \emph{restrict\_and\_coarsen} and \emph{update\_with\_coarse\_grid\_correction}.
\begin{listing}
\begin{minted}{python}
def restrict_and_coarsen(A_2h, x_2h, R_h, cycle):
    cycle = apply(R_h, cycle)
    return coarsening(A_2h, x_2h, cycle)

def update_with_coarse_grid_correction(relaxation_factor_index, P_2h, state):
    cycle = coarse_grid_correction(P_2h, state)
    return update(relaxation_factor_index, terminals.no_partitioning, cycle)
\end{minted}
	\caption{State Transition -- Inter-Grid Operations}
	\label{code:grammar:inter-grid-operations}
\end{listing}
In the case of the latter, we have included the prefix \emph{update\_with} to make its implementation distinguishable from the respective state transition function.
Finally, in contrast to the functions implemented so far, which can be applied on multiple levels, on the coarsest level, the only possible operation is the application of the coarse-grid solver, which corresponds to the Productions~\eqref{prod:coarse-grid-solver} and~\eqref{prod:coarse-grid-solver-correction}.
Here, Production~\eqref{prod:coarse-grid-solver} refers to the construction of the coarse problem, similar to the coarsening step defined in Production~\eqref{prod:coarsening}, while Production~\eqref{prod:coarse-grid-solver-correction} performs the actual correction step based on the exact solution obtained on the coarsest grid. 
However, note that only a single production is available for the variable $\ps{c_{16h}}$ in Algorithm~\ref{table:multigrid-grammar}, which means that the two productions are always applied in succession.
We can, therefore, combine the complete process of updating the current approximate solution based on the coarse-grid solver in a single Python function \emph{update\_with\_coarse\_grid\_solver}, whose implementation is shown in Listing~\ref{code:grammar:coarse-grid-solver}.
\begin{listing}
\begin{minted}{python}
def update_with_coarse_grid_solver(relaxation_factor_index, P_2h, CGS_2h, R_h, cycle):
    cycle = apply(R_h, cycle)
    cycle = apply(CGS_2h, cycle)
    cycle = apply(P_2h, cycle)
    return update(relaxation_factor_index, terminals.no_partitioning, cycle)
\end{minted}
	\caption{State Transition -- Coarse-Grid Solver}
	\label{code:grammar:coarse-grid-solver}
\end{listing}
For this purpose, similar to the \emph{restrict\_and\_coarsen} function, we first restrict the correction term of a given \texttt{Cycle} object.
The resulting error equation is then solved directly, which is denoted by the application of the coarse-grid solver.
As a final step, the obtained solution is transferred to the next higher level to correct the current approximate solution.
%TODO mention that EvoStencils includes additional functionality not mentioned here
\subsection{Genetic Programming in DEAP}
\label{sec:evostencils-part1:productions}
Finally, to generate the actual grammar, we need to assemble the respective subexpressions for each of the productions contained in Algorithm~\ref{table:multigrid-grammar} based on the terminal and state transition function implementation presented in the last section.
As already mentioned at the beginning of this chapter, we utilize the genetic programming (GP) module of the DEAP framework~\cite{rainville2012deap}. 
In principle, DEAP only offers support for untyped and strongly-typed tree-based GP and, therefore, does not allow implementing grammar-guided GP (G3P) directly.
However, as we have already discussed in Section~\ref{sec:gggp-representation}, G3P can be considered as a special variant of strongly-typed GP, whereby each variable that is placed on the left-hand side of a production encodes a unique type.
Before we discuss how the productions in Algorithm~\ref{table:multigrid-grammar} can be mapped to unique types, we need to introduce the relevant constructs already implemented in the framework.
In DEAP, the main data structure to represent a typed GP system is the class \emph{PrimitiveSetTyped}, which defines how a program can be constructed based on a set of \emph{Primitive} objects.
As each operation must adhere to these rules, it is ensured that only individuals fulfilling the specified type constraints are generated.
To demonstrate how a grammar can be implemented in the form of a \emph{PrimitiveSetTyped}, we consider the example grammar from Section~\ref{sec:gggp-representation}, whose productions are given by
\begin{equation*}
	\begin{split}
		\ps S \; \; \bnfpo & \; \; \ps E \\
		\ps E \; \; \bnfpo & \; \; \text{if} \; \ps B \; \text{then} \; \ps E \; \text{else} \; \ps E \; | \; \ps A \\
		\ps A \; \; \bnfpo & \; \; -\ps A \; | \; (\ps A + \ps A) \; | \; (\ps A - \ps A) \; | \\
		& \; \; (\ps A \cdot \ps A) \; | \; (\ps A / \ps A) \; | \ps{A}^{\ps{A}} \; | \; x \; | \; y \\  
		\ps B \; \; \bnfpo & \; \;  \neg \ps B \; | \; (\ps B \wedge \ps B) \; | \; (\ps B \vee \ps B) \; | \; u \; | \; v.
	\end{split}
\end{equation*}
Listing~\ref{code:grammar:pset-example} shows the corresponding implementation of the function \emph{generate\_grammar}, which generates a \emph{PrimitiveSetTyped} object.
\begin{listing}
\begin{minted}{python}
from deap import gp
import operator

def generate_grammar(x, y, u, v):
    E = type("E", (object,), {})
    A = type("A", (object,), {})
    B = type("B", (object,), {})
    pset = gp.PrimitiveSetTyped("main", [], ret_type=E)

    # <E> := if <B> then <E> else <E> | <A>
    pset.addPrimitive(lambda b, e1, e2: e1 if b else e2, [B, E, E], E, name="if_then_else")
    pset.addPrimitive(lambda _: _, [A], E, name="id")
    # <A> :=
    pset.addPrimitive(operator.neg, [A], A)         # -<A> |
    pset.addPrimitive(operator.add, [A, A], A)      # <A> + <A> |
    pset.addPrimitive(operator.sub, [A, A], A)      # <A> - <A> |
    pset.addPrimitive(operator.mul, [A, A], A)      # <A> * <A> |
    pset.addPrimitive(operator.truediv, [A, A], A)  # <A> / <A> |
    pset.addPrimitive(operator.pow, [A, A], A)      # <A>^<A>
    # <B> :=
    pset.addPrimitive(operator.not_, [B], B)        # not <B> |
    pset.addPrimitive(operator.and_, [B, B], B)     # <B> and <B> |
    pset.addPrimitive(operator.or_, [B, B], B)      # <B> or <B>

    def add_argument(pset, arg, type_list):
        symbolic = True
        for t in type_list:
            pset._add(gp.Terminal(arg, symbolic, t))
            pset.terms_count += 1
        pset.arguments.append(arg)

    add_argument(pset, x, [E, A])
    add_argument(pset, y, [E, A])
    add_argument(pset, u, [B])
    add_argument(pset, v, [B])
    return pset

pset = generate_grammar("x", "y", "u", "v")
expr = gp.genGrow(pset, 3, 5)
tree = gp.PrimitiveTree(expr)
f = gp.compile(tree, pset)
print(tree)
print(f(42, 7, False, True))
\end{minted}
	\caption{Example Grammar Generation with DEAP}
	\label{code:grammar:pset-example}
\end{listing}
Here, the first step is to create a unique type for each symbol that is contained in the set of variables $V = \left\{\ps{E}, \ps{A}, \ps{B} \right\}$, which can be accomplished using Python's builtin \emph{type} function as shown in line~4--6 of Listing~\ref{code:grammar:pset-example}.
Next, a new \emph{PrimitiveSetTyped} object is created in line~7, in which we set the return type to $\ps{E}$, similar to the choice of the start variable within the grammar.
In line~9--21, we then proceed with defining each production in the form of a \emph{Primitive} objective, which consists of a function, a list of input types, and an output type.
Counterintuitively, the output and not the input type defines which variable is placed on the left-hand side of each production.
To understand this contradiction, we need to revisit the process of tree initialization in G3P.
As we have seen in Section~\ref{sec:gggp-initialization}, a new derivation tree is generated starting with the variable $\ps{S}$ by recursively selecting productions for each leaf node of the tree that corresponds to a variable.
Similarly, to generate a tree based on a given \emph{PrimitiveSetTyped}, a \emph{Primitive} is picked randomly from those whose output type matches the specified return type.
After extending the tree accordingly, the process is continued with each of the input types of the chosen \emph{Primitive}.
Therefore, terminating this process at a certain point within the tree requires choosing a \emph{Primitive} with an empty list of input types, which corresponds to a production whose right-hand side does not contain any variables.
In DEAP, such a \emph{Primitive} is called a \emph{Terminal}, which should not be confused with the terminals introduced in Section~\ref{sec:gggp-representation} that simply refer to each non-variable symbol of a grammar.
To add either a \emph{Primitive} or \emph{Terminal} object to a given \emph{PrimitiveSetTyped}, the \emph{\_add} method can be used.
Additionally, the methods \emph{addPrimitive} and \emph{addTerminal} are also available, which are more convenient to use.
In order to make the productions of our example grammar available to the created \emph{PrimitiveSetTyped} object, we hence include a \emph{Primitive} for each of them, using the aforementioned function and the previously defined types.
Note that the production
\begin{equation*}
	\ps E \; \; \bnfpo \; \; \ps A
\end{equation*}
is implemented based on the identity function defined in line~10 of Listing~\ref{code:grammar:pset-example}.
After adding a \emph{Primitive} object for each of the grammar's productions, we need to take care of the four symbols $x$, $y$, $u$, and $v$.
We can represent these symbols as objects of the \emph{Terminal} class, which corresponds to a \emph{Primitive} object without any arguments, and thus an empty list of input types.
Optionally, a \emph{Terminal} object may also refer to a Python symbol, which can be triggered by setting the \emph{symbolic} argument accordingly.
If we assume that each of the four symbols represents an argument to a Python function generated by our grammar, we additionally have to append it to the list of arguments of the corresponding \emph{PrimitiveSetTyped} object.
For this purpose, we define the local function \emph{add\_argument} in line~22--27, which is then utilized in line~28--31 to include each of the four symbols as an argument. 
Note that for each of the two symbols $x$ and $y$, we create two \emph{Terminal} objects that only differ in their return type.
As a consequence, our implementation includes two additional productions
\begin{equation*}
	\ps E \; \; \bnfpo \; \; \ps x \bnfor \ps y,
\end{equation*}
which are not part of the original grammar.
The reason for this adaption is that DEAP's implementation of the \emph{grow} operator, as described in Section~\ref{sec:gggp-initialization}, expects the availability of at least one \emph{Terminal} and \emph{Primitive} object for each type within a \emph{PrimitiveSetTyped}.
In the given case, we only need to include additional \emph{Terminal} objects for the type $\ps{E}$, as the condition is already fulfilled for all other types.
While the example considered here does not require significant adjustments to the structure of the grammar as the type $\ps{E}$ can already be converted to $\ps{E}$ using the identity function, in general, this is not the case.
In particular, according to Algorithm~\ref{table:multigrid-grammar}, for the majority of the variables of our multigrid grammar, only non-terminal productions, i.e., productions that generate strings with at least one variable, are available.
To handle grammars that violate the requirement of having at least one terminal production for each of its variables, we need to adapt DEAP's implementation of the \emph{grow} operator.
The details of this adaptation will be discussed in Section~\ref{sec:evostencils-part1:evolutionary-program-synthesis}, where we will present the implementation of our evolutionary program synthesis method.
Finally, after constructing a \emph{PrimitiveSetTyped} object that corresponds to our example grammar, we can construct a random derivation tree using the \emph{genGrow} function, which corresponds to the aforementioned \emph{grow} initialization operator.
Based on the resulting tree, a function object is then generated using DEAP's \emph{compile} function, which can be executed similarly to any other Python function by providing a value for each of its arguments.
%Similar to \emph{Primitive} objects, a \emph{Terminal} can be added to an existing \emph{PrimitiveSetTyped} using the \emph{addTerminal} method.

\subsection{Variable Encoding}
After introducing the functionality of DEAP's GP module, we can proceed with the actual implementation of the productions of our multigrid grammar, as defined in Algorithm~\ref{table:multigrid-grammar}.
However, before we define a function for constructing the corresponding \emph{PrimitiveSetTyped}, we need to consider the unique structure of our grammar.
As we have discussed in Section~\ref{sec:multigrid-grammar}, Algorithm~\ref{table:multigrid-grammar} is obtained by repeating the same productions on each level for the corresponding set of variables.
In Listing~\ref{code:grammar:terminals}, we have already implemented a data structure that includes all required terminals for a particular level of the discretization hierarchy.
Therefore, instead of initializing a \emph{PrimitiveSetTyped} with the complete set of productions, we can instead add the respective \emph{Primitive} objects on a per-level basis.
Furthermore, note that Algorithm~\ref{table:multigrid-grammar} is defined for a specific number of coarsening steps.
We can thus implement a function that iteratively constructs a multigrid grammar independent of the total number of coarsening steps.
This allows us to deploy the same function for the generation of multigrid grammars defined on different discretization hierarchies, for instance, with a different number of coarsening steps.
As in the above example, we first need to generate a unique type for each variable that is defined on a certain level of our grammar.
While in Listing~\ref{code:grammar:pset-example}, we have encoded each variable with an actual Python type, at this point, we have to introduce an additional constraint that prevents us from pursuing the same approach.
In order to store the state of a program in the form of a serialized binary format, Python provides the \emph{Pickle} module.
This module enables the serialization of arbitrary Python objects, which can then later be restored.
Object serialization has multiple purposes in our implementation.
First of all, we want to be able to store the current state of our evolutionary algorithm such that, for instance, in case of a hardware failure, we can continue at the same position.
Furthermore, parallelizing certain parts of our implementation on a distributed computing system using the message-passing interface (MPI) requires us to transfer arbitrary Python objects via a communication network.
The Pickle module provides a simple and portable solution to this problem.
Unfortunately, Pickle enforces a number of additional constraints on the serializability of a Python object.
In particular, Python types can only be serialized if they are defined at the top-level domain of a module.
Since our goal is to create a \emph{PrimitiveSetTyped} object adapted to the properties of a grammar that is only provided during program execution, we need to create the respective types dynamically.
Unfortunately, dynamic types can not be defined at the top-level domain of the corresponding Python module, which renders their serialization impossible.
To resolve this issue, we propose a custom \emph{Type} class which is defined in Listing~\ref{code:grammar:typing}.
\begin{listing}
\begin{minted}{python}
class Type:
    def __init__(self, identifier, guard=False):
        self.identifier = identifier
        self.guard = guard

    def __eq__(self, other):
        if isinstance(other, type(self)):
            return self.identifier == other.identifier and self.guard == other.guard
        else:
            return False

    def __hash__(self):
        return hash((self.identifier, self.guard))
\end{minted}
	\caption{Variable Encoding}
	\label{code:grammar:typing}
\end{listing}
Since the instances of this class are regular Python objects, a Pickle-based serialization can be achieved without any further adaption.
To distinguish different \emph{Type} objects within a \emph{PrimitiveSetTyped}, the \emph{identifier} attribute is provided.
In addition, the class also includes a \emph{guard} attribute in the form of a boolean variable, whose purpose will be discussed later in this section. 
To enable type equivalence checking in an automatic manner, we provide a custom \emph{\_\_eq\_\_} method.
We first check whether the other object is also an instance of the \emph{Type} class, as we want to be able the compare a given \emph{Type} object to any other Python object.
If this condition is fulfilled, we then proceed to check if both the \emph{identifier} and \emph{guard} attributes of the \emph{Type} objects are equal.
In order to correctly store \emph{Type} objects in a Python dictionary, we implement the \emph{\_\_hash\_\_} method, where we utilize Python's builtin \emph{hash} function to generate a unique value based on the content of the two attributes. 
Finally, to ensure that the type checking is performed correctly, we adapt the original implementation of the \emph{PrimitiveSetTyped} class such that types are compared using the equality operator instead of Python's built-in \emph{issubclass} function.
The resulting implementation can be found in Section~\ref{appendix:gp} of the appendix.

Based on this tailored type representation, we can now implement a data structure that incorporates the type of each variable defined on a given level.
The implementation of this data structure in the form of the \emph{Types} class is shown in Listing~\ref{code:grammar:types}.
\begin{listing}
\begin{minted}{python}
class Types:
    @staticmethod
    def _init_type(identifier, types, type_attribute=None, guard=False):
        if type_attribute is None:
            type_attribute = identifier
        if types is None:
            return Type(identifier, guard)
        else:
            return getattr(types, type_attribute)

    def __init__(self, depth, previous_types=None):
        # Fine-grid types
        gen_id = lambda base: f"{base}_{depth}"
        self.S_h = self._init_type(gen_id("S"), previous_types, "S_2h")
        self.S_guard_h = self._init_type(gen_id("S_guard"), previous_types, "S_guard_2h", guard=True)
        self.C_h = self._init_type(gen_id("C"), previous_types, "C_2h")
        self.C_guard_h = self._init_type(gen_id("C_guard"), previous_types, "C_guard_2h", guard=True)
        self.x_h = self._init_type(gen_id("x"), previous_types, "x_2h")
        self.A_h = self._init_type(gen_id("A"), previous_types, "A_2h")
        self.B_h = self._init_type(gen_id("A"), previous_types, "B_2h")
        self.R_h = Type(f"R_{2 ** depth}h")

        # Coarse-grid types
        gen_id = lambda base: f"{base}_{depth+1}h"
        self.S_2h = Type(gen_id("S"))
        self.S_guard_2h = Type(gen_id("S_guard"), guard=True)
        self.C_2h = Type(gen_id("C"))
        self.C_guard_2h = Type(gen_id("C_guard"), guard=True)
        self.x_2h = Type(gen_id("x"))
        self.A_2h = Type(gen_id("A"))
        self.B_2h = Type(gen_id("B"))
        self.P_2h = Type(gen_id("P"))
        self.CGS_2h = Type(gen_id("CGC"))

        # General types
        self.Partitioning = self._init_type("Partitioning", previous_types)
        self.RelaxationFactorIndex = self._init_type("RelaxationFactorIndex", previous_types)
\end{minted}
  	\caption{Types Data Structure}
  	\label{code:grammar:types}
\end{listing}
%TODO include line numbers at the respective text positions
To initialize the type that corresponds to a certain variable of the grammar, we can either create a new \emph{Type} object based on an \emph{identifier} or retrieve the respective object if it is already contained in the \emph{Types} object that refers to the next higher level in the discretization hierarchy.
The resulting initialization procedure is implemented in the \emph{\_init\_type} method, which can be found in line~2--9 of Listing~\ref{code:grammar:types}.
By applying this function, we obtain a \emph{Type} object with unique \emph{identifier} for each variable defined on two subsequent levels of the discretization hierarchy, whereby each attribute with a subscript $h$ and $2h$ corresponds to a variable on the current and next coarser level, respectively.
Consequently, if we provide an existing \emph{Types} object for the initialization of its predecessor on the next lower level, the coarse-grid types correspond to the fine-grid types of this object.
We, therefore, only have to generate a unique type for each variable once, which can then be transferred to a lower level in the form of a reference to the respective \emph{Type} object.
For the sake of simplicity, we identify each \emph{Type} object with a Python \emph{string}. 
The resulting generation of each \emph{Type} object and its assignment to the respective attribute is performed in line~12--31.
In contrast to the types that need to be defined on each level, both the \emph{Partitioning} and \emph{RelaxationFactorIndex} attribute is level-independent, and thus only a single type needs to be created for each of them, which is shown in line~33 and~34.

\subsection{Productions}
With the implementation of a type system that allows expressing each production as a mapping between input and output types, we are now finally at the point where we can bring everything together.
Similar to the example shown in Listing~\ref{code:grammar:pset-example}, we first have to create a \emph{PrimitiveSetTyped} to which we can then add the respective \emph{Terminal} and \emph{Primitive} objects in an iterative manner.
The resulting implementation is shown in Listing~\ref{code:grammar:init-grammar}.
\begin{listing}
\begin{minted}{python}
import numpy as np
from evostencils.ir import partitioning
from evostencils.grammar.gp import PrimitiveSetTyped

def init_grammar(generator, x_h, b_h, max_level, samples, coarsest=False):
    A_h = generator.get_operator(max_level)
    A_2h = generator.get_operator(max_level - 1)
    restriction_operators, prolongation_operators = generator.get_inter_grid_operators(max_level, max_level - 1)
    CGS_2h = CoarseGridSolver("CGS", A_2h)
    relaxation_factor_interval = np.linspace(0.1, 1.9, samples)
    partitionings = [partitioning.RedBlack]

    terminals = Terminals(x_h, A_h, A_2h, restriction_operators, prolongation_operators, CGS_2h, relaxation_factor_interval, partitionings)
    types = Types(0)

    pset = PrimitiveSetTyped("main", [], types.S_h)
    pset.addTerminal((x_h, b_h), types.S_guard_h, 'initial_state')

    # Relaxation factors
    for i in range(0, samples):
        pset.addTerminal(i, types.RelaxationFactorIndex)

    # Partitionings
    pset.addTerminal(terminals.no_partitioning, types.Partitioning, terminals.no_partitioning.get_name())
    for p in terminals.partitionings:
        pset.addTerminal(p, types.Partitioning, p.get_name())

    add_level(pset, terminals, types, max_level, 0, samples, coarsest=coarsest)
    return pset, terminals, types
\end{minted}
	\caption{Grammar Initialization}
	\label{code:grammar:init-grammar}
\end{listing}
In line~5--12, we first construct the respective \emph{Terminals} and \emph{Types} data structures on the finest grid of the given discretization hierarchy.
Based on these data structures, we proceed with the creation of the \emph{PrimitiveSetTyped} object in line~13.
We then include the initial state as a \emph{tuple} consisting of the respective \emph{Approximation} and \emph{RightHandSide} object and all level-independent terminals in line~14--21.
For the sake of simplicity, we assume that it is possible to generate all required \emph{Approximation}, \emph{RightHandSide} and \emph{Operator} objects using a predefined \emph{Generator} class, which retrieves the required information automatically from a domain-specific representation of the problem.
All other \emph{Terminal} and \emph{Primitive} objects are then included in line~22 using the \emph{add\_level} function shown in Listing~\ref{code:grammar:add-level}.

At this point, we now need to come back to the definition of the \emph{Type} class, which includes an additional \emph{guard} attribute.
The purpose of this attribute is to enforce additional constraints that are not present in the productions in Algorithm~\ref{table:multigrid-grammar}.
Consider, for instance, the derivation
\begin{equation*}
\ps{S} \Rightarrow \ps{s_h} \Rightarrow (x_h^0, b_h, \bnfes, \bnfes),
\end{equation*}
which is valid according Algorithm~\ref{table:multigrid-grammar}.
However, the resulting method consists only of a single statement, which returns the initial approximate solutions.
It is therefore important to define the minimal requirements a sequence of operations needs to fulfill to qualify as a multigrid method. 
While the optimal amount of smoothing and coarse-grid corrections might differ for each case, we consider it essential to apply the coarse-grid solver at least once to obtain an accurate approximation of the solution of the error equation on the coarsest grid.
We can enforce this requirement with the mentioned \emph{guard} attribute of the \emph{Type} class.
As each derivation of our grammar ends with the generation of the initial state $\left(x_h^0, b_h\right)$ on the finest level, we have to prevent the application of the respective production until the coarse-grid solver has been utilized at least once.
For this purpose, we set the output type of the respective \emph{Terminal} object in line~14 of Listing~\ref{code:grammar:init-grammar} to the \emph{S\_guard\_h} attribute of the respective \emph{Types} data structure.
We then define all remaining productions in such a way that all possible derivations include the function \emph{update\_with\_coarse\_grid\_solver}, as shown in Listing~\ref{code:grammar:coarse-grid-solver}, at least once.
Listing~\ref{code:grammar:add-level} demonstrates how this can be achieved for each level within a discretization hierarchy of arbitrary depth.
\begin{listing}
\begin{minted}{python}
def add_level(pset, terminals, types, depth, coarsest=False):
    add_terminals(pset, terminals, types, depth, coarsest)

    def add_primitive(pset, f, fixed_types, input_types, output_types, name):
        for t1, t2 in zip(input_types, output_types):
            pset.addPrimitive(f, fixed_types + [t1], t2, name)

    # Productions
    add_primitive(pset, residual, [], [types.S_h, types.S_guard_h], [types.C_h, types.C_guard_h], f"residual_{depth}")
    add_primitive(pset, jacobi, [types.RelaxationFactorIndex, types.Partitioning], [types.C_h, types.C_guard_h], [types.S_h, types.S_guard_h], f"jacobi_{depth}")
    # Include additional smoothers here
    if not coarsest:
        add_primitive(pset, update_with_coarse_grid_correction, [types.RelaxationFactorIndex, types.P_2h], [types.S_2h, types.S_guard_2h], [types.S_h, types.S_guard_h], f"update_with_coarse_grid_correction_{depth}")
        add_primitive(pset, restrict_and_coarsen, [types.A_2h, types.x_2h, types.R_h], [types.C_h, types.C_guard_h], [types.C_2h, types.C_guard_2h], f"restrict_and_coarsen_{depth}")
    else:
        # Add transition C_guard_h <- S_h to enable derivation termination
        add_primitive(pset, update_with_coarse_grid_solver, [types.RelaxationFactorIndex, types.P_2h, types.CGS_2h, types.R_h], [types.C_h, types.C_guard_h], [types.S_h, types.S_h], f'update_with_coarse_grid_solver_{depth}')
        pset.addTerminal(terminals.CGS_2h, types.CGS_2h, f'CGS_{depth+1}')
\end{minted}
	\caption{Addition of Terminals and Primitives per Level}
	\label{code:grammar:add-level}
\end{listing}
In line~2, we first make use of the function \emph{add\_terminals}, whose implementation is shown in Listing~\ref{code:grammar:add-terminals}, to add all terminals that are required on the respective level to the given \emph{PrimitiveSetTyped}.
\begin{listing}
\begin{minted}{python}
def add_terminals(pset, terminals, types, depth, coarsest=False):
    if not coarsest:
        x_2h = Approximation(terminals.coarse_grid)
        pset.addTerminal(x_2h, types.x_2h, f'zero_{depth + 1}')
        pset.addTerminal(terminals.A_2h, types.A_2h, f'A_{depth + 1}')
    for P_2h in terminals.prolongation_operators:
        pset.addTerminal(P_2h, types.P_2h, P_2h.name)
    for R_h in terminals.restriction_operators:
        pset.addTerminal(R_h, types.R_h, R_h.name)
\end{minted}
	\caption{Addition of Terminals per Level}
	\label{code:grammar:add-terminals}
\end{listing}
To facilitate the process of generating  \emph{Primitive} objects, we implement the helper function \emph{add\_primitive} in line~3--5.
This function can generate multiple instances of the same \emph{Primitive} using different input and output types, which we utilize to include each type as a guarded and unguarded version.
Now note that each derivation starts with \emph{S\_h} and ends with the type \emph{S\_guard\_h}.
According to the productions in Algorithm~\ref{table:multigrid-grammar}, only four unique primitives are required on each level, except for the definition of different smoothers.
For each of the \emph{Primitive} objects generated in line~7--12 whose list of input types contains only unguarded types, the output type is also unguarded.
Consequently, so far, with none of the generated \emph{Primitive} objects, a transition from an unguarded to a guarded type is possible.
We can thus enforce that the coarse-grid solver is applied at least once within each derivation by including the necessary transition in line~15.
Since this transition represents the only way to reach the final state, the application of the \emph{update\_with\_coarse\_grid\_solver} function is a necessary requirement for the termination of each derivation.
Note that we only want to add this production on the lowest level, where it replaces the respective productions for the coarsening and coarse-grid correction steps.
We, therefore, have to check whether we have arrived at the coarsest grid of the discretization hierarchy, which is done in line~10 of Listing~\ref{code:grammar:add-level}.
Since we generate the respective \emph{Primitive} additionally with an unguarded input type, similar to most other productions, the coarse-grid solver can be applied multiple times within each derivation.

Finally, using the functions \emph{init\_grammar} and \emph{add\_level}, we can define the Python function \emph{generate\_grammar}, whose implementation is shown in Listing~\ref{code:grammar:generate-grammar}. 
This function generates a \emph{PrimitiveSetTyped} for discretization hierarchies of variable depth.
\begin{listing}
\begin{minted}{python}
import numpy as np

def generate_grammar(generator, max_level, depth, samples=37):
    coarsest = False
    if depth == 1:
        coarsest = True
    x_h = generator.get_approximation(max_level)
    b_h = generator.get_rhs(max_level)
    pset, terminals, types = init_grammar(x_h, b_h, max_level, samples, coarsest=coarsest)

    # Iteratively add productions per level
    for i in range(1, depth):
        level = max_level - i
        x_h = generator.get_approximation(level)
        A_h = terminals.A_2h
        if i == depth - 1:
            coarsest = True
        A_2h = generator.get_operator(level - 1)
        restriction_operators, prolongation_operators = generator.get_inter_grid_operators(level, level - 1)
        CGS_2h = CoarseGridSolver("CGS", A_2h)
        relaxation_factor_interval = np.linspace(0.1, 1.9, samples)
        terminals = Terminals(x_h, A_h, A_2h, restriction_operators, prolongation_operators, CGS_2h, relaxation_factor_interval, partitionings)
        previous_types = types
        types = Types(i, previous_types)
        add_level(pset, terminals, types, i, samples, coarsest=coarsest)

    return pset
\end{minted}
	\caption{Grammar Generation}
	\label{code:grammar:generate-grammar}
\end{listing}
%TODO maybe clarify the meaning of max_level in comparison to earlier sections
We first create an initial \emph{PrimitiveSetTyped} to which we add all required \emph{Terminal} and \emph{Primitive} objects on the finest grid using the previously defined function \emph{init\_grammar} in line~8.
In case the given discretization hierarchy consists of at least three levels, we iteratively extend our \emph{PrimitiveSetTyped} using the \emph{add\_level} function.
For this purpose, we only have to create the respective \emph{Terminals} and \emph{Types} data structures, as shown in line~20 and~22, which are then passed to \emph{add\_level} in line~23 to generate all required \emph{Primitive} and \emph{Terminal} objects.
After we have traversed all levels, the \emph{PrimitiveSetTyped} object returned at the end of the function precisely replicates the functionality of Algorithm~\ref{table:multigrid-grammar} formulated on the given hierarchy of discretizations.

\section{Evolutionary Program Synthesis}
\label{sec:evostencils-part1:evolutionary-program-synthesis}
In order to provide an accessible interface to EvoStencils' functionality, we implement all steps from the construction of the grammar, based on the given PDE-based problem, to the synthesis of an optimal program in the \emph{Optimizer} class, whose basic implementation is shown in Listing~\ref{code:optimization:optimizer}.
\begin{listing}
\begin{minted}{python}
class Optimizer:
    def __init__(self, program_generator):
        self.program_generator = program_generator
        self.pset = None
        self.toolbox = None
        self.params = None
        self.individual_cache = {}

    def run(self, params, use_random_search=False):
        self.params = params
        self.pset = generate_grammar(self.program_generator, params.max_level, params.depth, params.relaxation_factor_samples)
        min_level = params.max_level - params.depth
        self.program_generator.initialize_code_generation(min_level, params.max_level)
        self.init_toolbox(self.pset, params.tree_min_height, params.tree_max_height, params.multi_objective)
        if params.multi_objective:
            self.init_multi_objective_selection(self.pset, self.evaluate)
        else:
            self.init_single_objective_selection(self.pset, self.evaluate)
        final_population = self.evolutionary_search(params, use_random_search)
        return final_population
\end{minted}
	\caption{Optimizer Class}
	\label{code:optimization:optimizer}
\end{listing}
In addition to the generated \emph{PrimitiveSetTyped} object, this class relies on two main components, a \emph{Toolbox} object for the execution of the evolutionary program synthesis and a \emph{ProgramGenerator} object, which enables the generation and evaluation of solver implementations while also providing an interface to all required information about the underlying problem.
Based on the functionality of these two components, we can define the \emph{run} method of the \emph{Optimizer} class.
The implementation of this class is shown in line~8--19 of Listing~\ref{code:optimization:optimizer}, while we additionally define a data structure for storing its parameters in Listing~\ref{code:optimization:parameters}.
\begin{listing}
\begin{minted}{python}
class Parameters:
    def __init__(self, max_level, depth, generations, mu_, lambda_):
        self.max_level = max_level
        self.depth = depth
        self.generations = generations
        self.mu_ = mu_
        self.lambda_ = lambda_
        self.initial_population_size = 4 * mu_
        self.crossover_probability = 0.7
        self.relaxation_factor_samples = 37
        self.tree_min_height = 10
        self.tree_max_height = 50
        self.multi_objective = True
\end{minted}
	\caption{Parameters of the Optimizer Class}
	\label{code:optimization:parameters}
\end{listing}
This method first generates a \emph{PrimitiveSetTyped} and then, based on the initialized \emph{Toolbox} object, performs an evolutionary search that yields a population of individuals representing efficient multigrid methods for the given problem.
Note that within this procedure, the evaluation of each generated solver relies on an \emph{evaluate} method, which returns suitable metrics for assessing its quality.

After translating our class of multigrid grammars introduced in Section~\ref{sec:multigrid-grammar} into a strongly-typed GP system, we can now make use of the functionality of DEAP's GP module to implement the evolutionary search operators introduced in Section~\ref{sec:gggp}.
However, before we can define the individual operators, we first need to implement a data structure for storing each derivation tree generated based on the rules of our grammar.
For this purpose, DEAP provides the \emph{PrimitiveTree} class, which represents each tree as a list of \emph{Primitive} and \emph{Terminal} objects stored in depth-first order.
Each entry of this list corresponds to the choice of a particular production in the form of a \emph{Primitive} or \emph{Terminal} object, while its return type refers to the respective variable on the left-hand side.
Therefore, each entry in this list refers to a certain variable and its children within the derivation tree.
After obtaining a suitable data structure for the representation of each derivation tree generated by our grammar, the next step is the randomized creation of an initial population.
As we have already discussed in Section~\ref{sec:gggp-initialization}, the \emph{full} initialization operator~\cite{koza1994genetic,poli2008field} is only applicable in cases where a grammar allows the termination of unfinished branches at any given point in the derivation by invoking a terminal production.
However, according to Algorithm~\ref{table:multigrid-grammar}, not all variables of our grammar fulfill this condition.
For instance, the variable $\ps{c_h}$ always yields a residual expression containing $\ps{s_h}$.
In contrast, the \emph{grow} operator is a suitable initialization operator in the given case, as it supports the generation of trees with branches of variable length. 
The DEAP framework implements the generation of random trees in the \emph{generate} function, which constructs a tree based on the provided \emph{PrimitiveSetTyped} object while enforcing that the length of each of its branches satisfies a certain depth limit.
If a certain branch satisfies this limit, this function always tries to pick a \emph{Terminal} object among the available productions to finish the growth at the respective branch. 
Note that these productions are given in the form of a list of \emph{Primitive} and \emph{Terminal} objects whose output types match with the input type of the node from which the current branch should be grown further.
Unfortunately, the specified limit is applied in a strict sense, and thus, the \emph{generate} function fails if a \emph{Terminal} object is erroneously expected to be available for a certain type.
To resolve this issue, we, therefore, need to adapt DEAP's implementation of the \emph{generate} function, such that the specified limit is only applied whenever possible.
In order to circumvent cases where the constraint of choosing only \emph{Terminal} objects can not be satisfied, we simply ignore it in the current step, which means that we extend the given node with a random \emph{Primitive}.
We then proceed with the current branch while choosing a \emph{Terminal} object whenever possible until the growth ends.
As an additional option, we also allow the insertion of an existing branch into the generated tree, which corresponds to the \emph{subtree insertion} operator introduced in Section~\ref{sec:gggp-mutation-and-recombination}.
Based on the \emph{generate} function, we can finally implement the \emph{grow} operator in the form of the function \emph{genGrow}, which chooses the maximum depth of a tree randomly from an interval defined by the parameters \emph{min\_height} and \emph{max\_height}.
The resulting implementation is shown in Listing~\ref{code:gp:grow}, while the \emph{generate} function can be found in Section~\ref{appendix:gp} of the appendix.
\begin{listing}
\begin{minted}{python}
def genGrow(pset, min_height, max_height, type_=None, size_limit=150):
    def condition(height, depth):
        return depth < height
    result = generate(pset, min_height, max_height, condition, type_)
    while len(result) > size_limit:
        # Repeat the process until the generated individual 
        # does fulfill the specified limits
        result = generate(pset, min_height, max_height, condition, type_)
    return result
\end{minted}
	\caption{Grow Operator for Tree-Based Genetic Programming}
	\label{code:gp:grow}
\end{listing}
In order to prevent the generation of excessively large expressions, the \emph{genGrow} function additionally enforces a size limit, which is set to 150 by default.
Therefore, if a generated tree does not fulfill this constraint, we repeat the process, again starting from scratch, until the requirement is met.

\subsection{Evolutionary Operators}
Based on our implementation of the \emph{generate} function and DEAP's built-in functionality, we can now define all remaining components of an evolutionary program synthesis method.
To facilitate the utilization of different operator choices within the context of an evolutionary algorithm, DEAP provides the \emph{Toolbox} container for the creation of unified interfaces.
Listing~\ref{code:gp:toolbox} demonstrates how we generate such a container for the given case.
\begin{listing}
\begin{minted}{python}
from deap import base, gp, creator
from evostencils.gp import genGrow, mutInsert

def init_toolbox(self, pset, min_height, max_height, multi_objective=True):
    # Define an Individual together with its Fitness
    if multi_objective:
        creator.create("Fitness", base.Fitness, weights=(-1.0, -1.0))
    else:
        creator.create("Fitness", base.Fitness, weights=(-1.0,))
    creator.create("Individual", gp.PrimitiveTree, fitness=creator.Fitness)

    self.toolbox = deap.base.Toolbox()
    # Population initialization
    self.toolbox.register("generate_tree", genGrow, pset=pset, min_height=min_height, max_height=max_height)
    self.toolbox.register("generate_individual", tools.initIterate, creator.Individual, self.toolbox.generate_tree)
    self.toolbox.register("population", tools.initRepeat, list, self.toolbox.generate_individual)
    
    # Crossover and mutation
    self.toolbox.register("mate", gp.cxOnePoint)
    self.toolbox.register("mutate", mutateSubtree, pset_=pset, min_height=1, max_height=max_height / 2)
\end{minted}
	\caption{Toolbox Initialization}
	\label{code:gp:toolbox}
\end{listing}
First of all, we define a data structure for the fitness of each individual based on DEAP's built-in \emph{Fitness} class.
Here the \emph{weights} attribute determines the number of objectives, while a negative value means that the respective objective is to be minimized.
As we will later come back to the definition of the optimization objectives, at this point, the reader can just assume that either a single or a multi-objective minimization is performed.
Next, we represent each \emph{Individual} as a \emph{PrimitiveTree}, whose \emph{fitness} attribute is set to the previously defined data structure.
To facilitate the implementation of both the \emph{Fitness} and \emph{Individual} class, in line~5--9, we make use of DEAP's \emph{creator} module.
For a detailed description of this module's functionality, the reader is referred to the documentation of the DEAP framework\footnote{DEAP Creator Module: \url{https://deap.readthedocs.io/en/master/api/creator.html}}.
Based on the definition of the \emph{Fitness} and \emph{Individual} class, we can then create the actual \emph{Toolbox} object together with the evolutionary operators, which is shown in line~10-17 of Listing~\ref{code:gp:toolbox}.
To generate a random individual, we apply our custom \emph{genGrow} function, as described above.
A complete population is then generated by utilizing DEAP's \emph{initIterate} and \emph{initRepeat} functions\footnote{DEAP Evolutionary Tools: \url{https://deap.readthedocs.io/en/master/api/tools.html}}.
Next, we implement suitable crossover and mutation operators for our tree-based G3P method, as defined in Section~\ref{sec:gggp-mutation-and-recombination}.
For the former, we apply subtree crossover for which an implementation is already provided by DEAP in the form of the function \emph{cxOnePoint}.
To implement subtree replacement and insertion within a single function, we, again, utilize our custom \emph{generate} function, which already provides the required functionality.
The resulting implementation of the function \emph{mutateSubtree} is shown in Listing~\ref{code:gp:mutate}.
\begin{listing}[ht]
\begin{minted}{python}
import random

def mutateSubtree(individual, min_height, max_height, pset):
    index = random.randrange(len(individual))
    node = individual[index]
    slice_ = individual.searchSubtree(index)
    def condition(height, depth):
        return depth < height
    subtree = None
    if random.random() < 0.5:
        subtree = individual[slice_]
    new_subtree = generate(pset, min_height, max_height, condition, node.ret, subtree)
    individual[slice_] = new_subtree
    return individual,
\end{minted}
	\caption{Subtree Mutation Operator}
	\label{code:gp:mutate}
\end{listing}
Similar to the \emph{genGrow} function, \emph{mutateSubtree} includes an interval from which the maximum depth of the randomly generated subtree is chosen.
To prevent the size of the inserted subtree from exceeding that of the original tree, an interval smaller than that of the \emph{genGrow} function should be chosen.
Now before we can implement the actual evolutionary search method, we need to address the problem of selecting promising individuals for recombination and mutation.
Therefore, we first need to find a way to accurately evaluate the quality of a multigrid method in an automatic manner.

\subsection{Fitness Evaluation and Selection}
\label{sec:fitness-evaluation-and-selection}
In Section~\ref{sec:grammar-based-algorithm-generation},
we have already demonstrated that it is possible to translate a computational graph in the form of Figure~\ref{fig:example-three-grid-method-computational-graph} into the algorithmic representations of a multigrid method.
In Section~\ref{sec:intermediate-representation}, we have introduced an intermediate representation (IR) for multigrid methods that stores each computational graph as a hierarchical composition of \emph{Cycle} objects, as demonstrated by the example shown in Listing~\ref{code:ir:example.py}.
Consequently, similar to the discussion in Section~\ref{sec:grammar-based-algorithm-generation}, we can translate each \emph{Cycle} object generated by the respective grammar to an algorithmic representation of the corresponding multigrid method.
This transformation is achieved by traversing the given hierarchy of IR objects in a recursive bottom-up manner.
However, we still need to figure out a way to evaluate each multigrid method based on this representation, which requires us to define one or multiple metrics for assessing their quality.
For this purpose, we first have to revisit the general definition of an iterative method, as provided in Section~\ref{sec:basic-iterative-methods}.
In principle, an iterative method for solving the linear system $A \bm{x} = \bm{b}$ is characterized by the repeated application of a number of well-defined computational steps, whereby each step has the purpose of computing an improved approximation of this system's solution based on the previous one.
If we execute an iterative solver on a particular computing system, usually the goal is to approximate the solution of the target problem up to a certain accuracy in as little time as possible, which we define as the \emph{solving time} $T_{\varepsilon}$, where $\varepsilon$ is the desired reduction of the initial error. 
If we assume that an iterative method performs the same sequence of computations in each of its iterations, the solving time can be defined as
\begin{equation}
	T_{\varepsilon} = t \cdot n_{\varepsilon},
	\label{eq:solving-time-basic}
\end{equation} 
where $t$ is the execution time of each iteration of the method and $n_{\varepsilon}$ the number of iterations required to achieve an error reduction of $\varepsilon$.
Note that in contrast to $n_{\varepsilon}$, the execution time $t$ solely depends on the amount of computation performed within each iteration and is thus independent of $\varepsilon$.
Similar to Section~\ref{sec:basic-iterative-methods}, we next define the \emph{convergence factor} of an iterative method as the limit of the sequence
\begin{equation*}
	\rho = \lim \limits_{n \to  \infty} \Bigg( \underbrace{ \max \limits_{x^{(0)} \in \mathbb{R}^{n}} \frac{\norm{\bm{x}^{(n)} - \bm{x}^{*}}}{\norm{\bm{x}^{(0)} - \bm{x}^{*}}} }_{\varepsilon_n} \Bigg)^{\frac{1}{n}},
\end{equation*}
where $\varepsilon_n$  is the minimum error reduction after $n$ iterations of the method.
Similarly, we can define an iteration-dependent convergence factor
\begin{equation}
	\rho_n = \left( \varepsilon_n \right)^{\frac{1}{n}}.
	\label{eq:iteration-dependent-convergence-factor}
\end{equation} 
Now assume that our goal is to achieve an error reduction of $\varepsilon$.
We can set $\varepsilon_n = \varepsilon$ and apply the natural logarithm to both sides of Equation~\eqref{eq:iteration-dependent-convergence-factor} which yields
\begin{equation*}
	n = \frac{\ln \varepsilon}{\ln \rho_n},
\end{equation*} 
where $n$ now represents the required number of iterations to achieve an error reduction of $\varepsilon$.
For a particular convergence factor $\rho_n$, the number of iterations $n$ can now be considered as a function of $\varepsilon$. 
We can thus insert the resulting term into Equation~\eqref{eq:solving-time-basic} and obtain
\begin{equation*}
	T_{\varepsilon} = t \cdot n_{\varepsilon} = t \cdot \frac{\ln \varepsilon}{\ln \rho_{n_{\varepsilon}}}.
\end{equation*} 
If we additionally assume that $\rho_{n_{\varepsilon}} \approx \rho$ for a sufficiently high number of iterations, we can further simplify this equation to
\begin{equation}
	T_{\varepsilon} = t \cdot n_{\varepsilon} = t \cdot \frac{\ln \varepsilon}{\ln \rho}.
	\label{eq:solving-time}
\end{equation} 
Note that since the negative value of the natural logarithm decreases between zero and one, a smaller convergence factor leads to a faster solving time.
For example, achieving an error reduction of $\varepsilon = 10^{-6}$ with a convergence factor of $\rho = 0.5$ would require 20 iterations, while with a convergence factor of $\rho = 0.1$, the same accuracy is achieved in only six iterations.
Consequently, there are two possibilities to improve upon an existing solver: Reducing its execution time $t$ or convergence factor $\rho$. 
However, as we have discussed in Section~\ref{sec:multigrid-cycles}, the choice of each method usually represents a compromise between a lower number of computations and hence a faster execution time and a faster convergence, which usually requires more computations.

While we have now derived suitable metrics for assessing the quality of an iterative method, we have not yet answered how to choose one or multiple optimization objectives based on this information.
In addition to achieving a fast solving time, the main goal of automated algorithm design is \emph{generalization}, which means that we aim to find multigrid methods that are robust in solving problems with similar characteristics.
Theoretically, we could simply execute a solver on each individual problem instance of interest and then measure its solving time.
However, due to the high computational demands of many PDE-based problems, this approach is often infeasible.
One way to mitigate this issue is the application of \emph{predictive models} to obtain an estimation of each objective without needing to execute the solver.
While this approach is usually faster to execute, its application to non-standard multigrid cycles, as the one shown in Figure~\ref{fig:non-traditional-multigrid-cycle}, is not well researched.
Due to the entirely different nature of the two performance-guiding metrics $t$ and $\rho$, we need to consider them separately.
First of all, note that the execution time $t$ of each step of an iterative method depends entirely on the amount of computation and its execution efficiency on the given computing platform.
In contrast, the convergence factor $\rho$ is linked to the mathematical properties of the solver, for instance, its capability to quickly reduce certain error components, as we have discussed in Section~\ref{sec:multigrid-methods}. 
While the availability of performance models for modern computer architectures, as the roofline~\cite{williams2009roofline} and execution-cache memory model~\cite{hager2016exploring}, enables predictions in an automated and deterministic manner, analyzing the mathematical properties of a solver is often difficult.
Local Fourier analysis (LFA)~\cite{wienands2004practical} represents a promising approach for predicting the convergence behavior of an iterative solver in a resolution-independent manner.
Even though LFA has been applied successfully to numerous applications~\cite{rodrigo2017validity}, only recently software tools that automate its execution have become available~\cite{rittich2018extending,kahl2020automated}.
In~\cite{schmitt2020constructing,hoefer2020comparing}, we have performed a number of experiments to test whether the library LFA Lab\footnote{LFA Lab: \url{https://github.com/hrittich/lfa-lab}} is suitable for the automated convergence estimation of grammar-generated multigrid methods.
While LFA Lab yields reliable predictions for many model problems and common multigrid cycles, unfortunately, we could not obtain the same degree of accuracy for non-standard cycles and more complicated problems, such as systems of PDEs.

A second option to decrease the evaluation cost is to estimate the quality of a given solver by measuring its characteristics on a number of proxy problems that are faster to solve than the original problem while possessing similar mathematical properties.
One of the main features of multigrid is that, if properly constructed, these methods can achieve $h$-independent convergence, which means that we can apply the same solver to a larger instance of the same problem without slowing down its convergence~\cite{trottenberg2000multigrid}.
Therefore, in case our goal is to solve a certain PDE discretized on a grid with step size $h$, we can instead consider an instance of the same PDE on a similar grid with lower resolution and thus a larger spacing of the grid points.
If a given multigrid method is able to achieve $h$-independent convergence on a sequence of increasingly finer-resolved instances of the same PDE, there is a high probability that the method also generalizes to other problem instances.
However, the execution of this approach requires us to generate efficient implementations for each multigrid method considered, which can then be executed to obtain the relevant quality metrics on each proxy problem.
For this purpose, we need to automate the process of generating a solver implementation based on its algorithmic representation.
Recently, code generation techniques that grant us this ability have become available~\cite{kostler2020code,schmitt2018automating}.
ExaStencils\footnote{ExaStencils: \url{https://www.exastencils.fau.de}} is a code generation framework, implemented in the Scala programming language, that has been specifically designed for the automatic generation of scalable multigrid solver implementations on modern parallel computing hardware~\cite{lengauer2014exastencils,lengauer2020exastencils}.
It enables the formulation of multigrid methods in a discretization-independent domain-specific language (DSL) called ExaSlang~\cite{schmitt2014exaslang,schmitt2016systems}, which is almost indistinguishable from the textbook-like description of an algorithm.
Based on this DSL specification, the framework is able to generate highly-optimized C++ code for different problem sizes and computer architectures.
To illustrate this approach, Listing~\ref{code:exaslang-example} shows an ExaSlang implementation of the three-grid method formulated in Algorithm~\ref{alg:example-three-grid-method}.
% \begin{lstlisting}[language=ExaSlang3]
% Function VCycle@finest {
%   r = b - A * x
%   x@coarser = 0.0
%   // Restriction
%   b@coarser = R * r
%   r@coarser = b@coarser - A@coarser * x@coarser
%   x@(coarser-1) = 0.0
%   b@(coarser-1) = R@coarser * r@coarser
%   // Apply the Coarse-Grid Solver
%   VCycle@(coarser-1)
%   x@coarser += P@(coarser-1) * x@(coarser-1)
%   // Smoothing with damped Jacobi
%   x@coarser += 0.6 * diag_inv(A@coarser) * ( b@coarser - A@coarser * x@coarser )
%   // Prolongation and Coarse-Grid Correction
%   x += P@coarser * x@coarser
% }
% \end{lstlisting}
\begin{listing}
\begin{minted}{scala}
Function VCycle@finest {
  r = b - A * x
  x@(coarser) = 0.0
  // Restriction and coarsening
  b@(coarser) = R * r
  r@(coarser) = b@(coarser) - A@(coarser) * x@(coarser)
  x@(coarser-1) = 0.0
  b@(coarser-1) = R@(coarser) * r@(coarser)
  // Apply the coarse-grid solver
  VCycle@(coarser-1)
  x@(coarser) += P@(coarser-1) * x@(coarser-1)
  // Smoothing with damped Jacobi
  x@(coarser) += 0.6 * diag_inv(A@(coarser)) * (b@(coarser) - A@(coarser) * x@(coarser))
  // Prolongation and coarse-grid correction
  x += P@(coarser) * x@(coarser)
}
\end{minted}
	\caption{Three-Grid Example from Algorithm~\ref{alg:example-three-grid-method} in ExaSlang}
	\label{code:exaslang-example}
\end{listing}
As this example demonstrates, ExaSlang enables the formulation of the same operations on different levels of a discretization hierarchy based on so-called \emph{level declarations}\footnote{In case no level declaration is provided, ExaStencils assumes that the operation is applied on the current level.}.
Similar to the specification of the grid spacing as a subscript in Algorithm~\ref{alg:example-three-grid-method}, we can utilize this notation to formulate each operation within a multigrid method relative to the \emph{finest} grid.
Since the ExaSlang specification of a multigrid method is thus independent of the actual discretization, the ExaStencils compiler is able to generate implementations of the same solver for different grid sizes.
Furthermore, as ExaSlang supports the formulation of each multigrid operation in a high-level mathematics-like syntax, we can apply the same solver to different problems only by changing the definitions of the individual operators, such as the system, prolongation, and restriction operator.
Therefore, ExaStencils is ideally suited for the automatic generation and evaluation of multigrid methods within our evolutionary program synthesis approach.
While ExaStencils is implemented in the Scala programming language and can thus not be accessed directly within Python code, it provides simple configuration files, which allow us to adapt certain problem characteristics, such as the grid spacing and the number of coarsening steps.
Furthermore, additional code generation options, such as compiler-based performance optimizations, can be en- or disabled.
To encapsulate the usage of the ExaStencils compiler, EvoStencils provides a \emph{ProgramGenerator} class\footnote{Due to the vast and complicated implementation of this class, it is omitted in this thesis, while its functionality is described informally.}.
In addition to the automatic adaption of the respective configuration files and the actual execution of the code generation process, this class also includes functionality for extracting the required information about a given problem based on an existing ExaSlang specification. 
Therefore, in case this specification contains all the necessary information for generating the corresponding multigrid grammar, it allows us to integrate EvoStencils directly into ExaStencils' solver generation workflow.
Moreover, this approach enables the application of our evolutionary program synthesis method to many of the problems already available within ExaStencils.
Based on the functionality provided by the \emph{ProgramGenerator} class, we can implement an evaluation function in the form of the \emph{evaluate} method of the previously mentioned \emph{Optimizer} class.
The implementation of this method is shown in Listing~\ref{code:optimization:evaluate}.
\begin{listing}
\begin{minted}{python}
from deap import gp

# Implemented as a method of the Optimizer class
def evaluate(self, individual, pset, min_level, max_level, evaluation_samples=3, pde_parameter_values=None):
    if pde_parameter_values is None:
        pde_parameter_values = {}

    key = str(individual)
    if key in self.individual_cache:
        return self.individual_cache[key]
    expression, _ = gp.compile(individual, pset)
    solving_time, convergence_factor, number_of_iterations = self.program_generator.generate_and_evaluate(expression, min_level, max_level, evaluation_samples, pde_parameter_values)

    fitness = convergence_factor, solving_time / number_of_iterations
    self.individual_cache[key] = fitness
    return fitness
\end{minted}
	\caption{Optimizer Class -- Evaluate Method}
	\label{code:optimization:evaluate}
\end{listing}
It provides a high-level interface to the evaluation of an arbitrary individual generated based on the provided \emph{PrimitiveSetTyped} object.
In order to prevent the repeated evaluation of structurally-equal individuals, the \emph{Optimizer} class implements a caching mechanism that keeps track of all previously evaluated individuals.
Therefore, in line~6, we first check whether the given individual has already been evaluated before, in which case we simply return the cached objective function value in line~7.
Otherwise, we utilize the \emph{generate\_and\_evaluate} method of the \emph{ProgramGenerator} class in line~10 to generate a C++ implementation of the corresponding solver using the ExaStencils compiler, which is then executed on a proxy problem.
As a result, we obtain three metrics, the total solving time $T_\varepsilon$, the convergence factor $\rho$, and the number of iterations $n_\varepsilon$, whereby the desired error reduction $\varepsilon$ is usually determined by the given ExaStencils specification of the problem.
Since in general, the solution of the given PDE is not known in advance, we approximate the convergence factor using the formula
\begin{equation}\label{eq:asymptotic_convergence_factor}
	\rho \approx \left(\prod_{i=1}^{n}\tilde{\rho}_i \right)^{1/n},
\end{equation} where $n$ is the number of iterations until convergence and
\begin{equation}
	\tilde{\rho}_{i} = \frac{\norm{b_h - A_h x_{h}^{(i)} }}{\norm{b_h - A_h x^{(i-1)}_{h}}}
\end{equation}
is the L2-norm of the residual reduction in every iteration of the solver~\cite{trottenberg2000multigrid}.
To determine the execution time per iteration $t$ of the method, we simply divide the total solving time by the number of iterations.
To reduce potential hardware-based variations in the execution of each solver, $T_\varepsilon$ is obtained as an average over multiple evaluations of the same individual.
Note that, at this point, we assume that a multi-objective optimization is performed, as the returned \emph{Fitness} object includes two objectives.
However, it is also possible to implement a similar function that returns a single-objective fitness value, for instance, in the form of the total solving time.
Finally, the method enables the adaption of certain parameters of the underlying PDE for the solver evaluation.
For this purpose, a mapping from each parameter to the respective value can be provided in the form of a Python \emph{dictionary}.

After successfully obtaining a \emph{Fitness} object for each individual, we can complete our collection of evolutionary operators with the definition of suitable selection and elitism procedures.
As we have discussed in Section~\ref{sec:gggp-evaluation-and-selection}, in contrast to initialization, mutation, and crossover, the selection of individuals within an evolutionary algorithm is independent of their representation as a data structure but solely depends on the definition of their fitness.
Therefore, different selection operators have been proposed for single and multi-objective evolutionary algorithms, of which many are already available within the DEAP framework.  
Similar to the other evolutionary operators in Listing~\ref{code:gp:toolbox}, we can add a certain selection operator to a given \emph{Toolbox} object using its \emph{register} method.
Listing~\ref{code:optimization:selection} shows the example initialization of a \emph{Toolbox} object with the widely-used NSGA-II selection operator~\cite{deb2002fast}. 
\begin{listing}
\begin{minted}{python}
from deap import tools

def init_multi_objective_selection(self, pset, evaluation_function):
    self.toolbox.register('evaluate', evaluation_function, pset=pset)
    self.toolbox.register("select", tools.selTournamentDCD)
    self.toolbox.register("elitism", tools.selNSGA2)
\end{minted}
	\caption{Toolbox Initialization with the NSGA-II Selection Operator}
	\label{code:optimization:selection}
\end{listing}
In a similar manner, we can initialize a \emph{Toolbox} object with other selection operators and evaluation functions based on the definition of the \emph{Fitness} class and the functionality of the \emph{evaluate} method of our \emph{Optimizer} class.

\subsection{Search Algorithm}
\label{sec:evostencils-part1:search-algorithm}
We have now finally assembled all components that are required to implement an evolutionary search method for the grammar-based optimization of multigrid methods.
Since each component has been registered to the respective \emph{Toolbox} object, we can utilize the resulting interface to implement a search method that is independent of the actual implementation of each operator as a method of the \emph{Optimizer} class. 
For this purpose, we store the previously assembled \emph{Toolbox} object in the respective attribute of this class.
The resulting implementation of the \emph{evolutionary\_search} method can be found in Listing~\ref{code:optimization:evolutionary_search}. 
\begin{listing}
\begin{minted}{python}
def evolutionary_search(self, params, use_random_search=False):
    # Generate and evaluate initial population
    population = self.toolbox.population(n=params.initial_population_size)
    invalid_ind = [ind for ind in population]
    self.evaluate_individuals(invalid_ind)
    # Select mu individuals from the initial population
    population = self.toolbox.select(population, params.mu_)
    for gen in range(1, params.generations + 1):
        if use_random_search:
            # For random search simply generate individuals randomly
            offspring = self.toolbox.population(n=params.lambda_)
        else:
            # Select lambda parents for mutation and crossover
            selected = self.toolbox.select(population, params.lambda_)
            parents = [self.toolbox.clone(ind) for ind in selected]
            offspring = self.create_offspring(parents, params.crossover_probability)
        # Evaluate new individuals
        invalid_ind = [ind for ind in offspring if not ind.fitness.valid]
        self.evaluate_individuals(invalid_ind)
        # Select new population from the combined set of parents
        # and offspring (elitism)
        population = self.toolbox.elitism(population + offspring, params.mu_)

    return population
\end{minted}
	\caption{Optimizer Class -- Evolutionary Search Method}
	\label{code:optimization:evolutionary_search}
\end{listing}
Note that each part of this implementation corresponds to a particular step in our general description of an evolutionary program synthesis method in Algorithm~\ref{alg:genetic-programming}.
In accordance with the common notation to describe evolutionary algorithms~\cite{back1997handbook} based on the two parameters $\mu$ and $\lambda$, our method implements a $\left(\mu + \lambda \right)$ strategy.
Therefore, in each generation, $\lambda$ new individuals are created based on $\mu$ parent individuals.
The new population is then selected from the combined set of parent and child individuals using the defined elitism operator.
For the sake of simplicity, we implement the evaluation of a list of individuals and the creation of offspring using mutation and crossover in two separate functions, which are shown in Listing~\ref{code:optimization:evaluate-and-create-offspring}.
\begin{listing}
\begin{minted}{python}
import random

def evaluate_individuals(self, invalid_ind):
    fitnesses = self.toolbox.map(self.toolbox.evaluate, invalid_ind)
    for ind, fit in zip(invalid_ind, fitnesses):
        ind.fitness.values = fit

def create_offspring(self, parents, crossover_probability):
    offspring = []
    for ind1, ind2 in zip(parents[::2], parents[1::2]):
        child1, child2 = ind1, ind2
        key1, key2 = None, None
        count = 0
        operator_choice = random.random()
        while key1 or key2 in self.individual_cache and count < 10:
            if operator_choice < crossover_probability:
                child1, child2 = self.toolbox.mate(ind1, ind2)
            else:
                child1, = self.toolbox.mutate(ind1)
                child2, = self.toolbox.mutate(ind2)
            key1, key2 = str(child1), str(child2)
            count += 1
        del child1.fitness.values, child2.fitness.values
        offspring.extend([child1, child2])
    return offspring
\end{minted}
	\caption{Auxiliary Functions for Creating and Evaluating Offspring}
	\label{code:optimization:evaluate-and-create-offspring}
\end{listing}
To enforce the creation of novel individuals, which have not already been discovered in previous generations, we check whether individuals created through mutation or crossover are already present in the cache.
If this condition is fulfilled for at least one of the two children created from each pair of parent individuals, we repeat the application of the respective mutation or crossover operator until two suitable children are obtained.
To initiate the search with a sufficiently diverse population, the number of initially-generated individuals can be set higher than $\mu$, which is controlled by the parameter \emph{initial\_population\_size}.
In line~3, we then randomly generate an initial population, which is evaluated in line~5.
As a next step, $\mu$ individuals are selected for the first generation using the registered \emph{elitism} operator in line~5.
The actual search is then performed for a predefined number of generations.
As a fallback solution, we additionally implement a simple random search in line~11, where we randomly generate $\lambda$ individuals in each generation.
Within the evolutionary algorithm, we first select $\lambda$ individuals for crossover and mutation, using the \emph{select} operator in line~14.
As mutation represents a modifying operation, we first create an identical copy of each selected parent individual in line~15, based on which then either crossover or mutation is applied using the previously-mentioned \emph{create\_offspring} method in line~16.
Finally, in line~18--22, the resulting newly created individuals are evaluated, and a new elitist population of size $\mu$ is selected from the combined set of parent and child individuals.
Note that this step is identical to the random search variant of our implementation.

With the definition of this method, we have now completed our implementation for the automated grammar-based design of multigrid methods for solving PDE-based problems.
Before we evaluate the effectiveness of our approach on a number of benchmark problems, we will discuss a few important extensions of this basic implementation in the next chapter.
In particular, we will present an extension of the evolutionary search method shown in Listing~\ref{code:optimization:evolutionary_search} that enables the systematic generalization of a population of multigrid methods to a sequence of increasingly-difficult instances of the same problem.
Furthermore, in order to accelerate the evaluation of the large number of individuals required to obtain competitive solvers for many PDEs, we demonstrate how our implementation can be parallelized on multi-node systems using the message-passing interface (MPI).
%TODO include again if systems of PDEs are covered in the next section
%Finally, we will also briefly discuss how the intermediate representation presented in Section~\ref{sec:intermediate-representation} can be extended to support the generation of multigrid methods for systems of PDEs.

%% file: contents/evostencils_part2.tex
In the following, we will discuss two extensions of the evolutionary program synthesis framework described in the last chapter, which are crucial for the effective application of our approach to many applications - \emph{generalization} and \emph{parallelization}.  
First of all, to discover multigrid methods that are capable of solving different instances of a PDE, we are concerned with their generalizability.
For this purpose, we will derive an adapted version of our original evolutionary algorithm that utilizes the special properties of multigrid methods.
Furthermore, since the successful application of this method requires the code generation-based evaluation of a large number of individuals, we will present a distributed parallelization scheme, which enables the execution of our implementation on current multi-node computing systems.
In the final chapter of this thesis, we will then demonstrate how the implementation described here can be successfully applied to different PDEs, yielding multigrid-based solvers that are competitive with hand-optimized methods. 

\section{Generalization}
\label{sec:generalization}
As we have briefly discussed in Section~\ref{sec:fitness-evaluation-and-selection}, if properly constructed, the error reduction capabilities of a multigrid method are independent of the discretization width $h$, which is usually described with the term $h$-independent convergence.
Therefore, the same method can often be successfully applied to different systems that arise from similar discretizations of the same PDE.
We have already introduced the idea of evaluating each multigrid method on a number of proxy applications whose properties are similar to the problem that we actually aim to solve.
The motivation for this idea is that, while we are usually interested in solving a problem instance of a specific size, the evaluation of each solver on this instance requires an excessive amount of computational resources.
In many cases, it is possible to construct such a set of proxy applications by discretizing the same PDE with varying step sizes $h$.
If we are able to design a solver that achieves $h$-independent convergence, it can be expected to solve each of these proxy applications using the same number of operations per grid point.
In addition to this requirement, we want to identify the method that leads to the fastest solving time $T_\varepsilon$ for the target problem with a discretization width of $h$.
The main challenge is thus to identify this method within the evolved population while performing the majority of evaluations on smaller problem instances with a discretization width $H$ greater than $h$.
Note that since evolutionary algorithms are usually not guaranteed to find the global optimum, we restrict ourselves to identifying the optimum among the individuals discovered while evolving the population.
%First of all, recall that our evolutionary algorithm, whose implementation can be found in Listing~\ref{code:optimization:evolutionary_search}, performs a search by evolving a population of individuals for a specified number of generations.
In each generation, new individuals are created by applying mutation and crossover.
However, only a limited number of individuals, chosen from the combined set of child and parent individuals, are allowed to enter the new population.
Whether an individual is accepted for the next generation solely depends on its fitness, as defined by one or multiple objectives.
Therefore, to achieve generalization, it is crucial that we define the fitness of an individual in a way that maximizes the probability that the optimum, according to our criterion $T_{\varepsilon}$, is contained in the final population.
We thus have to prevent the eviction of this individual from the population at any point during the execution of our algorithm.

\subsection{Objective Function Definition}
\label{sec:generalization:objective-function-definition}
As we have shown in Section~\ref{sec:fitness-evaluation-and-selection} the solving time $T_{\varepsilon}$ is given by the formula
\begin{equation*}
	T_{\varepsilon} = t \cdot \frac{\ln \varepsilon}{\ln \rho},
\end{equation*}
where $t$ is the execution time per iteration and $\rho$ the (asymptotic) convergence factor of the iterative solver.
Therefore, there are two ways to define the fitness of an individual based on this metric:
\begin{enumerate}
	\item Single-objective: $T_{\varepsilon}$
	\item Multi-objective:  $t$ and $\rho$
\end{enumerate}
The main difference here is that while a single-objective evaluation always returns a single individual as the optimum, a multi-objective evaluation instead identifies a set of Pareto-optimal individuals, i.e., individuals that do not \emph{dominate} each other, which means that they are unable to achieve a better value in both objectives.
Since in the given case, an improvement in either of the two objectives, $t$ and $\rho$, necessarily leads to a faster solving time $T_{\varepsilon}$, the single-objective optimum is always contained in the Pareto-front obtained from a multi-objective evaluation of the same individuals.
Therefore, the main question regarding generalization is whether the individual with the fastest solving time $T_{\varepsilon}$ for our target problem will be consistently selected as an optimum when it is evaluated on smaller instances of the same problem.
First of all, note that the solving time $T_{\varepsilon}$, as a single objective, does not necessarily lead to the same ranking of individuals for each problem instance.
While the convergence factor $\rho$ of a functioning multigrid method executed as an iterative solver is expected to be constant, its execution time per iteration $t$ is drastically affected by hardware effects.
For instance, if the memory requirement for a certain problem size exceeds the capacity of the cache, the execution time is expected to increase substantially compared to a problem instance that still fits into the cache\footnote{This is only true for memory-bound computations. A property that is, however, fulfilled for the majority of stencil operations.}.
As a consequence, the execution time per iteration of a multigrid method that achieves the fastest solving time for a small problem might drastically increase for a larger problem instance, which means that the solver might no longer be optimal.
%On the downside, if we consider the solver that achieves the fastest solving time for a certain problem instance.
%While this method might no longer be optimal with respect to its solving time for problems of smaller size, there is a high probability that it is still contained in the Pareto-front obtained through a multi-objective evaluation.
%In general, faster convergence is either achieved by applying a higher number of smoothing or coarse-grid correction steps.

Now consider an example of two different non-dominating multigrid V-cycles.
The first cycle uses two smoothing steps per level leading to a convergence factor of 0.15 and an execution time per iteration of one millisecond on a grid with step size $h$, while the second one achieves a convergence factor of 0.1 and an execution time of 1.5 milliseconds using three smoothing steps per level.
If we now execute both methods on a smaller grid with step size $2h$, it is very unlikely that the first cycle will converge faster than the second one.
Likewise, three smoothing steps per level will also lead to a higher execution time per iteration for a smaller problem.
As a consequence, the dominance relation between both methods is preserved.
In contrast, assume that the second method achieves a slightly faster solving time on the larger problem.
Since it is impossible to predict the order of magnitude of change in the value of $t$ for both methods without actually executing them, we can not be sure whether this is still the case for a smaller problem.
While treating the design of a generalizable multigrid method as a multi-objective optimization problem increases the probability that the final population contains the fastest solver, this approach still has limitations.
If we consider different choices for each smoothing and coarse-grid correction step, our assumption that a certain sequence of operations also leads to a faster execution time for a smaller problem size is no longer a certainty, as some operations might possess a different computational complexity.
A simple example of such an operation is line smoothing which becomes significantly more expensive when it is applied to larger problems~\cite{trottenberg2000multigrid}.
One caveat for this issue, which has already been mentioned in Section~\ref{sec:fitness-evaluation-and-selection}, is to predict the execution time of a solver with a performance model~\cite{williams2009roofline,hager2016exploring}.
However, since in this work, we only consider operations whose complexity is independent of the problem size, i.e., pointwise smoothers and block smoothers with blocks of fixed size, the measured execution time per iteration provides a sufficient prediction for larger instances of the same problem.

\subsection{Generalization Procedure}
\label{sec:generalization-procedure}
Based on the observations made in the last section, we can now formulate an extension of our evolutionary algorithm for the systematic generalization of multigrid methods to a given problem class, which is summarized in Algorithm~\ref{alg:generalization-procedure}.
\begin{algorithm}
	\caption{Generalization Procedure}
	\label{alg:generalization-procedure}
	\begin{algorithmic} % The number tells where the line numbering should start
		\State \textbf{Construct} the grammar $G_0$ for the initial problem
		\State \textbf{Initialize} the population $P_0$ based on $G_0$
		\State \textbf{Evaluate} $P_0$ on the initial problem with respect to $t$ and $\rho$
		\For{$i := 0, \dots, n-1$}
		\If{$i > 0$ and $i \mod m = 0$}
		\State $j := i / m$ 
		\State Increase the problem size
		\State \textbf{Construct} the corresponding grammar $G_j$
		\State  \textbf{Adapt} the current population $P_i$ to $G_j$
		\State \textbf{Evaluate} $P_i$ on the new problem with respect to $t$ and $\rho$
		\EndIf
		\State \textbf{Generate} new solutions $C_i$ based on $P_i$ and $G_j$
		\State \textbf{Evaluate} $C_i$ on the current problem with respect to $t$ and $\rho$
		\State \textbf{Select} $P_{i+1}$ from $C_i \cup P_i$
		\EndFor
		\State \textbf{Construct} the grammar $G$ for the target problem
		\State  \textbf{Adapt} the current population $P_n$ to $G$
		\State \textbf{Identify the best solver} by evaluating $P_{n}$ on the target problem
	\end{algorithmic}
\end{algorithm}
The first step of this procedure is the choice of an initial problem size, which should be small enough to enable the fast evaluation of a large number of randomly generated and thus often inefficient solvers.
As the search progresses and the average quality of the individuals in the population improves, the problem size can then be iteratively increased toward the target size.
While each problem size adaption increases the required time to evaluate each solver, it also improves the accuracy of evaluation with respect to both objectives.
As discussed in Section~\ref{sec:generalization:objective-function-definition}, evaluating a multigrid-based solver on a sequence of increasingly larger instances of the same problem allows us to assess whether the choice of the grid spacing $h$ affects its convergence.
Also, note that if there is only a small difference between the execution times of two non-dominating solvers, even slight hardware-based fluctuations in the measurements might perturb the outcome of an evaluation.
Considering a larger instance of the same problem reduces the relative magnitude of these fluctuations compared to the overall evaluation time.
At the end of Algorithm~\ref{alg:generalization-procedure}, we obtain a population that has evolved against a sequence of problem instances that iteratively approaches the size of the target instance.
We, therefore, identify the fastest solver by only considering those individuals contained in the first non-dominated front of this population, i.e., the subset in which none of the individuals is dominated by any of those present in the population.

Finally, note that so far, we have only considered increasing the problem size while evolving the population.
However, in certain cases, a PDE contains additional parameters which need to be adapted accordingly.
One prominent example that we consider in this thesis is the indefinite Helmholtz equation, as given by
\begin{equation}
	-(\nabla^{2} + k^{2})u = f,
\end{equation}
where $\nabla^{2}$ is the Laplace operator, $k$ the \emph{wavenumber} and $f$ the source term.
In general, the difficulty of solving this problem increases with the value of $k$.
However, many applications require the discretization width $h$ to fulfill an accuracy requirement, such as $h k \leq 0.625$. 
As a consequence, in order to solve this problem on a coarser grid, we also need to adapt the wavenumber accordingly, which results in a sequence of problem instances not only increasing in size but also in difficulty.
In Chapter~\ref{chapter:experiments}, we will demonstrate that our generalization procedure can cope with this challenge, yielding efficient multigrid methods for Helmholtz problems of varying size and difficulty.

\subsection{Implementation}
After we have now both motivated and described a procedure for the grammar-based design of generalizable multigrid methods, the remaining step is its successful implementation within the EvoStencils framework.
First of all, note that the individual evolutionary operators, i.e., initialization, mutation, crossover, and selection, are all implemented based on the created \emph{PrimitiveSetTyped} object. 
Their application is thus independent of the underlying problem size, and we only have to attach each operator to the respective \emph{Toolbox} object, as shown in Listing~\ref{code:gp:toolbox}.
We can generate a grammar for grid hierarchies with varying step size $h$ but the same number of coarsening steps by utilizing the \emph{generate\_grammar} function defined in Listing~\ref{code:grammar:generate-grammar}.
For this purpose, we only need to provide a different value for the argument \emph{max\_level}.
However, an essential issue that has not been addressed yet is how we can adapt the current population to a new grammar.
For this purpose, we need to return to the original formulation of our multigrid grammar, whose productions are shown in Algorithm~\ref{table:multigrid-grammar}.
Note that the level of each variable and terminal, unless it is level-independent, is denoted by its subscript.
Each of these terms is an expression whose value depends on $h$, i.e., the spacing of the finest grid.
In other words, if we want to apply a multigrid method whose derivation tree has been generated based on a discretization hierarchy with the step size $h$ to a different one of similar depth but with a step size of $H$, we only have to replace each occurrence of $h$ with the value of $H$.
Now recall that in our implementation, each derivation tree is represented as a \emph{PrimitiveTree} object of DEAP's GP module, which internally stores a list of its nodes in depth-first order.
Listing~\ref{code:gp:primitive} contains a minimal implementation of the \emph{Primitive} and \emph{Terminal} class within DEAP.
\begin{listing}
\begin{minted}{python}
class Primitive(object):
    def __init__(self, name, args, ret):
        self.name = name # name as a string
        self.arity = len(args)
        self.args = args # list of argument types
        self.ret = ret # return type

class Terminal(object):
    def __init__(self, terminal, symbolic, ret):
        self.ret = ret # return type
        self.value = terminal
        self.name = str(terminal) # name as a string
        self.conv_fct = str if symbolic else repr

    @property
    def arity(self):
        return 0
\end{minted}
	\caption{Primitive and Terminal Class in DEAP}
	\label{code:gp:primitive}
\end{listing}
According to this implementation, each \emph{Primitive} and \emph{Terminal} object is identified by three attributes - its \emph{name}, argument types (\emph{args}), and return type (\emph{ret}).
Note that a \emph{Terminal} object does not possess any input types, and if the \emph{terminal} argument of its \emph{\_\_init\_\_} method is provided as a string, the \emph{terminal} and \emph{name} attributes are identical.
However, neither of the classes incorporates the necessary information for its compilation to executable Python code.
Therefore, in addition to the given \emph{PrimitiveTree}, the respective \emph{PrimitiveSetTyped} object needs to be provided, which then enables constructing a sequence of function applications according to the order of the tree nodes.  
As a consequence, in case two different \emph{PrimitiveSetTyped} objects are structurally equal and employ the same name string for each of their \emph{Primitive} and \emph{Terminal} objects, they enable the compilation of the same \emph{PrimitiveTree} objects.
Now note that we generate the \emph{name} of each \emph{Primitive} and \emph{Terminal} in Listing~\ref{code:grammar:add-level} and~\ref{code:grammar:add-terminals} using the same problem-independent naming convention.
Since the subscript included in each name specifies the current level within the discretization hierarchy only symbolically, each \emph{PrimitiveTree} generated based on a \emph{PrimitiveSetTyped} object returned by our \emph{generate\_grammar} function can be compiled by any other \emph{PrimitiveSetTyped} object that has been obtained using the same value for the \emph{depth} argument.
While this provides us with a way to compile a given \emph{PrimitiveTree} after replacing our previous \emph{PrimitiveSetTyped} object with a new one, we have not yet addressed the question of how this affects the generation of new individuals through mutation and crossover.
In the case of both operators, the creation of new individuals is subject to the type constraints specified within the \emph{Primitive} and \emph{Terminal} objects of the parent individuals.
Therefore, whenever a modification is to be applied at a specific position within a \emph{PrimitiveTree}, the types of the respective \emph{Primitive} and \emph{Terminal} objects need to match.
If we consider the initialization of each type in the respective method of the \emph{Types} class defined in Listing~\ref{code:grammar:types}, we can see that each \emph{Type} object is created as a function of the \emph{depth} argument.
It is independent of the details of the discretization hierarchy based on which a \emph{PrimitiveSetTyped} object is constructed.
Therefore, exchanging a \emph{PrimitiveSetTyped} with one of similar \emph{depth} does not affect the semantics of each \emph{Type} object, which means that the same mutation and crossover operators can be applied without further adaption.
We can conclude the previous discussion by stating that while each \emph{PrimitiveTree} contains the computational structure of the corresponding multigrid method, all required information for its interpretation with respect to a particular discretization hierarchy is contained in the supplied \emph{PrimitiveSetTyped} object.
Therefore, whenever we want to increase the problem size within our evolutionary program synthesis procedure, we only have to generate a new \emph{PrimitiveSetTyped} object based on which we then update all operations registered at the current \emph{Toolbox} object.
At the same time, all individuals in the population remain unchanged and only need to be reevaluated on the updated problem instance.
The resulting steps can then be implemented in the form of an additional method of the \emph{Optimizer} class, which is shown in Listing~\ref{code:optimization:adapt-problem-size}.
\begin{listing}
\begin{minted}{python}
def adapt_problem_size(max_level, params, multi_objective=True):
    self.pset = generate_grammar(self.program_generator, max_level, params.depth, params.relaxation_factor_samples)
    # Initialize the code generator
    # according to the specified maximum level
    self.program_generator.initialize_code_generation(max_level - params.depth, max_level)
    self.init_toolbox(self.pset, params.tree_min_height, params.tree_max_height, multi_objective)
    if multi_objective:
        self.init_multi_objective_selection(self.pset, self.evaluate)
    else:
        self.init_single_objective_selection(self.pset, self.evaluate)
\end{minted}
	\caption{Optimizer Class -- Problem Size Adaption}
	\label{code:optimization:adapt-problem-size}
\end{listing}
This implementation leverages the \emph{initialize\_code\_generation} method of the provided \emph{CodeGenerator} object to adapt the respective ExaStencils configuration files according to the provided \emph{max\_level} argument, which defines the maximum level of the discretization hierarchy.
What now remains is the integration of this method into the implementation of our evolutionary algorithm defined in Listing~\ref{code:optimization:evolutionary_search}.
However, we postpone this step until we have described the next major extension of our basic implementation, which is the distributed parallelization of our approach using the message-passing interface (MPI).

\section{Distributed Parallelization}
\label{sec:distributed-parallelization}
As we have already discussed in Section~\ref{sec:search-space-estimation}, the size of the search space spanned by our family of grammars makes it infeasible to evaluate every single multigrid method that can be generated based on their productions.
Even though the utilization of search heuristics in the form of an evolutionary algorithm allows us to reduce the number of considered individuals significantly, we still have to execute a large number of solvers on different problem instances.
Furthermore, note that for the evaluation of each solver, we need to utilize the ExaStencils framework to automatically generate a C++ implementation, which then has to be compiled into an executable program.
Both steps induce a significant overhead and hence further increase the evaluation time.
A common approach to accelerate the computationally intensive parts of an evolutionary algorithm is to distribute their execution to several compute nodes such that each computational step can be performed in parallel.
An overview of different approaches for the distributed parallelization of evolutionary algorithms can be found in~\cite{gong2015distributed}.
In principle, we can distinguish between approaches that are behaviorally equivalent to a sequential evolutionary algorithm and those that do not fulfill this property, usually in order to achieve better scalability.
Unfortunately, which approach leads to the best outcome can not be answered in general.
Before deciding on a specific parallelization method, we hence first need to investigate which parts of our evolutionary search method have to be parallelized to achieve good scalability.

\subsection{Empirical Execution Time Analysis}
\label{sec:execution-time-analysis}
According to Algorithm~\ref{alg:genetic-programming}, each step of our evolutionary algorithm consists of the following four operations:
\begin{enumerate}
	\item Parent Selection
	\item Child Creation
	\item Child Evaluation
	\item Population Selection
\end{enumerate}
In order to estimate the expected speedup of a parallel implementation of each of these operations, we first determine their fraction of the algorithm's total execution time.
As a representative example, we consider Poisson's equation on the unit square $\left[0,1\right]^2$ with Dirichlet boundary conditions. 
We discretize this equation using a uniform grid with step size $h = 1/2^{11}$ and the common five-point stencil
\begin{equation*}
	\nabla^2_h = 
	\frac{1}{h^2} \begin{bmatrix}
		0 & 1 & 0\\
		1 & -4 & 1 \\
		0 & 1 & 0  
	\end{bmatrix},
\end{equation*} 
which leads to a system of linear equations with $4\,190\,209$ unknowns.
The complete specification of this test problem can be found in Section~\ref{sec:poisson-equation}.
In order to obtain representative measurements for each of the four steps of our evolutionary algorithm, we execute it in the form of the \emph{run} method defined in Listing~\ref{code:optimization:optimizer} for a total number of 250 generations on a single socket of the \emph{Meggie} compute cluster of the Erlangen National High-Performance Computing Center (NHR)\footnote{As an exception, the child evaluation step is performed in parallel on multiple sockets of the same type. Since the order of evaluation does not change the behavior of our algorithm, this decision does not affect the measurements.}.
Note that we do not apply the generalization procedure introduced in Section~\ref{sec:generalization} here, which means that the problem size is kept constant throughout the execution of our evolutionary algorithm.
In each generation, we select $\lambda = 256$ individuals from the current population, based on which we create $\lambda$ children.
A new population of $\mu = 256$ individuals is then obtained from the combined set of $512$ individuals using the NSGA-II non-dominated sorting procedure described in~\cite{deb2002fast}.
We measure the average time required for each of the four steps over all generations, which is shown in Table~\ref{table:evolutionary-search-profiling}.
\begin{table}
	\caption{Average time required for each step performed within one generation of the evolutionary algorithm.}
	\label{table:evolutionary-search-profiling}
	\centering
	\begin{tabular}{l c}
		\toprule
		Step & Average Time \\
		\midrule
		Parent Selection & 0.68 ms \\
		\midrule
		Child Creation  & 0.32 s \\
		\midrule
		Child Evaluation  & 3.31 h \\
		\midrule
		Population Selection  & 0.20 s \\
		\bottomrule
	\end{tabular}
\end{table}
According to these measurements, the overall execution time of our implementation is heavily dominated by the evaluation step, which is reflected in the fact that the combined execution times of all other steps do not even account for one percent of the overall time.
Consequently, we can drastically reduce the execution time of our implementation by performing the evaluation of multiple individuals in parallel, while the parallelization of any other step will only result in a negligible speedup. 
Since we have to evaluate at least a single individual per compute node, the maximum achievable speedup is equal to $\lambda$, i.e., the number of children created in each generation of our method.
Finally, we need to answer the question of whether the previous statements still hold when we apply our evolutionary program synthesis method to other PDE-based problems.
Therefore, note that the two-dimensional Poisson equation represents a common test problem that is well known to be efficiently solvable by multigrid.
In fact, Poisson's equation has motivated the development of multigrid methods by Fedorenko~\cite{fedorenko1962relaxation}, Brandt~\cite{brandt1977multi} et al., and hence these methods achieve the highest possible degree of efficiency in solving it.
For the majority of other PDE-based problems of similar size, we can expect a further increase in the execution time, which leads to an even larger relative time consumption of the evaluation step.
Therefore, we can safely assume that our observation for the given test problem can be safely carried over to other PDE-based problems of similar or greater size and difficulty.

\subsection{Parallelization Method}
Based on the previous discussion, we can now derive a suitable parallelization scheme for our evolutionary algorithm.
However, while we have already estimated the impact of parallel execution of each of its steps, we have not yet discussed how to parallelize its individual operations on a given number of processing units. 
In general, if each operation within a sequence of computations can be performed independently, which means that it is not affected by the result of any other operation, the sequence is trivially parallelizable.
As this condition is fulfilled for step~2--3 of our evolutionary algorithm, these operations can be performed in a fully parallel manner.
In contrast, in both selection steps of our method, each processing unit needs to access the complete population.
We can thus distinguish two fundamentally different ways to parallelize them on a multi-processor system:
\begin{enumerate}
	\item Duplicate the population on each processing unit.
	\item Split the population into subpopulations and perform the selection on each of them independently while allowing periodic migration between certain subpopulations. This approach is usually described with the term \emph{island-based} evolutionary algorithm.
\end{enumerate}   
The first approach achieves behavioral equivalence to a sequential evolutionary algorithm at the cost that both the memory and computational requirements increase with the population size $\mu$, which restricts its applicability to only medium-sized populations.
In contrast, depending on the amount of migration between the individual subpopulations, which require a certain amount of communication, island-based models can yield higher scalability since all operations are performed on completely independent subpopulations.
On the downside, an island-based approach comprises the risk of selecting a higher percentage of inferior individuals, as only individuals in the respective subpopulation are considered for selection, which might lead to slower convergence compared to its sequential counterpart.
Considering the relatively low cost of selection even compared to the evaluation of a single individual\footnote{If we divide the total evaluation time per generation in Table~\ref{table:evolutionary-search-profiling} by the number of children, we obtain an average evaluation time of 47 seconds per individual.}, we can conclude that a duplication of the whole population is feasible for most experiments performed on small to medium-sized compute clusters.
It is, therefore, the chosen method of parallelization within our implementation.
While this decision theoretically limits scalability, for all experiments considered in this thesis, which do not employ populations larger than 256 individuals, a duplication of the complete population is feasible.
However, an island-based parallelization can be considered a viable extension for potential future applications that require us to execute our algorithm with a significantly larger population.

\subsection{Implementation}
After deriving a suitable parallelization approach, we can now proceed with its implementation as an extension of our previously defined \emph{Optimizer} class shown in Listing~\ref{code:optimization:optimizer}.
For this purpose, we utilize the message-passing interface (MPI), which is available in the form of the Python package \emph{mpi4py}.
While MPI only defines interfaces for the C and Fortran programming languages, mpi4py provides an additional layer of abstractions that enables the exchange of arbitrary Python objects between different processes using Python's Pickle library.
As we have already discussed in Section~\ref{sec:evostencils-part1:productions}, Pickle provides a unified way to serialize and deserialize objects based on a portable binary format, which can then be transmitted using the core functionality provided by MPI.
Since MPI represents the de facto standard for distributed computing and is thus supported on the majority of high-performance computing systems, its usage facilitates the portability of our implementation.
Furthermore, the availability of highly-optimized MPI implementations, which are often developed in cooperation with hardware manufacturers, enables communication between different processors in a highly-efficient manner.
As a first step towards the parallelization of our evolutionary program synthesis method with mpi4py, we extend the previously defined \emph{Optimizer} class found in Listing~\ref{code:optimization:optimizer-mpi} by providing an interface to all required MPI operations. 
\begin{listing}
\begin{minted}{python}
def merge_lists(lst):
    return [item for sublist in lst for item in sublist]

class Optimizer:
    def __init__(self, program_generator, mpi_comm=None, nprocs=1, rank=0):
        self.program_generator = program_generator
        self.pset = None
        self.toolbox = None
        self.params = None
        self.individual_cache = {}
        self.mpi_comm = mpi_comm
        self.nprocs = nprocs
        self.rank = rank

    def is_root(self):
        return self.rank == 0

    def allgather(self, data):
        if self.mpi_comm is None:
            return data
        else:
            return merge_lists(self.mpi_comm.allgather(data))
\end{minted}
	\caption{Optimizer Class -- MPI Extension}
	\label{code:optimization:optimizer-mpi}
\end{listing}
For this purpose, we add the MPI communicator object, the number of processes, and the process rank as additional arguments to the initialization method.
Note that all MPI operations can then be performed solely based on this information and the respective communicator object.
In order to exchange individuals between the processes, we utilize the \emph{allgather} operation, which first collects a list of objects from all processes that are then distributed to each individual process.
As in our case, each of these objects corresponds to a list of individuals, we additionally employ the \emph{merge\_lists} function to merge all sublists into a single list, which then contains the complete set of individuals collected from all processes.
Note that in case only a single process exists, we simply return the passed object without modification, which results in a unified interface for the sequential and parallel execution of our evolutionary algorithm.
The resulting operation is then implemented in Lines~15--19 of Listing~\ref{code:optimization:optimizer-mpi}. 

As a final step, we can now utilize this interface to implement a parallel version of the generalization procedure described in Algorithm~\ref{alg:generalization-procedure} as an extension of the basic implementation of our evolutionary algorithm.
The resulting Python implementation is shown in Listing~\ref{code:optimization:evolutionary-search-mpi}.
\begin{listing}
\begin{minted}{python}
import random

def evolutionary_search(self, params, use_random_search=False, generalization_interval=None):
    if generalization_interval is None:
        generalization_interval = params.generations
    # Generate and evaluate the initial population
    n = params.initial_population_size / self.nprocs
    population = self.toolbox.population(n)
    invalid_ind = [ind for ind in population]
    self.evaluate_individuals(invalid_ind)
    population = self.allgather(population)
    # Select mu individuals from the initial population
    population = self.toolbox.select(population, params.mu_)
    max_level = params.max_level
    for gen in range(1, params.generations + 1):
        if gen > 1 and gen % generalization_interval == 0:
            # Increase the problem size
            max_level += 1
            self.adapt_problem_size(max_level, params, params.multi_objective)
            # Reevaluate individuals
            invalid_ind = [ind for ind in population]
            self.evaluate_individuals(invalid_ind)
        lambda_ = params.lambda_ / self.nprocs
        if use_random_search:
            # For random search simply generate individuals randomly
            offspring = self.toolbox.population(n=lambda_)
        else:
            # Select lambda parents for mutation and crossover
            selected = self.toolbox.select(population, lambda_)
            parents = [self.toolbox.clone(ind) for ind in selected]
            offspring = self.create_offspring(parents, params.crossover_probability)
        # Evaluate new individuals
        invalid_ind = [ind for ind in offspring if not ind.fitness.valid]
        self.evaluate_individuals(invalid_ind)
        offspring = self.allgather(offspring)
        # Select new population from the combined set of parents
        # and offspring (elitism)
        population = self.toolbox.elitism(population + offspring, params.mu_)

    return population
\end{minted}
	\caption{Evolutionary Search Method with Generalization and Parallelization}
	\label{code:optimization:evolutionary-search-mpi}
\end{listing}
Since most of the required functionality is already implemented within the \emph{adapt\_problem\_size} and \emph{allgather} methods, only a few adaptions of our original implementation are required.
First of all, we include an additional argument \emph{generalization\_interval}, which corresponds to the parameter $m$ in Algorithm~\ref{alg:generalization-procedure}.
Based on the value of this argument, we iteratively increase the problem size after a specified number of generations, which is implemented in line~15--21.
Finally, we parallelize the individual steps of our evolutionary algorithm using the previously defined MPI interface.
According to the MPI programming model, if the number of processes is larger than one, each process executes its own instance of the same program.
We thus only have to implement the required synchronization points between the processes.
Similar to our original implementation, the first step within our evolutionary program synthesis method is the generation and evaluation of an initial population.
In line~6--9, each process first generates and evaluates its respective fraction of the initial population, which is then combined and distributed using the \emph{allgather} method in line~10.
In a similar fashion, each process generates its fraction of the offspring by first selecting the respective number of parents in line~28, based on which new individuals are created in line~30 using crossover and mutation.
Again, after each process has finished the evaluation of its local individuals, they are combined using the \emph{allgather} method in line~34.
Since after this operation, each process has an exact copy of all newly created individuals, the subsequent elitist selection in line~37 consistently yields the same population, based on which the algorithm proceeds until the maximum number of generations has been reached.
With the implementation of a scalable evolutionary program synthesis method that can leverage the compute capabilities of current multi-node systems, we can now evaluate the effectiveness of our approach for the grammar-based design of multigrid methods in a number of experiments, where we consider two common PDE-based model problems, Poisson's equation, and a linear elastic boundary value problem.
As a final evaluation step, we then assess the efficiency and generalizability of the multigrid methods obtained with the generalization procedure described in Algorithm~\ref{alg:generalization-procedure} on a difficult benchmark problem, the indefinite Helmholtz equation with large wavenumbers.

% TODO Maybe include this section again
%\section{Systems of Partial Differential Equations}

%% file: contents/experiments.tex
As a first step in the evaluation of our evolutionary program synthesis method, we consider two PDE-based model problems, Poisson's equation and a linear elastic boundary value problem, which can already be solved efficiently by applying common multigrid cycles iteratively.
Here our goal is to demonstrate that our approach is able to reliably find functioning multigrid cycles in a number of randomized experiments.
Furthermore, since classical multigrid cycles already provide a strong baseline for these problems, we can investigate whether the methods designed with our approach are able to achieve a similar degree of efficiency.
In addition, by considering both two- and three-dimensional problems as well as a system of PDEs, we can demonstrate that our implementation is able to handle PDEs of different types.
\section{Multigrid Cycles for Solving Common Partial Differential Equations}
\label{sec:experiments-part1}
The goal of this section is to evaluate the effectiveness of our evolutionary program synthesis method for the automated design of multigrid cycles that are used as an iterative method for solving a discretized PDE.
Therefore, the problem instances considered in this section are chosen in a way that facilitates their efficient solution by multigrid.
This is reflected in the fact that common multigrid cycles, as those described in Section~\ref{sec:multigrid-cycles}, are able to quickly converge to the correct solution of each of the resulting systems of linear equations.
We begin this section by introducing the considered problem instances and their mathematical formulation.
At this point, we want to emphasize that all results presented in this section have originally been published in~\cite{schmitt2021evostencils}.
However, this thesis complements this work with an additional analysis of the multigrid solvers discovered by our evolutionary algorithm.
\subsection{Problem Formulation}
\subsubsection{Poisson's Equation}
\label{sec:poisson-equation}
Poisson's equation is an elliptic PDE that occurs in the study of many physical phenomena~\cite{folland2020introduction} and is defined as
\begin{equation}
	\begin{split}
		-\nabla^{2} u & = f \quad \text{in} \; \Omega \\
		u & = g \quad \text{on} \; \partial \Omega.
	\end{split}
	\label{eq:poisson}
\end{equation}
In our experimental evaluation, we consider two different instances of Poisson's equation with Dirichlet boundary conditions, which are summarized in Table~\ref{table:poisson-problems}.
\begin{table}
	\caption{Considered instances of Poisson's equation.}
	\begin{tabular}{r l l}
		\toprule
		Problem & 2D Poisson & 3D Poisson \\
		\midrule
		$\Omega = $ & $ (0, 1)^2$ & $(0, 1)^3$ \\
		\midrule
		$f(\bm{x}) = $ & $\pi^2 \cos(\pi x) - 4 \pi^2 \sin(2 \pi y)$ & $x^2 - 0.5 y^2 - 0.5 z^2$ \\
		\midrule
		$g(\bm{x}) = $ & $\cos(\pi x) - \sin(\pi y)$ & $0$ \\
		\bottomrule
	\end{tabular}
	\label{table:poisson-problems}
\end{table}
Note that in Section~\ref{sec:execution-time-analysis}, we have employed the same two-dimensional instance of Poisson's equation to estimate the relative cost of each operation within our evolutionary algorithm.
We discretize the Laplace operator $\nabla^{2}$ with finite differences on a uniform cartesian grid with step size $h = 1/2^{l_{max}}$, which yields the five-point stencil
\begin{equation*}
	\nabla^2_h = 
	\frac{1}{h^2} \begin{bmatrix}
		0 & 1 & 0\\
		1 & -4 & 1 \\
		0 & 1 & 0  
	\end{bmatrix},
\end{equation*}
in two dimensions, and the seven-point stencil
\begin{equation*}
\nabla^2_h = 
\frac{1}{h^2} \begin{bmatrix}
	\begin{bmatrix}
	0 & 0 & 0 \\
	0 & 1 & 0 \\
	0 & 0 & 0
	\end{bmatrix}
	&		
	\begin{bmatrix}
	0 & 1 & 0 \\
	1 & -6 & 1 \\
	0 & 1 & 0 
	\end{bmatrix} &
	\begin{bmatrix}
	0 & 0 & 0 \\
	0 & 1 & 0 \\
	0 & 0 & 0
\end{bmatrix}
\end{bmatrix}
\end{equation*} in three dimensions.
We choose a maximum level of $l_{max} = 11$ in 2D and $l_{max} = 7$ in 3D, which results in systems of linear equations consisting of $4\,190\,209$ and $2\,048\,383$ unknowns, respectively.

\subsubsection{Linear Elasticity}
\begin{figure}
	\centering
	\includegraphics[width=0.5\textwidth]{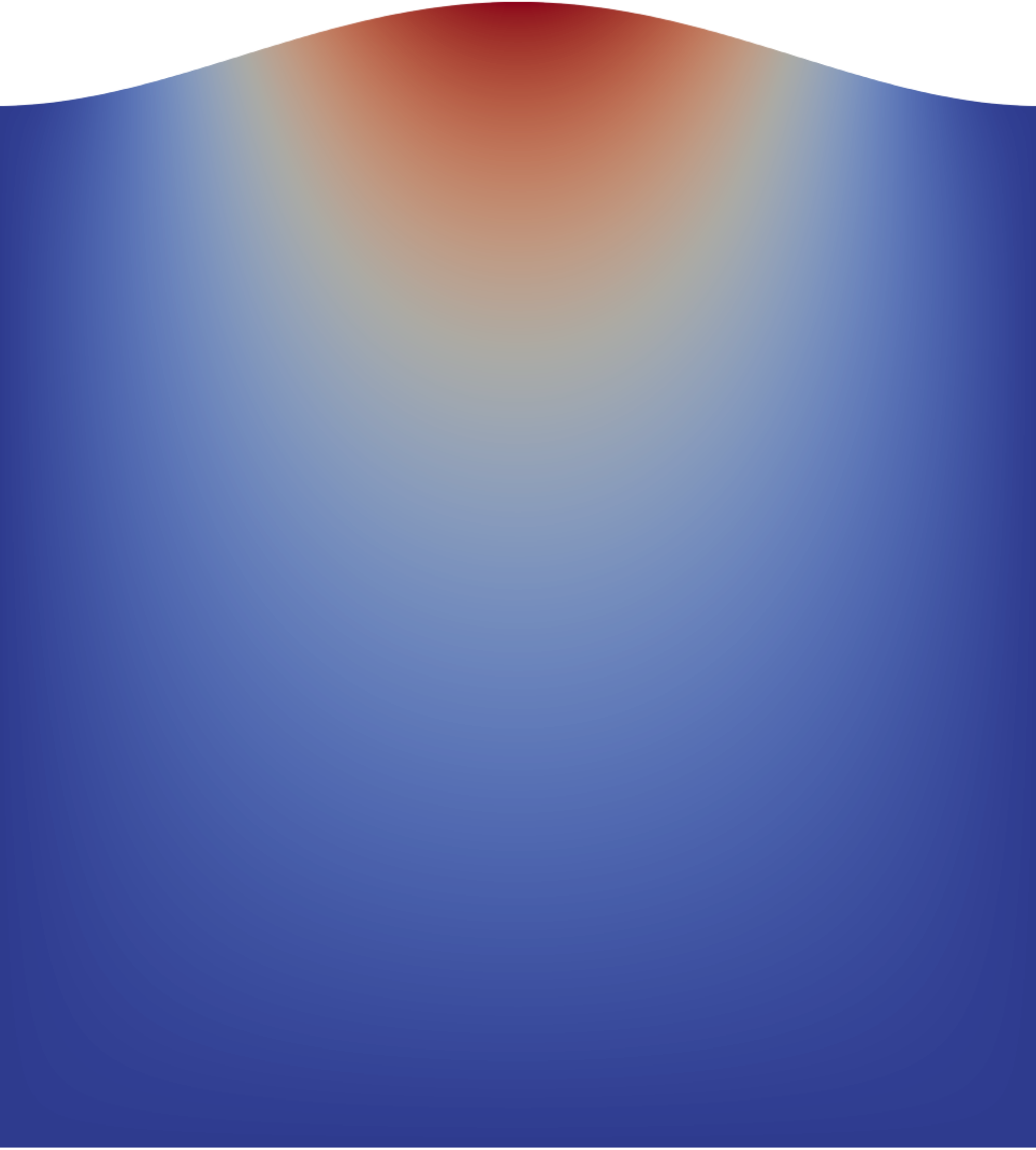}
	\caption[Visualization of the considered linear elastic boundary value problem]{Visualization of the considered linear elastic boundary value problem. A two-dimensional rectangular body undergoes an elastic deformation into the y-direction.}
	\label{fig:visualization-linear-elasticity}
\end{figure}

Linear elasticity is a fundamental branch of solid mechanics with numerous applications in engineering and material science~\cite{holzapfel2001nonlinear}.
It is derived from the more general theory of nonlinear continuum mechanics by assuming a linear relationship between stress and strain during elastic deformation.
We consider a two-dimensional linear elastic boundary value problem given by the system of PDEs
\begin{equation}
	\begin{split}
		(\alpha + \beta) \cdot (\frac{\partial^2}{\partial x^2} u + \frac{\partial^2}{\partial x \partial y} v) + \alpha \nabla^2 u & = 0 \quad \text{in} \; \Omega \\
		(\alpha + \beta) \cdot (\frac{\partial^2}{\partial x \partial y} u + \frac{\partial^2}{\partial y^2} v) + \alpha \nabla^2 v & = 0 \quad \text{in} \; \Omega \\
		u = 0 \quad \text{and} \quad v & = g \quad \text{on} \; \partial \Omega 
		\label{eq:linear-elasticity}
	\end{split}
\end{equation}
where $\Omega = (0,1)^2$, $\alpha = 195$, $\beta = 130$ and
\begin{equation*}
	g(x,y) = 0.4 \, (1 - x) \, x y \, \sin(\pi x).
\end{equation*}
From a physical point of view, this system represents a two-dimensional rectangular body that undergoes an elastic deformation into the y-direction, as it can be seen in Figure~\ref{fig:visualization-linear-elasticity}.
We discretize Equation \eqref{eq:linear-elasticity} with finite differences on a cartesian grid using a step size $h = 1/2^{10}$ such that $l_{max} = 10$, which yields a system of linear equations $A \bm{u} = \bm{f}$ with 
\begin{equation*}
	A =
	\begin{pmatrix}
		(\alpha + \beta) \frac{\partial^2}{\partial x^2} + \alpha \nabla^2 & (\alpha + \beta) \frac{\partial^2}{\partial x \partial y} \\
		(\alpha + \beta) \frac{\partial^2}{\partial x \partial y} & (\alpha + \beta) \frac{\partial^2}{\partial y^2} +  \alpha \nabla^2
	\end{pmatrix},
\end{equation*}
\begin{equation*}
	\bm{u} = \begin{pmatrix}
		u \\ v
	\end{pmatrix}, \quad
	\bm{f} =
	\begin{pmatrix}
		f_{u} \\ f_{v}
	\end{pmatrix} =
	\begin{pmatrix}
		0 \\ 0
	\end{pmatrix},
\end{equation*}
whereby the differential operators $\nabla^2$, $\frac{\partial^2}{\partial x^2}$, $\frac{\partial^2}{\partial y^2}$ and $\frac{\partial^2}{\partial x \partial y}$ are approximated by their discrete counterparts
\begin{equation*}
	\left(\nabla^2 u\right)_{i,j} = 
	\frac{1}{h^2} \begin{bmatrix}
		0 & 1 & 0\\
		1 & -4 & 1 \\
		0 & 1 & 0  
	\end{bmatrix},
\end{equation*}
\begin{equation*}
	\left(\frac{\partial^2}{\partial x^2} u\right)_{i,j} =
		\frac{1}{h^2} \begin{bmatrix}
		0 & 0 & 0\\
		1 & -2 & 1 \\
		0 & 0 & 0  
	\end{bmatrix},
\end{equation*}
\begin{equation*}
	\left(\frac{\partial^2}{\partial y^2} u\right)_{i,j} =
	\frac{1}{h^2} \begin{bmatrix}
		0 & 1 & 0\\
		0 & -2 & 0 \\
		0 & 1 & 0  
	\end{bmatrix},
\end{equation*}
\begin{equation*}
	\left(\frac{\partial^2}{\partial x \partial y} u\right)_{i,j} = 
	\frac{1}{4 h^2} \begin{bmatrix}
		-1 & 0 & 1\\
		0 & 0 & 0 \\
		1 & 0 & -1  
	\end{bmatrix}.
\end{equation*}
Similar to the above case, we employ a uniform cartesian grid with a step size $h = 1/2^{l_{max}}$ with $l_{max} = 10$, such that the resulting system of linear equations contains $2\,093\,058$ unknowns.
\subsection{Solver Configuration}
\label{sec:experiments1-multigrid-configuration}
To design an efficient multigrid method, we consider each of the given problems on a grid hierarchy consisting of five discretization levels 
\begin{equation*}
    l \in \left[l_{max} - 4, l_{max}\right],
\end{equation*} where the grid spacing on each level is given by the formula $h = 1/2^l$.
We then obtain the respective operator on each level by applying the same discretization method as on the finest grid.
Therefore, the resulting grammar is structurally similar to the one shown in Algorithm~\ref{table:multigrid-grammar}.
Within each grammar production, we consider the following components:
\begin{description}
	\item[\textbf{Smoothers}:] Decoupled / Collective Jacobi and red-black Gauss-Seidel (RB-GS), block Jacobi with rectangular blocks up to a maximum number of six terms~\cite{trottenberg2000multigrid}.
	\item[\textbf{Restriction}:] Full-weighting restriction.
	\item[\textbf{Prolongation}:] Bilinear interpolation.
	\item[\textbf{Relaxation factors}:] $\omega \in \left( 0.1 + 0.05i \right)_{i = 0}^{36} = \left(0.1, 0.15, 0.2, \dots, 1.9 \right)$
	\item[\textbf{Coarse-grid solver}:] Conjugate gradient method, in case $l = l_{max} - 4$.
\end{description}
Here we generate block Jacobi smoothers by defining a splitting $A = L + D + U$ where $D$ is a block diagonal matrix of the form
\begin{equation*}
	D = 
		\begin{pmatrix}A_{11}&0&\cdots &0\\
			0&A_{22}&\cdots &0\\
			\vdots &\vdots &\ddots &\vdots \\0&0&\cdots &A_{mm}\end{pmatrix},
\end{equation*}
where each matrix $A_{ij}$ corresponds to a set of adjacent grid points contained in the respective rectangular block, as it has been discussed in Section~\ref{subsec:block-smoothing}.
A more detailed treatment of block smoothers can be found in~\cite{trottenberg2000multigrid}.
For each smoothing and coarse-grid correction step, the relaxation factor $\omega$ can be chosen from the above interval.
As a baseline for assessing the efficiency and generalizability of the designed multigrid solvers, we consider a number of common multigrid cycles with RB-GS smoothing and optimized relaxation factors.
We hence formulate these methods on the same five-grid hierarchy using the same restriction, prolongation operators, and coarse-grid solver.
In each case, we consider the corresponding linear system as solved when the initial residual has been reduced by a factor of $10^{-12}$.

\subsection{Experimental Settings and Evaluation Platform}
\label{sec:optimization-settings}
After specifying the operator and parameter choices considered within the construction of each multigrid solver, we next describe the settings under which we perform each experiment.
Here we utilize the EvoStencils framework, whose implementation has been described in detail in Chapter~\ref{chapter:evostencils-1} and~\ref{chapter:evostencils-2}.
The goal of each experiment is to evolve a set of non-dominated individuals according to the two objectives convergence factor $\rho$ and execution time per iteration $t$, as described in Section~\ref{sec:fitness-evaluation-and-selection}, which are evaluated by applying each multigrid cycle as an iterative solver to the respective test problem.
The resulting individuals are then subject to a subsequent evaluation and comparison with the available reference methods. 
Table~\ref{table:gp-parameters} gives an overview of the algorithmic configuration used within each experiment.
\begin{table}
	\centering
	\caption{Summary of the genetic programming configuration parameters.}
	\label{table:gp-parameters}
	\begin{tabular}{l c}
		\toprule
		Parameter & Value \\
		\midrule 
		Evolutionary algorithm type & $(\mu + \lambda)$ \\
		\midrule
		Objectives & $t, \rho$ \\
		\midrule
		Number of generations & 250 \\
		\midrule
		Initial population size & 2048 \\
		\midrule
		$\lambda$ & 256 \\
		\midrule
		$\mu$ & 256 \\
		\midrule
		Number of MPI processes & 64 \\
		\midrule
		Non-dominated sorting procedure & \cite{deb2002fast} \\ 
		\midrule
		Selection operator & \cite{deb2002fast} \\ 
		\midrule
		Crossover operator & Subtree crossover \\
		\midrule
		Crossover probability & $2/3$ \\
		\midrule
		Mutation operator & Random subtree insertion \\
		\midrule 
		Probability to mutate a terminal symbol & $1/3$ \\
		\bottomrule
	\end{tabular}
\end{table}
In each run, we perform an evolutionary search for 250 generations starting with a randomly generated population of 2048 individuals.
In each generation, we create new individuals by first selecting $\lambda = 256$ candidates from the current population.
We then apply mutation and crossover to each pair of selected candidates to create two child individuals, whereby the crossover probability is set to $2/3$, and in the case of mutation, we choose a terminal symbol with a probability of $1/3$.
As a mutation operator, we employ subtree insertion whenever possible. Otherwise, the respective subtree is completely replaced by a randomly-created one, as described in Section~\ref{sec:gggp-mutation-and-recombination}.
The resulting individuals are then evaluated according to the two objectives by generating a parallel C++ solver implementation using the ExaStencils framework, which is applied to the respective problem as described above.
Hereby we distribute the evaluation of all 256 individuals to 64 MPI processes, such that each process is responsible for the evaluation of four individuals.
The resulting fitness values are distributed to all 64 processes, such that they possess an identical copy of each child individual together with its fitness value, as it has been described in Chapter~\ref{sec:distributed-parallelization}.
Finally, we select $\mu = 256$ individuals as a population for the next generation from the combined set of parent and child individuals using the NSGA-II non-dominated sorting procedure~\cite{deb2002fast}.

As an evaluation platform for running each experiment, we employ 32 nodes of the Meggie compute cluster of the Erlangen National High-Performance Computing Center (NHR), where each node of the system consists of two sockets with ten physical CPU cores.
Each process is thus pinned and executed on a dedicated socket. 
For the evaluation of each solver, we employ a thread-based parallelization using ten OpenMP threads, whereby each tread is pinned to a distinct physical compute core on the respective socket.
For the parallelization of each nested loop within the generated solver implementation, we employ static scheduling based on the outer-loop range, such that each thread processes a consecutive chunk of iterations.
To generate a thread-parallel executable of each solver, we employ the GCC 9.3.0 compiler with the -O3 optimization level.
To reduce statistical variations between individual evaluation runs, we execute each solver three times and then compute the average for both objectives.

\subsection{Analysis of the Evolutionary Algorithm}
\label{sec:experiments1-algorithm-behavior-analysis}
\begin{figure}[h]
	\centering
	\begin{subfigure}[b]{0.49\textwidth}
		\centering
		\includegraphics[width=\textwidth]{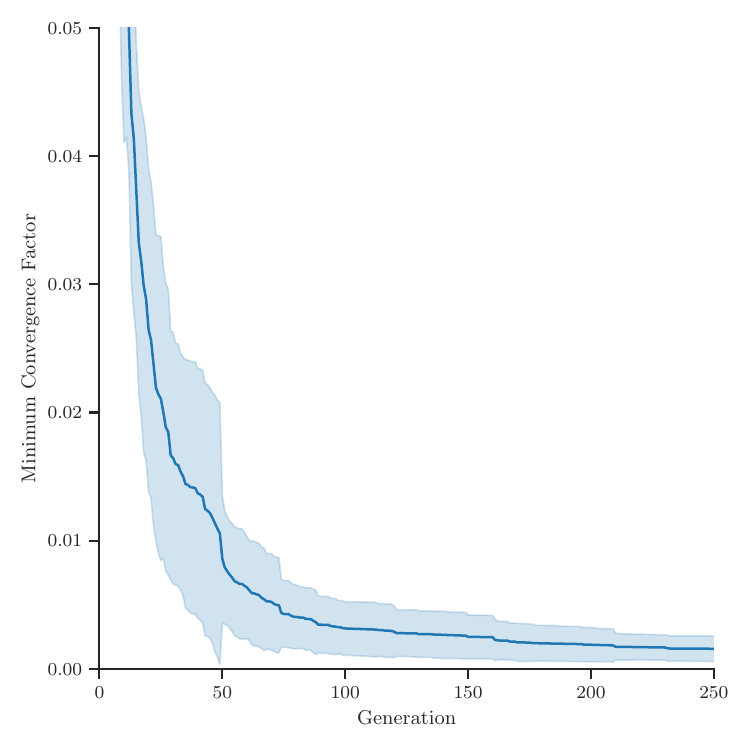}
		\caption{Minimum Convergence Factor}
		\label{fig:poisson-2D-minimum-convergence-factor}
	\end{subfigure}
	\hfill
	\begin{subfigure}[b]{0.49\textwidth}
		\centering
		\includegraphics[width=\textwidth]{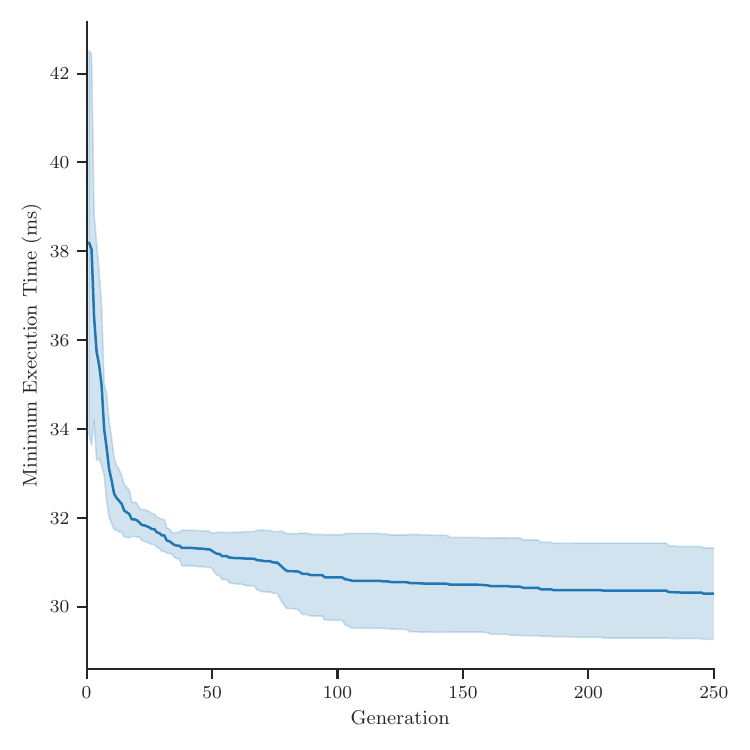}
		\caption{Minimum Execution Time per Iteration}
		\label{fig:poisson-2D-minimum-execution-time}
	\end{subfigure}
	\caption[2D Poisson -- Mean and standard deviation of the minimum objective function values]{2D Poisson -- Mean and standard deviation of the minimum objective function values of all evolutionary algorithm runs.}
	\label{fig:poisson-2D-minimum-objectives}
\end{figure}

\begin{figure}
	\centering
	\begin{subfigure}[b]{0.49\textwidth}
		\centering
		\includegraphics[width=\textwidth]{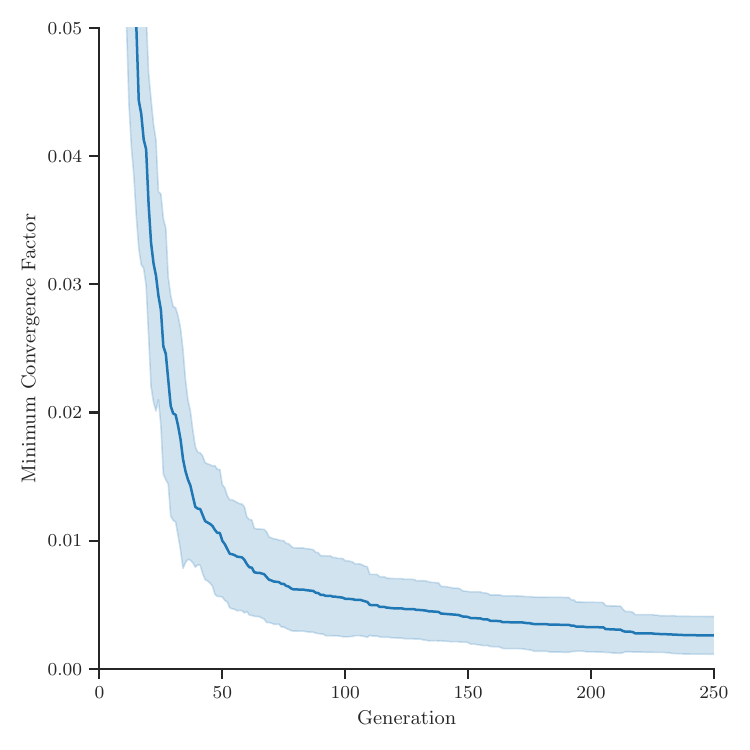}
		\caption{Minimum Convergence Factor}
		\label{fig:poisson-3D-minimum-convergence-factor}
	\end{subfigure}
	\hfill
	\begin{subfigure}[b]{0.49\textwidth}
		\centering
		\includegraphics[width=\textwidth]{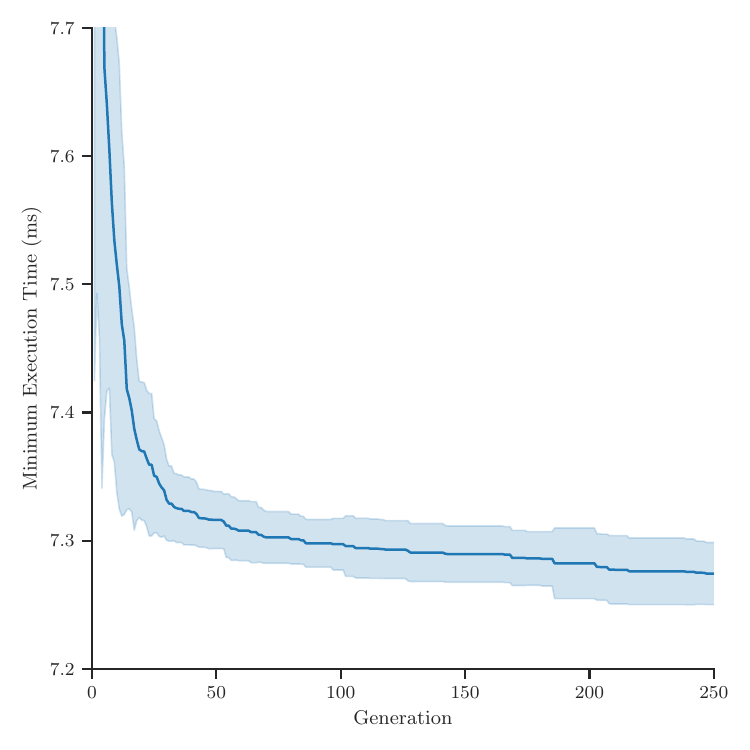}
		\caption{Minimum Execution Time per Iteration}
		\label{fig:poisson-3D-minimum-execution-time}
	\end{subfigure}
	\caption[3D Poisson -- Mean and standard deviation of the minimum objective function values]{3D Poisson -- Mean and standard deviation of the minimum objective function values of all evolutionary algorithm runs.}
	\label{fig:poisson-3D-minimum-objectives}
\end{figure}

\begin{figure}
	\centering
	\begin{subfigure}[b]{0.49\textwidth}
		\centering
		\includegraphics[width=\textwidth]{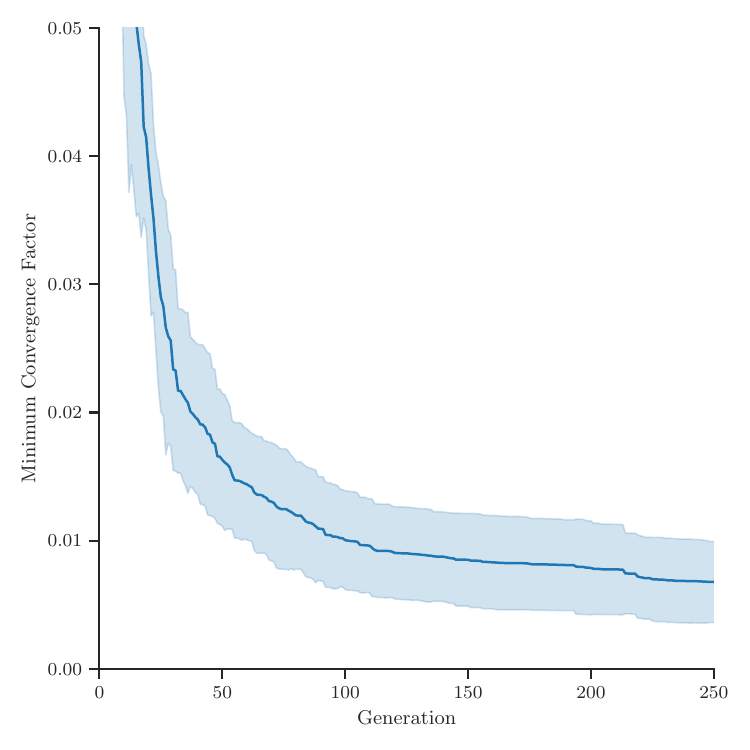}
		\caption{Minimum Convergence Factor}
		\label{fig:linear-elasticity-2D-minimum-convergence-factor}
	\end{subfigure}
	\hfill
	\begin{subfigure}[b]{0.49\textwidth}
		\centering
		\includegraphics[width=\textwidth]{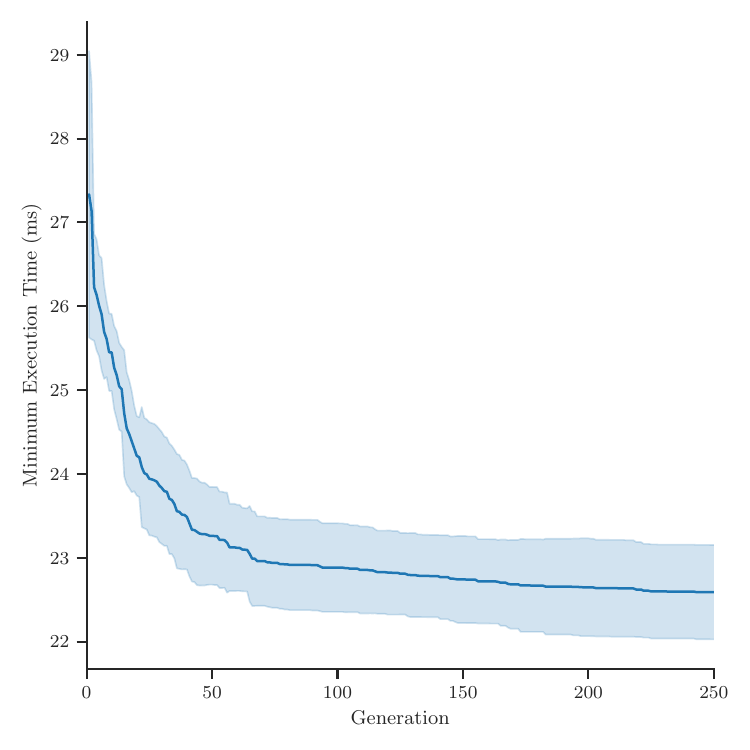}
		\caption{Minimum Execution Time per Iteration}
		\label{fig:linear-elasticity-2D-minimum-execution-time}
	\end{subfigure}
	\caption[2D Linear Elasticity -- Mean and standard deviation of the minimum objective function values]{2D Linear Elasticity -- Mean and standard deviation of the minimum objective function values of all evolutionary algorithm runs.}
	\label{fig:linear-elasticity-2D-minimum-objectives}
\end{figure}

As a first step towards a quantitative evaluation of our evolutionary program synthesis method, we assess whether our algorithm is able to effectively find good solutions with respect to our two optimization objectives.
For this purpose, we measure the minimum of each of the two objectives within the population throughout each of the ten experiments for all three test problems.
As a result, Figure~\ref{fig:poisson-2D-minimum-objectives},~\ref{fig:poisson-3D-minimum-objectives} and~\ref{fig:linear-elasticity-2D-minimum-objectives} shows the mean and standard deviation of the current optimum of both objectives over all experiments.
First of all, we can conclude that in all three cases, our algorithm is able to quickly reduce the minimum value of both objectives within the population.
If we compare the slope of the two objectives for the three different problems, we can see that in the case of the first objective $\rho$, significantly more generations are required to achieve the same degree of reduction as for the second objective $t$.
However, the algorithm still significantly improves upon the current minimum of $\rho$ beyond the first 50 generations, which can be seen in Figure~\ref{fig:poisson-2D-minimum-convergence-factor},~\ref{fig:poisson-3D-minimum-convergence-factor} and~\ref{fig:linear-elasticity-2D-minimum-convergence-factor}.
For the second objective $t$, the majority of improvement happens within the first 50--70 generations, as it can be seen in Figure~\ref{fig:poisson-2D-minimum-execution-time},~\ref{fig:poisson-3D-minimum-execution-time} and~\ref{fig:linear-elasticity-2D-minimum-execution-time}.
Note that while the convergence factor is constant for each execution of the same solver, $t$ is obtained by measuring its execution time on the respective compute node.
Therefore, due to manufacturing and temperature-based variations, we can expect a certain degree of fluctuations when measuring the execution time of the same solver on different compute nodes during consecutive runs.
Consequently, even though the algorithm is able to reduce the second objective faster, this effect induces a larger deviation between the individual experiments and thus an overall higher standard deviation.

As a next step, we can assess the difficulty of the underlying problem by considering the absolute value of the minimum convergence factor attained in each of the three different cases.
While the execution time per iteration $t$ is solely determined by the computational complexity of each solver and the properties of the given computer architecture, a smaller convergence factor indicates that the underlying problem is easier to solve.
Poisson's equation represents an often-studied model problem whose strong ellipticity enables its easy solution by multigrid~\cite{trottenberg2000multigrid}.
As a consequence, our method consistently discovers multigrid methods that achieve fast convergence, with minimum convergence factors of less than 0.005, for both the two and three-dimensional instances of this equation, which can be seen in Figure~\ref{fig:poisson-2D-minimum-convergence-factor} and~\ref{fig:poisson-3D-minimum-convergence-factor}.
In the case of the linear elastic boundary value problem, both the mean and standard deviation of the first objective remains higher throughout each experiment.
However, on average, our method is still able to discover multigrid methods that achieve a convergence factor of 0.01 or less, leading to extremely quick convergence.
In summary, we can conclude that our evolutionary program synthesis method consistently yields multigrid cycles that represent satisfactory minima for both objectives in all three test cases considered.

\begin{figure}
\centering
	\includegraphics[scale=0.725]{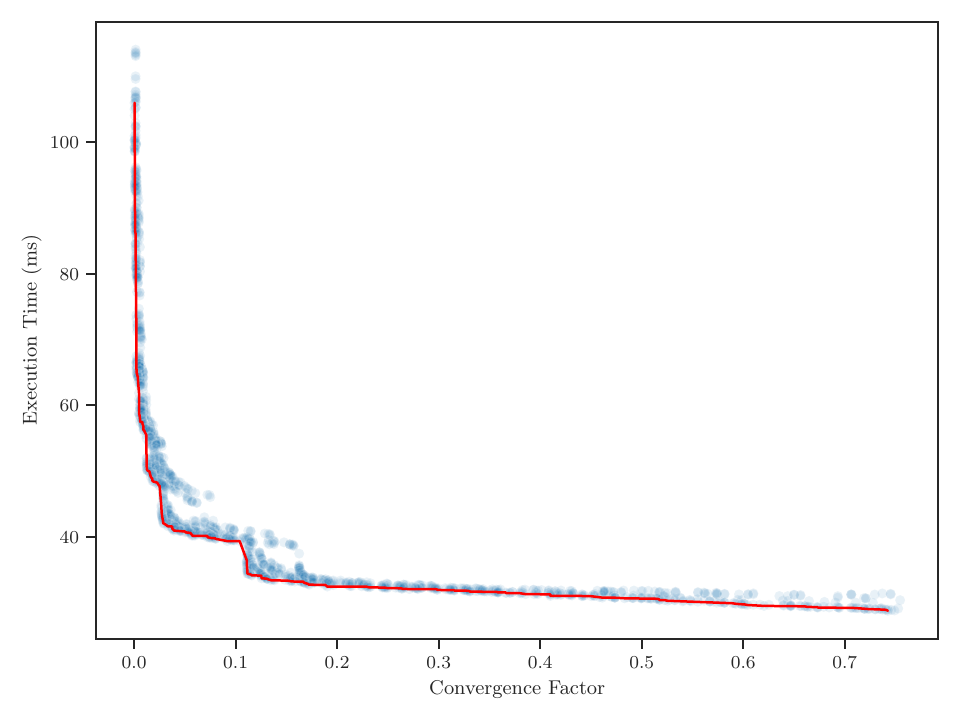}
	\caption[2D Poisson -- Distribution of non-dominated individuals at the end of all ten experiments]{2D Poisson -- Distribution of non-dominated individuals at the end of all ten experiments. The red line denotes the combined Pareto front.}
	\label{fig:pareto-front-2D-poisson}
\end{figure}

\begin{figure}
\centering
\includegraphics[scale=0.725]{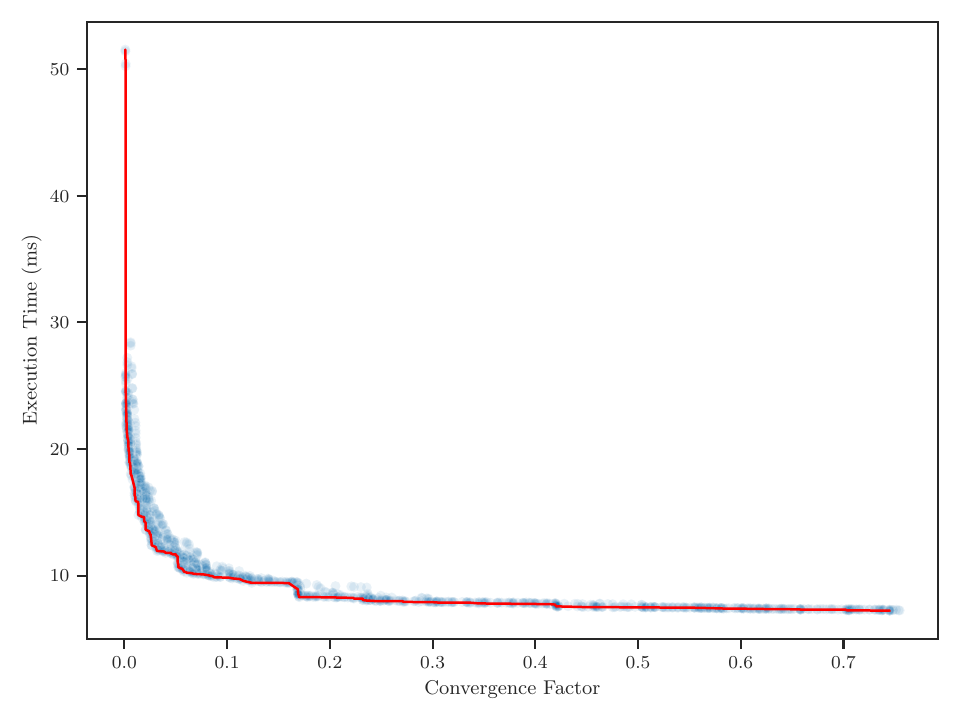}
	\caption[3D Poisson -- Distribution of non-dominated individuals at the end of all ten experiments]{3D Poisson -- Distribution of non-dominated individuals at the end of all ten experiments. The red line denotes the combined Pareto front.}
	\label{fig:pareto-front-3D-poisson}
\end{figure}

\begin{figure}
\centering
\includegraphics[scale=0.725]{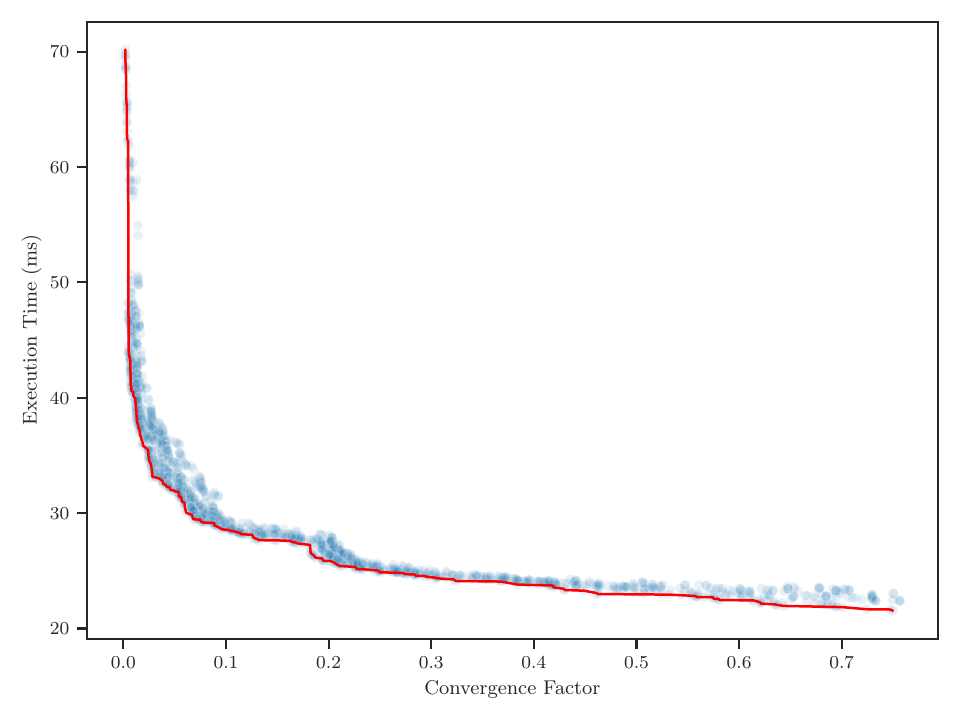}
	\caption[2D Linear Elasticity -- Distribution of non-dominated individuals at the end of all ten experiments]{2D Linear Elasticity -- Distribution of non-dominated individuals at the end of all ten experiments. The red line denotes the combined Pareto front.}
	\label{fig:pareto-front-2D-linear-elasticity}
\end{figure}

To further analyze the behavior of our multi-objective evolutionary algorithm, we consider the distribution of non-dominated individuals at the end of all experiments, which is shown in Figure~\ref{fig:pareto-front-2D-poisson},~\ref{fig:pareto-front-3D-poisson} and~\ref{fig:pareto-front-2D-linear-elasticity}.
Here the red line denotes the combined front, while the color density of the distribution indicates where most of the individuals are located at the end of all experiments. 
In all three cases, the majority of non-dominated individuals obtained within each experiment can be found close to the combined front.
While in Figure~\ref{fig:pareto-front-2D-poisson}, the number of individuals that are distinctly located outside of the front is slightly higher than in the other two cases, their distance to the red line is still comparably small compared to the complete objective space.
Again we can attribute part of this effect to variations in the execution time on different compute nodes.
In the majority of the objective space, but especially in its center, the individuals are evenly distributed alongside the front.
Here only the left-upper part of Figure~\ref{fig:pareto-front-3D-poisson} and~\ref{fig:pareto-front-2D-linear-elasticity} represents a noteworthy exception, where our approach struggles to discover the same individuals within each experiment.
An explanation for this effect is that the individuals located in this part of the space are characterized by an extremely fast convergence.
Since the convergence of a multigrid method can primarily be accelerated with the addition of smoothing and coarse-grid correction steps, discovering the same non-dominated individuals with fast convergence requires us to construct the same large expressions starting from different random initializations.
Recently, the impaired capabilities of NSGA-II to evolve large non-dominated expressions have been analyzed in~\cite{liu2022evolvability}, which similarly explains why our implementation struggles to discover the same fast-converging solvers consistently.

\subsection{Evaluation of the Evolved Multigrid Methods}
Finally, in order to investigate whether our evolutionary program synthesis approach yields multigrid methods that are competitive with well-known multigrid cycles, we consider two different multi-core CPU architectures for evaluation: Intel Xeon E5-2630v4 (Broadwell) and Intel Xeon 2660v2 (Ivy Bridge).
In both cases, each compute node consists of two sockets with 20 physical CPU cores and a cache-coherent NUMA architecture.
In order to assess the solving speed of the discovered methods, we consider two different problem sizes for each of the three test cases.
While in the first case, the problem size is identical to the one employed within our evolutionary algorithm, in the second case, we consider a larger problem instance by doubling the number of grid points in each dimension.
To measure the solving speed of each method, we execute it on a dedicated node using a thread-based parallelization with 20 OpenMP threads, whereby we pin each thread to a unique physical core and employ the same parallelization approach as described in Section~\ref{sec:optimization-settings}.
Even though generalizability is not the focus of this section, it is still worthwhile to investigate whether the methods found with our approach are capable of solving larger instances of the same problem.
For comparison, we consider a number of different V-cycles with at most four RB-GS pre- and post-smoothing steps.
With the exception of full-multigrid methods (FMG), which require a different formulation than the one considered here, these cycles are known to be the fastest multigrid-based solvers for Poisson's equation~\cite{trottenberg2000multigrid}.
As we have already investigated in~\cite{schmitt2020constructing}, the same is true for the linear elastic boundary value problem considered here.
To achieve a fair comparison, we determine the optimum relaxation factor $\omega$ for each test case from the same interval considered within our experiments, which leads to $\omega = 1.15$ for the two-dimensional Poisson equation and $\omega = 1.25$ both for the three-dimensional Poisson equation and the linear elastic boundary value problem.
Table~\ref{table:poisson-2D-reference-methods},~\ref{table:poisson-3D-reference-methods} and~\ref{table:linear-elasticity-2D-reference-methods} contain the resulting solving times and the number of iterations required to achieve the desired defect reduction for each test case.
\begin{table}
	\caption[2D Poisson -- Number of iterations and solving times of the reference methods]{2D Poisson -- Measured number of iterations and solving times of the reference methods on 20 cores and two sockets.}
	\label{table:poisson-2D-reference-methods}
	\centering
	\begin{tabular}{l c c c c c c}
		\toprule
		& \multicolumn{2}{c}{Iterations} & \multicolumn{2}{c}{Broadwell (ms)} & \multicolumn{2}{c}{Ivy Bridge (ms)} \\
		\cmidrule(r){2-3} \cmidrule(r){4-5} \cmidrule(r){6-7}
		$l_{max}$ & $11$& $12$ & $11$ & $12$ & $11$ & $12$\\
		\midrule
		V(1, 0) & 21 & 21 & 969 & 2810 & 879 & 2652 \\
		\midrule
		V(1, 1) & 9 & 9 & 461 & 1359 & 411 & 1287 \\
		\midrule
		V(2, 1) & 7 & 7 & 377 & 1137 & 334 & 1087\\
		\midrule
		V(2, 2) & 6 & 6 & 344 & 1056 & 302 & 1007 \\
		\midrule
		V(3, 2) & 6 & 6 & 378 & 1160 & 324 & 1112 \\
		\midrule
		V(3, 3) & 6 & 6 & 397 & 1255 & 344 & 1201 \\
		\midrule
		V(4, 3) & 6 & 6 & 425 & 1350 & 366 & 1306 \\
		\midrule
		V(4, 4) & 6 & 6 & 448 & 1449 & 383 & 1409\\
		\bottomrule
	\end{tabular}
	%FMG for max_level = 12: 530.3038500000001 ms, 3 Iterations
\end{table}
\begin{table}
	\caption[3D Poisson -- Number of iterations and solving times of the reference methods]{3D Poisson -- Measured number of iterations and solving times of the reference methods on 20 cores and two sockets.}
	\label{table:poisson-3D-reference-methods}
	\centering
	\begin{tabular}{l c c c c c c}
		\toprule
		& \multicolumn{2}{c}{Iterations} & \multicolumn{2}{c}{Broadwell (ms)} & \multicolumn{2}{c}{Ivy Bridge (ms)} \\
		\cmidrule(r){2-3} \cmidrule(r){4-5} \cmidrule(r){6-7}
		$l_{max}$ & $7$& $8$ & $7$ & $8$ & $7$ & $8$\\
		\midrule
		V(1, 0) & 29 & 30 & 121.3 &1221 & 134.6 & 1470 \\
		\midrule
		V(1, 1) & 13 & 13 & 70.8 & 682 & 79.9 & 838 \\
		\midrule
		V(2, 1) & 9 & 9 & 59.0 & 582 & 66.2 & 708 \\
		\midrule
		V(2, 2) & 7 & 7 & 54.6 & 531 & 65.4 & 654 \\
		\midrule
		V(3, 2) & 7 & 7 & 61.9 & 610 & 74.6 & 757 \\
		\midrule
		V(3, 3) & 7 & 7 & 72.6 & 690 & 86.6 & 857 \\
		\midrule
		V(4, 3) & 7 & 6 & 77.9 & 656 & 87.3 & 825 \\
		\midrule
		V(4, 4) & 6 & 6 & 73.2 & 725 & 82.5 & 906 \\
		\bottomrule
	\end{tabular}
\end{table}
\begin{table}
	\caption[2D Linear Elasticity -- Number of iterations and solving times of the reference methods]{2D Linear Elasticity -- Measured number of iterations and solving times of the reference methods on 20 cores and two sockets.}
	\label{table:linear-elasticity-2D-reference-methods}
	\centering
	\begin{tabular}{l c c c c c c}
		\toprule
		& \multicolumn{2}{c}{Iterations} & \multicolumn{2}{c}{Broadwell (ms)} & \multicolumn{2}{c}{Ivy Bridge (ms)} \\
		\cmidrule(r){2-3} \cmidrule(r){4-5} \cmidrule(r){6-7}
		$l_{max}$ & $10$& $11$ & $10$ & $11$ & $10$ & $11$\\
		\midrule
		V(1, 0) & 32 & 31 & 872 & 4306 & 828 & 4128 \\
		\midrule
		V(1, 1) & 15 & 15 & 439 & 2118 & 418 & 2075\\
		\midrule
		V(2, 1) & 10 & 10 & 318 & 1529 & 312 & 1529 \\
		\midrule
		V(2, 2) & 9 & 9 & 314 & 1449 & 316 & 1476 \\
		\midrule
		V(3, 2) & 8 & 8 & 297 & 1368 & 304 & 1388 \\
		\midrule
		V(3, 3) & 7 & 7 & 283 & 1247 & 288 & 1288 \\
		\midrule
		V(4, 3) & 7 & 7 & 293 & 1320 & 313 & 1397 \\
		\midrule
		V(4, 4) & 7 & 7 & 311 & 1378 & 334 & 1471 \\
		\bottomrule
	\end{tabular}
\end{table}
Here, for instance, the abbreviation V(2,1) denotes a V-cycle with two pre- and one post-smoothing step of RB-GS.
First of all, as we can expect from a functioning multigrid method, the number of iterations stays constant for both problem sizes, only with the exception of the V(1,0)-cycle, where we can observe a slight increase for the linear elastic boundary value problem.
Overall, we can conclude that while the V(2,2)-cycle represents the fastest solver for both cases of Poisson's equation, the V(3,3) cycle leads to the fastest solving time for the linear elastic boundary value problem.

Finally, we evaluate the multigrid methods obtained with our evolutionary program synthesis approach in each of the ten experiments under the same conditions.
As the number of non-dominated individuals within the population at the end of each experiment is unrestricted and can thus be too large for a direct evaluation, we heuristically identify the 50 most promising individuals by sorting them according to the metric
\begin{equation*}
	T_{\varepsilon} = t \cdot \frac{\log(\varepsilon)}{\log(\rho)},
\end{equation*}
where $\varepsilon = 10^{-12}$ is the desired defect reduction factor and $\rho$ and $t$ the attained objective function values.
The resulting multigrid methods are then executed as a solver for each test problem on a Broadwell compute node consisting of two sockets with 20 CPU cores.

After we have identified the multigrid method that leads to the fastest solving time in each case, we additionally evaluate it on the larger problem instance using the same settings as for the evaluation of the reference methods.
Table~\ref{table:poisson-2D-evolved-methods},~\ref{table:poisson-3D-evolved-methods} and~\ref{table:linear-elasticity-2D-evolved-methods} contain the resulting measured solving times for each case.
\begin{table}
	\caption[2D Poisson -- Number of iterations and solving times of the evolved methods]{2D Poisson -- Measured number of iterations and solving times of the evolved multigrid methods on 20 cores and two sockets.}
	\label{table:poisson-2D-evolved-methods}
	\centering
	\begin{tabular}{l c c c c c c}
		\toprule
		& \multicolumn{2}{c}{Iterations} & \multicolumn{2}{c}{Broadwell (ms)} & \multicolumn{2}{c}{Ivy Bridge (ms)} \\
		\cmidrule(r){2-3} \cmidrule(r){4-5} \cmidrule(r){6-7}
		$l_{max}$ & $11$& $12$ & $11$ & $12$ & $11$ & $12$\\
		\midrule
		ES-1 & 5 & 5 & 338 & 1064 & 304 & 1055\\
		\midrule
		ES-2 & 6 & 6 & 371 & 1163 & 330 & 1133 \\
		\midrule
		ES-3 & 5 & 5 & 311 & 988 & 279 & 976 \\
		\midrule
		ES-4 & 6 & 6 & 380 & 1188 & 338 & 1153 \\
		\midrule
		ES-5 & 5 & 5 & 312 & 978 & 279 & 963 \\
		\midrule
		ES-6 & 5 & 5 & 349 & 1123 & 309 & 1106 \\
		\midrule
		ES-7 & 6 & 6 & 354 & 1096 & 320 & 1068 \\
		\midrule
		ES-8 & 6 & 6 & 347 & 1081 & 310 & 1056 \\
		\midrule
		ES-9 & 6 & 6 & 353 & 1079 & 313 & 1045 \\
		\midrule
		ES-10 & 5 & 5 & 310 & 960 & 275 & 934 \\
		\bottomrule
	\end{tabular}
\end{table}
\begin{table}
	\caption[3D Poisson -- Number of iterations and solving times of the evolved methods]{3D Poisson -- Measured number of iterations and solving times of the evolved multigrid methods on 20 cores and two sockets.}
	\label{table:poisson-3D-evolved-methods}
	\centering
	\begin{tabular}{l c c c c c c}
		\toprule
		& \multicolumn{2}{c}{Iterations} & \multicolumn{2}{c}{Broadwell (ms)} & \multicolumn{2}{c}{Ivy Bridge (ms)} \\
		\cmidrule(r){2-3} \cmidrule(r){4-5} \cmidrule(r){6-7}
		$l_{max}$ & $7$& $8$ & $7$ & $8$ & $7$ & $8$\\
		\midrule
		ES-1 & 10 & 11 & 55.3 & 577 & 70.0 & 704\\
		\midrule
		ES-2 & 8 & 9 & 57.2 & 578 & 64.3 & 716 \\
		\midrule
		ES-3 & 8 & 9 & 59.0 & 671 & 65.3 & 824 \\
		\midrule
		ES-4 & 8 & 9 & 54.6 & 576 & 62.7 & 710 \\
		\midrule
		ES-5 & 8 & 10 & 54.6 & 641 & 60.9 & 789 \\
		\midrule
		ES-6 & 9 & 10 & 59.4 & 716 & 67.1 & 891 \\
		\midrule
		ES-7 & 6 & 8 & 56.2 & 702 & 70.9 & 880 \\
		\midrule
		ES-8 & 5 & 5 & 56.7 & 589 & 74.0 & 724 \\
		\midrule
		ES-9 & 10 & 10 & 61.0 & 568 & 66.3 & 681 \\
		\midrule
		ES-10 & 10 & 11 & 55.4 & 581 & 61.3 & 705 \\
		\bottomrule
	\end{tabular}
\end{table}
\begin{table}
	\caption[2D Linear Elasticity -- Number of iterations and solving times of the evolved methods]{2D Linear Elasticity -- Measured number of iterations and solving times of the evolved multigrid methods on 20 cores and two sockets.}
	\label{table:linear-elasticity-2D-evolved-methods}
	\centering
	\begin{tabular}{l c c c c c c}
		\toprule
		& \multicolumn{2}{c}{Iterations} & \multicolumn{2}{c}{Broadwell (ms)} & \multicolumn{2}{c}{Ivy Bridge (ms)} \\
		\cmidrule(r){2-3} \cmidrule(r){4-5} \cmidrule(r){6-7}
		$l_{max}$ & $10$& $11$ & $10$ & $11$ & $10$ & $11$\\
		\midrule
		ES-1 & 6 & 6 & 234 & 1117 & 235 & 1137 \\
		\midrule
		ES-2 & 6 & 6 & 216 & 1033 & 211 & 1035 \\
		\midrule
		ES-3 & 7 & 7 & 258 & 1225 & 259 & 1231 \\
		\midrule
		ES-4 & 6 & 6 & 226 & 1077 & 219 & 1093 \\
		\midrule
		ES-5 & 6 & 6 & 235 & 1121 & 229 & 1139 \\
		\midrule
		ES-6 & 6 & 6 & 220 & 1083 & 213 & 1093 \\
		\midrule
		ES-7 & 7 & 7 & 238 & 1191 & 236 & 1186 \\
		\midrule
		ES-8 & 6 & 6 & 217 & 1037 & 223 & 1039 \\
		\midrule
		ES-9 & 6 & 6 & 224 & 1039 & 222 & 1058 \\
		\midrule
		ES-10 & 7 & 7 & 243 & 1188 & 238 & 1188 \\
		\bottomrule
	\end{tabular}
\end{table}
In general, we can conclude that in all three cases, our evolutionary algorithm was consistently able to discover well-functioning multigrid methods, leading to fast solving times for both problem sizes.
Furthermore, in the case of the two-dimensional Poisson equation and linear elasticity, our approach yields multigrid methods that achieve an even higher error reduction efficiency than the best reference cycle.
Here the fastest discovered method for the two-dimensional Poisson equation, ES-10, leads to a $9 \%$ solving-time improvement compared to the V(2,2)-cycle on both architectures, while for linear elasticity the ES-2 method achieves an even larger speedup of 17--27 \% compared to the V(3,3)-cycle.
%Interestingly, in contrast to the other two cases, for the two-dimensional Poisson equation, all solvers achieve faster execution times on the older Ivy Bridge computer architecture.
%Further investigating this phenomenon would require us to perform an in-depth analysis of the generated code.
%However, since our focus is a comparison of the relative performance of different multigrid methods on the same architecture, we have not yet performed such an analysis.
While in the case of the three-dimensional Poisson, the methods discovered with our approach still represent competitive solvers, they are not able to achieve the same degree of efficiency in solving both problem sizes.
In particular, with the exception of the ES-8 and ES-9 methods, we observe a slight increase in the number of iterations for the larger instance of this problem, which leads to worse solving times compared to the reference cycles.
This effect indicates that not all multigrid methods obtained for a particular instance of this test problem can be generalized to larger instances without further adaption.
In Section~\ref{sec:generalization}, we have already addressed this issue by proposing a multigrid-specific generalization scheme, whose effectiveness will be investigated in the next section of this chapter.   
%we have already presented an adapted version of our multi-objective evolutionary algorithm that aims to overcome these limitations by iteratively increasing the problem size during the search.
%In the next section of this chapter, we will, thus, investigate the effectiveness of this approach on the indefinite Helmholtz equation, a problem of substantially higher difficulty than those considered within this section.

To conclude our experimental analysis, Figure~\ref{fig:evolved-methods-graphical-representation:2D} and~\ref{fig:evolved-methods-graphical-representation:3D} contain a graphical representation of the discovered multigrid method that achieves the fastest solving time for the larger problem instance of two and three-dimensional Poisson equation, respectively. 
\begin{figure}
	%2D Poisson
 \captionsetup[subfigure]{justification=centering}
		\begin{subfigure}[b]{\columnwidth}
			\scalebox{0.725}{%
			\begin{tikzpicture}
				\node   (h) at (-0.75, 4){$h$};
				\node   (2h) at (-0.75, 3){$2h$};
				\node   (4h) at (-0.75, 2){$4h$};
				\node   (8h) at (-0.75, 1){$8h$};
				\node   (16h) at (-0.75, 0){$16h$};
				\node	(a) at (0,4) [draw, fill=lightred, circle,inner sep=0pt,minimum size=5mm] {\tiny $1.15$};
				\node	(b) at (0.5,3) [draw, circle,inner sep=0pt,minimum size=5mm] {\phantom{\tiny $1.00$}};
				\node	(c) at (1,2) [draw, circle,fill=lightred,inner sep=0pt,minimum size=5mm] {\tiny $0.80$};
				\node	(d) at (1.5,1) [draw, circle,inner sep=0pt,minimum size=5mm] {\phantom{\tiny $1.00$}};
				\node	(e) at (2,0) [draw, circle,fill=black, inner sep=0pt,minimum size=5mm] {\phantom{\tiny $1.00$}};
				\node	(f) at (2.5,1) [draw, circle,  fill=lightred,inner sep=0pt,minimum size=5mm] {\tiny $1.90$};
				\node	(g) at (3,2) [draw, circle,fill=lightred,inner sep=0pt,minimum size=5mm] {\tiny $1.40$};
				\node	(h) at (4,2) [draw, circle,fill=lightred,inner sep=0pt,minimum size=5mm] {\tiny $1.65$};
				\node	(i) at (5,2) [draw, circle,fill=lightred,inner sep=0pt,minimum size=5mm] {\tiny $1.45$};
				\node	(j) at (6,2) [draw, circle,fill=lightred,inner sep=0pt,minimum size=5mm] {\tiny $1.05$};
				\node	(k) at (6.5,3) [draw, circle,fill=lightred,inner sep=0pt,minimum size=5mm] {\tiny $1.30$};
				\node	(l) at (7.5,3) [draw, circle,fill=lightred,inner sep=0pt,minimum size=5mm] {\tiny $0.55$};
				\node	(m) at (8,4) [draw, circle,fill=lightred,inner sep=0pt,minimum size=5mm] {\tiny $1.05$};
				\node	(n) at (9,4) [draw, circle,fill=lightred,inner sep=0pt,minimum size=5mm] {\tiny $1.05$};
				\node	(o) at (9.5,3) [draw, circle, inner sep=0pt,minimum size=5mm] {\phantom{\tiny $1.00$}};
				\node	(p) at (10,2) [draw, circle, inner sep=0pt,minimum size=5mm] {\phantom{\tiny $1.00$}};
				\node	(q) at (10.5,1) [draw, circle, fill=lightred, inner sep=0pt,minimum size=5mm] {\tiny $0.90$};
				\node	(r) at (11.5,1) [draw, circle, fill=lightred, inner sep=0pt,minimum size=5mm] {\tiny $1.40$};
				\node	(s) at (12,2) [draw, circle, fill=lightred, inner sep=0pt,minimum size=5mm] {\tiny $0.95$};
				\node	(t) at (13,2) [draw, circle, fill=lightred, inner sep=0pt,minimum size=5mm] {\tiny $1.25$};
				\node	(u) at (14,2) [draw, circle, fill=lightred, inner sep=0pt,minimum size=5mm] {\tiny $1.05$};
				\node	(v) at (14.5,3) [draw, circle, fill=lightred, inner sep=0pt,minimum size=5mm] {\tiny $1.00$};
				\node	(w) at (15,4) [draw, circle, fill=lightred, inner sep=0pt,minimum size=5mm] {\tiny $1.05$};
				\draw 
				(a) edge[->] (b) 
				(b) edge[->] (c)
				(c) edge[->] (d)
				(d) edge[->] (e)   
				(e) edge[->] node[near end,left] {\tiny 0.80} (f)
				(f) edge[->] node[near end,left] {\tiny 1.20} (g)
				(g) edge[->] (h)
				(h) edge[->] (i)
				(i) edge[->] (j) 
				(j) edge[->] node[near end,left] {\tiny 0.90} (k)
				(k) edge[->] (l)
				(l) edge[->] node[near end,left] {\tiny 1.15} (m)   
				(m) edge[->] (n)
				(n) edge[->] (o)
				(o) edge[->] (p)
				(p) edge[->] (q)
				(q) edge[->] (r)
				(r) edge[->] node[near end,left] {\tiny 1.20} (s)
				(s) edge[->] (t)
				(t) edge[->] (u)
				(u) edge[->] node[near end,left] {\tiny 0.90} (v)
				(v) edge[->] node[near end,left] {\tiny 1.30} (w)
				;
			\end{tikzpicture}
		}
	\caption{2D Poisson: ES-10}
	\label{fig:evolved-methods-graphical-representation:2D}
	\end{subfigure}
  \par\bigskip
	\begin{subfigure}[b]{\columnwidth}
	\scalebox{0.725}{%
		\begin{tikzpicture}
			\node   (h) at (-0.75, 4){$h$};
			\node   (2h) at (-0.75, 3){$2h$};
			\node   (4h) at (-0.75, 2){$4h$};
			\node   (8h) at (-0.75, 1){$8h$};
			\node   (16h) at (-0.75, 0){$16h$};
			\node	(a) at (0,4) [draw, circle,inner sep=0pt,minimum size=5mm] {\phantom{\tiny $1.00$}};
			\node	(b) at (0.5,3) [draw, circle, fill=lightred, inner sep=0pt,minimum size=5mm] {\tiny $1.30$};
			\node	(c) at (1,2) [draw, circle,inner sep=0pt,minimum size=5mm] {\phantom{\tiny $1.00$}};
			\node	(d) at (1.5,1) [draw, circle,fill=lightred, inner sep=0pt,minimum size=5mm] {\tiny $0.85$};
			\node	(e) at (2,0) [draw, circle,fill=black, inner sep=0pt,minimum size=5mm] {\phantom{\tiny $1.00$}};
			\node	(f) at (2.5,1) [draw, circle,inner sep=0pt,minimum size=5mm] {\phantom{\tiny $1.00$}};
			\node	(g) at (3,2) [draw, circle,fill=lightblue,inner sep=0pt,minimum size=5mm, label=above:{\tiny $(2,1,1)$}] {\tiny 1.00}; %TODO add block info
			\node	(h) at (4,2) [draw, circle,fill=lightblue,inner sep=0pt,minimum size=5mm] {\tiny $1.5$};
			\node	(i) at (5,2) [draw, circle,fill=lightred,inner sep=0pt,minimum size=5mm] {\tiny $1.85$};
			\node	(j) at (6,2) [draw, circle,fill=lightred,inner sep=0pt,minimum size=5mm] {\tiny $0.75$};
			\node	(k) at (7,2) [draw, circle,fill=lightblue,inner sep=0pt,minimum size=5mm] {\tiny $0.60$};
			\node	(l) at (8,2) [draw, circle,fill=lightred,inner sep=0pt,minimum size=5mm] {\tiny $0.85$};
			\node	(m) at (8.5,1) [draw, circle,fill=lightred,inner sep=0pt,minimum size=5mm] {\tiny $0.85$};
			\node	(n) at (9,0) [draw, circle,fill=black,inner sep=0pt,minimum size=5mm] {\phantom{\tiny $1.00$}};
			\node	(o) at (9.5,1) [draw, circle,inner sep=0pt,minimum size=5mm] {\phantom{\tiny $1.00$}};
			\node	(p) at (10,2) [draw, circle,  fill=lightred, inner sep=0pt,minimum size=5mm] {\tiny $1.55$};
			\node	(q) at (10.5,3) [draw, circle, fill=lightred, inner sep=0pt,minimum size=5mm] {\tiny $1.30$};
			\node	(r) at (11,4) [draw, circle, fill=lightred, inner sep=0pt,minimum size=5mm] {\tiny $1.05$};
			\node	(s) at (12,4) [draw, circle, fill=lightred, inner sep=0pt,minimum size=5mm] {\tiny $1.25$};
			\draw 
			(a) edge[->] (b) 
			(b) edge[->] (c)
			(c) edge[->] (d)
			(d) edge[->] (e)   
			(e) edge[->] node[near end,left] {\tiny 0.80} (f)
			(f) edge[->] node[near end,left] {\tiny 1.20} (g)
			(g) edge[->] (h)
			(h) edge[->] (i)
			(i) edge[->] (j) 
			(j) edge[->] (k)
			(k) edge[->] (l)
			(l) edge[->] (m)   
			(m) edge[->] (n)
			(n) edge[->] node[near end,left] {\tiny 0.85}(o)
			(o) edge[->] node[near end,left] {\tiny 1.00}(p)
			(p) edge[->] node[near end,left] {\tiny 1.30}(q)
			(q) edge[->] node[near end,left] {\tiny 0.90}(r)
			(r) edge[->] (s)
			;
		\end{tikzpicture}
	}
	\caption{3D Poisson: ES-9}
	\label{fig:evolved-methods-graphical-representation:3D}
	\end{subfigure}
	\caption[Poisson -- Algorithmic structure of the discovered multigrid methods]{Algorithmic structure of the discovered multigrid methods for Poisson's Equation. The color of the node denotes the type of operation. Black: Coarse-grid solver, Blue: Block Jacobi smoothing, Red: Red-black Gauss-Seidel smoothing, White: No operation. The relaxation factor of each smoothing step is included in each node, while for coarse-grid correction, it is attached to the respective edge. For block smoothers, the dimension of the block is specified on top of the respective node. If no block size is specified, pointwise smoothing is applied.}
	\label{fig:evolved-methods-graphical-representation}
\end{figure}
The first observation that can be made from investigating the computational structure of these methods is that none of them can be clearly characterized as a V-, F- or W-cycle.
It is, therefore, impossible to formulate them within the framework of classical multigrid cycles, as shown in Algorithm~\ref{alg:multigrid-cycle}.
While, for instance, the first part of Figure~\ref{fig:evolved-methods-graphical-representation:2D} can be characterized as a V-cycle, the method proceeds with an additional coarse-grid correction step that is based on a purely smoothing-based error reduction on lower levels.
Figure~\ref{fig:evolved-methods-graphical-representation:3D} starts off in a similar fashion but then applies an additional three-grid V-cycle on the third level with a step size of $4h$, before it transfers the computed correction back to the finest grid.
An analysis performed in the recent work by Avnat and Yavneh already suggests that employing cycles of non-classical structure can lead to a significantly higher degree of efficiency in solving certain PDEs~\cite{avnat2022recursive}. 
However, in addition to their non-classical structure, both multigrid methods discovered with our approach employ a different number of smoothing steps on each level using a wide range of different relaxation factors. 
Moreover, significantly more smoothing is performed on certain levels.
In particular, smoothing on the third level ($4h$) seems to be exceptionally effective in reducing the most significant error components, and hence, the number of smoothing steps on this level is larger than on any other level.
While both methods predominantly employ RB-GS as a smoother, Figure~\ref{fig:evolved-methods-graphical-representation:3D} also includes individual pointwise and block Jacobi steps.
Finally, the even more complicated computational structure of the methods discovered for the linear elastic boundary prevents us from including their graphical representations similar to Figure~\ref{fig:evolved-methods-graphical-representation}.
However, if we consider the method ES-2, which leads to the fastest solving time for the larger instance of the linear elastic boundary value problem with $l_{max} = 11$, we can observe a number of structural similarities with the ones shown in Figure~\ref{fig:evolved-methods-graphical-representation}.
In particular, ES-2 applies the coarse-grid solver only once within its computations and, therefore, resembles a V-cycle but also includes additional smoothing-based coarse-grid correction steps with a varying amount of smoothing on each level.
Furthermore, with the exception of a single Jacobi step on the second-coarsest level, it employs exclusively RB-GS smoothing.

\section{Multigrid-Based Preconditioners for the Indefinite Helmholtz Equation}
\begin{figure}
		\includegraphics[scale=0.725]{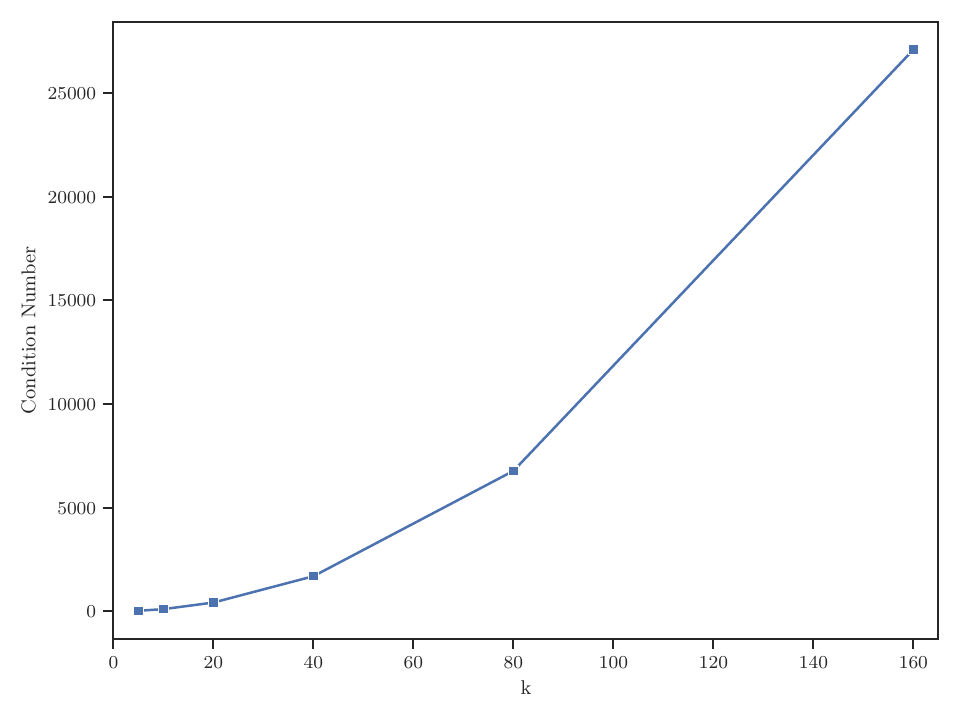}
		\caption[Condition number resulting from a finite-difference discretization of the two-dimensional Helmholtz equation]{Condition number of the system matrix resulting from a finite-difference discretization of the two-dimensional Helmholtz equation with $kh = 0.625$.}
		\label{fig:condition-number-helmholtz}
\end{figure}

While within the last section, we could already demonstrate that our evolutionary program synthesis approach leads to the automated design of multigrid methods with competitive performance compared to common multigrid cycles, none of the PDEs considered there is particularly challenging to solve.   
Therefore, in order to evaluate whether our approach can achieve any advantages compared to state-of-the-art methods in the case of a problem that is difficult to solve numerically, we consider the indefinite Helmholtz equation.
At this point, we would like to point out that all results presented in this section have been originally published in~\cite{schmitt2022evolving}.
The Helmholtz equation, whose basic form is given as 
\begin{equation}
	-\nabla ^{2}u - k^{2}u = f,
	\label{eq:helmholtz}
\end{equation} 
is a famously-difficult test problem for the application of numerical methods that also has practical relevance for many real-world applications, such as~\cite{versteeg1994marmousi,martin2006marmousi2,billette20052004,gray1995migration}.
The main difficulty in solving this equation is that the system of linear equations resulting from its discretization becomes indefinite and highly ill-conditioned for large values of the wavenumber $k$~\cite{ernst2012difficult}.
For instance, a discretization with ten grid points per wavelength\footnote{Usually the relation of the wavelength $\lambda$ to the wavenumber $k$ is defined as $k = 2 \pi/\lambda$.}, as it is common in geophysical applications~\cite{erlangga2006multigrid}, requires us to fulfill a second-order accuracy requirement of $kh = 0.625$.
Figure~\ref{fig:condition-number-helmholtz} shows the condition number of the system matrix $A$ that results from a discretization of Equation~\eqref{eq:helmholtz} for different values of the wavenumber $k$.
The condition number gives us a measure of how much an initial approximation error is magnified by the application of the operator $A$.
As Figure~\ref{fig:condition-number-helmholtz} illustrates, the condition number of $A$ increases dramatically for larger values of $k$, which leads to an extreme accumulation of numerical errors. 
The necessity to handle these errors makes the design of an efficient or even functioning numerical method for the indefinite Helmholtz equation outstandingly difficult.
In particular, the use of classical multigrid methods as an iterative solver for indefinite Helmholtz problems is often infeasible~\cite{ernst2012difficult}.
One multigrid-based approach that mitigates this problem is the application of the method as a preconditioner instead of applying it to the discretized Helmholtz system directly.
In general, preconditioning has the purpose of modifying a given system of linear equations by applying a so-called preconditioning matrix $M$ to obtain a new system that is easier to solve.
For instance, right-preconditioning the system $A \bm{x} = \bm{b}$ with the matrix $M$ results in
\begin{equation}
	A M^{-1} \bm{y} = \bm{b},
	\label{eq:right-preconditioning}
\end{equation}
with $M \bm{x} = \bm{y}$. 
The main consequence of this formulation is that we now have to solve an additional system of linear equations in the form of $M \bm{x} = \bm{y}$ whenever the operator $A$ is applied within our original solution method.
If we consider the extreme choice of $M = A$, this means that Equation~\eqref{eq:right-preconditioning} reduces to $\bm{y} = \bm{b}$, which, however, means that the system $M \bm{x} = \bm{y}$ becomes just as ill-conditioned as the original one.
Therefore, the choice of the preconditioning matrix $M$ usually represents a compromise between the accurate approximation of the original matrix $A$ and the ease of solvability of $M \bm{x} = \bm{y}$~\cite{benzi2002preconditioning}.
Since its invention by Erlangga et al.~\cite{erlangga2004preconditioner}, the choice of $M$ as a complex-shifted version of the original operator, which leads to
\begin{equation*}
	M = -\nabla ^{2} - (k^{2} + i \varepsilon),
\end{equation*}
has proven its effectiveness for solving indefinite Helmholtz problems numerically~\cite{erlangga2008advances,cocquet2017shift,umetami2009multigrid,cools2013analysis}.
%TODO include all suitable citations
As it has been shown in~\cite{cocquet2017shift}, a shift $\varepsilon \approx \mathcal{O}(k^2)$ enables the efficient inversion of $M$ by multigrid, whereas the choice of a smaller shift increases the preconditioning effectiveness.
Since the resulting linear system is often complex-symmetric but non-hermitian, it is commonly solved using a suitable Krylov subspace method, such as GMRES and BiCGSTAB~\cite{saad2003iterative}.
\subsection{Problem Formulation}
After introducing the indefinite Helmholtz equation and the difficulties, its solution involves, we can now proceed with defining a representative instance of this equation for the evaluation of our grammar-based approach for the automated design of multigrid methods.
For this purpose, we consider the two-dimensional Helmholtz equation on a unit square with Dirichlet boundary conditions at the top and bottom and Robin radiation conditions at the left and right, as given by
\begin{equation}
	\label{eq:helmholtz-test-problem}
	\begin{split}
		(-\nabla ^{2} - k^{2}) u & = f \quad \text{in} \; \left( 0, 1 \right)^2 \\
		u & = 0 \quad \text{on} \; \left( 0, 1 \right) \times \{0\}, \left( 0, 1 \right) \times \{1\} \\
		\partial_{\mathbf{n}} u - iku & = 0 \quad \text{on} \; \{0\} \times \left( 0, 1 \right), \{1\} \times \left( 0, 1 \right) \\
		f(x, y) & = \delta(x - 0.5, y - 0.5),
	\end{split}
\end{equation}
where $\delta(\bm{x})$ is the Dirac delta function.
We discretize this equation on a uniform Cartesian grid using the classical five-point stencil
\begin{equation*}
	A_h = \frac{1}{h^2} \begin{bmatrix}
		& -1 & \\
		-1 & 4 - (k h)^2 & -1 \\
		& -1 &  
	\end{bmatrix}.
\end{equation*}
In addition, $\delta(\bm{x})$ is approximated with a second-order Zenger correction~\cite{koestler2004extrapolation}.
The spacing $h$ of the grid is chosen to fulfill the second-order accuracy requirement $h k = 0.625$ as described above.
Finally, we apply the shifted Laplace operator
\begin{equation*}
	M = -\nabla^{2} - (k^{2} + 0.5 i k^{2}),
\end{equation*}
as a preconditioner, following the suggestion in~\cite{erlangga2008advances}.
This operator is discretized similarly to the original one, utilizing a five-point stencil:
\begin{equation*}
	M_h = \frac{1}{h^2} \begin{bmatrix}
		& -1 & \\
		-1 & 4 - (1.0 + 0.5i)(k h)^2 & -1 \\
		& -1 &  
	\end{bmatrix}.
\end{equation*}
\subsection{Solver Configuration}
\label{sec:solver-configuration-helmholtz}
As a result of the above formulation of our test problem, we obtain two systems of linear equations
\begin{equation}
	A_h M_h^{-1} \bm{y}_h = \bm{b}_h,
	\label{eq:helmholtz-test-problem-discretized-and-preconditioned}
\end{equation}
where $\bm{b}_h$ contains the values of $\delta(\bm{x})$ at each grid point, and
\begin{equation}
	M_h \bm{x}_h = \bm{y}_h,
	\label{eq:helmholtz-test-problem-preconditioning}
\end{equation}
where $\bm{x}_h$ represents the approximate solution of Equation~\eqref{eq:helmholtz-test-problem}.
While for each of these two systems of linear equations, a functioning solver is needed, the focus of our experimental evaluation is the design of an efficient multigrid method for the approximate solution of Equation~\eqref{eq:helmholtz-test-problem-preconditioning}.
Here, to limit the cost of preconditioning, we assume that the application of a single multigrid cycle is sufficient to compute a reasonable approximation for $M_{h}^{-1}$, as proposed in~\cite{erlangga2008advances}.
After designing a suitable multigrid-based preconditioner, Equation~\eqref{eq:helmholtz-test-problem-discretized-and-preconditioned} is solved using the biconjugate gradient stabilized method (BiCGSTAB)~\cite{saad2003iterative}.
The resulting iterative solution scheme is summarized in Algorithm~\ref{alg:preconditioned-bicgstab}, where we omit the grid spacing $h$ and parentheses in superscripts for simplicity.
\begin{algorithm}
	\caption{Right-Preconditioned BiCGSTAB}
	\label{alg:preconditioned-bicgstab}
	\begin{algorithmic}[1] % The number tells where the line numbering should start
			\State $\bm{r}^0 = \bm{b} - A \bm{x}^0$
			\State $\bm{\hat{r}}^0 = \bm{r}^0$
			\State $\alpha_0 = \beta_0 = \rho_0 = \omega_0 = 1$
			\State $\bm{p}^0 = \bm{q}^0 = \bm{0}$
			\For{$i := 1, \dots, n$}
			\State $\rho_i = \bm{\hat{r}}^{i-1} \cdot \bm{r}^{i-1}$
			\State $\beta_i = \frac{\rho_i }{\rho_{i-1} }\frac{\alpha_{i-1}}{ \omega_{i-1}}$
			\State $\bm{p}^i = \bm{r}^{i-1} + \beta_i (\bm{p}^{i-1} - \omega_{i-1} \bm{q}^{i-1})$
			\State Solve $M \bm{x}^i = \bm{p}^i$
			\State $\bm{q}^i = A \bm{x}^i$
			\State $\alpha_i = \rho_i / (\bm{\hat{r}}^{i-1} \cdot \bm{q}^i)$
			\State $\bm{h}^i = \bm{x}^{i-1} + \alpha_i \bm{x}^i$	
			\State $\bm{s}^i = \bm{r}^{i-1} - \alpha_i \bm{r}^i$
			\State Solve $M \bm{x}^i = \bm{s}^i$
			\State $\bm{t}^i = A \bm{x}^i$
			\State $\omega_i = (\bm{t}^i \cdot \bm{s}^i) / (\bm{t}^i \cdot \bm{t}^i)$
			\State $\bm{x}^i = \bm{h}^i + \omega_i \bm{x}^i$
			\State $\bm{r}^i = \bm{s}^i - \omega_i \bm{t}^i$
			\If{$\norm{\bm{r}^i}/\norm{\bm{r}^0} < \epsilon$}
			\Return $\bm{x}^i$
			\EndIf
			\EndFor
		\end{algorithmic}
\end{algorithm}
In each step of this iterative scheme, it is necessary to compute an approximate solution for two systems of linear equations in the form of $M \bm{x}^i = \bm{y}^i$, %, in line 9 and 14, 
which in both cases is achieved by applying a single multigrid cycle.
To obtain an efficient method for this task, we consider the class of five-grid methods that are defined on a hierarchy of discretizations with $h = 1/2^{l}$ on each level 
\begin{equation*}
    l \in \left[l_{max} - 4,l_{max}\right]. 
\end{equation*}
Similar to Section~\ref{sec:experiments-part1}, we then consider the following components for constructing a multigrid method:
\begin{description}
	\item[\textbf{Smoothers}:] Pointwise and block Jacobi with rectangular blocks up to a maximum number of six terms, red-black Gauss-Seidel (RB-GS)
	\item[\textbf{Restriction}:] Full-weighting restriction
	\item[\textbf{Prolongation}:] Bilinear interpolation
	\item[\textbf{Relaxation factors}:] $\left( 0.1 + 0.05i \right)_{i = 0}^{36} = \left(0.1, 0.15, 0.2, \dots, 1.9 \right)$
	\item[\textbf{Coarse-grid solver}:] BiCGSTAB for $l = l_{max} - 4$
\end{description}
The productions of the resulting multigrid grammar are similar to Algorithm~\ref{table:multigrid-grammar}. 
However, each occurrence of a system matrix $A_H$ and right-hand side $\ps{b_H}$ on a level with step size $H$ is replaced by the respective preconditioning matrix $M_H$ and right-hand side $\ps{y_H}$~\cite{schmitt2022evolving}.
Similar to Section~\ref{sec:experiments1-multigrid-configuration}, block Jacobi smoothers are generated based on rectangular blocks of grid points, while the relaxation factor $\omega$ of each smoothing and coarse-grid correction step is chosen from the above interval.
To assess the efficiency and generalizability of the multigrid preconditioners designed with our grammar-based approach, we consider the set of classical multigrid cycles that can be constructed based on the same components.
To ensure that each method uses the optimal smoother and relaxation factor, we evaluate each possible combination on the largest problem size for which convergence can be achieved.
Due to the ill-conditioning of the indefinite Helmholtz equation, we consider an approximate solution to be sufficient when the initial residual has been reduced by a factor of $10^{-7}$ for $k \leq 160$ and $10^{-6}$ for all larger wavenumbers.

\subsection{Experimental Settings and Evaluation Platform}
The ease of solvability of the PDEs considered in Section~\ref{sec:experiments-part1} allowed us to keep the problem size constant throughout each experiment.
In the case of the indefinite Helmholtz equation, this strategy is infeasible due to a number of reasons.
First of all, doubling the wavenumber requires us to use twice as many grid points in each dimension because of the requirement $kh = 0.625$.
This means that we end up with a system of linear equations that is not only significantly worse conditioned but also has four times the number of unknowns.
As a consequence, solving Helmholtz problems with a large wavenumber becomes tremendously expensive, which makes the evaluation of a large number of different preconditioners infeasible.
Furthermore, our goal is to design multigrid methods that can be applied to a wide range of different problem instances, and hence, generalizability is one of our main concerns.
In Section~\ref{sec:generalization-procedure}, we have already proposed a systematic procedure for the generalization of a population of individuals to multiple instances of the same problem.
This method aims to evolve a population of generalizable multigrid methods by iteratively increasing the size of the test problem considered within the evaluation of each solver after a certain number of generations $m$.
While the formulation of our problem based on the wavenumber-dependent grid spacing $h$ allows us to construct problem instances of larger size and difficulty in a straightforward manner, it is unclear after how many generations this operation should be performed.
In the experiments performed in Section~\ref{sec:experiments-part1}, we could observe that the majority of improvement in both objectives is achieved within the first 50 generations of our evolutionary algorithm.
Therefore, setting the generalization interval $m$ to 50 represents a reasonable compromise between allowing the population to adapt to modified conditions and preventing it from overfitting to the characteristics of a particular problem instance.
To initiate the execution of our evolutionary algorithm, we choose $k = 80$ as a problem instance, which leads to a maximum level $l_{max} = 7$ and a system of linear equations consisting of 16129 unknowns. 
Therefore, the cost of evaluation at the beginning is significantly lower compared to later generations.
We then execute Algorithm~\ref{alg:generalization-procedure} with a population size of 128 and a total number of 150 generations, whereby after 50 and 100 generations, we increase the value of wavenumber to 160 and 320, respectively.
Table~\ref{table:gp-parameters-helmholtz} gives an overview of the algorithmic parameters used in each experiment, which are mostly similar to the ones employed in Section~\ref{sec:experiments-part1}.
\begin{table}
	\centering
	\caption{Summary of genetic programming configuration parameters.}
	\label{table:gp-parameters-helmholtz}
	\begin{tabular}{l c}
		\toprule
		Parameter & Value \\
		\midrule 
		Evolutionary algorithm type & $(\mu + \lambda)$ \\
		\midrule
		Objectives & $t, n$ \\
		\midrule
		Number of generations & 150 \\
		\midrule
		$k$ & 80, 160, 320 \\
		\midrule
		Generalization interval & 50 \\
		\midrule
		Initial population size & 1024 \\
		\midrule
		$\lambda$ & 128 \\
		\midrule
		$\mu$ & 128 \\
		\midrule
		Number of MPI processes & 64 \\
		\midrule
		Non-dominated sorting procedure & \cite{deb2002fast} \\ 
		\midrule
		Selection operator & \cite{deb2002fast} \\ 
		\midrule
		Crossover operator & Subtree crossover \\
		\midrule
		Crossover probability & $2/3$ \\
		\midrule
		Mutation operator & Random subtree insertion \\
		\midrule 
		Probability to mutate a terminal symbol & $1/3$ \\
		\bottomrule
	\end{tabular}
\end{table}
In order to evaluate each multigrid method, we measure the number of iterations $n$ that Algorithm~\ref{alg:preconditioned-bicgstab} requires until the desired defect reduction is achieved, together with its execution time per iteration $t$.
These two metrics then serve as the two objectives for our multi-objective evolutionary algorithm.
Note that in previous formulations of our method, we have employed the solver's convergence factor $\rho$ as a first objective.
However, due to the difficulty of the problem considered here, we expect the average number of iterations to be significantly larger than in our previous experiments. 
The small potential difference in this metric between two preconditioners relative to their total number of iterations thus means that the absolute values of their convergence factors are hardly distinguishable, which has the potential to lead to the accumulation of numerical errors in the computation of population-density metrics, such as the crowding distance~\cite{deb2002fast}.
Using the number of iterations instead as the first objective avoids this issue, while, due to the ill-conditioning of the considered problem, even a small improvement in the accuracy of the preconditioner will lead to a measurable decrease in the value of this metric.

Similar to Section~\ref{sec:experiments-part1}, we select individuals for crossover and mutation using a binary tournament selection based on their dominance relation and crowding distance.
New individuals are created using the subtree crossover operator with a probability of 2/3 and random subtree insertion with a probability of 1/3, whereby in the latter case, we insert the generated subtree whenever possible within the selected branch and otherwise replace it completely, as described in Section~\ref{sec:gggp-mutation-and-recombination}.
To evaluate the effectiveness of our evolutionary program synthesis approach, we perform a total number of ten experiments with a random initialization, each of which is executed on the SuperMUC-NG cluster\footnote{SuperMUC-NG: \url{https://doku.lrz.de/display/PUBLIC/SuperMUC-NG}} of the Leibniz Supercomputing Center (LRZ).
Even though the population size, as well as the number of generations, is smaller than within our experimental evaluation in Section~\ref{sec:experiments-part1}, the time required to evaluate each solver drastically increases for larger wavenumbers, leading to an overall higher computational cost for each experiment.
To parallelize our evolutionary algorithm, we utilize the message-passing interface (MPI) in the form of the \emph{mpi4py} library as described in Section~\ref{sec:distributed-parallelization}. 
We, therefore, create 64 MPI processes that are distributed to eight compute nodes of the system, where each process is exclusively executed on one of the eight available islands.
Each process then creates and evaluates two new individuals per generation, which are executed on the respective island using an OpenMP-based shared-memory parallelization with 12 threads.

\subsection{Evaluation of the Evolved Multigrid Methods}
To assess whether our evolutionary program synthesis method can discover competitive multigrid-based preconditioners for the indefinite Helmholtz equation, we first need to investigate the preconditioning efficiency of classical multigrid cycles.
As in previous evaluations, we consider the time required to achieve a certain error reduction as a measure of the efficiency of a solver.
However, since in the given case, we apply the same iterative solver to each problem instance and only vary the multigrid cycle used for preconditioning, this metric instead relates to the efficiency of the latter.
For instance, a multigrid cycle that is able to compute an improved approximation of the solution of the preconditioning system, as given by Equation~\eqref{eq:helmholtz-test-problem-preconditioning}, reduces the number of iterations required for solving the preconditioned system, which might come at the prize of a higher computational cost.
The second quality metric that we aim to investigate in this section is the capability of each multigrid preconditioner to generalize to different instances of the same problem.
Therefore, we consider three different wavenumbers $k = 160, 320$, and $640$ for the evaluation of each method, which we perform on a single compute node of SuperMUC-NG using an OpenMP-based shared-memory parallelization  with 48 threads, such that each thread is executed on a dedicated CPU core.
As multigrid-based preconditioners, we consider different V-, W-, and F-cycles with up to five pre- and post-smoothing steps based on the configuration space described in Section~\ref{sec:solver-configuration-helmholtz}.
Similar to previous works on multigrid-based preconditioners for the indefinite Helmholtz equation~\cite{erlangga2006multigrid,erlangga2008advances,cocquet2017shift}, the same smoother and relaxation factor is applied for a fixed number of steps on each level.
To find out which combination of smoother and relaxation factor works best for each cycle, we have first tested each of them on the largest and thus worst conditioned problem with $k = 640$.
However, since none of the available configuration options led to a converging solver for this case, we instead consider the next smaller problem size with a wavenumber of $k = 320$ to find out which pair of smoother and relaxation factor leads to the fastest solving time for each multigrid cycle.
In general, according to our experimental evaluation RB-GS represents the single most efficient smoother for the given case, while Jacobi-type smoothers only lead to a converging solver for wavenumbers of at most $k = 160$, however, with less efficiency compared to RB-GS.
Table~\ref{table:reference-methods-helmholtz} shows the resulting average solving times for each multigrid cycle with RB-GS smoothing and an optimum relaxation factor for $k = 320$, where the omission of a number means that convergence could not be achieved within 20,000 iterations. 
\begin{table}
	\caption[2D Helmholtz -- Number of iterations and solving time of the reference methods]{Reference methods -- Optimum relaxation factors $\omega$ for $k = 320$, number of iterations, and average time required for solving a problem with the particular wavenumber ($k$).}
	\label{table:reference-methods-helmholtz}
	\centering
	\begin{tabular}{l l l l l l}
		\toprule
		& $\omega$ & \multicolumn{2}{c}{Iterations} & \multicolumn{2}{c}{Solving Time (s)} \\
		\cmidrule(r){3-4} \cmidrule(r){5-6}
		$k$ & & $160$ & $320$ & $160$ & $320$ \\
		\midrule
		V(0, 1) & $1.25$ & $2078$ & $6297$ & $6.38$ & $35.11$ \\
		\midrule
		V(1, 1) & $0.6$ & $1880$ & $6297$ & $7.66$ & $44.27$ \\
		\midrule
		V(2, 1) & $0.6$ & $-$ & $5532$ & $-$ & $47.0$ \\
		\midrule
		V(2, 2) & $0.5$ & $1627$ & $5115$ & $9.93$ & $50.54$  \\
		\midrule
		V(3, 3) & $0.4$ & $1753$ & $5168$ & $13.97$ & $76.00$ \\
		\midrule
		F(0, 1) & $1.15$ & $1467$ & $4028$ & $8.15$ & $42.87$  \\
		\midrule
		F(1, 1) & $0.75$ & $1546$ & $3988$ & $11.21$ & $54.51$ \\
		\midrule
		F(2, 1) & $0.55$ & $1146$ & $3934$ & $10.87$ & $67.62$ \\
		\midrule
		F(2, 2) & $0.65$ & $1060$ & $3213$ & $13.92$ & $65.06$ \\
		\midrule
		F(3, 3) & $0.45$ & $1085$ & $3464$ & $18.88$ & $92.97$ \\
		\midrule
		W(0, 1) & $0.75$ & $1265$ & $4215$ & $8.67$ & $72.08$ \\
		\midrule
		W(1, 1) & $0.8$ & $1208$ & $3570$ & $13.08$ & $76.22$ \\
		\midrule
		W(2, 1) & $0.6$ & $1313$ & $3074$ & $17.71$ & $79.67$ \\
		\midrule
		W(2, 2) & $0.5$ & $1069$ & $3376$ & $17.14$ & $101.6$ \\
		\midrule
		W(3, 3) & $0.45$ & $942$ & $2976$ & $19.65$ & $117.8$ \\
		\bottomrule
	\end{tabular}
\end{table}
Since the use of more than three pre- and post-smoothing steps did not yield any additional benefits, we have omitted all multigrid cycles with a higher number of smoothing steps.
To cut down the variability between individual measurements, each value contained in the table represents the average of 50 solver executions.
As expected, a higher number of smoothing and coarse-grid correction steps increases the efficiency of the preconditioner and hence leads to a lower number of iterations until convergence.
However, if we consider the configuration that leads to the fastest solving time, which is the V(0,1), we can see that a single coarse-grid correction followed by a single step of RB-GS smoothing is the best choice for all wavenumbers considered.
After determining the best possible configuration for each classical multigrid cycle, we can finally investigate whether the multigrid methods designed with our evolutionary program synthesis approach are able to achieve the same degree of efficiency and generalizability in preconditioning indefinite Helmholtz problems.

\begin{figure}
\centering
	\includegraphics[scale=0.725]{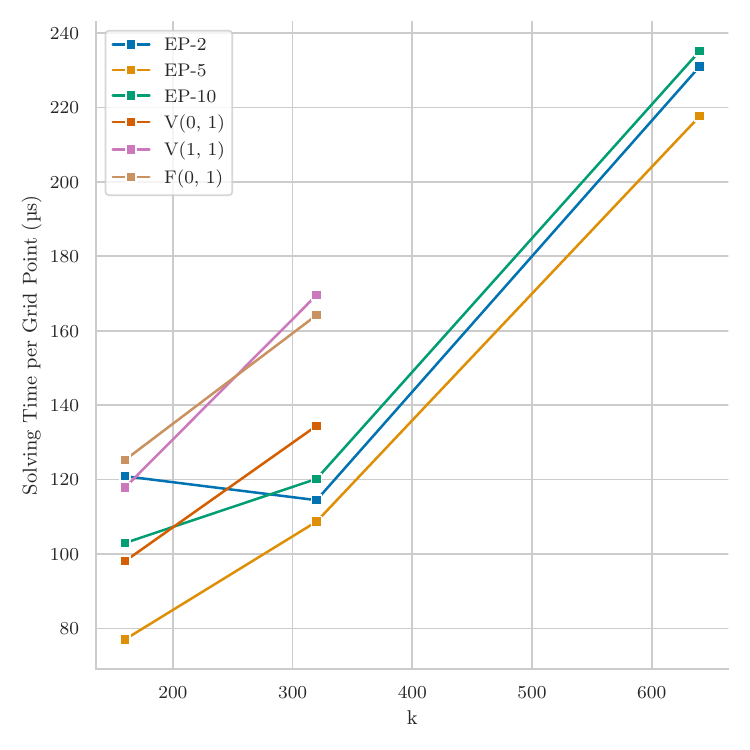}
	\caption[2D Helmholtz -- Solving time comparison of the best preconditioners]{Solving time comparison of the best preconditioners according to the product of both objectives for different wavenumbers ($k$).}
	\label{fig:helmholtz-solving-time-comparison}
\end{figure}

In order to identify the most promising multigrid-based preconditioners without needing to perform an excessive amount of evaluations, we sort the population at the end of each experiment according to the product of both objectives.
Note that while this metric corresponds to the solving time achieved on a single island, we intend to execute each method on a full node of the system, which can lead to a different outcome.
We then evaluate the ten best multigrid methods by employing them as a preconditioner for Equation~\eqref{eq:helmholtz-test-problem-preconditioning} with $k = 640$ while we measure the solving time of the resulting preconditioned iterative method.
In case none of the methods considered is able to reduce the initial residual by the required factor within 20,000 iterations, we repeat the evaluation with the next smaller wavenumber $k = 320$.
The best-performing method is then evaluated on the full set of wavenumbers, i.e., $k = 160, 320, 640$.
Table~\ref{table:evolved-solvers} shows the required number of iterations and solving time when using the respective multigrid method as a preconditioner for Algorithm~\ref{alg:preconditioned-bicgstab}.
\begin{table}
	\caption[2D Helmholtz -- Number of iterations and solving time of the best evolved methods]{Best evolved preconditioners according to the product of both objectives -- Number of iterations and the average time required for solving a problem with the particular wavenumber ($k$).}
	\label{table:evolved-solvers}
	\centering
	\begin{tabular}{l l l l l l l }
		\toprule
		& \multicolumn{3}{c}{Iterations} & \multicolumn{3}{c}{Solving Time (s)} \\
		\cmidrule(l){2-4} \cmidrule(l){5-7}
		$k$ & $160$ & $320$ & $640$ & $160$ & $320$ & $640$ \\
		\midrule
		EP-1 & $1178$ & $3399$ & $-$ & $6.29$ & $28.07$ & $-$ \\
		\midrule
		EP-2 & $795$ & $2160$ & $8449$ & $7.86$ & $29.89$ & $241.7$\\
		\midrule
		EP-3 & $933$ & $2827$ & $11143$ & $6.08$ & $27.58$ & $257.8$ \\
		\midrule
		EP-4 & $637$ & $2509$ & $7901$ & $7.17$& $41.04$ & $268.2$ \\
		\midrule
		EP-5 & $539$ & $1838$ & $7765$ & $5.01$ & $28.39$ & $227.7$ \\
		\midrule
		EP-6 & $941$ & $2103$ & $-$ & $9.58$ & $30.76$ & $-$ \\
		\midrule
		EP-7 & $955$ & $2701$ & $-$  & $6.45$& $27.84$ & $-$ \\
		\midrule
		EP-8 & $945$ & $2870$ & $10839$ & $7.24$ & $33.02$ & $276.9$ \\
		\midrule
		EP-9 & $3436$ & $3872$ & $-$  & $15.15$ & $27.51$ & $-$ \\
		\midrule
		EP-10 & $586$ & $1881$ & $8855$ & $6.70$ & $31.39$ & $246.1$ \\
		\bottomrule
	\end{tabular}
\end{table}
With the exception of EP-4, the multigrid-based preconditioners discovered with our evolutionary algorithm all achieve faster solving times than the best reference method V(0,1) for $k = 320$, while the same methods also work efficiently for $k = 160$ in the majority of cases.
If we consider the required number of iterations as a second metric, we can conclude that the methods discovered with our approach achieve at least a similar degree of preconditioning effectiveness as most W-cycles, however, with a significantly higher degree of computational efficiency.
Furthermore, even though the largest wavenumber considered within our evolutionary algorithm is 320, our approach yields multigrid-based preconditioners that lead to a converging solver for $k = 640$ in six out of ten experiments, while all classical multigrid cycles fail to achieve convergence in this case.
Figure~\ref{fig:helmholtz-solving-time-comparison} shows a direct comparison of the preconditioners from both categories that lead to the fastest solving time for large wavenumbers.
To enhance the interpretability of the results, we have normalized each solving time by the number of grid points.
Since the system of linear equations resulting from a discretization of Equation~\ref{eq:helmholtz-test-problem} becomes increasingly ill-conditioned for larger values of the wavenumber $k$, the relative cost to solve for each unknown grows accordingly.
All three discovered multigrid methods included in Figure~\ref{fig:helmholtz-solving-time-comparison} achieve a faster solving time than the best reference method V(0,1) for $k = 320$ while also staying competitive for smaller wavenumbers.
In particular, the overall best preconditioner EP-5 leads to a consistent solving-time improvement of 20 \% compared to the V(0,1) cycle.

\begin{figure}
\centering
	\includegraphics[scale=0.725]{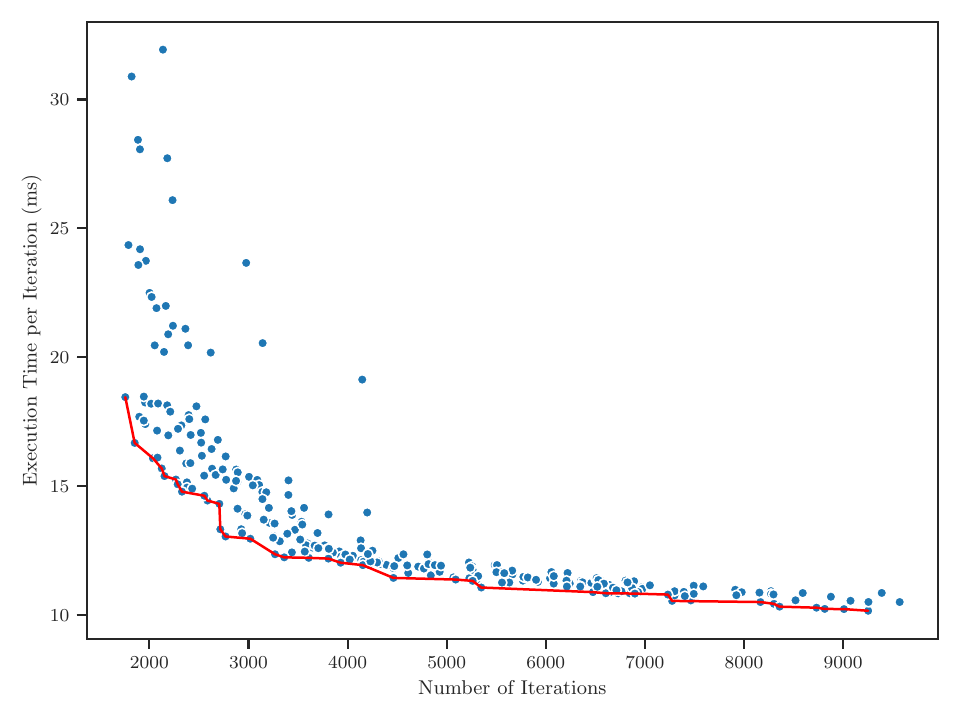}
	\caption[2D Helmholtz -- Distribution of non-dominated individuals at the end of all ten experiments]{Distribution of non-dominated individuals at the end of all ten experiments for $k = 320$. The red line denotes the combined front.}
	\label{fig:pareto-front-helmholtz}
\end{figure}

As a second evaluation step, we address the question of why our evolutionary algorithm failed to consistently discover multigrid-based preconditioners that achieve the same degree of efficiency and generalizability in each experiment.
Similar to Section~\ref{sec:experiments1-algorithm-behavior-analysis}, we consider the space of non-dominated individuals at the end of each experiment, which is shown in Figure~\ref{fig:pareto-front-helmholtz}.
First of all, note that the number of non-dominated individuals is smaller than in the experiments performed in Section~\ref{sec:experiments-part1}, which can be attributed to the smaller population size, which is only half as large as in our previous experiments.
Furthermore, since we iteratively increase the problem size and difficulty during the execution of our evolutionary algorithm, it is forced to adapt the current population to the characteristics of a new problem instance within only 50 generations, which increases the difficulty of evolving the same Pareto-optimal individuals.
Despite this fact, our method achieves a high degree of consistency in finding individuals that are located in the right half of the objective function space.
In contrast, in the left half of Figure~\ref{fig:pareto-front-helmholtz}, which corresponds to those individuals that lead to increasingly effective but more costly preconditioners, the algorithm is not able to find the same non-dominated individuals in each experiment.
Note that all multigrid-based preconditioners considered in our final evaluation shown in Table~\ref{table:evolved-solvers} correspond to individuals located here.
The larger distance of certain individuals in Figure~\ref{fig:pareto-front-helmholtz} to the combined front thus explains why we could not achieve the same degree of efficiency in each run.
To reduce the number of iterations through preconditioning, we need to improve the accuracy of the approximate solution computed for Equation~\eqref{eq:helmholtz-test-problem-preconditioning}.
Since this requires us to perform more smoothing and coarse-grid correction steps, the size of the corresponding derivation trees grows accordingly.
In order to obtain a particular derivation tree, the same sequence of productions needs to be discovered in each individual run of our algorithm, starting from a random initialization.
As we have already investigated in Section~\ref{sec:experiments1-algorithm-behavior-analysis}, NSGA-II struggles to achieve this goal consistently~\cite{liu2022evolvability}, which leads to a suboptimal outcome in the leftmost part of the objective function space in certain experiments.
\subsection{Analysis of the Discovered Algorithmic Features}
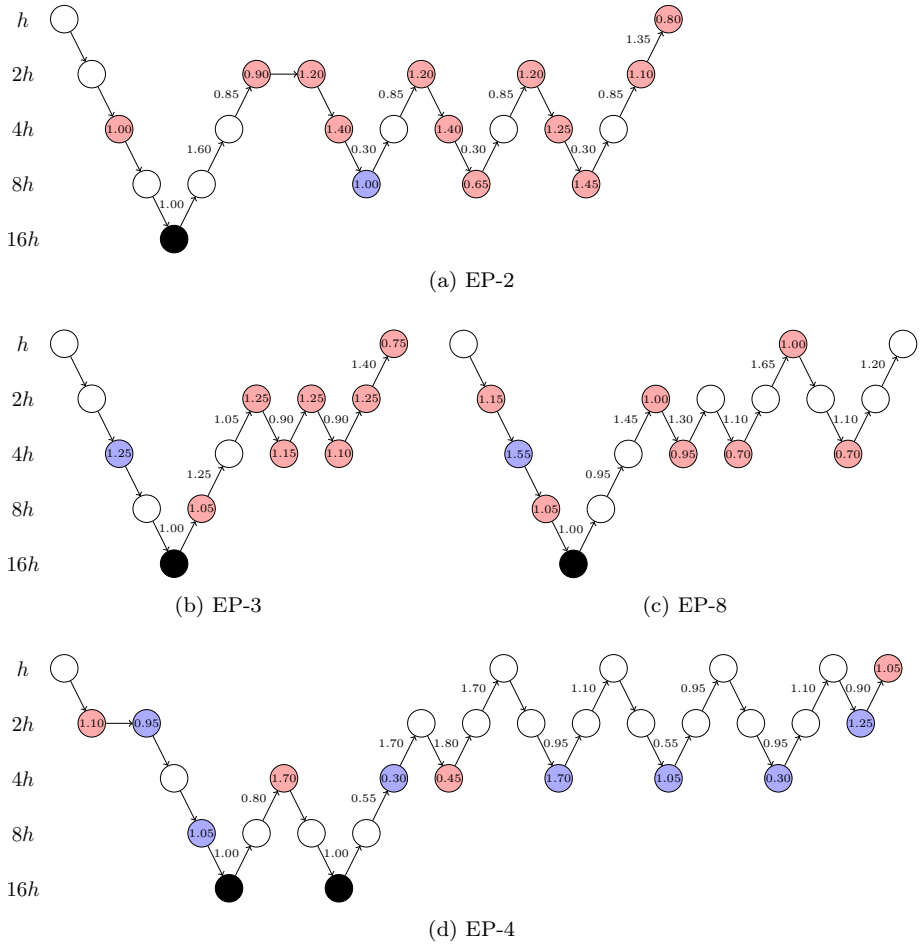
\begin{figure}
\captionsetup[subfigure]{justification=centering}
\begin{subfigure}[b]{\columnwidth}
		\scalebox{0.725}{%
			\begin{tikzpicture}
			\node   (h) at (-0.75, 4){$h$};
			\node   (2h) at (-0.75, 3){$2h$};
			\node   (4h) at (-0.75, 2){$4h$};
			\node   (8h) at (-0.75, 1){$8h$};
			\node   (16h) at (-0.75, 0){$16h$};
			\node	(a) at (0,4) [draw, circle,inner sep=0pt,minimum size=5mm] {\phantom{\tiny 1.00}};
			\node	(b) at (0.5,3) [draw, circle,inner sep=0pt,minimum size=5mm] {\phantom{\tiny 1.00}};
			\node	(c) at (1,2) [draw, circle,fill=lightred,inner sep=0pt,minimum size=5mm] {\tiny 1.00};
			\node	(d) at (1.5,1) [draw, circle,inner sep=0pt,minimum size=5mm] {\phantom{\tiny 1.00}};
			\node	(e) at (2,0) [draw, circle,fill=black, inner sep=0pt,minimum size=5mm] {\phantom{\tiny 1.00}};
			\node	(f) at (2.5,1) [draw, circle,inner sep=0pt,minimum size=5mm] {\phantom{\tiny 1.00}};
			\node	(g) at (3,2) [draw, circle,inner sep=0pt,minimum size=5mm] {\phantom{\tiny 1.00}};
			\node	(h) at (3.5,3) [draw, circle,fill=lightred,inner sep=0pt,minimum size=5mm] {\tiny 0.90};
			\node	(i) at (4.5,3) [draw, circle,fill=lightred,inner sep=0pt,minimum size=5mm] {\tiny 1.20};
			\node	(j) at (5,2) [draw, circle,fill=lightred, inner sep=0pt,minimum size=5mm] {\tiny 1.40};
			\node	(k) at (5.5,1) [draw, circle,fill=lightblue,inner sep=0pt,minimum size=5mm] {\tiny 1.00};
			\node	(l) at (6,2) [draw, circle,inner sep=0pt,minimum size=5mm] {\phantom{\tiny 1.00}};
			\node	(m) at (6.5,3) [draw, circle, fill=lightred, inner sep=0pt,minimum size=5mm] {\tiny 1.20};
			\node	(n) at (7,2) [draw, circle,fill=lightred, inner sep=0pt,minimum size=5mm] {\tiny 1.40};
			\node	(o) at (7.5,1) [draw, circle, fill=lightred, inner sep=0pt,minimum size=5mm] {\tiny 0.65};
			\node	(p) at (8,2) [draw, circle, inner sep=0pt,minimum size=5mm] {\phantom{\tiny 1.00}};
			\node	(q) at (8.5,3) [draw, circle, fill=lightred, inner sep=0pt,minimum size=5mm] {\tiny 1.20};
			\node	(r) at (9,2) [draw, circle, fill=lightred, inner sep=0pt,minimum size=5mm] {\tiny 1.25};
			\node	(s) at (9.5,1) [draw, circle, fill=lightred, inner sep=0pt,minimum size=5mm] {\tiny 1.45};
			\node	(t) at (10,2) [draw, circle, inner sep=0pt,minimum size=5mm] {\phantom{\tiny 1.00}};
			\node	(u) at (10.5,3) [draw, circle, fill=lightred, inner sep=0pt,minimum size=5mm] {\tiny 1.10};
			\node	(v) at (11,4) [draw, circle, fill=lightred, inner sep=0pt,minimum size=5mm] {\tiny 0.80};
			%\node	(w) at (11,4) [draw, circle, fill=lightred, inner sep=0pt,minimum size=5mm] {};
			\draw 
			(a) edge[->] (b) 
			(b) edge[->] (c)
			(c) edge[->] (d)
			(d) edge[->] (e)   
			(e) edge[->] node[near end,left] {\tiny 1.00} (f)
			(f) edge[->] node[near end,left] {\tiny 1.60} (g)
			(g) edge[->] node[near end,left] {\tiny 0.85} (h)
			(h) edge[->] (i)
			(i) edge[->] (j) 
			(j) edge[->] (k)
			(k) edge[->] node[near end,left] {\tiny 0.30} (l)
			(l) edge[->] node[near end,left] {\tiny 0.85} (m)   
			(m) edge[->] (n)
			(n) edge[->] (o)
			(o) edge[->] node[near end,left] {\tiny 0.30} (p)
			(p) edge[->] node[near end,left] {\tiny 0.85} (q)
			(q) edge[->] (r)
			(r) edge[->] (s)
			(s) edge[->] node[near end,left] {\tiny 0.30} (t)
			(t) edge[->] node[near end,left] {\tiny 0.85} (u)
			(u) edge[->] node[near end,left] {\tiny 1.35} (v)
			%(v) edge[->] (w)
			;
			\end{tikzpicture}
		}
		\caption{EP-2}
		\label{fig:ep-2}
	\end{subfigure}
  \par\bigskip
	\begin{subfigure}[b]{0.465\columnwidth}
		\scalebox{0.725}{%
			\begin{tikzpicture}
			\node   (h) at (-0.75, 4){$h$};
			\node   (2h) at (-0.75, 3){$2h$};
			\node   (4h) at (-0.75, 2){$4h$};
			\node   (8h) at (-0.75, 1){$8h$};
			\node   (16h) at (-0.75, 0){$16h$};
			\node	(a) at (0,4) [draw, circle,inner sep=0pt,minimum size=5mm] {\phantom{\tiny 1.00}};
			\node	(b) at (0.5,3) [draw, circle,inner sep=0pt,minimum size=5mm] {\phantom{\tiny 1.00}};
			\node	(c) at (1,2) [draw, circle,fill=lightblue,inner sep=0pt,minimum size=5mm] {\tiny 1.25};
			\node	(d) at (1.5,1) [draw, circle,inner sep=0pt,minimum size=5mm] {\phantom{\tiny 1.00}};
			\node	(e) at (2,0) [draw, circle,fill=black, inner sep=0pt,minimum size=5mm] {\phantom{\tiny 1.00}};
			\node	(f) at (2.5,1) [draw, circle,fill=lightred,inner sep=0pt,minimum size=5mm] {\tiny 1.05};
			\node	(g) at (3,2) [draw, circle,inner sep=0pt,minimum size=5mm] {\phantom{\tiny 1.00}};
			\node	(h) at (3.5,3) [draw, circle,fill=lightred,inner sep=0pt,minimum size=5mm] {\tiny 1.25};
			\node	(i) at (4,2) [draw, circle,fill=lightred,inner sep=0pt,minimum size=5mm] {\tiny 1.15};
			\node	(j) at (4.5,3) [draw, circle,fill=lightred, inner sep=0pt,minimum size=5mm] {\tiny 1.25};
			\node	(k) at (5,2) [draw, circle,fill=lightred,inner sep=0pt,minimum size=5mm] {\tiny 1.10};
			\node	(l) at (5.5,3) [draw, circle, fill=lightred, inner sep=0pt,minimum size=5mm] {\tiny 1.25};
			\node	(m) at (6,4) [draw, circle,fill=lightred, inner sep=0pt,minimum size=5mm] {\tiny 0.75};
			%\node	(w) at (11,4) [draw, circle, fill=lightred, inner sep=0pt,minimum size=5mm] {};
			\draw 
			(a) edge[->] (b) 
			(b) edge[->] (c)
			(c) edge[->] (d)
			(d) edge[->] (e)   
			(e) edge[->] node[near end,left] {\tiny 1.00} (f)
			(f) edge[->] node[near end,left] {\tiny 1.25} (g)
			(g) edge[->] node[near end,left] {\tiny 1.05} (h)
			(h) edge[->] (i)
			(i) edge[->] node[near end,left] {\tiny 0.90} (j) 
			(j) edge[->] (k)
			(k) edge[->] node[near end,left] {\tiny 0.90} (l)
			(l) edge[->] node[near end,left] {\tiny 1.40} (m)   
			;
			\end{tikzpicture}
		}
		\caption{EP-3}
		\label{fig:ep-3}
	\end{subfigure}
     \begin{subfigure}[b]{0.5\columnwidth}
		\scalebox{0.725}{%
			\begin{tikzpicture}
			\node	(a) at (0,4) [draw, circle,inner sep=0pt,minimum size=5mm] {\phantom{\tiny 1.00}};
			\node	(b) at (0.5,3) [draw, circle,fill=lightred,inner sep=0pt,minimum size=5mm] {\tiny 1.15};
			\node	(c) at (1,2) [draw,circle,fill=lightblue, inner sep=0pt,minimum size=5mm] {\tiny 1.55};
			\node	(d) at (1.5,1) [draw,circle,fill=lightred, inner sep=0pt,minimum size=5mm] {\tiny 1.05};
			\node	(e) at (2,0) [draw, circle,fill=black, inner sep=0pt,minimum size=5mm] {\phantom{\tiny 1.00}};
			\node	(f) at (2.5,1) [draw, circle, inner sep=0pt,minimum size=5mm] {\phantom{\tiny 1.00}};
			\node	(g) at (3,2) [draw, circle, inner sep=0pt,minimum size=5mm] {\phantom{\tiny 1.00}};
			\node	(h) at (3.5,3) [draw, circle,fill=lightred, inner sep=0pt,minimum size=5mm] {\tiny 1.00};
			\node	(i) at (4,2) [draw, circle,fill=lightred, inner sep=0pt,minimum size=5mm] {\tiny 0.95};
			\node	(j) at (4.5,3) [draw, circle, inner sep=0pt,minimum size=5mm] {\phantom{\tiny 1.00}};
			\node	(k) at (5,2) [draw, circle,fill=lightred, inner sep=0pt,minimum size=5mm] {\tiny 0.70};
			\node	(l) at (5.5,3) [draw, circle, inner sep=0pt,minimum size=5mm] {\phantom{\tiny 1.00}};
			\node	(m) at (6,4) [draw, circle,fill=lightred, inner sep=0pt,minimum size=5mm] {\tiny 1.00};
			\node	(n) at (6.5,3) [draw, circle, inner sep=0pt,minimum size=5mm] {\phantom{\tiny 1.00}};
			\node	(o) at (7,2) [draw, circle, fill=lightred,inner sep=0pt,minimum size=5mm] {\tiny 0.70};
			\node	(p) at (7.5,3) [draw, circle, inner sep=0pt,minimum size=5mm] {\phantom{\tiny 1.00}};
			\node	(q) at (8,4) [draw, circle, inner sep=0pt,minimum size=5mm] {\phantom{\tiny 1.00}};
			\draw 
			(a) edge[->] (b) 
			(b) edge[->] (c)
			(c) edge[->] (d)
			(d) edge[->] (e)   
			(e) edge[->] node[near end,left] {\tiny 1.00}  (f)
			(f) edge[->] node[near end,left] {\tiny 0.95} (g)
			(g) edge[->] node[near end,left] {\tiny 1.45} (h)
			(h) edge[->] (i)
			(i) edge[->] node[near end,left] {\tiny 1.30} (j) 
			(j) edge[->] (k)
			(k) edge[->] node[near end,left] {\tiny 1.10} (l)
			(l) edge[->] node[near end,left] {\tiny 1.65} (m)   
			(m) edge[->] (n)
			(n) edge[->] (o)
			(o) edge[->] node[near end,left] {\tiny 1.10} (p)
			(p) edge[->] node[near end,left] {\tiny 1.20} (q)
			;
			\end{tikzpicture}
		}
		\caption{EP-8}
		\label{fig:ep-8}
	\end{subfigure}
    \par\bigskip
  	\begin{subfigure}[b]{\columnwidth}
			\scalebox{0.725}{%
			\begin{tikzpicture}
			\node   (h) at (-0.75, 4){$h$};
			\node   (2h) at (-0.75, 3){$2h$};
			\node   (4h) at (-0.75, 2){$4h$};
			\node   (8h) at (-0.75, 1){$8h$};
			\node   (16h) at (-0.75, 0){$16h$};
			\node	(a) at (0,4) [draw, circle,inner sep=0pt,minimum size=5mm] {\phantom{\tiny 1.00}};
			\node	(b) at (0.5,3) [draw, circle,fill=lightred,inner sep=0pt,minimum size=5mm] {\tiny 1.10};
			\node	(c) at (1.5,3) [draw, circle,fill=lightblue,inner sep=0pt,minimum size=5mm] {\tiny 0.95};
			\node	(d) at (2,2) [draw, circle,inner sep=0pt,minimum size=5mm] {\phantom{\tiny 1.00}};
			\node	(e) at (2.5,1) [draw, circle,fill=lightblue,inner sep=0pt,minimum size=5mm] {\tiny 1.05};
			\node	(f) at (3,0) [draw, fill=black,circle,inner sep=0pt,minimum size=5mm] {\phantom{\tiny 1.00}};
			\node	(g) at (3.5,1) [draw, circle,inner sep=0pt,minimum size=5mm] {\phantom{\tiny 1.00}};
			\node	(h) at (4,2) [draw, circle,fill=lightred,inner sep=0pt,minimum size=5mm] {\tiny 1.70};
			\node	(i) at (4.5,1) [draw, circle,inner sep=0pt,minimum size=5mm] {\phantom{\tiny 1.00}};
			\node	(j) at (5,0) [draw, circle,fill=black,inner sep=0pt,minimum size=5mm] {\phantom{\tiny 1.00}};
			\node	(k) at (5.5,1) [draw, circle,inner sep=0pt,minimum size=5mm] {\phantom{\tiny 1.00}};
			\node	(l) at (6,2) [draw, circle,fill=lightblue,inner sep=0pt,minimum size=5mm] {\tiny 0.30};
			\node	(m) at (6.5,3) [draw, circle, inner sep=0pt,minimum size=5mm] {\phantom{\tiny 1.00}};
			\node	(n) at (7,2) [draw, circle,fill=lightred, inner sep=0pt,minimum size=5mm] {\tiny 0.45};
			\node	(o) at (7.5,3) [draw, circle, inner sep=0pt,minimum size=5mm] {\phantom{\tiny 1.00}};
			\node	(p) at (8,4) [draw, circle, inner sep=0pt,minimum size=5mm] {\phantom{\tiny 1.00}};
			\node	(q) at (8.5,3) [draw, circle, inner sep=0pt,minimum size=5mm] {\phantom{\tiny 1.00}};
			\node	(r) at (9,2) [draw, circle, fill=lightblue, inner sep=0pt,minimum size=5mm] {\tiny 1.70};
			\node	(s) at (9.5,3) [draw, circle, inner sep=0pt,minimum size=5mm] {\phantom{\tiny 1.00}};
			\node	(t) at (10,4) [draw, circle, inner sep=0pt,minimum size=5mm] {\phantom{\tiny 1.00}};
			\node	(u) at (10.5,3) [draw, circle, inner sep=0pt,minimum size=5mm] {\phantom{\tiny 1.00}};
			\node	(v) at (11,2) [draw, circle, fill=lightblue, inner sep=0pt,minimum size=5mm] {\tiny 1.05};
			\node	(w) at (11.5,3) [draw, circle, inner sep=0pt,minimum size=5mm] {\phantom{\tiny 1.00}};
   			\node	(x) at (12,4) [draw, circle, inner sep=0pt,minimum size=5mm] {\phantom{\tiny 1.00}};
      		\node	(y) at (12.5,3) [draw, circle, inner sep=0pt,minimum size=5mm] {\phantom{\tiny 1.00}};
			\node	(z) at (13,2) [draw, circle, fill=lightblue, inner sep=0pt,minimum size=5mm] {\tiny 0.30};
			\node	(aa) at (13.5,3) [draw, circle, inner sep=0pt,minimum size=5mm] {\phantom{\tiny 1.00}};
   			\node	(bb) at (14,4) [draw, circle, inner sep=0pt,minimum size=5mm] {\phantom{\tiny 1.00}};
            \node	(cc) at (14.5,3) [draw, circle, fill=lightblue, inner sep=0pt,minimum size=5mm] {\tiny 1.25};
   			\node	(dd) at (15,4) [draw, circle, fill=lightred, inner sep=0pt,minimum size=5mm] {\tiny 1.05};
			\draw 
			(a) edge[->] (b) 
			(b) edge[->] (c)
			(c) edge[->] (d)
			(d) edge[->] (e)   
			(e) edge[->] (f)
			(f) edge[->] node[near end,left] {\tiny 1.00} (g)
			(g) edge[->] node[near end,left] {\tiny 0.80} (h)
			(h) edge[->] (i)
			(i) edge[->] (j) 
			(j) edge[->] node[near end,left] {\tiny 1.00} (k)
			(k) edge[->] node[near end,left] {\tiny 0.55} (l)
			(l) edge[->] node[near end,left] {\tiny 1.70} (m)   
			(m) edge[->] (n)
			(n) edge[->] node[near end,left] {\tiny 1.80} (o)
			(o) edge[->] node[near end,left] {\tiny 1.70} (p)
			(p) edge[->] (q)
			(q) edge[->] (r)
			(r) edge[->] node[near end,left] {\tiny 0.95} (s)
			(s) edge[->] node[near end,left] {\tiny 1.10} (t)
			(t) edge[->] (u)
			(u) edge[->] (v)
			(v) edge[->] node[near end,left] {\tiny 0.55} (w)
   			(w) edge[->] node[near end,left] {\tiny 0.95} (x)
      		(x) edge[->] (y)
			(y) edge[->] (z)
			(z) edge[->] node[near end,left] {\tiny 0.95} (aa)
   			(aa) edge[->] node[near end,left] {\tiny 1.10} (bb)
      		(bb) edge[->] (cc)
   			(cc) edge[->] node[near end,left] {\tiny 0.90} (dd)
			;
			\end{tikzpicture}
		}
		\caption{EP-4}
		\label{fig:ep-4}
	\end{subfigure}
	\caption[2D Helmholtz -- Algorithmic structure of the discovered multigrid preconditioners EP-2, EP-3, EP-8 and EP-4]{Algorithmic structure of the discovered multigrid preconditioners EP-2, EP-3, EP-8 and EP-4. The color of the node denotes the type of operation. Black: Coarse-grid solver, Blue: Pointwise Jacobi smoothing, Red: Red-black Gauss-Seidel smoothing, White: No operation. The relaxation factor of each smoothing step is included in each node, while for coarse-grid correction, it is attached to the respective edge.}
	\label{fig:structure-evolved-preconditioners1}
\end{figure}

\begin{figure}
\captionsetup[subfigure]{justification=centering}
	\begin{subfigure}[b]{\columnwidth}
		\scalebox{0.725}{%
			\begin{tikzpicture}
			\node   (h) at (-0.75, 4){$h$};
			\node   (2h) at (-0.75, 3){$2h$};
			\node   (4h) at (-0.75, 2){$4h$};
			\node   (8h) at (-0.75, 1){$8h$};
			\node   (16h) at (-0.75, 0){$16h$};
			\node	(a) at (0,4) [draw, circle,inner sep=0pt,minimum size=5mm] {\phantom{\tiny 1.00}};
			\node	(b) at (0.5,3) [draw, circle,fill=lightred,inner sep=0pt,minimum size=5mm] {\tiny 1.30};
			\node	(c) at (1,2) [draw,circle,fill=lightblue, inner sep=0pt,minimum size=5mm] {\tiny 1.45};
			\node	(d) at (1.5,1) [draw,circle,fill=lightblue, inner sep=0pt,minimum size=5mm] {\tiny 1.35};
			\node	(e) at (2,0) [draw, circle,fill=black, inner sep=0pt,minimum size=5mm] {\phantom{\tiny 1.00}};
			\node	(f) at (2.5,1) [draw, circle, inner sep=0pt,minimum size=5mm] {\phantom{\tiny 1.00}};
			\node	(g) at (3,2) [draw, circle, inner sep=0pt,minimum size=5mm] {\phantom{\tiny 1.00}};
			\node	(h) at (3.5,3) [draw, circle, inner sep=0pt,minimum size=5mm] {\phantom{\tiny 1.00}};
			\node	(i) at (4,4) [draw, circle,fill=lightred, inner sep=0pt,minimum size=5mm] {\tiny 1.20};
			\node	(j) at (4.5,3) [draw, circle, inner sep=0pt,minimum size=5mm] {\phantom{\tiny 1.00}};
			\node	(k) at (5,2) [draw, circle,fill=lightred, inner sep=0pt,minimum size=5mm] {\tiny 1.30};
			\node	(l) at (5.5,3) [draw, circle, inner sep=0pt,minimum size=5mm] {\phantom{\tiny 1.00}};
			\node	(m) at (6,2) [draw, circle,fill=lightblue, inner sep=0pt,minimum size=5mm] {\tiny 1.80};
			\node	(n) at (7,2) [draw, circle,fill=lightblue, inner sep=0pt,minimum size=5mm] {\tiny 0.55};
			\node	(o) at (7.5,3) [draw, circle, inner sep=0pt,minimum size=5mm] {\phantom{\tiny 1.00}};
			\node	(p) at (8,2) [draw, circle, fill=lightred, inner sep=0pt,minimum size=5mm] {\tiny 1.80};
			\node	(q) at (8.5,1) [draw, circle, fill=lightred, inner sep=0pt,minimum size=5mm] {\tiny 1.45};
			\node	(r) at (9,2) [draw, circle, inner sep=0pt,minimum size=5mm] {\phantom{\tiny 1.00}};
			\node	(s) at (9.5,3) [draw, circle, inner sep=0pt,minimum size=5mm] {\phantom{\tiny 1.00}};
			\node	(t) at (10,2) [draw, circle, fill=lightblue, inner sep=0pt,minimum size=5mm] {\tiny 0.95};
			\node	(u) at (11,2) [draw, circle, fill=lightred, inner sep=0pt,minimum size=5mm] {\tiny 0.70};
			\node	(v) at (11.5,3) [draw, circle, fill=lightred, inner sep=0pt,minimum size=5mm] {\tiny 1.00};
			\node	(w) at (12,4) [draw, circle, fill=lightred, inner sep=0pt,minimum size=5mm] {\tiny 0.90};
			\draw 
			(a) edge[->] (b) 
			(b) edge[->] (c)
			(c) edge[->] (d)
			(d) edge[->] (e)   
			(e) edge[->] node[near end,left] {\tiny 1.00}  (f)
			(f) edge[->] node[near end,left] {\tiny 1.00} (g)
			(g) edge[->] node[near end,left] {\tiny 0.80} (h)
			(h) edge[->] node[near end,left] {\tiny 1.60} (i)
			(i) edge[->] (j) 
			(j) edge[->] (k)
			(k) edge[->] node[near end,left] {\tiny 0.85} (l)
			(l) edge[->] (m)   
			(m) edge[->] (n)
			(n) edge[->] node[near end,left] {\tiny 0.85} (o)
			(o) edge[->] (p)
			(p) edge[->] (q)
			(q) edge[->] node[near end,left] {\tiny 1.05} (r)
			(r) edge[->] node[near end,left] {\tiny 1.10} (s)
			(s) edge[->] (t)
			(t) edge[->] (u)
			(u) edge[->] node[near end,left] {\tiny 1.05} (v)
			(v) edge[->] node[near end,left] {\tiny 1.80} (w)
			;
			\end{tikzpicture}
		}
		\caption{EP-5}
		\label{fig:ep-5}
	\end{subfigure}
  \par\bigskip
 	\begin{subfigure}[b]{\columnwidth}
			\scalebox{0.725}{%
			\begin{tikzpicture}
			\node   (h) at (-0.75, 4){$h$};
			\node   (2h) at (-0.75, 3){$2h$};
			\node   (4h) at (-0.75, 2){$4h$};
			\node   (8h) at (-0.75, 1){$8h$};
			\node   (16h) at (-0.75, 0){$16h$};
			\node	(a) at (0,4) [draw, circle,inner sep=0pt,minimum size=5mm] {\phantom{\tiny 1.00}};
			\node	(b) at (0.5,3) [draw, circle,inner sep=0pt,minimum size=5mm] {\phantom{\tiny 1.00}};
			\node	(c) at (1,2) [draw, circle,inner sep=0pt,minimum size=5mm] {\phantom{\tiny 1.00}};
			\node	(d) at (1.5,1) [draw, circle,inner sep=0pt,minimum size=5mm] {\phantom{\tiny 1.00}};
			\node	(e) at (2,0) [draw, circle,fill=black, inner sep=0pt,minimum size=5mm] {\phantom{\tiny 1.00}};
			\node	(f) at (2.5,1) [draw, circle,inner sep=0pt,minimum size=5mm] {\phantom{\tiny 1.00}};
			\node	(g) at (3,2) [draw, circle,fill=lightred,inner sep=0pt,minimum size=5mm] {\tiny 0.55};
			\node	(h) at (3.5,3) [draw, circle,fill=lightred,inner sep=0pt,minimum size=5mm] {\tiny 1.40};
			\node	(i) at (4,2) [draw, circle,fill=lightred,inner sep=0pt,minimum size=5mm] {\tiny 1.15};
			\node	(j) at (4.5,3) [draw, circle,fill=lightred, inner sep=0pt,minimum size=5mm] {\tiny 1.10};
			\node	(k) at (5,2) [draw, circle,fill=lightblue,inner sep=0pt,minimum size=5mm] {\tiny 1.50};
			\node	(l) at (6,2) [draw, circle,fill=lightred,inner sep=0pt,minimum size=5mm] {\tiny 0.20};
			\node	(m) at (7,2) [draw, circle, fill=lightred, inner sep=0pt,minimum size=5mm] {\tiny 0.55};
			\node	(n) at (7.5,3) [draw, circle,fill=lightred, inner sep=0pt,minimum size=5mm] {\tiny 1.40};
			\node	(o) at (8,2) [draw, circle, fill=lightred, inner sep=0pt,minimum size=5mm] {\tiny 0.65};
			\node	(p) at (8.5,3) [draw, circle,fill=lightred, inner sep=0pt,minimum size=5mm] {\tiny 1.10};
			\node	(q) at (9,2) [draw, circle, fill=lightblue, inner sep=0pt,minimum size=5mm] {\tiny 0.25};
			\node	(r) at (10,2) [draw, circle, fill=lightred, inner sep=0pt,minimum size=5mm] {\tiny 1.05};
			\node	(s) at (10.5,3) [draw, circle, inner sep=0pt,minimum size=5mm] {\phantom{\tiny 1.00}};
			\node	(t) at (11,4) [draw, circle, inner sep=0pt,minimum size=5mm] {\phantom{\tiny 1.00}};
			\node	(u) at (11.5,3) [draw, circle, inner sep=0pt,minimum size=5mm] {\phantom{\tiny 1.00}};
			\node	(v) at (12,2) [draw, circle, fill=lightred, inner sep=0pt,minimum size=5mm] {\tiny 0.65};
			\node	(w) at (12.5,3) [draw, circle, fill=lightred, inner sep=0pt,minimum size=5mm] {\tiny 1.25};
   			\node	(x) at (13,4) [draw, circle, fill=lightred, inner sep=0pt,minimum size=5mm] {\tiny 0.80};
			\draw 
			(a) edge[->] (b) 
			(b) edge[->] (c)
			(c) edge[->] (d)
			(d) edge[->] (e)   
			(e) edge[->] node[near end,left] {\tiny 1.00} (f)
			(f) edge[->] node[near end,left] {\tiny 0.70} (g)
			(g) edge[->] node[near end,left] {\tiny 1.10} (h)
			(h) edge[->] (i)
			(i) edge[->] node[near end,left] {\tiny 0.85} (j) 
			(j) edge[->] (k)
			(k) edge[->] (l)
			(l) edge[->] (m)   
			(m) edge[->] node[near end,left] {\tiny 1.00} (n)
			(n) edge[->] (o)
			(o) edge[->] node[near end,left] {\tiny 1.40} (p)
			(p) edge[->] (q)
			(q) edge[->] (r)
			(r) edge[->] node[near end,left] {\tiny 1.70} (s)
			(s) edge[->] node[near end,left] {\tiny 1.65} (t)
			(t) edge[->] (u)
			(u) edge[->] (v)
			(v) edge[->] node[near end,left] {\tiny 1.40} (w)
   			(w) edge[->] node[near end,left] {\tiny 0.80} (x)
			;
			\end{tikzpicture}
		}
		\caption{EP-10}
		\label{fig:ep-10}
	\end{subfigure}
	\caption[2D Helmholtz -- Algorithmic structure of the discovered multigrid preconditioners EP-5 and EP-10]{Algorithmic structure of the discovered multigrid preconditioners EP-5 and EP-10. %The color of the node denotes the type of operation. Black: Coarse-grid solver, Blue: Pointwise Jacobi smoothing, Red: Red-black Gauss-Seidel smoothing, White: No operation. The relaxation factor of each smoothing step is included in each node, while for coarse-grid correction, it is attached to the respective edge.
 }
	\label{fig:structure-evolved-preconditioners2}
\end{figure}
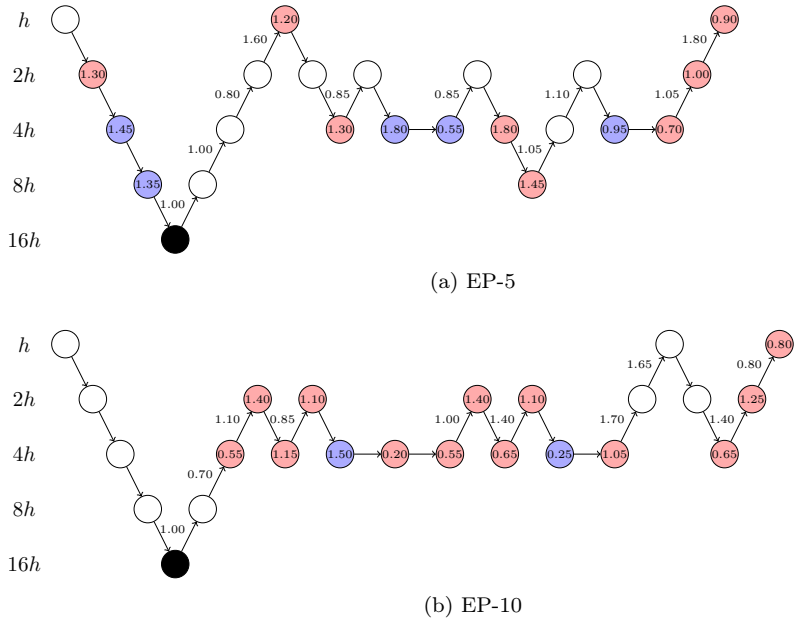

While the results shown in Table~\ref{table:evolved-solvers} indicate that our automatically-designed preconditioners yield faster solvers than classical multigrid cycles for a wide range of different wavenumbers, we have not investigated how these methods are able to achieve this feat.
For this purpose, we first need to gain an understanding of the algorithmic features of these methods and how they compare to those of classical multigrid cycles.
Figure~\ref{fig:structure-evolved-preconditioners1} and~\ref{fig:structure-evolved-preconditioners2} contain the graphical representations of those automatically-designed preconditioners that yield a converging solver for a wavenumber of $k = 640$.
Note that we are particularly interested in this case since all V-, F-, and W-cycles fail to achieve convergence when applied to this problem instance.
As illustrated by their graphical representations, the employed sequence of multigrid operations varies significantly between the individual methods, and none of them can easily be characterized as one of the classical multigrid cycles.
However, we can still identify a number of common characteristics.
While the search space also includes different block Jacobi smoothers, all methods employ exclusively pointwise Jacobi and RB-GS smoothing, whereby, with the exception of EP-3, the latter is used predominantly.
If we consider the total amount of smoothing within each method, on average, the Jacobi method is used in $1/3$ of the smoothing steps.
However, in contrast to classical multigrid cycles, which usually utilize the same smoother with a fixed number of steps on each level, our evolutionary algorithm can alter each operation independently, yielding varying smoothers and relaxation factors in each individual step.
If we examine the computational structure of each of the six multigrid-based preconditioners more closely, we can see that many coarse-grid correction steps originate from one or multiple smoothing operations on intermediate levels, often using a different combination of relaxation factors.
With the exception of EP-4, all methods apply the coarse-grid solver only once but use multiple smoothing-based coarse-grid corrections on intermediate levels, whereby often individual smoothing steps are completely skipped.
Considering that EP-5 represents the overall most efficient preconditioner for the considered range of wavenumbers, a combination of different smoothers and relaxation factors within these correction steps seems to be the most effective approach.
Furthermore, in contrast to classical multigrid cycles, the number of smoothing and coarse-grid correction steps is unevenly distributed over the range of levels.
Table~\ref{table:operations-per-level} shows the average number of smoothing steps, coarse-grid correction steps, and the number of updates of the approximate solution\footnote{Note that this includes both smoothing and coarse-grid correction steps} on every level.
\begin{table}
	\caption[2D Helmholtz -- Operation statistics]{Statistics about the operations that are performed by the multigrid-based preconditioners yielding a converging solver for $k = 640$.}
	\label{table:operations-per-level}
	\centering
	\begin{tabular}{l c c c c c c c}
		\toprule
  		& \multicolumn{3}{c}{Smoothing} & \multicolumn{3}{c}{Coarse-Grid Correction} & \\
		\cmidrule(l){2-4} \cmidrule(l){5-7}
		Level & Mean & Min & Max & Mean & Min & Max & Mean Updates \\
        \midrule
        $h$ & $1.17$ & $1$ & $2$ & $2.17$ & $1$ & $5$ & $3.33$\\
        \midrule
        $2h$ & $3.33$ & $2$ & $5$ & $4.50$ & $3$ & $6$ & $7.83$ \\
        \midrule
        $4h$ & $5.50$ & $3$ & $9$ & $1.83$ & $1$ & $4$ & $7.33$ \\
        \midrule
        $8h$ & $1.33$ & $0$ & $3$ & $1.17$ & $1$ & $2$ & $2.50$ \\
		\bottomrule
	\end{tabular}
\end{table}
As this analysis reveals, the number of operations on the intermediate levels $2h$ and $4h$ is significantly larger, which means that these operations are most cost-effective for accelerating the convergence of the preconditioned iterative solver.
Here it is especially interesting that instead of applying the combined error reduction of all smoothing steps in a single coarse-grid correction, as it is usually performed in multigrid cycles, the correction is split into multiple steps.
As a consequence, the average number of coarse-grid corrections is roughly equal to the number of smoothing steps.
To gain a better understanding of the consequences of this strategy, consider the EP-3 method, which uses the lowest total number of operations and, therefore, leads to the highest number of iterations.
If we compare the solving time and the number of iterations this method achieves according to Table~\ref{table:evolved-solvers}, with the characteristics of the V-, F- and W-cycles shown in Table~\ref{table:reference-methods-helmholtz}, we can see that it converges as fast as the W(3,3)-cycle while requiring 69 - 75 \% less time per iteration.
While especially EP-2, EP-5 and EP-10 utilize a higher number of coarse-grid correction and smoothing steps to achieve faster convergence, the simplicity of the EP-3 method demonstrates that the addition of only two smoothing-based correction steps on intermediate levels using the right combination of smoothers and relaxation factors yields a remarkably efficient preconditioner for the indefinite Helmholtz equation.
The EP-4 method represents another noteworthy exception. 
In contrast to the other multigrid-based preconditioners, which are more similar to V-cycles, as they all apply the coarse-grid solver only once, this method combines a four-grid W-cycle with a series of smoothing-based coarse-grid corrections.
Each of these corrections is structurally similar to a three-grid V-cycle that does not use any smoothing while replacing the coarse-grid solver with a single Jacobi step.
In general, the distinct sequences of operations shown in Figure~\ref{fig:structure-evolved-preconditioners1} and~\ref{fig:structure-evolved-preconditioners2} demonstrate that our evolutionary program synthesis approach can discover multigrid methods with novel algorithmic strategies that have not been investigated in the literature before.
Finally, since the use of varying relaxation factors is integral to all six considered multigrid-based preconditioners, it is important to investigate this feature in more detail. 
Figure~\ref{fig:histograms} shows the distribution of relaxation factors over all RB-GS, Jacobi, and coarse-grid correction steps within Figure~\ref{fig:structure-evolved-preconditioners1} and~\ref{fig:structure-evolved-preconditioners2} in the form of a histogram.
\begin{figure}
\captionsetup[subfigure]{justification=centering}
	\centering
\begin{subfigure}[b]{0.495\columnwidth}
    \includegraphics[width=\columnwidth]{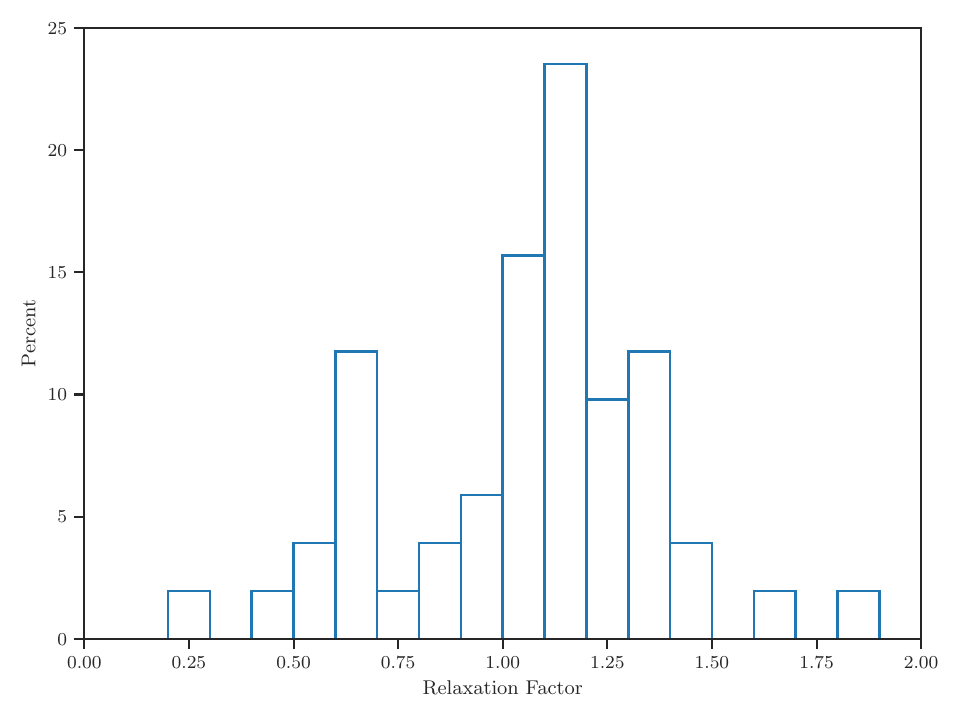}
	\caption{Red-black Gauss-Seidel}
	\label{fig:histogram-rbgs}
\end{subfigure}
	\centering
\begin{subfigure}[b]{0.495\columnwidth}
	\centering
    \includegraphics[width=\columnwidth]{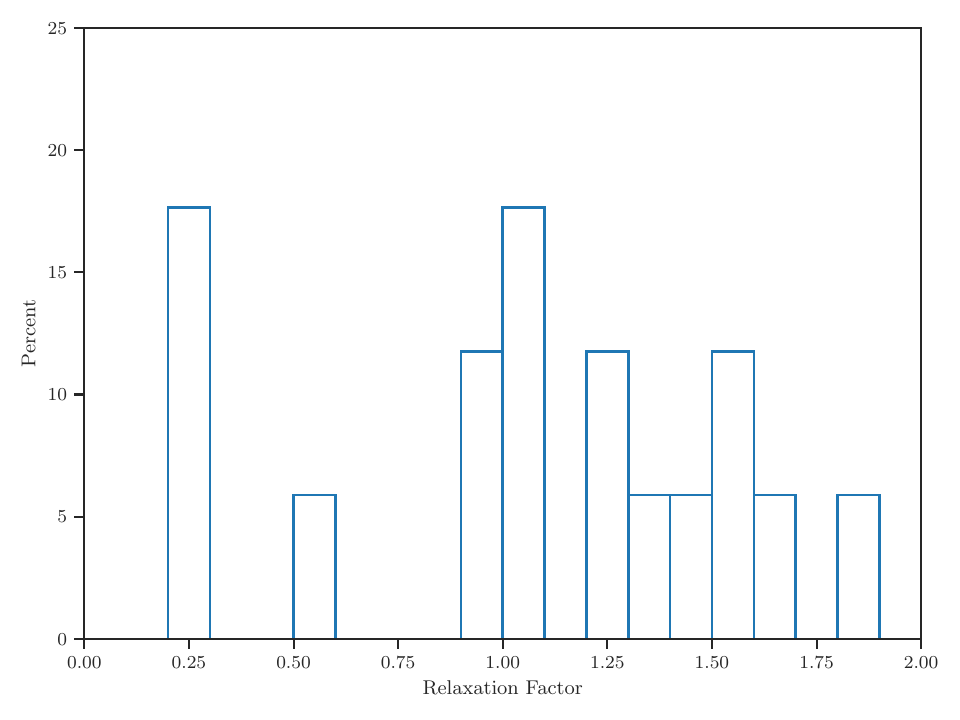}
	\caption{Jacobi}
	\label{fig:histogram-jacobi}
\end{subfigure}
\par\bigskip
	\centering
\begin{subfigure}[b]{0.495\columnwidth}
    \includegraphics[width=\columnwidth]{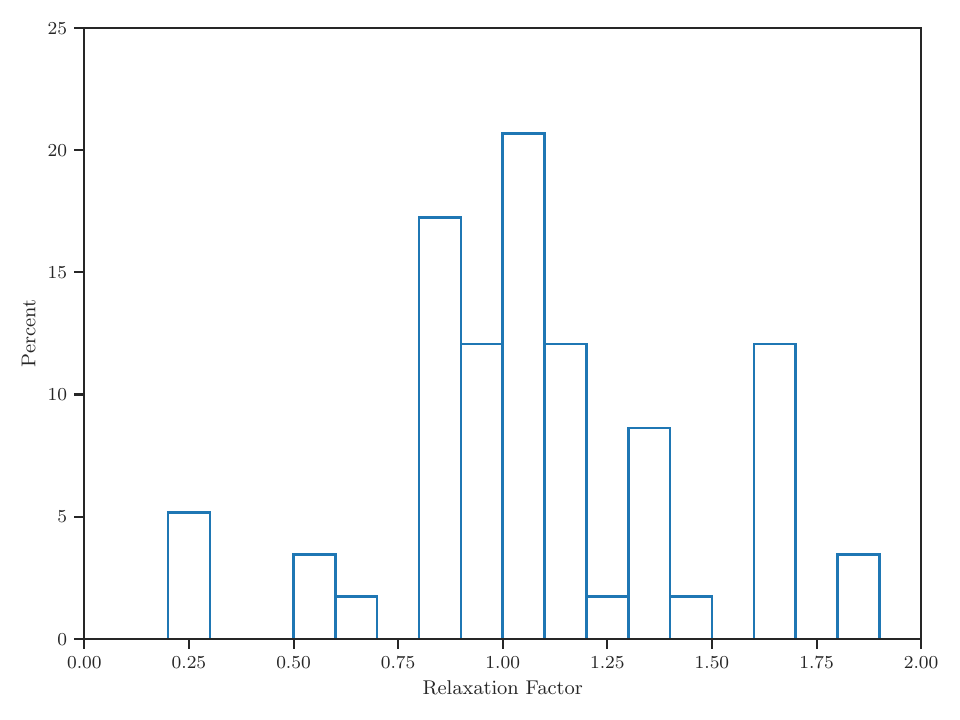}
	\caption{Coarse-grid correction}
	\label{fig:histogram-cgc}
\end{subfigure}
    \caption[2D Helmholtz -- Histogram of the relaxation factors]{Histogram of the relaxation factor distributions of the six multigrid-based preconditioners that yield a converging solver for $k = 640$.}
    \label{fig:histograms}
\end{figure}
First of all, note that both for RB-GS and the coarse-grid correction, most relaxation factors are concentrated around the value one, whereby the former slightly favors higher values and thus more overrelaxation.
However, in both cases, also distinct peaks comprising roughly 12 \% of the values can be found outside this range, which are located around 0.65 in Figure~\ref{fig:histogram-rbgs} and 1.65 in Figure~\ref{fig:histogram-cgc}.
In the case of Jacobi smoothing, which is shown in Figure~\ref{fig:histogram-jacobi}, the relaxation factor values are more scattered. 
While as in the case of RB-GS, overrelaxation has a higher prevalence, a drastic underrelaxation with $\omega \in \left[0.25, 0.30\right]$ is also used in 17.5 \% of the smoothing steps.
In general, this analysis underlines our observation that a combination of varying relaxation factors in consecutive smoothing and coarse-grid correction steps has the potential to improve the effectiveness of a multigrid method, as illustrated by the faster speed of convergence attained with our automatically-designed preconditioners compared to most classical multigrid cycles. 
However, note that especially in the case of the Jacobi method, which is employed in only $1/3$ of the smoothing steps, the distribution shown in Figure~\ref{fig:histogram-jacobi} comprises a certain degree of statistical uncertainty.

%% file: contents/related_work.tex
\begin{figure}
	\includegraphics[width=\textwidth]{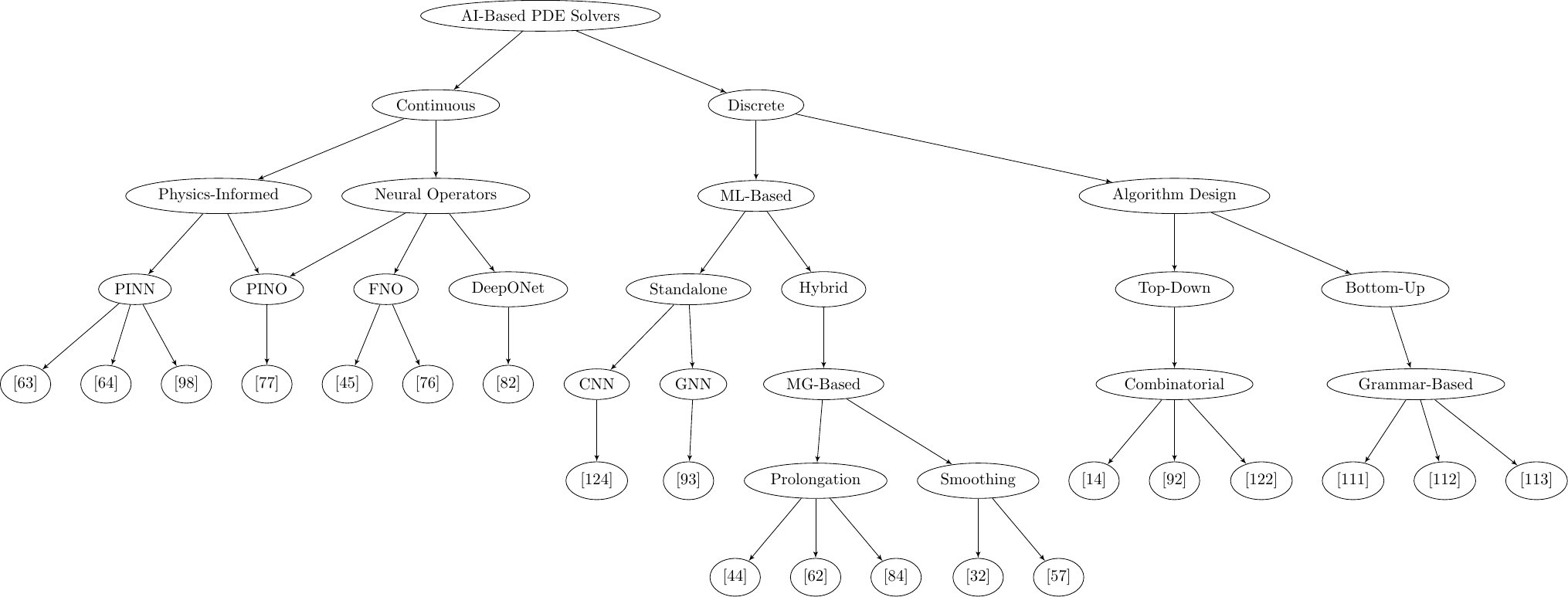}
	\caption[Overview of artificial intelligence-based methods for solving partial differential equations]{Overview of artificial intelligence-based methods for solving partial differential equations. Each inner node represents a different class of methods introduced in at least one research paper.}
	\label{fig:overview-ai-based-methods}
\end{figure}

%\section{Comparison with Related Work}
In the last chapter, we could demonstrate that our evolutionary program synthesis approach leads to the discovery of multigrid methods with novel algorithmic features.
By applying the generalization procedure presented in  Section~\ref{sec:generalization}, we were able to evolve multigrid methods that yield a high degree of efficiency in solving different instances of the same discretized PDE, one of them being intractable using classical multigrid methods.
However, since the application of artificial intelligence (AI) and automated algorithm design to the solution of PDEs includes a wide range of different methods, it is important to classify the approach presented here within this quickly growing field.
% In contrast to other %TODO provide more explanation
% We can classify AI-based methods for solving PDEs according to different criteria.
Figure~\ref{fig:overview-ai-based-methods} gives an overview of the state of AI-based methods for solving PDEs at the time this thesis was published.
First of all, we can distinguish methods that aim to solve a PDE directly in the continuous domain and those that require it to be formulated as a discrete problem, usually obtained by applying a specific discretization method.
The currently most popular\footnote{Here, we consider the number of citations as a popularity metric.} methods based on machine learning (ML), physics-informed neural networks~\cite{karniadakis2021physics,raissi2019physics,kharazmi2019variational,kharazmi2021hp} and neural operators~\cite{li2020fourier,guibas2021efficient,lu2021learning,li2021physics} fall into the former category.
This class of methods aims to approximate the function that represents the solution or the operator of a given PDE by exploiting the fact that neural networks can act as universal function approximators when given enough data~\cite{hornik1989multilayer}.
Physics-informed methods try to improve over purely data-driven methods by directly incorporating the physical constraints of the given PDE into the learning process.
Neural operators aim to achieve a higher degree of generalization by learning a representation of the operator of a PDE instead of approximating its solution.
Recently, the usability of ML-based methods has been tremendously improved through the availability of easy-to-use and well-maintained implementations, such as DeepXDE~\cite{lu2021deepxde} and NVIDIA Modulus~\cite{hennigh2021nvidia}. 
Instead of directly targeting a PDE in the continuous domain, the second branch of AI-based methods operates on its discrete version, obtained after applying a suitable discretization method.
Therefore, these methods can either act as a direct replacement for classical numerical solvers or operate in combination with them.
An early example of the former is the neural network-based PDE solver proposed by Lagaris et al.~\cite{lagaris1998artificial} but also more recent approaches based on convolutional~\cite{thuerey2020deep} and graph neural networks~\cite{pfaff2020learning} fall into this category.
In contrast, AI-based methods that work in combination with existing solvers do not try to replace the method as a whole but rather aim to enhance it, for instance, by adding or replacing certain steps of the method or by finding an optimal configuration for each of its parameters and design options.
Multigrid methods are an often-considered target since these methods have the potential to achieve an optimal asymptotic complexity while possessing a large number of configuration options and complex interactions between each of their components.
A first step towards the automated design of multigrid methods has been the work by Oosterlee and Wienands~\cite{oosterlee2003genetic}, which uses a genetic algorithm to optimize the choice of each multigrid component.
Similarly, Thekale et al.~\cite{thekale2010optimizing} aim to optimize the number of multigrid cycles within a full-multigrid method (FMG) using a branch-and-bound approach, while Brown et al.~\cite{brown2021tuning} optimize the algorithmic parameters of a two-grid method by solving a minimax problem obtained from an LFA-based analysis.
All these works have in common that they are based on the classical formulation of multigrid methods, as shown in Algorithm~\ref{alg:multigrid-cycle}.
Since the resulting optimization capabilities are still restricted to a set of global parameters, these approaches can be classified as top-down algorithm design or algorithm configuration methods.
In contrast, as we have shown in this work, expressing each possible computational step of a multigrid method as a separate production within a context-free grammar allows treating the task of designing an optimal multigrid method as a program synthesis problem.
%Therefore, according to our classification of algorithm design methods in Chapter~%TODO insert ref her
Therefore, together with this thesis, the papers~\cite{schmitt2020constructing,schmitt2021evostencils,schmitt2022evolving} can be considered the first implementation of a bottom-up approach for the automated design of multigrid methods.
However, while our approach offers the flexibility to construct arbitrary sequences of multigrid operations on a given hierarchy of discretizations, similar to~\cite{oosterlee2003genetic,thekale2010optimizing,brown2021tuning}, we consider the internal structure of each individual operation as immutable.
Recently, ML-based approaches have been utilized to enhance or replace certain operations within a multigrid method.
A first example is the works by Katrutsa et al.~\cite{katrutsa2020black}, Greenfeld et al.~\cite{greenfeld2019learning}, and Luz et at.~\cite{luz2020learning}, which utilizes ML to discover optimized prolongation operators.
Huang et al.~\cite{huang2021learning}, and Fanaskov~\cite{fanaskov2021neural} apply a similar approach to the optimization of smoothers.
The works by Taghibakhshi et al.~\cite{taghibakhshi2021optimization}, and Markidis~\cite{markidis2021old} go even one step further.
In the former, the authors replace the coarsening step within an algebraic multigrid method altogether with an ML system, while in~\cite{markidis2021old}, the author proposes to employ a physics-informed neural network as a coarse-grid solver.
Finally, Hsieh et al.~\cite{hsieh2019learning} take a different direction by enhancing the approximation obtained in each step of an iterative solver with an additional neural network-based correction.
Since these approaches all focus on the optimization of the individual components of a multigrid method but consider the method's algorithmic structure as immutable, they can be considered as a complementary approach to the algorithm design methods discussed earlier in this section.
However, a combination of both paradigms has not yet been considered and could be a promising research direction for the future.
One possibility to achieve this goal would be to incorporate learned operators or even ML-based methods that act as a replacement for certain solver components into the search space of a top-down or bottom-up algorithm design method.
%We will discuss this possibility together with other potential extensions of our approach in the following final section of this thesis.

%% file: contents/conclusion.tex
\paragraph{Conclusion}
We have started this thesis with an introduction to automated algorithm design and its successful application to different domains.
Since our main goal has been to harness the potential of evolutionary program synthesis for the automated design of efficient multigrid methods, the question of whether we have achieved our original goals remains to be answered.
After establishing a theoretical foundation about multigrid methods, formal languages, and evolutionary program synthesis, in Chapter~\ref{chapter:multigrid-formal-language}, we have derived a novel context-free grammar, which enables us to construct sequences of multigrid operations that can not be obtained with the classical formulation of these methods.
Building on this foundation, we have presented a prototypical implementation of our evolutionary program synthesis tool, called \emph{EvoStencils}, which automates the design and implementation of efficient and generalizable multigrid-based PDE solvers by leveraging the capabilities of the evolutionary computation library DEAP~\cite{rainville2012deap} and the ExaStencils~\cite{lengauer2020exastencils} code generation framework.
In Chapter~\ref{chapter:experiments}, we could then finally demonstrate that this approach yields efficient and generalizable multigrid methods for different PDEs, whereby using the example of the indefinite Helmholtz equation, our automatically-designed solvers were able to achieve super-human performance in an outstandingly-difficult benchmark problem for the application of numerical PDE solvers\footnote{For this result, the corresponding paper~\cite{schmitt2022evolving} has been awarded the 2022 Humies Gold Award for Human-Competitive Results (\url{https://www.human-competitive.org/awards})}.
While we believe that this lays the foundation for the utilization of automated algorithm design methods within the domain of numerical PDE solvers, the multigrid methods considered in this work represent only a tiny fraction of this vast research area.
%We, therefore, believe that there exists a multitude of promising extensions of this work, some of which we briefly want to discuss in the following.
A promising extension of this approach would thus be its application to other multigrid variants or even other classes of numerical solvers.
While this thesis focuses on classical geometric multigrid (GMG) methods, since the invention of these methods, several other variants tailored to different use cases have been developed.
One example is the full-approximation scheme (FAS), which can be considered a non-linear version of multigrid~\cite{trottenberg2000multigrid,briggs2000multigrid}.
While the formulation of an FAS method requires replacing both the smoother and coarse-grid solver with non-linear variants that are, for instance, based on Newton's or Picard's method, the application of our program synthesis approach requires only a minimal amount of adaption.
To illustrate this, Section~\ref{appendix:fas} shows the necessary adaptions of our original state transition functions and grammar productions for generating FAS-style multigrid methods. 
As it can be seen there, apart from adjusting the respective productions, only the two functions \textsc{coarsening} and \textsc{cgc} need to be changed, while, additionally, a single new state transition function \textsc{cgs} for the application of the coarse-grid solver needs to be provided. 
Note that since the application of the operator $A_h$ no longer represents a linear operation in the case of FAS, we instead denote it as a function application.
Another promising idea would be an extension of our approach to full-multigrid methods.
While in Section~\ref{sec:experiments-part1}, we have demonstrated that the multigrid methods designed with our approach can solve different PDEs faster than common cycles, in many linear cases, including Poisson's equation, FMG potentially represents the fastest and most efficient solver available~\cite{trottenberg2000multigrid}.
Since FMG is based on applying multigrid cycles on different levels of a given discretization hierarchy, one possibility would be to apply our evolutionary program synthesis method to each of these individual cycles.
While this approach would result in an even larger search space than the ones considered in this work, it could yield methods that achieve an even higher degree of efficiency in solving PDE-based problems.
Finally, another class of multigrid methods not yet considered in this work is algebraic multigrid (AMG).
In contrast to GMG, AMG methods derive their operations directly from an algebraic formulation of the system of linear equations in the form of a (usually sparse) matrix.
While this makes these methods fundamentally different from the multigrid methods considered in this work, the actual algorithmic structure of AMG methods is similar to its geometric counterpart.
We could, therefore, formulate a grammar structurally similar to the one shown in Algorithm~\ref{table:multigrid-grammar} that replaces each individual operation applied within their productions with their algebraic equivalent, which would allow us to utilize the same program synthesis approach for the automated design of AMG methods.
Apart from considering different multigrid variants, another interesting direction would be to incorporate additional components and even completely different classes of solvers into the algorithm design space.
As Figure~\ref{fig:overview-ai-based-methods} illustrates, recently, a whole new class of AI-based solvers has been developed, such as data-driven~\cite{thuerey2020deep} and physics-informed~\cite{karniadakis2021physics} surrogate models.
While these methods show promise in individual domains, it is not yet clear how these methods can be integrated or combined with established numerical solvers.
Automated algorithm design offers an attractive solution to this problem by integrating these methods into the corresponding search space, for instance, in the form of the class of multigrid grammars introduced in this thesis.
For instance, an idea proposed in~\cite{markidis2021old} is to utilize physics-informed neural networks (PINNs) as a coarse-grid solver within multigrid methods.
A second complementary extension would be to target the current limitations in the evaluation accuracy of predictive models, which have been briefly discussed in Section~\ref{sec:fitness-evaluation-and-selection}.
While local Fourier analysis (LFA) has been successfully applied to different applications~\cite{rodrigo2017validity}, the missing availability of broadly-tested open-source tools for its automated use currently limits its applicability to the approach presented in this thesis~\cite{schmitt2020constructing}.
In case this situation might not change in the future, an alternative would be to use a statistical model to predict the quality of a certain multigrid method based on a history of samples.
One such approach, which has proven to be successful in other domains of automated algorithm design, such as automated machine learning (AutoML)~\cite{kotthoff2019auto,snoek2012practical}, is Bayesian optimization (BO)~\cite{frazier2018tutorial}.
While, in contrast to evolutionary algorithms, there is usually a limit to the number of parameters that can be optimized using BO, recent work in the field of AutoML demonstrates that these methods can also be scaled to grammar-based search spaces~\cite{schrodi2022towards}.
Finally, if we think in broader terms and consider the complete field of scientific computing, automation does not need to stop at the solver level but could include the full simulation design space.
For instance, in many cases, the choice of a suitable solver is not the main issue, but the considered system consists of different physical domains that each need to be simulated separately~\cite{gomes2018co}.
Only the coupling of these submodules then yields an understanding of the dynamic behavior of the complete system.
Similar to the recent success of automated algorithm design methods in the domain of data science, where systems exhibit a similarly high degree of complexity, the development of a unified approach for the automated design of simulation-based systems could yield great benefits in the future.

%% file: contents/appendix.tex
\section{Intermediate Representation}
\label{appendix:ir}
\begin{listing}[!htb]
\begin{minted}{python}
class InterGridOperator(Operator):

    def __init__(self, name, grid, fine_grid, coarse_grid, stencil_generator):
        self._fine_grid = fine_grid
        self._coarse_grid = coarse_grid
        super().__init__(name, grid, stencil_generator)

    @property
    def fine_grid(self):
        return self._fine_grid

    @property
    def coarse_grid(self):
        return self._coarse_grid
\end{minted}
	\caption{IR -- Inter-Grid Operator}
	\label{code:ir:inter-grid-operator}
\end{listing}
\begin{listing}[!htb]
\begin{minted}{python}
from operator import mul
from functools import reduce

class Restriction(InterGridOperator):
    def __init__(self, name, fine_grid, coarse_grid, stencil_generator=None):
        super().__init__(name, coarse_grid, fine_grid, coarse_grid, stencil_generator)
        tmp1 = reduce(mul, fine_grid.size)
        tmp2 = reduce(mul, coarse_grid.size)
        self._shape = (tmp2, tmp1)

    @property
    def fine_grid(self):
        return self._fine_grid

    @property
    def coarse_grid(self):
        return self._coarse_grid
\end{minted}
	\caption{IR -- Restriction}
	\label{code:ir:restriction}
\end{listing}
\begin{listing}[!htb]
\begin{minted}{python}
from operator import mul
from functools import reduce

class Prolongation(InterGridOperator):
    def __init__(self, name, fine_grid, coarse_grid, stencil_generator=None):
        super().__init__(name, fine_grid, fine_grid, coarse_grid, stencil_generator)
        tmp1 = reduce(mul, fine_grid.size)
        tmp2 = reduce(mul, coarse_grid.size)
        self._shape = (tmp1, tmp2)

    @property
    def fine_grid(self):
        return self._fine_grid

    @property
    def coarse_grid(self):
        return self._coarse_grid
\end{minted}
	\caption{IR -- Prolongation}
	\label{code:ir:prolongation}
\end{listing}
\begin{listing}[!htb]
\begin{minted}{python}
import evostencils.stencils as stencils

class Diagonal(UnaryExpression):
    def generate_stencil(self):
        return stencils.multiple.diagonal(self.operand.generate_stencil())

class BlockDiagonal(UnaryExpression):
    def __init__(self, operand, block_size):
        self._block_size = block_size
        super().__init__(operand)

    def generate_stencil(self):
        operand_stencil = self.operand.generate_stencil()
        return stencils.multiple.block_diagonal(operand_stencil, self.block_size)

    @property
    def block_size(self):
        return self._block_size
\end{minted}
	\caption{IR -- Diagonal and Block-Diagonal}
	\label{code:ir:diagonal}
\end{listing}
\begin{listing}[!htb]
\begin{minted}{python}
import evostencils.stencils as stencils

class Multiplication(BinaryExpression):

    def __init__(self, operand1, operand2):
        assert operand1.shape[1] == operand2.shape[0], "Operand shapes are not aligned"
        self._shape = (operand1.shape[0], operand2.shape[1])
        super().__init__(operand1, operand2)

    @property
    def grid(self):
        return self.operand1.grid

    @property
    def shape(self):
        return self._shape
\end{minted}
	\caption{IR -- Operator Application}
	\label{code:ir:multiplication}
\end{listing}
\clearpage
\section{Genetic Programming}
\label{appendix:gp}
\begin{listing}[!htb]
\begin{minted}{python}
from deap import gp

class PrimitiveSetTyped(gp.PrimitiveSetTyped):
    def _add(self, prim):
        def addType(dict_, ret_type):
            if ret_type not in dict_:
                new_list = []
                for type_, list_ in dict_.items():
                    # Use __equal__ instead of issubclass
                    if ret_type == type_:
                        for item in list_:
                            if item not in new_list:
                                new_list.append(item)
                dict_[ret_type] = new_list

        addType(self.primitives, prim.ret)
        addType(self.terminals, prim.ret)

        self.mapping[prim.name] = prim
        if isinstance(prim, gp.Primitive):
            for type_ in prim.args:
                addType(self.primitives, type_)
                addType(self.terminals, type_)
            dict_ = self.primitives
        else:
            dict_ = self.terminals

        for type_ in dict_:
            if type_ == prim.ret:
                dict_[type_].append(prim)

\end{minted}
	\caption{PrimitiveSetTyped}
	\label{code:gp:primitive-set-typed}
\end{listing}
\begin{listing}[!htb]
\begin{minted}{python}
import random
from inspect import isclass

def generate(pset, min_height, max_height, condition, return_type=None, subtree=None):
    type_ = return_type
    if type_ is None:
        type_ = pset.ret
    expression = []
    height = random.randint(min_height, max_height)
    stack = [(0, type_)]
    max_depth = 0
    subtree_inserted = False
    if subtree is None:
        subtree_inserted = True
    while len(stack) != 0:
        depth, type_ = stack.pop()
        if not subtree_inserted and type_ == return_type and len(expression) > 0:
            expression.extend(subtree)
            subtree_inserted = True
            continue
        max_depth = max(max_depth, depth)
        terminals_available = len(pset.terminals[type_]) > 0
        if condition(height, depth):
            nodes = pset.terminals[type_] + pset.primitives[type_]
        else:
            if terminals_available:
                nodes = pset.terminals[type_]
            else:
                nodes = pset.primitives[type_]
        if len(nodes) == 0:
            raise RuntimeError(f"Neither terminal nor primitive available for {type_}")
        choice = random.choice(nodes)
        if choice.arity > 0:
            for arg in reversed(choice.args):
                stack.append((depth + 1, arg))
        else:
            if isclass(choice):
                choice = choice()
        expression.append(choice)
    return expression
\end{minted}
	\caption{Adapted Implementation of DEAP's Generate Function}
	\label{code:gp:generate}
\end{listing}
\clearpage
\section{Full-Approximation Scheme}
\label{appendix:fas}
\begin{bnf*}
\bnfprod{$s_{h}$} {
	\bnfts{\textnormal{\textsc{update}}}(\bnfts{$\omega$}, \bnfsp \bnfes, \bnfsp \bnfts{\textnormal{\textsc{cgc}}}(\bnfts{$I_{2h}^{h}$}, \bnfts{$I_{h}^{2h}$}, \bnfsp \bnfpn{$s_{2h}$}))
 }
 \\
\bnfprod{$s_{2h}$} {
    \bnfts{\textnormal{\textsc{cgs}}}(\bnfts{$I_{4h}^{2h}$}, \bnfsp \bnfts{$A_{4h}$}, \bnfsp \bnfts{$I_{2h}^{4h}$}, \bnfpn{$s_{2h}$})
}
\\
\bnfprod{$c_{2h}$} {
	\bnfts{\textnormal{\textsc{coarsening}}}(\bnfts{$A_{2h}$}, \bnfsp \bnfts{$x^0_{2h}$}, \bnfsp, \bnfts{$I_{h}^{2h}$}, \bnfsp \bnfts{\textnormal{\textsc{apply}}}(\bnfts{$I_h^{2h}$}, \bnfsp \bnfpn{$c_h$}))
}
\end{bnf*}
\begin{table}[!htb]
	%\caption{State transition functions adapted for FAS}
	%\label{table:fas-grammar-semantics}
	\begin{algorithmic}
	\Function{coarsening}{$A_{2h}$, $x_{2h}^0$, $I_{h}^{2h}$, ($x_h$, $b_{h}$, $c_{2h}$, $Z_{h/2}$)}
	\State $x_{2h} \gets x_{2h}^0$ 
	\State $b_{2h} \gets c_{2h} + A_{2h} \left( I_{h}^{2h} x_h \right) $
	\State $c_{2h} \gets b_{2h} - A_{2h} \left( x_{2h} \right)$ 
	\State $Z_h \gets$ ($x_{h}$, $b_{h}$, $\lambda$, $Z_{h/2}$)
	\State return ($x_{2h}$, $b_{2h}$, $c_{2h}$, $Z_h$)
	\EndFunction
	\State
	\Function{cgc}{$I_{2h}^{h}$, $I_{h}^{2h}$, $(x_{2h}, b_{2h}, \lambda, Z_{h})$}
	\State ($x_h$, $f_{h}$, $c_h$, $Z_{h/2}$) $\gets Z_{h}$
	\State $c_h \gets I_{2h}^{h} \cdot (x_{2h} - I_{h}^{2h} x_h)$
	\State return ($x_h$, $f_{h}$, $c_h$, $Z_{h/2}$)
	\EndFunction
    \State
	\Function{cgs}{$I_{2h}^{h}$, $A_{2h}$, $I_{h}^{2h}$, ($x_h$, $b_{h}$, $c_{h}$, $Z_{h/2}$)}
	\State $x_{h} \gets I_{2h}^{h} \left( A_{2h}^{-1}\left( I_{h}^{2h} c_h \right) - I_{h}^{2h} x_h \right)$ 
	\State return ($x_{h}$, $b_{h}$, $c_{h}$, $Z_{h/2}$)
	\EndFunction
	\end{algorithmic}
\end{table}